\documentclass[12pt]{article}
\usepackage{amscd,amsfonts,amsmath,amssymb}
\usepackage[french]{babel}
\font\tenmsb=msbm10
\font\sevenmsb=msbm7
\font\fivemsb=msbm5

\newfam\msbfam
\textfont\msbfam=\tenmsb
\scriptfont\msbfam=\sevenmsb
\scriptscriptfont\msbfam=\fivemsb
\def\Bbb#1{{\fam\msbfam #1}}

\input xy
\xyoption{all}
\usepackage[all]{xy}

\makeatletter
\ifnum\@ptsize=0  \fi \ifnum\@ptsize=2
 \fi 

\newcommand\sF{{\cal F}}
\newcommand\sG{{\cal G}}

\newcommand\sO{{\cal O}}
\newcommand\sS{{\cal S}}

\newcommand\sX{{\cal X}}
\newcommand\sC{{\cal C}}

\newcommand\bZ{{\Bbb Z}}
\newcommand\bC{{\Bbb C}}
\newcommand\bQ{{\Bbb Q}}
\newcommand\bN{{\Bbb N}}

\newcommand\bP{{\Bbb P}}

\newtheorem{theorem}{Th\' eor\`eme}[section]
\newtheorem{lemma}[theorem]{Lemme}
\newtheorem{corollary}[theorem]{Corollaire}
\newtheorem{proposition}[theorem]{Proposition}
\newtheorem{question}[theorem]{Question}
\newtheorem{re}[theorem]{Remarque}

\newtheorem{definition}[theorem]{Definition}
\newtheorem{conjecture}[theorem]{Conjecture}

\newtheorem{example}[theorem]{Exemple}
\newtheorem{notation}[theorem]{Notation}

\begin{document}

\title {ORBIFOLDES G\'EOM\'ETRIQUES SP\'ECIALES ET CLASSIFICATION BIM\'EROMORPHE DES VARI\'ET\' ES K\" AHL\'ERIENNES COMPACTES}

\author{Fr\'ed\'eric Campana}

\maketitle

\
\tableofcontents
\


\section{INTRODUCTION}

\subsection{Abstract}

This is the ``geometric-orbifold" version of [Ca01/04]. We define the bimeromorphic {\it category} of geometric orbifolds. These interpolate between (compact K\" ahler) manifolds and such manifolds with logarithmic structure, and may be considered as ``virtual" ramified covers of the underlying manifold. These geometric orbifolds are here considered as fully geometric objects, and thus come naturally equipped with the usual invariants of varieties: morphisms and bimeromorphic maps, differential forms, fundamental groups and universal covers, Kobayashi pseudometric, fields of definition and rational points. The general expectation is that their geometry is qualitatively the same as that of manifolds with similar invariants. The most elementary of such geometric properties are established here, by direct adaptation of the arguments of [Ca01].

The motivation is that the natural frame for the theory of bimeromorphic classification of compact K\" ahler (and complex projective) manifolds without orbifold structure unavoidably seems to be the category of geometric orbifolds, as shown here (and in [Ca01] for manifolds) by the fonctorial two-step decomposition of arbitray manifolds: first by the ``core" fibration with fibres special and orbifold base of general type. Then {\it special} orbifolds are canonically (but conditionally) decomposed in towers of orbifolds with fibres having either $\kappa_+=-\infty$ or $\kappa=0$. An orbifold is {\it special} if it does not map ``stably" onto a (positive-dimensional) orbifold of general type, while having $\kappa_+=-\infty$ means that it maps only onto orbifolds  having $\kappa=-\infty$, and is expected to mean rationally connected in the orbifold category. In a slightly different context, and for seemingly different reasons, the Log minimal model program, also considers ``pairs" because most proofs naturally work, inductively on the dimension, only after the adjunction of a ``boundary". 

Moreover, fibrations enjoy in the bimeromorphic category of geometric orbifolds extension (or ``additivity") properties {\it not satisfied} in the category of varieties without orbifold structure, permitting to express invariants of the total space as the extension (or ``sum") of those of the generic fibre and of the base. For example, the natural sequence of fundamental groups always becomes exact  in the orbifold category. Also the total space of a fibration is special if so are the generic orbifold fibre and the orbifold base. In fact, geometric orbifolds were initially introduced precisely to remedy this last lack of ``additivity". 

 This makes this category of geometric orbifolds suitable to lift properties from orbifolds having either $\kappa_+=-\infty$ or $\kappa=0$ to those which are special. And even leads to expect that {\it specialness} is the exact geometric characterisation of some important properties (such as potential density or vanishing of the Kobayashi pseudometric).

There are still many open basic problems reated to the bimeromorphic equivalence in this orbifold category which need to be studied (such as its extension to the log-canonical case).

\subsection{Introduction}

L'objectif du texte est en priorit\' e de {\it d\'efinir} et d'\' etablir les propri\'et\'es de base de nouveaux objets (les ``orbifoldes g\'eom\'etriques") qui semblent \^etre essentiels pour la compr\' ehension de la structure birationnelle des vari\'et\'es projectives ou K\" ahl\'eriennes compactes, et en donnent une vue synth\'etique globale tr\`es simple. Les d\'emonstrations donn\'ees reposent cependant sur les techniques usuelles de la g\'eom\'etrie alg\'ebrique/analytique. De nombreuses questions ou conjectures \`a leur sujet sont \'egalement formul\'ees.

Le pr\'esent texte est la suite de [Ca01/04]. Dans ce texte \'etait introduite la notion de {\it base orbifolde} $(Y\vert\Delta(f))$ d'une fibration $f:X\to Y$, avec $X$ compacte K\" ahler, cette base orbifolde \'etant vue comme un rev\^etement ramifi\'e ``virtuel" de $Y$ \'eliminant virtuellement les fibres multiples en codimension $1$ de $f$. Ces fibres multiples forment l'obstruction principale \`a exprimer les invariants g\'eom\'etriques fondamentaux de $X$ comme extension (ou ``somme") de ceux de la fibre g\'en\'erique $X_y$, et de la base $Y$ de $f$. On y avait introduit, de plus, une nouvelle classe de vari\'et\'es, dites {\it sp\'eciales}. Ce sont, par d\'efinition, celles n'ayant pas de fibration sur une base orbifolde de type g\'en\'eral. Un $X$ arbitraire y \'etait scind\'e, \`a l'aide d'une unique fibration fonctorielle (son ``coeur"), en ses parties antith\'etiques: sp\'eciale (les fibres) et de type g\'en\'eral (la base orbifolde).  On y avait d\'ecompos\'e\footnote{Conditionnellement en une version orbifolde $C^{orb}_{n,m}$ de la conjecture $C_{n,m}$ d'Iitaka.} fonctoriellement toute vari\'et\'e sp\'eciale en tours de fibrations \`a fibres {\it orbifoldes} ayant soit $\kappa_+=-\infty$ (version conjecturale de la connexit\'e rationnelle), soit $\kappa=0$. La g\'eom\'etrie sp\'eciale apparaissant ainsi comme la combinaison au sens orbifolde de ces deux g\'eom\'etries classiques.

Cette d\'ecomposition (et aussi la th\'eorie des mod\`eles minimaux o\`u l'adjonction d'un ``bord" joue un r\^ole crucial dans les d\'emonstrations, bas\'ees sur une r\' ecurrence sur la dimension) semble indiquer que le cadre naturel de la th\'eorie de la classification des vari\'et\'es (projectives complexes ou K\" ahl\'eriennes compactes) est la {\it cat\'egorie} bim\'eromorphe des orbifoldes g\'eom\'etriques. Cette cat\'egorie restant, dans une premi\`ere \'etape, \`a d\'efinir.

L'objectif principal du pr\'esent texte est de le faire (du moins dans le cas ``lisse"), et d'en d\'evelopper les propri\'et\'es de base les plus simples (celles dont la d\'emonstration s'adapte plus ou moins directement du cas des vari\'et\'es \`a celle des orbifoldes). Le point de vue adopt\'e est que ces objets, qui interpolent entre les vari\'et\'es sans structure orbifolde et les vari\'et\'es logarithmiques, sont des objets g\'eom\'etriques \`a part enti\`ere, et sont donc \'equip\'es de tous les attributs des vari\'et\'es usuelles: morphismes, transformations bim\'eromorphes, formes diff\'erentielles, groupe fondamental et rev\^etement universel, corps de d\'efinition, points rationnels. Nous esp\'erons d\'evelopper des aspects plus profonds ult\'erieurement. En particulier, seules les orbifoldes lisses sont \'etudi\'ees ici. Il semble indispensable de consid\'erer un cadre plus g\'en\'eral (celui des orbifoldes ``klt" ou ``l.c") pour les d\'eveloppements ult\'erieurs.

Notons que les fibrations poss\`edent bien, dans cette cat\'egorie, les propri\'et\'es attendues d'additivit\'e qui font d\'efaut dans celle des vari\'et\'es sans structure orbifolde, et qui sont \`a l'origine de leur introduction. Par exemple, la suite des groupes fondamentaux devient exacte dans cette cat\'egorie. De m\^eme, l'espace total d'une fibration \`a fibres et base orbifoldes sp\'eciales est sp\'eciale. Cette cat\'egorie devrait ainsi permettre, par {\it d\'evissage}, de ``relever" \`a la ``g\'eom\'etrie sp\'eciale" certaines des propri\'et\'es attendues des ``g\'eom\'etries $\kappa_+=-\infty$ et $\kappa=0$". Et m\^eme de conjecturer que certaines propri\'et\'es importantes (densit\'e potentielle et nullit\'e de la pseudom\'etrique de Kobayashi) caract\'erisent exactement les orbifoldes sp\'eciales.

Remarquons que les invariants birationnels des orbifoldes logarithmiques lisses fournissent de nouveaux invariants pour les vari\'et\'es quasi-projectives lisses.

Le contenu du texte est le suivant: le \S2 introduit la cat\'egorie des orbifoldes, les morphismes \'etant d\'efinis de $3$ fa\c cons diff\'erentes (\'equivalentes): multiplicit\'es (voie g\'eom\'etrique; voir d\'efinition \ref{morphorb}, qui est probablement la contribution principale du pr\'esent texte), pr\'eservation des faisceaux de formes diff\'erentielles, et disques testants. 
Le \S3 d\'efinit la base orbifolde ``stable"  d'une fibration $f:(X\vert\Delta)\to Y$, et calcule la base orbifolde ``stable" d'une compos\'ee. 
Dans le \S4, on \'etablit l'invariance bim\'eromorphe de la dimension de Kodaira de la base orbifolde stable d'une fibration, r\'esultat utilis\'e constamment dans la suite.  
Le \S5 d\'efinit les courbes rationnelles orbifoldes, et pose la question de savoir si leurs propri\'et\'es sont analogues \`a celles du cas non orbifolde. On montre que c'est bien le cas lorsque l'orbifolde g\'eom\'etrique consid\'er\'ee est un ``quotient global" au sens de \ref{dqg}. La consid\'eration du champ alg\'ebrique associ\'e devrait permettre de traiter le cas g\'en\'eral avec des arguments similaires si les multiplicit\'es sont enti\`eres. 
Le \S6 contient des rappels extraits de [Ca01] sur l'additivit\'e de la dimension de Kodaira dans le cadre orbifolde. C'est le r\'esultat technique central du pr\'esent texte. 
Le \S8 contient des pr\'eliminaires techniques sur les fibrations de type g\'en\'eral, n\'ecessaires pour la construction du coeur, effectu\'ee au \S9, apr\`es un expos\'e au \S7 des propri\'et\'es et exemples de base des orbifoldes sp\'eciales. 
Au \S10, on d\'ecompose (conditionnellement en $C_{n,m}^{orb})$ fonctoriellement le coeur en tour de fibrations \`a fibres orbifoldes ayant soit $\kappa_+=-\infty$, soit $\kappa=0$. Ce d\'evissage permet de formuler des conditions sous lesquelles ``relever" aux orbifoldes sp\'eciales les propri\'et\'es attendues des orbifoldes ayant soit $\kappa_+=-\infty$, soit $\kappa=0$. 
Dans le \S11, on d\'efinit et \'etudie le groupe fondamental et le rev\^etement universel, et on construit la $\Gamma$-r\'eduction dans le cadre orbifolde (avec morphismes au sens divisible). 
Dans le \S12, on \'enonce des conjectures qui \'etendent directement au cadre orbifolde celles de [Ca01], et r\'esultent, pour la plupart d'entre elles, du ``rel\`evement" du \S 10. 

De nombreux probl\`emes de base de la g\'eom\'etrie bim\'eromorphe dans le contexte orbifolde restent \`a \'etudier (dont l'extension au cas log-canonique).

 Le projet d'\' etendre aux orbifoldes g\'eom\'etriques les r\'esultats de [Ca01] y \'etait  d\'ej\`a propos\'e. Il a \'et\'e aussit\^ot  abord\'e par S. Lu dans $math.AG/0211029$. Une notion de morphisme orbifolde y est \'evoqu\'ee en termes de ``relevant multiplicities". Certaines des constructions de [Ca01] sont ainsi directement utilis\'ees, mais sans interpr\'etation g\'eom\'etrique, indispensable pour les applications et d\'eveloppements ult\'erieurs.

 Je voudrais, par ailleurs, remercier D. Greb et K. Jabbusch pour m'avoir signal\'e une erreur et des impr\'ecisions dans la version initiale du pr\'esent texte, ainsi que A. Levin pour une tr\`es int\'eressante discussion, et des r\'ef\'erences (dont [B 87]), sur les aspects hyperbolique et arithm\'etique (voir le \S \ref{sha}).

\section{LA CAT\' EGORIE DES ORBIFOLDES G\'EOM\'ETRIQUES}

\subsection{Diviseurs orbifoldes}

\begin{definition}\label{orbgeom} Soit $X$ un espace analytique complexe normal (d\'enombrable \`a l'infini; $X$ sera, de plus, suppos\'e compact et connexe dans la suite de ce texte, \`a l'exception du pr\'esent \S2 essentiellement) . On note $W(X)$ l'ensemble des diviseurs de Weil irr\'eductibles de $X$.

Une {\bf multiplicit\'e orbifolde} sur $X$ est une application $m: W(X)\to (\bQ\cup +\infty):=\overline{\Bbb Q^+}$ telle que $m(D)\geq 1$ pour tout $D\in W(X)$, et telle que $m(D)=1$ pour localement presque tout $D\in W(X)$ (ie: pour tout compact $K$ de $X$, $m(D)=1$ pour tous les $D\in W(X)$ rencontrant $K$, sauf un nombre fini d'entre eux). Une telle multiplicit\'e orbifolde est dite {\bf enti\`ere} si elle prend ses valeurs dans $(\bN\cup +\infty)$.

Une multiplicit\'e orbifolde $m$ sera dite {\bf finie} si elle ne prend pas la valeur $+\infty$.

Le {\bf diviseur orbifolde} associ\'e est $\Delta:=\sum_{D\in W(X)}(1-\frac{1}{m(D)}).D$. C'est une somme localement finie d\'efinissant un $\bQ$-diviseur effectif sur $X$. (Par convention, $1/(+\infty):=0)$. On notera $m_{\Delta}: W(X)\to (\bQ\cup +\infty)$ la multiplicit\'e orbifolde d\'efinissant $\Delta$.

Le {\bf support} de $\Delta$, not\'e $supp(\Delta)$, ou $\lceil \Delta\rceil$, est la r\'eunion des $D\in W(X)$ tels que $m(D)>1$.

Une {\bf orbifolde g\'eom\'etrique}\footnote{Le terme ``g\'eom\'etrique" provient de ce que l'on ne conserve de la structure orbifolde que le support (un diviseur) et les multiplicit\'es, mais que l'on ne se donne pas d'action locale d'un groupe. En ce sens, une ``orbifolde g\'eom\'etrique" peut \^etre vue comme la trace g\'eom\'etrique en codimension $1$ d'une orbifolde au sens de la terminologie existante.} est un couple, not\'e $(X\vert\Delta)$, de la forme pr\'ec\'edente. Une orbifolde g\'eom\'etrique $(X\vert\Delta)$ sera dite {\bf finie} si sa multiplicit\'e $m_{\Delta}$ est finie, et {\bf enti\`ere} si sa multiplicit\'e l'est.

On dira que $(X\vert\Delta)$ est {\bf lisse} si $X$ est lisse et si $supp(\Delta)$ est un diviseur \`a croisements normaux.
\end{definition}

\begin{re}\label{rmk1} 1. Lorsque $m$ ne prend que les valeurs $1$ et $+\infty$, le diviseur $\Delta$ est entier, r\'eduit, de multiplicit\'e $1$. On dira que l'orbifolde g\'eom\'etrique $(X\vert\Delta)$ est {\bf ouverte}, ou: {\bf logarithmique}. On notera $(X\vert D)$ une telle orbifolde g\'eom\'etrique. Le cas g\'en\'eral interpole donc entre les cas propre (o\`u $\Delta=0)$ et logarithmique.

Toute vari\'et\'e lisse et quasi-projective $U$ s'\'ecrit: $U=X-D$, avec $(X\vert D)$ lisse. Les propri\'et\'es de $U$ sont celles de $(X\vert D)$ ne d\'ependant que de $X-D$. Ce sont en fait les propri\'et\'es ``birationnelles" de $(X\vert D)$ au sens d\'efini ci-dessous. Nous construirons ainsi de nouveaux invariants des vari\'et\'es quasi-projectives lisses.

On \'ecrira aussi $\Delta=\sum_{j\in J} (1-\frac{1}{m_j}).D_j$ si $J$ est un sous-ensemble localement fini de $W(X)$ contenant tous les $D$ tels que $m(D)>1$, avec $m_j=m(D_j),\forall j\in J$. 

2. On peut consid\'erer, plus g\'eneralement, des fonctions de multiplicit\'e \`a valeurs r\'eelles dans $ \lbrace[1,+\infty[\cup \lbrace +\infty\rbrace\rbrace$. Les d\'efinitions et r\'esultats ci-dessous s'appliquent, ainsi que leurs d\'emonstrations, avec des modifications mineures, \`a ce cadre \'elargi (susceptible d'autres applications). 

3. Si $\Delta, \Delta'$ sont des diviseurs orbifoldes sur $X$, de fonctions de multiplicit\'e $m,m'$ respectivement, on dira que $\Delta\geq \Delta'$ si $m\geq m'$. (Ceci signifie en effet que $(\Delta-\Delta')$ est un diviseur effectif). Si ces deux orbifoldes g\'eom\'etriques sont enti\`eres, on dira que $\Delta'$ {\bf divise} $ \Delta$ (not\'e $m_{\Delta'}(D) \vert m_{\Delta'}(D)$ divise $m_{\Delta}(D), \forall D\in W(X)$ (avec la convention: $n$ divise $+\infty, \forall n>0,$ entier).

4. On d\'efinit de mani\`ere \'evidente le produit de deux orbifoldes g\'eom\'etriques. Ce produit est donc lisse (resp. entier, fini) si les facteurs le sont. 

5. Si $\Delta,\Delta'$ sont deux (ou m\^eme une famille d') orbifoldes g\'eom\'etriques sur $X$ on d\'efinit de mani\`ere \'evidente $sup\{\Delta,\Delta'\}$ et $inf\{\Delta,\Delta'\}$.

6. D. Abramovich a introduit, dans {[Abr 07]}, le terme de {\bf constellation} et pour son raffinement toroidal, qui lui est d\^u, celui de {\bf firmament}. La notion de ``constellation" consiste en la donn\'ee d'un syst\`eme (compatible) de diviseurs orbifoldes g\'eom\'etriques sur toutes les modifications propres de $X$, consid\'er\'ees simultan\'ement. Nous tentons ici de consid\'erer les propri\'et\'es bim\'eromorphes (d\'efinies de mani\`ere ad\'equate) des ``orbifoldes g\'eom\'etriques" individuelles, sans inclure la totalit\'e de la constellation associ\'ee dans la donn\'ee initiale. La diff\'erence entre les deux notions semble cependant inessentielle. 

D'ailleurs, dans le cas consid\'er\'e ici, o\`u $X$ est projective lisse et  le support de $\Delta$ \`a croisement normaux, la notion d'orbifolde g\'eom\'etrique est un cas particulier de celle de ``champ alg\'ebrique" lisse de Deligne-Mumford, et le terme est donc compatible (au niveau des objets, sinon des morphismes d\'efinis ci-dessous) avec les terminologies ant\'erieures.  

L' expression ``orbifolde g\'eom\'etrique" sera parfois abr\'eg\'ee en ``orbifolde" dans la suite. On parlera ainsi de ``morphismes orbifoldes, fibre ou base orbifolde" d'une fibration. Ceci ne devrait pas cr\'eer de confusion avec la terminologie ant\'erieure, puisque les ``orbifoldes" d\'ej\`a existantes ne sont jamais consid\'er\'ees dans le pr\'esent texte. \end{re}

\subsection{Invariants: dimension canonique, groupe fondamental, points entiers.}\label{invorb}

Le {\bf principe} que nous voudrions illustrer dans le texte qui suit est le suivant:

1. Les orbifoldes g\'eom\'etriques (lisses) sont des objets g\'eom\'etriques \`a part enti\`ere, au m\^eme titre que les vari\'et\'es complexes (projectives ou compactes): on peut leur attribuer en particulier les invariants g\'eom\'eriques d\'efinis ci-dessous, ainsi que des notions de morphismes et d'\'equivalence bim\'eromorphe.  

2. Leurs propri\'et\'es g\'eom\'etriques sont les m\^emes que celles des vari\'et\'es (sans structure orbifolde g\'eom\'etrique) ayant des invariants analogues.

3. Nombre de ces propri\'et\'es sont \'etablies en adaptant (en g\'en\'eral sans difficult\'es majeures) les d\'emonstrations \'etablissant celles des vari\'et\'es sans structure orbifolde g\'eom\'etrique. Pour certaines propri\' et\'es cependant, l'adaptation semble requ\'erir des id\'ees nouvelles. Des exemples de propri\'et\'es pour lesquelles l'adaptation n'est cependant pas imm\'ediate sont: celles des courbes rationnelles orbifoldes, ou celles (de finitude ou d'ab\'elianit\'e) du groupe fondamental des orbifoldes g\'eom\'etriques qui sont soit  Fano, soit \`a fibr\'e canonique trivial.

4. L'origine des orbifoldes g\'eom\'etriques (l'\'elimination virtuelle des fibres multiples en codimension $1)$ conduit naturellement \`a les consid\'erer comme des objets g\'eom\'etriques. Le LMMP (programme des Log-mod\`eles minimaux), qui consid\`ere ces m\^emes objets pour des raisons apparemment diff\'erentes (la possibilit\'e de faire des r\'ecurrences sur la dimension par extraction de diviseurs \`a partir de ``centres Log-canoniques"), ne les munit au contraire que d'un seul invariant: le fibr\'e canonique $K_X+\Delta$, consid\'er\'e pour ses propri\'et\'es num\'eriques d'intersection.

$\square$ Si $(X\vert\Delta)$ est une orbifolde g\'eom\'etrique, on d\'efinit (voir [Ca04], les d\'efinitions d\'etaill\'ees sont aussi donn\'ees  ci-dessous, dans les chapitres correspondants ):

1. Son fibr\'e canonique $K_{X\vert\Delta}:=K_X+\Delta$: c'est un $\Bbb Q$-diviseur sur $X$.

2. Sa dimension canonique (ou de Kodaira\footnote{Habituellement appel\'ee ``dimension de Kodaira", elle est en fait introduite pour les surfaces dans le livre de Shafarevitch et al. sur la classification des surfaces projectives, et en g\'en\'eral par S. Iitaka et B. Moishezon (dont nous suivons la terminologie). K. Kodaira ne l'a utilis\'ee qu'en 1975, alors qu'elle \'etait d\'ej\`a d'usage courant.}): $\kappa(X\vert\Delta):=\kappa(X,K_X+\Delta)\geq \kappa(X)\in \lbrace-\infty,0,1, ...,dim(X)\rbrace$, lorsque $X$ est compacte et irr\'eductible. 

On dira que $(X\vert\Delta)$ est {\bf de type g\'en\'eral} si $\kappa(X\vert\Delta)=dim(X)>0$.

3. On d\'efinira aussi des faisceaux de formes (pluri)-diff\'erentielles et, plus g\'en\'eralement, de tenseurs holomorphes de tous les types sur une orbifolde g\'eom\'etrique lisse. (Voir \S\ref{fdiff} et \S \ref{tdiff} ci-dessous).

4. La notion de courbe $\Delta$-rationnelle (en trois versions).

$\square$ Lorsque $(X\vert\Delta)$ est, de plus, enti\`ere, on d\'efinit aussi:

5. Son groupe fondamental $\pi_1(X\vert\Delta)$ et son rev\^etement universel (voir [Ca04] et le \S\ref{gf} ci-dessous).

6. Sa pseudom\'etrique de Kobayashi $d_{X\vert\Delta}$ (voir [Ca04] et le \S\ref{pskob} ci-dessous).

7. Ses points entiers sur un corps de nombres, pour un mod\`ele donn\'e. (Voir [Ca05], et le \S\ref{arithm} ci-dessous).

8. La notion d'orbifolde g\'eom\'etrique peut naturellement \^etre d\'efinie en g\'eom\'etrie alg\'ebrique sur d'autres corps que $\bC$. Un cas int\'eressant est celui des corps de fonctions (m\'eromorphes sur une courbe projective d\'efinie sur  $\bC$ ou sur un corps fini). Voir [Ca01] ou le \S\ref{fonct} ci-dessous.

\subsection{Morphismes orbifoldes.}\label{mo}

La d\'efinition suivante est  la contribution principale du pr\'esent texte:

\begin{definition}\label{morphorb} Soit $f:Y\to X$ une application holomorphe entre espaces analytiques complexes normaux, et $\Delta_Y,\Delta_X$ des diviseurs orbifoldes sur $Y,X$ respectivement. On note  ici $m_Y,m_X$ les multiplicit\'es orbifoldes associ\'ees. On dit que $f$ induit un {\bf morphisme orbifolde} (not\'e alors $f:(Y\vert\Delta_Y)\to (X\vert\Delta_X))$ si:

1. $f(Y)$ n'est pas contenu dans $ \lceil\Delta_X\rceil$.

2. $X$ est $\bQ$-factorielle au sens alg\'ebrique (ie: tout diviseur de Weil irr\'eductible sur $X$ est $\bQ$-Cartier).

3. pour tout $D\in W(X)$, et tout $E\in W(Y)$, on a: $t.m_Y(E)\geq m_X(D)$, si $t>0$, o\`u: $t=t_{E,D}\in \bQ^+$ est
 tel que $f^*(D)=t_{E,D}.E+R$, avec $R$ un diviseur effectif de $Y$ ne contenant pas $E$.

\end{definition}

\begin{definition}\label{morphorb*}

Dans l'\'etude du groupe fondamental orbifolde (d\'efini seulement pour les orbifoldes g\'eom\'etriques enti\`eres, voir \S\ref{gf}), et aussi pour la notion de courbe $\Delta$-rationnelle (voir \S5), nous utiliserons une notion plus restrictive de morphisme orbifolde, celle de {\bf morphisme orbifolde divisible}: il s'agit d'un morphisme orbifolde dans le sens pr\'ec\'edent, mais satisfaisant la variante ``muliplicit\'es divisibles": 
$m_X(D)$ {\bf divise} $t.m_Y(E)$\footnote{Au lieu de: $t.m_Y(E)\geq m_X(D).$}, pour tous $D,E$ comme dans la condition 3 ci-dessus. 
(Par convention, $+\infty$ est multiple de tout entier $m>0$, et ne divise que lui-m\^eme).

Pour distinguer ces deux types de morphismes, nous appellerons ``morphismes non-classiques" ceux d\'efinis en \ref{morphorb} ci-dessus, et ``morphismes classiques", ou ``divisibles" ceux d\'efinis ici. 

On notera $Georb^Q$ la cat\'egorie des orbifoldes g\'eom\'etriques $\bQ$-factorielles munie des morphismes non-classiques, et $Georb^Z$ sa sous-cat\'egorie pleine constitu\'ee des orbifoldes enti\`eres. On notera en fin $Georb^{div}$ la cat\'egorie des orbifoldes g\'eom\'etriques $\bQ$-factorielles enti\`eres munie des morphismes ``classiques" (ou ``divisibles").

Pour indiquer qu'un morphisme orbifolde $f:(Y\vert\Delta_Y)\to (X\vert\Delta_X)$ appartient \`a l'une de ces trois cat\'egories, on le notera: $f:(Y\vert\Delta_Y)^*\to (X\vert\Delta_X)^*$, avec $*=div, Z, Q$ selon le cas.
 \end{definition}

\begin{re} \label{defmobim}

\

0. Les motivations pour la condition 3. de la d\'efinition \ref{morphorb} pr\'ec\'edente sont multiples, bien que peu \'evidentes {\it a priori}: tout d'abord, cette condition est exactement celle qui pr\'eserve les faisceaux de formes pluri-diff\'erentielles (voir proposition \ref{modiff}), ainsi que les morphismes du disque unit\'e pour les orbifoldes g\'eom\'etriques enti\`eres (voir proposition \ref{disctest}). 

Enfin, une motivation g\'eom\'etrique plus directe (dans le cas divisible) est la compatibilit\'e avec la composition des fibrations: soient $f:X\to Y$ et $g:Y\to Z$ deux fibrations entre vari\'et\'es complexes projectives normales. On suppose $f$ et $g$ \`a fibres \'equidimensionnelles pour simplifier.  

Soit $D$ un diviseur de Weil sur $Z$. Alors la multiplicit\'e de la fibre de $(g\circ f)$ au-dessus du point g\'en\'erique de $D$ est le $pgcd$, not\'e $\mu(D)$, des $t_j.m_j$ si $g^*(D)=\sum_jt_j.(E_j)$, o\`u $m_j=pgcd(s_{ij})$ , avec: $f^*(E_j)=\sum_is_{ij}.F_{ij}$. (Les entiers $t_j,s_{ij}$ sont bien d\'efinis au-dessus des points g\'en\'eriques de $D,E_j$, puisque $f,g$ sont \`a fibres \'equidimensionnelles, de sorte que tous les $E_j$ ont pour image $D$). Les $m_j$ sont donc simplement les multiplicit\'es usuelles des fibres de $f$ au-dessus des points g\'en\'eriques des $E_j$.

Si l'on d\'efinit une orbifolde g\'eom\'etrique $(Y\vert \Delta_f)$ sur $Y$ en posant: $m_{\Delta_f}(E_j):=m_j$ pour tout diviseur de Weil $E_j$ sur $Y$ comme ci-dessus, et une orbifolde g\'eom\'etrique $(Z\vert \Delta_{g\circ f})$ en posant $m_{\Delta_{g\circ f}}(D):=\mu(E)$, pour tout $E\subset Z$, alors: pour tout diviseur orbifolde $\Delta_Y$ sur $Y$, l'application $g:(Y\vert\Delta_Y)\to (Z\vert \Delta_{g\circ f})$ est un morphisme orbifolde (au sens de \ref{morphorb}) si et seulement si $\Delta_Y\geq \Delta_f$.

Cette situation sera \'etudi\'ee en d\'etail aux \S\ref{basmor} et \S\ref{comfib}.

1. Si $m_X(D)=1$ et si $X$ est factorielle, la condition 2. est vide (puisque $t>0$ est entier et $m_Y(E)\geq 1)$. Seuls les $D\subset \lceil\Delta_X\rceil$ et $E$ tels que $f(E)\subset  \lceil\Delta_X\rceil$ fournissent donc des conditions non vides. Si $X$ et $Y$ sont compacts, les conditions \`a v\'erifier sont donc en nombre fini.

2. Si $f:(Y\vert\Delta_Y)\to (X\vert\Delta_X)$ et $g:(Z\vert\Delta_Z)\to (Y\vert\Delta_Y)$ sont des morphismes, le compos\'e $f\circ g$ aussi. La d\'emonstration est imm\'ediate.

3. Si $\Delta_X=\Delta_Y=0$, toute application holomorphe $f:X\to Y$ est un morphisme orbifolde. Si  $f:(Y\vert\Delta_Y)\to (X\vert\Delta_X)$ est un morphisme orbifolde, et si $\Delta_Y^+\geq \Delta_Y$, alors  $f:(Y\vert\Delta_Y^+)\to (X\vert\Delta_X)$ est un morphisme orbifolde. De m\^eme: si $\Delta_X\geq \Delta^-_X$, alors  $f:(Y\vert\Delta_Y)\to (X\vert\Delta^-_{X})$ est un morphisme orbifolde.

4. Si $(X\vert\Delta_X)$ est logarithmique, alors  $f:(Y\vert\Delta_Y)\to (X\vert\Delta_X)$ est un morphisme orbifolde si et seulement si $\Delta_Y\geq f^{-1}(\Delta_X)$.

5. On dit que  $f:(Y\vert\Delta_Y)\to (X\vert\Delta_X)$ est un {\bf morphisme orbifolde bim\'eromorphe \'el\'ementaire} (ou une {\bf modification orbifolde \'el\'ementaire}) si c'est un morphisme orbifolde, si $f:Y\to X$ est birationnel, et si $f_*(\Delta_Y)=\Delta_X$. On suppose ici que les deux espaces $Y,X$ consid\'er\'es sont alg\'ebriquement $\bQ$-factoriels. Noter que les trois propri\'et\'es sont ind\'ependantes (deux d'entre elles n'impliquent pas la troisi\`eme). On dira que deux telles orbifoldes g\'eom\'etriques sont {\bf bim\'eromorphiquement \'equivalentes} s'il existe une chaine de morphismes orbifoldes bim\'eromorphes \'el\'ementaires les reliant. De mani\`ere \'equivalente: c'est la relation d'\'equivalence engendr\'ee par les morphismes orbifolde bim\'eromorphes \'el\'ementaires.

 Voir \ref{crem} et \ref{crem'} pour des exemples d'\'equivalences bim\'eromorphes non domin\'ees bim\'eromorphiquement par une troisi\`eme.

6. Si $f:(Y\vert\Delta_Y)\to (X\vert\Delta)$ et $f':(Y\vert\Delta_Y)\to (X\vert\Delta')$ sont des morphismes orbifoldes, ils se factorisent par $f^+:(Y\vert\Delta_Y)\to (X\vert\Delta^+)$, avec $\Delta^+:=sup\{\Delta,\Delta'\}$. De m\^eme,  si $f:(X\vert\Delta)\to (Y\vert\Delta_Y)$ et $f':(X\vert\Delta')\to (Y\vert\Delta_Y)$ sont des morphismes orbifoldes, ils se factorisent par $f^-:(X\vert\Delta^-)\to (Y\vert\Delta_Y)$, avec $\Delta^-:=inf\{\Delta,\Delta'\}$.

7. Si $f:Y\to (X\vert\Delta_X)$ est une application holomorphe propre et surjective avec $Y$ et $(X\vert\Delta_X)$ lisses, et $f^{-1}(\lceil\Delta_X\rceil)$ \`a croisement normaux, il existe un \'el\'ement minimum, not\'e $f^+(\Delta_X),$ parmi les diviseurs orbifoldes $\Delta_Y$ sur $Y$ tels que $f:(Y\vert\Delta_Y)\to (X\vert\Delta_X)$ soit un morphisme orbifolde. On l'appelle {\bf rel\`evement de $\Delta_X$ \`a $Y$}. Pour tout $E\in W(Y)$, la $f^+(\Delta_X)$-multiplicit\'e de $E$ est $m_Y(E):=max\lbrace 1, sup_D\{\frac{m(D)}{t_{E,D}}\} \rbrace$, $D\in W(X)$ tel que $f^*(D)=t_{E,D}.E+R$, avec $t_{E,D}>0$ et $R$ effectif, de support ne contenant pas $E$. 

Si l'on veut des orbifoldes enti\`eres (resp. et des morphismes divisibles), on doit remplacer $\lbrace \frac{m(D)}{t_{E,D}}\rbrace$ ci-dessus par $\lbrace\lceil \frac{m(D)}{t_{E,D}}\rceil\rbrace$ (resp. par: $ppcm_D\{\frac{m(D)}{pgcd\{m(D),t_{E,D}\}}\}$. 

On impose la condition $m_Y(E)\geq 1, \forall E\in W(Y)$, puisqu'elle n'est pas toujours satisfaite (par exemple si $f:Y\to X$ est un rev\^etement ramifi\'e de courbes projectives, si $\Delta_X=0$, et si $E$ est un point de ramification de $f)$.

On notera que si $g:Z\to Y$ est une seconde application holomorphe propre et surjective, alors $(f\circ g)^+(\Delta_X)\leq g^+(f^+(\Delta_X))$, mais que l'on n'a pas \'egalit\'e, en g\'en\'eral.

Par exemple: pour $f$, \'eclater  la surface $X$ en un point lisse du support (non vide, de dimension $1$) de $(\Delta_X)$, puis pour $g$, \'eclater l'intersection de la transform\'ee stricte du support de $\Delta_X$ avec le diviseur exceptionnel. Si la composante de $\Delta_X$ contenant le premier point \'eclat\'e est de multiplicit\'e $m\geq 2,$ la multiplicit\'e du diviseur exceptionnel du second \'eclatement dans $(f\circ g)^+(\Delta_X)$ (resp. dans $g^+(f^+(\Delta_X)))$ est $ m/2 $ (resp. $m)$.
\end{re}

\begin{definition} Soit $\Bbb D$ le disque unit\'e de $\Bbb C$, et $X$ lisse. On note $Hol(\Bbb D, (X\vert\Delta _X))$ l'ensemble, muni de la topologie de la convergence compacte, des applications holomorphes $h:\Bbb D\to X$ qui sont des morphismes orbidolde lorsque $\Bbb D$ (resp. $X)$ est muni du diviseur orbifolde vide (resp. $\Delta_X)$. (voir [C-W05]).
\end{definition}

\begin{proposition}\label{disctest} Soit $(Y\vert\Delta_Y), (X\vert\Delta_X)$ des orbifoldes g\'eom\'etriques, avec $Y,X$ lisses, et $f:Y\to X$ holomorphe.  Alors: 

L'application de composition: $f_*:Hol(\Bbb D, (Y\vert\Delta _Y))\to Hol(\Bbb D, X)$ d\'efinie par: $f_*(h):=f\circ h$ a son image contenue dans $Hol(\Bbb D, (X\vert\Delta _X))$ si $f:(Y\vert\Delta_Y)\to (X\vert\Delta_X)$ est un morphisme (d'orbifoldes g\'eom\'etriques). Si $(Y\vert\Delta_Y)$ est enti\`ere, la r\'eciproque est vraie (ie: $f$ est un morphisme orbifolde si 'image de $f_*$ est contenue dans $Hol(\Bbb D, (X\vert\Delta _X)))$. 
\end{proposition}

{\bf D\'emonstration:} Si $f$ est un morphisme d'orbifoldes, la remarque \ref{rmk1}(2)  pr\'ec\'edente montre que l'image de $f_*$ est contenue dans $Hol(\Bbb D, (X\vert\Delta _X))$. R\'eciproquement, avec les notations de \ref{morphorb}, supposons que $f^*(D)=t.E+R$, avec $t>0$. Supposons que $t.m_Y(E)<m_X(D)$. Soit $y\in E$ un point g\'en\'erique lisse, et $h:\Bbb D\to Y$ holomorphe telle que $h(0)=y$, et que $h^*(E)=m_Y(E).\lbrace 0\rbrace$, avec: $(h(\Bbb D)\cap \lceil \Delta_Y\rceil)=\lbrace 0\rbrace$. (On peut r\'ealiser ces conditions en restreignant suffisamment $\Bbb D$, puisque $(Y\vert\Delta_Y)$ est enti\`ere).  Alors: $(f\circ h)^*(D)=t.h^*(E)+h^*(R)=t.m_Y(E).\lbrace 0\rbrace$. Puisque $t.m_Y(E)<m_X(D)$, $f$ n'est donc pas un morphisme d'orbifoldes $\square$

\subsection{Faisceaux de formes diff\'erentielles sur les orbifoldes g\'eom\'etriques lisses.}\label{fdiff}

\

{\bf Notations:} Soit $(X\vert\Delta)$ une orbifolde g\'eom\'etrique {\bf lisse} (ie: $X$ est lisse, et $\lceil \Delta\rceil$ est \`a croisements normaux). On note $\Delta=\sum_{h\in H} (1-\frac{1}{m_h}).D_h=\sum_{h\in H} a_h.D_h$, $[a]=\lfloor a\rfloor$ \'etant la partie enti\`ere du r\'eel $a$, et $\lceil a\rceil:=-\lfloor (-a)\rfloor$ son ``arrondi sup\'erieur". Ici les multiplicit\'es $m_j$ sont donc, soit des rationnels, soit $+\infty$. 

Soit $q\geq 0$ un entier, et $\Omega^q_X(log \lceil \Delta\rceil)$ le faiceau des germes de formes diff\'erentielles \`a p\^oles logarithmiques dans $\lceil \Delta\rceil$. 

 Nous allons maintenant d\'efinir, pour une orbifolde g\'eom\'etrique {\it lisse} $(X\vert \Delta)$ les faisceaux $S^N(\Omega^q(X\vert\Delta))$, analogues des faisceaux $Sym^N(\Omega^q_X)$ lorsque $\Delta=0$. 

Localement, dans des coordonn\'ees locales $x=(x_1,...,x_p)$ {\bf  adapt\'ees \`a $\Delta$}, c'est-\`a-dire telles que $\Delta$ ait pour \'equation\footnote{symbolique: les $m_j$ sont les $\Delta$-multiplicit\'es des hyperplans de coordonn\'ees.}: $\Pi_{j=1}^{j=p} x_j^{(1-\frac{1}{m_j})}$, le faisceau $\Omega^q_X(log \lceil \Delta\rceil)$ admet comme $\cal O$$_X$-module les g\'en\'erateurs: $\frac{dx_J}{x_J}$, pour $J$ partie ordonn\'ee \`a $q$ \'el\'ements de $\lbrace 1,...,p\rbrace$. On a not\'e: $\frac{dx_J}{x_J}:=\wedge _{j\in J}\frac{dx_j}{x_j}$.

On note alors, pour tous les entiers non-n\'egatifs $q,N$: $S^N_q(X\vert\Delta):=S^{N}(\Omega^q(X\vert\Delta))$ le sous-faisceau analytique coh\'erent localement libre\footnote{C'est la propri\'et\'e cruciale, utilis\'ee constamment dans la suite. Elle permet de n\'egliger les sous-ensembles de codimension $2$ ou plus, et donc de ne faire intervenir que le lieu de codimension $1$ du diviseur orbifolde.} de $Sym^N(\Omega^q_X(log \lceil \Delta\rceil))$ engendr\'e par les \'el\'ements: $$\frac{dx_(J)}{x_(J)}:=x^{\lceil k/m\rceil}.\otimes_{l=1}^{l=N}\frac{dx_{_{J_l}}}{x_{_{J_l}}}=x^{\lceil -k.a\rceil}.\otimes_{l=1}^{l=N}dx_{_{J_l}}, (J)=(J_1,\dots,J_N)\},$$ d\'efinis comme suit:

1. Les $J_l$ sont les parties ordonn\'ees (croissantes) \`a $q$ \'el\'ements de $\lbrace 1,...,p\rbrace$. Les $N$-uplets $(J_1,\dots,J_N)$ consid\'er\'es sont croissants au sens large pour l'ordre lexicographique sur les parties ordonn\' ees \`a $q$ \'el\'ements de $\{1,\dots,p\}$. 

2. Pour tout $j=1,...,p$, on note $k_j$ le nombre d'occurrences de $j$ dans la suite $J_1,...,J_N$. (C'est-\`a-dire que $k_j=\sum_{l=1}^{l=N}k_{j,l}$, o\`u $k_{j,l}=1$ si $j\in J_l$, et $k_{j,l}=0$ sinon).

3. On a aussi not\'e $k/m$ le $p$-uplet $(k_1/m_1,...,k_p/m_p)$, et $\lceil k/m\rceil$ le $p$-uplet: $(\lceil k_1/m_1\rceil,...,\lceil k_p/m_p\rceil)$. Enfin: $x^{\lceil k/m\rceil}:=\Pi_{j=1}^{j=p}x_j^{\lceil k_j/m_j\rceil}$.

4. On a enfin not\'e $(-k.a)$ le $p$-uplet $(-k_1.a_1,...,-k_p.a_p)$, et $\lceil -k.a\rceil$ le $p$-uplet: $(\lceil -k_1.a_1\rceil,...,\lceil -k_p.a_p\rceil)$. Enfin: $x^{\lceil -k.a\rceil}:=\Pi_{j=1}^{j=p}x_j^{\lceil -k_j.a_j\rceil}$, avec: $a_j:=(1-\frac{1}{m_j})$.

 Cette d\'efinition est clairement ind\'ependante des cartes adapt\'ees locales utilis\'ees.

On v\'erifie alors ais\'ement que l'application naturelle $S^N_q(X\vert\Delta)\otimes S_M,q(X\vert\Delta)\to S_{N+M},q(X\vert\Delta)$ est bien d\'efinie (ie: prend bien ses valeurs dans le membre de droite).

\begin{re}\label{origdiff}

1. A nouveau, la d\'efinition pr\'ec\'edente s'applique avec changement mineur au cas o\`u les multiplicit\'es de $\Delta$ sont r\'eelles.

2. L'origine des faisceaux $S^N_q(X\vert\Delta)=S^N_q$ (d\'efinis dans [Ca01]) est la suivante, pour les orbifoldes g\'eom\'etriques enti\`eres: localement, dans la carte $x$ pr\'ec\'edente, $X\vert\Delta$ a un rev\^etement universel local $f:Y\to X$ donn\'e dans la carte $y=(y_1,\dots ,y_p)$ par: $f(y)=x=(x_1=:y_1^{m_1},\dots, x_p=:y_p^{m_p})$\footnote{Avec la convention: $x_i=y_i^{+\infty}:=exp(y_i),$ ou encore: $dx_i/x_i=dy_i$, si $m_i=+\infty$.}. 

Lorsque $N$ est suffisamment divisible, $f^*(S_{N,1})=Sym^N(\Omega_Y^1)$. 

Pour $N$ g\'en\'eral, $S^N_q$ est le plus grand sous-faisceau $\sF$ de $Sym^N(\Omega_X(log(\lceil\Delta\rceil))$ tel que $f^*(\sF)\subset Sym^N(\Omega_Y^q)$.

3. Ces faisceaux ont \'et\'e utilis\'es de mani\`ere cruciale dans le cadre de la th\'eorie de Nevanlinna dans [C-P05], avec $p=2,q=1$. 

4. En g\'en\'eral, $S^N_q(X\vert\Delta)$ contient, mais n'est pas \'egal \`a $Sym^{N}(\Omega^q(X\vert\Delta))$.
\end{re}

\begin{example}\label{exproj}

1. \label{K}Si $q=dim(X)$, $S^N_q(X\vert\Delta)=N.K_X+\sum_{h\in H}\{\lfloor(N(1-\frac{1}{m_h})\rfloor.D_h\}:=N.K_X+\lfloor N.\Delta\rfloor:=\lfloor N.(K_X+\Delta)\rfloor$. Le $\Bbb Q$-diviseur $K_X+\Delta$ est le fibr\'e canonique de $(X\vert\Delta)$.

2. Si $\Delta _X=0$, $S^N_q(X\vert\Delta)=Sym^N(\Omega^q_X),\forall N,q$.

3. Si  $\Delta=\lceil \Delta\rceil$, alors: $S^N_q(X\vert\Delta)=Sym^N(\Omega^q_X(log(\Delta))), \forall N,q$. 

On a bien s\^ur: $S^N_q(X\vert\Delta)\subset S^N_q(X\vert\Delta')$ si $\Delta\leq \Delta'$. 

Les faisceaux $S^N_q(X\vert\Delta)$ interpolent donc, en g\'en\'eral, entre $Sym^N(\Omega^q_X)$ et $Sym^N(\Omega^q_X(log(\lceil\Delta\rceil)))$.

4. Soit $(\bP^1\vert D)$ l'orbifolde g\'eom\'etrique logarithmique de dimension $1$, avec $D$ r\'eduit de support  $2$ points distincts de $\bP^1$ (par exemple $0$ et $\infty)$. Alors les faisceaux $S_{N,1}(\bP^1\vert D)$ sont triviaux, de rang $1$. On en d\'eduit que les faisceaux $S^N_q((\bP^1\vert D)^r)$ sont tous triviaux.
\end{example}

\begin{proposition}\label{modiff} Soient $(Y\vert\Delta_Y)$ et $(X\vert\Delta_X)$ des orbifoldes g\'eom\'etriques lisses, et $f:Y\to X$ holomorphe.  

1. Si $f$ induit un morphisme d'orbifoldes g\'eom\'etriques, alors $f^*(S^N_q(X\vert\Delta_X))\subset S^N_q(Y\vert\Delta_Y))$ pour tous $N,q$.

2. Si $f^*(S_{N,1}(X\vert\Delta_X))\subset S_{N,1}(Y\vert\Delta_Y))$ pour $N=ppcm_j(num(m_j))$, avec $\Delta_X=\sum_{j\in J} (1-\frac{1}{m_j}).D_j$, alors $f$ induit un morphisme d'orbifoldes g\'eom\'etriques $f:(Y\vert\Delta_Y)\to (X\vert\Delta_X)$\footnote{Cette propri\'et\'e ne subsiste donc pas en g\'en\'eral pour les morphismes divisibles entre orbifoldes enti\`eres et finies, morphismes consid\'er\'es dans la th\'eorie des champs alg\'ebriques de Deligne-Mumford.}. 

On a not\'e $num(m_j)$ le num\'erateur $u_j$ de $m_j=\frac{u_j}{v_j}$, si $u_j,v_j\in \bZ$ sont premiers entre eux.
\end{proposition}

{\bf D\'emonstration:} Les notations sont celles introduites ci-dessus. 

Pour 1, il suffit de montrer que, pour tous $N,q$,  $f^*(x^{\lceil k/m\rceil}.\otimes_{l=1}^{l=N}\frac{dx_{_{J_l}}}{x_{_{J_l}}})\in S^N_q(Y\vert\Delta_Y)$. Puisque ce dernier faisceau est localement libre, il suffit (Hartogs) de v\'erifier cettte inclusion en codimension un dans $Y$. Soit donc $E\in W(Y)$ et $b\in E$ un point g\'en\'erique. Soit $y=(y_1,...,y_n)$ des coordonn\'ees locales de $Y$ en $b$ telles que $E$ ait pour \'equation locale $y_1=0$. Au voisinage de $b$, on a donc: $f(y)=(y_1^{t_1}.f_1(y), y_1^{t_2}.f_2(y),...,y_1^{t_p}.f_p(y))$, avec $f_j(b)\neq 0$ pour $j\geq 1$. Notre hypoth\`ese est que $t_j.m'\geq m_j$, pour $j=1,...,p$, si $m':=m_{\Delta_Y}(E)$.

Donc: $f^*(\frac{dx_J}{x_J})=\frac{dy_1}{y_1}\wedge u_J$, $\forall J$, avec $u_J,\forall J,$ une $(q-1)$-forme holomorphe, si $\vert J\vert=q$, et si $t_j\geq1$ pour un $j\in J$ au moins. Par suite, \`a des termes holomorphes pr\`es: $f^*(x^{\lceil k/m\rceil}.\otimes_{l=1}^{l=N}\frac{dx_{_{J_l}}}{x_{_{J_l}}})=g(y).y_1^s.(\otimes_{l=1}^{l=k'}(\frac{dy_1}{y_1}\wedge u_{J_l}))\otimes_{l=k'+1}^{l=N} w_l$, avec $g$ et $w_l$ holomorphes, et $s:=\sum_{j=1,...,p}t_j\lceil k_j/m_j\rceil$.

On veut montrer que $s\geq \lceil k'/m'\rceil$, o\`u $k'$ est l'entier ci-dessus (nombre d'occurences de $\frac{dy_1}{y_1}$ dans la forme pr\'ec\'edente). On a \'evidemment: $k'\leq \sum_{j=1}^{j=p} k_j$. Par ailleurs: $t_j.m'\geq m_j, \forall j$, donc: $t_j/m_j\geq 1/m',\forall j$. Donc, pour tout $j$: $t_j.\lceil k_j/m_j\rceil\geq t_j.k_j/m_j\geq k_j/m'$. Donc: $s=\sum t_j.\lceil k_j/m_j\rceil\geq \sum k_j/m'\geq k'/m'$. Puisque $s$ est entier, on a donc bien: $s\geq \lceil k'/m'\rceil$.

Pour 2, soit $E, D,b,t,f$ comme ci-dessus. On a donc, en particulier: $x_1(f(y))=y_1^{t_1}$. Donc: $f^*(x_1^{\lceil N/m_1\rceil}.\frac{dx_1}{x_1}^{\otimes N})=t_1^N.y_1^{t_1\lceil N/m_1\rceil}.(\frac{dy_1}{y_1})^{\otimes N}$. On a, par hypoth\`ese: $t_1.\lceil N/m_1\rceil\geq N/m'$. Or, $num(m_1)$ divise $N=ppcm(num(m_j))$, et on a donc: $t_1.N/m_1=t_1.\lceil N/m_1\rceil\geq N/m'$. C'est pr\'ecis\'ement  (multipli\'ee par $N)$ l'in\'egalit\'e d\'efinissant les morphismes orbifoldes $\square$

De l'exemple \ref{K} on d\'eduit:

\begin{corollary} Si $f:(Y\vert\Delta_Y)\to (X\vert\Delta_X)$ est un morphisme orbifolde surjectif et g\'en\'eriquement fini entre orbifoldes g\'eom\'etriques lisses, on a: $(K_Y+\Delta_Y)\geq f^*(K_X+\Delta_X)$ (signifiant que la diff\'erence est $\Bbb Q$-effective). En particulier, $\kappa(Y\vert\Delta_Y)\geq \kappa(X\vert\Delta_X)$ si $Y$ est compacte et connexe.
\end{corollary}

\begin{re} Dans [Ca04, d\'efinition 2.40, p. 541], la notion de morphisme \'etait d\'efinie seulement \`a l'aide des fibr\'es canoniques. Cette notion est trop faible pour \'etudier la cat\'egorie des orbifoldes g\'eom\'etriques dans un cadre bim\'eromorphe.\end{re}

\subsection{Tenseurs holomorphes orbifoldes.}\label{tdiff}

\

On reprend les hypoth\`eses et notations de la section pr\'ec\'edente d\'efinissant les faisceaux $S^N_q(X\vert\Delta)$ sur $(X\vert\Delta)$ lisse. 
On va d\'efinir plus g\'en\'eralement, de mani\`ere enti\`erement similaire, les faisceaux $T^r_s(X\vert\Delta)$ de tenseurs holomorphes $r$-contravariants et $s$-covariants, muni des op\'erations de contraction et de tensorisation usuelles. Cette d\'efinition est motiv\'ee par une question de M. Pa\u un.

Dans des cordonn\'ees locales adapt\'ees $(x_1,...,x_n)$, c'est le faisceau $T^r_s(X\vert\Delta)$ localement libre engendr\'e comme $\cal O$$_X$-module par les:
$$t^u_v:= x^{\lceil(h-k).a\rceil}.\bigotimes_{j=1}^{j=s} dx_{v(j)} \bigotimes_{i=1}^{i=r} \partial/\partial x_{u(i)},$$

o\`u $u:[1,r]\to [1,p]$ et $v:[1,s]\to [1,p]$ sont des applications quelconques, le $p$-uplet $\lceil(h-k).a\rceil$ (\`a valeurs dans $\bZ$) \'etant d\'efini comme dans la section pr\'ec\'edente.

De la d\'efinition on d\'eduit les deux propri\'et\'e usuelles suivantes:

$\square$ On a une application de tensorisation naturelle: $$\otimes: T^r_s(X\vert\Delta)\bigotimes T^p_q(X\vert\Delta)\to T^{r+p}_{s+q}(X\vert\Delta),\forall p,q,r,s.$$

$\square$ On a une application de contraction naturelle: $$c: T^{r+p}_{s+p}(X\vert\Delta)\to T^{r}_{s}(X\vert\Delta),\forall p,r,s.$$

\begin{re} 

\

1. Les faisceaux $T^r(X\vert \Delta)$ des germes de $r$-champs de vecteurs tangents \`a $\Delta$ interviennent dans l'\'etude des d\'eformations de courbes $\Delta$-rationnelles au \S5.

2. On pourrait d\'efinir de m\^eme des espaces de jets sur $(X\vert\Delta)$ lisse, jets qui interviennent dans l'\'etude de l'hyperbolicit\'e.
\end{re}

\subsection{Dimension canonique d'un faisceau diff\'erentiel de rang $1$}\label{kL}

{\bf Notations} Soit $(X\vert\Delta)$ une orbifolde g\'eom\'etrique lisse, avec $X$ {\bf compacte} et connexe, et $N,q \geq 0$, entiers.

Soit $L_N\subset S^N_q(X\vert\Delta)$ un sous faisceau analytique coh\'erent de rang $1$. On notera $\overline{L_N}^{\Delta}:=\footnote{Si aucune confusion n'est possible sur $\Delta$}$ $\overline{L_N} \subset S^N_q(X\vert\Delta)$ sa saturation dans $S^N_q(X\vert\Delta)$, saturation qui est par d\'efinition le plus grand sous-faisceau analytique coh\'erent de rang un de $S^N_q(X\vert\Delta)$ contenant $L_N$. Ce faisceau est donc sans cotorsion, \'etant aussi d\'efini comme le noyau du morphisme compos\'e de l'injection de $L_N$ dans $S^N_q(X\vert \Delta)$ avec le quotient: $(S^N_q(X\vert \Delta)/L_N)/Torsion$.

Soit maintenant $L\subset \Omega^q_X$ un sous-faisceau analytique coh\'erent de rang $1$, et pour tout $N\geq 0$, soit $L_N\subset S^N_q(X\vert\Delta)$ l'image de $L^{\otimes N}$ dans $Sym^N(\Omega^q_X)\subset S^N_q(X\vert\Delta)$, et soit enfin $\overline{L_N}\subset S^N_q(X\vert\Delta)$ sa saturation dans $S^N_q(X\vert\Delta)$.

On note: $H^0(X\vert\Delta,L_N):=H^0(X, \overline{L_N}^{\Delta})$: c'est aussi le sous-espace de $H^0(X,S^N_q(X\vert\Delta))$ constitu\'e des sections dont l'image est contenue dans $L_N$ au point g\'en\'erique de $X$. On note enfin $p_N(X\vert\Delta,L):=h^0(X, \overline{L_N})$ sa dimension complexe.

On a des applications naturelles $\overline{L_N}\otimes\overline{L_{N'}}\to \overline{L_{N+N'}}$, qui induisent au niveau des sections une structure d'anneau gradu\'e sur $R(X\vert\Delta, L):=\oplus_{N\geq 0}H^0(X,\overline{L_N})$.

 On notera $\kappa(X\vert\Delta, L)$ le degr\'e de transcendance sur $\Bbb C$ de cet anneau, diminu\'e d'une unit\'e si ce degr\'e est au moins $1$; si ce degr\'e de transcendance est $0$, on notera $\kappa(X\vert\Delta,L)=-\infty$.

\begin{lemma}\label{kodL} Si $(X\vert\Delta)$ est lisse, on a, notant $\Phi_N(X)$ l'application m\'eromorphe induite par le syst\`eme lin\'eaire complet $\overline{L_N}$:
$$\kappa((X\vert\Delta),L)=\overline{lim}_{N\to +\infty}(\frac{log (p_N(X\vert\Delta,L)}{log N})=max_{N>0}(dim(\Phi_N(X)))$$\end{lemma}

{\bf D\'emonstration:} Ces assertions se d\'emontrent comme dans le cas d'un fibr\'e en droites (Voir, par exemple, [U75, theorem 5.10, p. 58]) $\square$

\

\begin{definition}On appelle $\kappa(X\vert\Delta, L)\in \lbrace -\infty,0,1,...,p:=dim(X)\rbrace$ {\bf la $\Delta$-dimension de $L$}.

La $\bC$-alg\`ebre $K(X\vert\Delta,L):=\oplus_{N\geq 0}H^0(X,(\overline{L_{N}})$ est appel\'ee {\bf la $\Delta$-alg\`ebre} de $L$. 

Lorsque $L=K_X$, on note simplement $K(X\vert\Delta)$ cette alg\`ebre, appel\'ee {\bf l'alg\`ebre canonique} de $(X\vert\Delta)$. 
\end{definition}

\begin{example} 

1. Si $q=dim(X)$, $\kappa(X\vert\Delta, K_X)=\kappa(X,K_X+\Delta)=\kappa(X\vert\Delta)$.

2. Soit $(X\vert\Delta):=(\bP^1\vert D)^r, r\geq 1$ l'exemple 4 de \ref{exproj}. Pour tout $q\geq 1$ et tout sous-faisceau $L$ de rang $1$ de $\Omega_X^q$, on a donc: $\kappa(X\vert\Delta, L)\leq 0$.

\end{example}

\subsection{Invariance bim\'eromorphe de la dimension canonique}

\begin{proposition}\label{mobimdiff} Soit $(Y\vert\Delta')$ et $(X\vert\Delta)$ des orbifoldes g\'eom\'etriques {\bf lisses}, $X,Y$ compactes et connexes. Soit $f:Y\to X$ une application bim\'eromorphe. On notera $N_0(\Delta)$ le plus petit commun multiple des $num(m_{\Delta}(D))$ pour les $D\in W(X)$ tel que $m_{\Delta}(D)$ soit fini.

On a \'equivalence entre les deux conditions suivantes:

1. $f:(Y\vert\Delta')\to (X\vert\Delta)$ est un morphisme orbifolde bim\'eromorphe \'el\'ementaire.\footnote{ Voir remarque \ref{defmobim}, (5) ci-dessus pour la d\'efinition.}

2.Il existe $N$, divisible par $N_0(\Delta)$ et $N_0(\Delta')$, tel que $f^*(S_{N,1}(X\vert\Delta))\subset S_{N,1}(Y\vert\Delta'))$, et $f_*(S_{N,1}(Y\vert\Delta'))\subset S_{N,1}(X\vert\Delta))$.
\end{proposition}

{\bf D\'emonstration:} Si la condition 1. est satisfaite, la premi\`ere (resp. seconde) des propri\'et\'es 2. est satisfaite par \ref{modiff} (resp. par le lemme d'Hartogs, et le fait que $S_{N,1}(X\vert\Delta))$ est localement libre).

Si la condition 2. est satisfaite, alors $f$ est un morphisme orbifolde par la premi\`ere condition 2. et \ref{modiff}. De plus, $f_*(\Delta')=\Delta$ par la seconde condition 2. En effet, cette condition garantit qu'en codimension $1$, les multiplicit\'es de $\Delta$ coincident avec celles de $\Delta'$ sur les transform\'es stricts des composantes du support de $\Delta$ $\square$

\begin{theorem}\label{invbir'} Soit $X, X'$ lisses compactes et connexes, et $u:(X'\vert\Delta')\to (X\vert\Delta)$ un morphisme orbifolde bim\'eromorphe \'el\'ementaire entre orbifoldes g\'eom\'etriques lisses. Soit $L\subset \Omega_X^q$ un sous-faisceau coh\'erent de rang $1$, et  $L':=u^*(L)\subset \Omega_{X'}^q$. On note, pour tout $N\geq 0$ $u^*:\overline{L_N}\to \overline{L'_N}$ le morphisme naturel de faisceaux. Alors:

1. $u^*:H^0(X\vert\Delta, L_N)\to H^0(X'\vert\Delta', L'_N)$ est un isomorphisme. 

2. $p_N(X\vert\Delta,L)=p_N(X'\vert\Delta',L')$.

2. $\kappa(X\vert\Delta, L)=\kappa(X'\vert\Delta', L')$.

\end{theorem}

{\bf D\'emonstration:} Soit $A\subset X$ le lieu (de codimension $2$ au moins) au-dessus duquel $u$ n'est pas un isomorphisme. Soit $r: H^0(X,\overline {L_N})\to H^0(X-A,\overline {L_N})$ la restriction. Par Hartogs, c'est un isomorphisme (d'espaces vectoriels complexes). Soit $r': H^0(X',\overline {L'_N})\to H^0(X-A,\overline {L_N})$ la restriction (compos\'ee avec $u)$. Elle est injective. Puisque $r'\circ u^*=r, u^*$ est un isomorphisme $\square$

\begin{re} En g\'en\'eral, si $L$ est localement libre, $\kappa(X,L)\leq \kappa(X\vert \Delta,L)$ et l'in\'egalit\'e peut-\^etre stricte, le membre de droite d\'ependant de $\Delta$.
\end{re}

\begin{example}\label{exproj'} Le cas des orbifoldes g\'eom\'etriques logarithmiques est simple: si $(X'\vert D')$  et $(X \vert D)$ sont des orbifoldes g\'eom\'etriques logarithmiques lisses (ie: si $X,X'$ le sont, et $D,D'$ r\'eduites \`a croisements normaux), ces deux orbifoldes g\'eom\'etriques logarithmiques sont bim\'eromorphes si et seulement s'il existe une suite $w_j=u_j'\circ (u_j^{-1}):X_{j-1}\dasharrow X_j$ d'\'equivalences bim\'eromorphes, pour $j=1,\dots,m$, telle que $(u_j^{-1}(D_j))=(u_j')^{-1}(D_{j-1})$, avec $(X_0\vert D_0)=(X\vert D)$, et $(X_m\vert D_m)=(X'\vert D')$ .

La condition``croisements normaux" n'est pas superflue: $X=\bP^2$, $D$ la r\'eunion de $3$ droites concourantes, et $X'=\bP^1\times \bP^1$, $D'$ la r\'eunion de $3$ fibres de la projection sur le premier facteur, le montrent: la dimension canonique n'est pas m\^eme pr\'eserv\'ee.

  Par exemple, les orbifoldes g\'eom\'etriques $(X\vert\Delta):=(\bP^1\vert D_1)^r, r\geq 1$ de l'exemple 4 de \ref{exproj}, et $(X'\vert\Delta'):=(\bP^r\vert D_r)$, o\`u $D_r$ est la r\'eunion r\'eduite et \`a croisements normaux de $(r+1)$ hyperplans projectifs (distincts) sont bim\'eromorphiquement \'equivalentes. On en d\'eduit que les $S^N_q(X'\vert\Delta')$ sont aussi tous triviaux, et $\kappa(X'\vert\Delta',L)\leq 0$ pour tout sous-faisceau $L$ de rang $1$ de $\Omega_{\bP^r}^q$. 
  \end{example}

 \begin{re} Si $U=X-D$ est une vari\'et\'e quasi-projective lisse (connexe), avec $(X\vert D)$ lisse, les propri\'et\'es de $(X\vert D)$ qui sont des invariants bim\'eromorphes sont donc des invariants de $U$. Tels sont donc $\kappa(X\vert D),H^0(X,S^N(\Omega^p(X\vert D))$, et la classique irr\'egularit\'e logarithmique $q(U):=H^0(X,\Omega^1(X\vert D))$.
\end{re}

\subsection{Fibration de Moishezon-Iitaka.}\label{moiit}

Soit $(X\vert\Delta)$ une orbifolde g\'eom\'etrique lisse telle que $\kappa(X\vert\Delta)\geq 0$. On suppose $X$ compacte et connexe. Pour tout entier $m>0$ tel que $m\Delta$ soit entier et $H^0(X,m(K_X+\Delta))$ soit non nul, on note $\Phi_m:X\dasharrow \bP_{N_m}$ l'application m\'eromorphe associ\'ee \`a ce syst\`eme lin\'eaire.

Les arguments classiques (voir [U75], par exemple), montrent que pour $m$ assez grand et divisible, cette application est \`a fibres connexes, ind\'ependantes \`a bim\'eromorphie pr\`es de $m$, et que ces fibres orbifoldes stables g\'en\'eriques ont $\kappa=0$, et que la base de cette fibration est de dimension $\kappa(X\vert\Delta)$.

On la notera $M_{(X\vert\Delta)}:(X\vert\Delta)\dasharrow M(X\vert\Delta)$, et on l'appellera la {\bf fibration de Moishezon-Iitaka}\footnote{Appel\'ee habituellement ``fibration d'Iitaka", bien qu'introduite ind\'ependamment par B. Moishezon au CIM de Nice (1970).} de $(X\vert\Delta)$. 

On d\'eduit de la d\'emonstration de [U75, theorem 5.10, p. 58] le:

\begin{theorem}\label{klo} Si $(X\vert \Delta)$ est lisse, $X$ compacte et connexe, avec $\kappa(X\vert \Delta)\geq 0$, et si $g:(Y\vert \Delta_Y)\to(X\vert \Delta)$ est un morphisme orbifolde bim\'eromorphe \'el\'ementaire tel que $M_Y:=M_{(X\vert \Delta)}\circ g: Y\to M(X\vert \Delta)$ soit holomorphe, alors la fibre g\'en\'erale\footnote{g\'en\'erique si $X$ est projective.} orbifolde $(Y_m\vert \Delta_m)$\footnote{Voir d\'efinitions \ref{fiborb} et \ref{dfo}.} de $M_Y$ est de dimension canonique nulle (i.e: $\kappa(Y_m\vert \Delta_m)=0)$.
\end{theorem}

La fibration $M$ peut aussi \^etre caract\'eris\'ee, avec les m\^emes arguments que dans le cas non-orbifolde, comme suit: toute fibration $g:(X\vert\Delta)\dasharrow Y$ dont les fibres orbifolde $F$ g\'en\'erales (voir d\'efinition \ref{dfo}) satisfont $\kappa(F)=0$ domine $M_{(X\vert\Delta_X)}$.

\subsection{Invariance \'etale}

\begin{definition}\label{defet} Soit $v:(Y'\vert\Delta')\to (Y\vert\Delta)$ un morphisme orbifolde entre orbifoldes g\'eom\'etriques lisses et enti\`eres de m\^eme dimension $q$. 

On dit que $v$ est {\bf \'etale en codimension $1$} s'il existe un sous-ensemble analytique ferm\'e $A$ de codimension au moins $2$ de $Y$  tel que, sur le compl\'ementaire de $A$ on ait l'\'egalit\'e: $v^*(N.(K_Y+\Delta))=N.(K_{Y'}+\Delta')),$ pour un entier $N$ tel que $N. \Delta$ soit entier.
\end{definition}

\begin{proposition}\label{etale}Soit $v:(Y'\vert\Delta')\to (Y\vert\Delta)$ un morphisme orbifolde divisible entre orbifoldes g\'eom\'etriques lisses et enti\`eres. 

1. On a \'equivalence entre les $3$ propri\'et\'es (a),(b),(c) suivantes:

(a) $v$ est \'etale en codimension $1$.

(b) Il existe un ferm\'e analytique $A\subset Y$ de codimension au moins $2$ tel que, au-dessus du compl\'ementaire de $A$ dans $Y$, $\Delta'$ {\bf divise exactement} $v^*(\Delta)$ (ie: pour tout $D'\in W(Y)$ rencontrant $(Y'-v^{-1}(A)))$, $r(D'):=\frac{m_{\Delta}(v(D'))}{m_{\Delta'}(D')}$ est un entier \'egal \`a l'ordre de ramification de $v$ au point g\'en\'erique de $D')$.

(c) Il existe un rev\^etement orbifolde \'etale (au sens du \S\ref{revun}) $u:(\bar Y\vert \bar\Delta)\to (Y\vert \Delta)$ fini (avec $\bar Y$ fini, normal et \`a singularit\'es quotient), et un morphisme orbifolde divisible bim\'eromorphe $\bar v: (Y'\vert \Delta')\to (\bar Y\vert \bar\Delta)$ tels que $v=u\circ \bar v$.

2. De plus, si $v$ est \'etale en codimension $1$, alors pour tout sous faisceau $L\subset \Omega^q_Y$ coh\'erent de rang $1$, $\kappa(Y\vert\Delta,L)=\kappa(Y'\vert\Delta',v^*(L))$. 

En particulier: $\kappa(Y\vert\Delta)=\kappa(Y'\vert\Delta')$.
\end{proposition}

{\bf D\'emonstration:} Assertion 1: localement en $b'\in Y'-A'$, $A':=v^{-1}(A)$, quitte \`a augmenter un peu $A$ de sorte que $\Delta\cap (Y-A)$ soit lisse, on a des coordonn\'ees $y'=(y'_1,\dots,y'_n)$, ainsi que des coordonn\'ees locales $y=(y_1,\dots,y_n)$ en $b:=v(b')$ telles que $v(y')=({y'_1}^r,y'_2,\dots, y'_n)$, et telles que les \'equations de $\Delta$ en $b$ (resp. $\Delta'$ en $b')$ soient: ${y_1}^{1-\frac{1}{m}}$ (resp. ${y'_1}^{1-\frac{1}{m}'})$. Un calcul imm\'ediat montre alors que les prori\'et\'es (a) et (c) sont simultan\'ement v\'erifi\'ees, ceci si et seulement si: $r.m'=m$. D'o\`u l'\'equivalence de (a) et (c). La propri\'et\'e (c) implique (a), en prenant $A:=S_{\Delta}$ le lieu singulier du support de $\Delta$. L'implication r\'eciproque s'obtient en prenant pour couple $u,\bar v$ la factorisation de Stein du morphisme $v$, et $\bar \Delta:=\bar v_*(\Delta')$, et en utilisant le fait que les rev\^etements orbifoldes \'etales sont d\'etermin\'es par leurs restrictions au-dessus de $X-S_{\Delta}$, qui est localement simplement connexe pr\`es de $S_{\Delta}$. Voir \S\ref{revun}.

Pour d\'emontrer l'assertion 2, on peut supposer que $u$ est un rev\^etement Galoisien de groupe (fini) $G$. Si $N$ est grand et suffisamment divisible, alors (par [U75,lemma 5.12]): 

$H^0(Y,\overline{L_N}^{\Delta})=H^0(\bar Y, u^*(\overline{L_N}^{\Delta}))^G=H^0(\bar Y\vert \bar\Delta, (\overline{u^*L})_N)^G$, puisque les faisceaux consid\'er\'es sont des sous-faisceaux des faisceaux localement libres $S^N_q(Y\vert\Delta)$ et $u^*(S^N_q(Y\vert\Delta))$  auquel on peut donc appliquer (sur $Y$, aux fonctions sym\'etriques des sections sur $\bar Y$) le prolongement de Hartogs, qui montre que $H^0(\bar Y, F)=H^0(\bar Y-\bar A, \bar F)$, si $\bar F=u^*(F)$ est un sous-faisceau satur\'e d'un sous-faisceau localement libre sur $\bar Y$ provenant aussi de $Y$, avec: $\bar A:=u^{-1}(A)$.

On a enfin, par l'\'egalit\'e pr\'ec\'edente:  $H^0(Y'\vert \Delta',v^*(L_N))^G\subset H^0(\bar Y\vert\bar \Delta,u^*(L_N)^G$, ce qui suffit \`a montrer, par [U75, theorem 5.13], que $\kappa(Y'\vert \Delta', v^*(L))\leq \kappa(Y\vert \Delta,L)$, l'in\'egalit\'e oppos\'ee \'etant \'evidente $\square$

\subsection{In\'egalit\'e de Bogomolov (version orbifolde)}

\

Le th\'eor\`eme suivant est d\^ u \`a F. Bogomolov ([Bo 78]) lorsque $D$ est vide. Sa d\'emonstration s'\'etend imm\'ediatement au cas logarithmique, en utilisant le fait, d\^ u \`a P. Deligne ([De 74]), que les sections de $\Omega_X^p(log D))$ sont $d$-ferm\'ees.\footnote{Cette version logarithmique de [Bo 78] a \'et\'e observ\'ee ind\'ependamment par S. Lu, \`a la suite de [Ca 01].}. Nous en donnons ci-dessous une version orbifolde g\'en\'erale.

\begin{theorem}\label{boglog} Soit $X$ une vari\'et\'e K\" ahl\'erienne compacte et connexe\footnote{Ou, plus g\'en\'eralement, bim\'eromorphe \`a une telle vari\'et\'e.}, $D$ un diviseur \`a croisements normaux sur $X$, et $L\subset \Omega_X^p(log D))$ un sous-faisceau coh\'erent de rang $1$. Alors:

(1) $\kappa(X,L)\leq p$.

(2) Si $\kappa(X,L)=p$, alors il existe une (unique) fibration m\'eromorphe $f:X\to Y$ telle que les saturations dans $\Omega_X^p(log D))$ de $L$ et de $f^*(K_Y)$ coincident. (En particulier, $dim(Y)=p$, et $f$ est la fibration d'Iitaka d\'efinie par le syst\`eme lin\'eaire des sections de $\overline{L_m}$, pour $m>0$ assez grand, la saturation de $L^{\otimes m}$ \'etant prise dans $S_{m,p}(X\vert D)=Sym^m(\Omega_X^p(log D)))$.
\end{theorem}

{\bf D\'emonstration:} Elle est identique \`a celle de [Bo 78] (dans le cas projectif), ou de [Ca 04] (dans le cas K\" ahler), en utilisant [De 74] au lieu de la classique fermeture des formes diff\'erentielles holomorphes de la th\'eorie de Hodge $\square$

\begin{corollary}\label{bogorb} Soit $X$ une vari\'et\'e K\" ahl\'erienne compacte et connexe, et  soit $\Delta$ un diviseur orbifolde sur $X$. Soit $L\subset \Omega_X^p$ un sous-faisceau coh\'erent de rang $1$. Alors:

(1) $\kappa(X\vert\Delta,L)\leq p$.

(2) Si $\kappa(X\vert\Delta,L)=p$, alors il existe une (unique, \`a \'equivalence bim\'eromorphe pr\`es) fibration m\'eromorphe $f:X\to Y$ telle que $L$ et $f^*(K_Y)$ coincident au point g\'en\'erique de $X$. (En particulier, $dim(Y)=p$, et $f$ est la fibration d'Iitaka d\'efinie par le syst\`eme lin\'eaire des sections de $\bar L_m$, pour $m>0$ assez grand, la saturation de $L^{\otimes m}$ \'etant prise dans $S_{m,p}(X\vert\Delta))$.
\end{corollary}

{\bf D\'emonstration:} Soit $D$ le support de $\Delta$. On a, bien s\^ur: $\kappa(X\vert\Delta,L)\leq \kappa(X/D,L)\leq p$. D'o\`u (1), par \ref{boglog} (1) ci-dessus, appliqu\'e \`a $L$ dans $\Omega_X^p(log D))$. Si on a \'egalit\'e, on a aussi: $p=\kappa(X\vert\Delta,L)\leq \kappa(X\vert D,L)\leq p$. D'o\`u (2), en appliquant de la m\^eme fa\c con \ref{boglog} (2) $\square$

\begin{definition}\label{fbog} Soit $(X\vert \Delta)$ une orbifolde g\'eom\'etrique lisse, avec $X\in \cal C$. Un sous-faisceau coh\'erent $L\subset \Omega_X^p$ de rang $1$, pour un $p>0$, est dit {\bf $\Delta$-Bogomolov} si $\kappa(X\vert \Delta,L)=p$. Cette propri\'et\'e ne d\'epend que de la saturation de $L$ dans $\Omega_X^p$, et de la classe d'\'equivalence bim\'eromorphe de $(X\vert \Delta)$.

On note $Bog(X\vert \Delta)$ l'ensemble (eventuellement vide) des faisceaux $\Delta$-Bogomolov satur\'es de $X$.
\end{definition}

\begin{question} L'ensemble $Bog(X\vert \Delta)$ est-il toujours fini? Cette propri\'et\'e g\'en\'eraliserait la finitude, due \`a Kobayashi-Ochiai, de l'ensemble des classes d'\'equivalence bim\'eromorphes d'applications m\'eromorphes connexes dominantes $X\dasharrow Y$ ($X$ fix\'ee, $Y$ variable) avec $Y$ de type g\'en\'eral.
\end{question}

\subsection{\'Equivalence bim\'eromorphe.}\label{equbim}

Rappelons (d\'efinition \ref{morphorb} (5)) que $f:(Y\vert\Delta_Y)\to (X\vert\Delta_X)$ est un {\bf morphisme orbifolde bim\'eromorphe \'el\'ementaire} si c'est un morphisme orbifolde, si $f:Y\to X$ est birationnel, et si $f_*(\Delta_Y)=\Delta_X$. On suppose ici que $Y,X$ sont $\bQ$-factorielles. Noter que les trois propri\'et\'es sont ind\'ependantes (deux d'entre elles n'impliquent pas la troisi\`eme).

\begin{definition} \label{defbim} On dira que deux orbifoldes g\'eom\'etriques {log-canoniques} sont {\bf bim\'eromorphiquement \'equivalentes} s'il existe une chaine de morphismes orbifoldes bim\'eromorphes \'el\'ementaires les reliant, et tels que chacunes des orbifoldes g\'eom\'etriques de cette chaine soit log-canonique. De mani\`ere \'equivalente: c'est la relation d\'equivalence engendr\'ee par les morphismes orbifolde bim\'eromorphes \'el\'ementaires entre orbifoldes g\'eom\'etriques log-canoniques.
\end{definition}

\begin{re} Il r\'esulte de la remarque \ref{icanlc} que la dimension canonique d'une orbifolde g\'eom\'etrique log-canonique est un invariant birationnel, au sens d\'efini ci-dessus.
\end{re}

L'\'equivalence bim\'eromorphe dans cette cat\'egorie pr\'esente (au moins) une diff\'erence notable avec le cas non orbifolde g\'eom\'etrique: l'exemple \ref{crem} ci-dessous montre que deux orbifoldes g\'eom\'etriques lisses (projectives) bim\'eromorphiquement \'equivalentes ne sont pas n\'ecessairement domin\'ees par une troisi\`eme qui leur est bim\'eromorphiquement \'equivalente. Cette diff\'erence complique consid\'erablement l'\'etude de l'\'equivalence bim\'eromorphe orbifolde.

\begin{example}\label{crem} {\bf (transformation de Cremona)} Soit $X:=\bP^2$, et $3$ points $a,b,c\in X$ non align\'es. On note $A,B,C$ les trois droites projectives de $X$ passant par $2$ de ces $3$ points. Soit $u:X'\to X$ l'\'eclat\'e de $X$ en les $3$ points $a,b,c$, et $A',B',C'$ les transform\'ees strictes des droites $A,B,C$ dans $X'$. 

On note enfin $E,F,G$ les trois diviseurs irr\'eductibles de $u$. Alors $X'$ admet une seconde contraction $v:X'\to Y$ sur $Y\cong \bP^2$, qui contracte $A',B',C'$ sur trois points $a',b',c'\in Y$, et $E,F,G$ sur trois droites projectives $E',F',G'$ de $Y\cong \bP^2$. ($v\circ u^{-1}:X\dasharrow Y$ n'est donc autre que la transformation de Cremona). 

On peut munir $X'$ des deux structures orbifoldes g\'eom\'etriques (logarithmiques): 

$\Delta:=(E+F+G)$ et $\Delta':=(A'+B'+C')$. 

Alors: $u:(X'\vert\Delta)\to X$ et $v:(X'\vert\Delta')\to Y$ sont des morphismes orbifoldes \'el\'ementaires, puisque le support de $\Delta$ (resp. $\Delta')$ est $u$-exceptionnel (resp. $v$-exceptionnel). 

De plus, $u:X'=(X'\vert 0)\to X$ et $v:X'\to Y$ sont des morphismes orbifoldes \'el\'ementaires, par d\'efinition m\^eme. 

Donc, l'identit\'e ensembliste $1_{X'}:(X'\vert\Delta)\to X'=(X'\vert 0)$ et $1_{X'}:(X'\vert\Delta')\to X'$ sont des morphismes orbifoldes birationnels (non \'el\'ementaires). 

Par suite, $(X'\vert\Delta)$ et $(X'\vert\Delta')$ sont birationnellement \'equivalentes.

Soit maintenant $\Delta":=sup(\Delta,\Delta')=(A'+B'+C'+E+F+G)$: c'est un diviseur orbifolde (logarithmique) sur $X'$. De plus, $\kappa(X'\vert\Delta")=\kappa(X/D)=0$, si $D:=(A+B+C)$. Donc $(X'\vert\Delta")$ n'est birationnellement \'equivalent ni \`a $(X'\vert\Delta)$ ni \`a $(X'\vert\Delta')$, puisque $\kappa(X'\vert\Delta)=\kappa(X\vert 0)=-\infty$. 

Si l'orbifolde g\'eom\'etrique lisse $(Z\vert\Delta_Z)$ domine (par des morphismes orbifoldes) $g:(Z\vert\Delta_Z)\to (X'\vert\Delta)$ et $g':(Z\vert\Delta_Z)\to (X'\vert\Delta')$, elle domine aussi $(X'\vert\Delta")$. En effet, on a, ensemblistement $g=g'$. Si $E'\in W(X')$, si $H\in W(Z)$, et si $g^*(E')=t.H+\dots$, avec $t>0$, on a: $t.m_Z(F)\geq sup\{m_{\Delta}(E'),m_{\Delta'}(E')\}$, par hypoth\`ese. Mais c'est justement la conclusion cherch\'ee.

Donc $\kappa(Z\vert\Delta_Z)\geq \kappa(X'\vert\Delta")=0$, et $(Z\vert\Delta_Z)$ n'est pas birationnelle \`a $(X'\vert\Delta)$. 
\end{example}

\begin{example}\label{crem'}{\bf Transformation \'el\'ementaire.} Le ph\'enom\`ene pr\'ec\'edent existe aussi pour les orbifoldes g\'eom\'etriques \`a multiplicit\'es (enti\`eres) finies.

Soit $B$ une courbe elliptique, et $b\in B$. Soit $p:X\to B$ une surface g\'eom\'etriquement r\'egl\'ee. On note $\Bbb P^1\cong F:=p^{-1}(b)$ la fibre de $p$ au-dessus de $b$. Soit $a\in F$ un point arbitraire.

Soit $t:X"\to X$ l'\'eclatement de $X$ en $a$, $F"$ la transform\'ee stricte de $F$ dans $X"$, et $E"$ le diviseur exceptionnel de $t$ au-dessus de $a$. Soit $t':X"\to X'$ la contraction de la $(-1)$-courbes $F"$ dans $X"$, et $E'\subset X'$ l'image de $E"$ dans $X'$ par $t'$. Soit $b':X'\to B$ la fibration telle que $p'\circ t'=p\circ t$. Et $E'=(p')^{-1}(b)$.

Soit $m>1$ un entier. On munit $X"$ de deux diviseurs orbifolde $\Delta_i$, avec $\Delta_1:=(1-\frac{1}{m}).[F"]$, et $\Delta_2:=(1-\frac{1}{m}).[E"].$ On a donc: 

$\Delta":=sup\{\Delta_1,\Delta_2\}=(1-\frac{1}{m}).[E"+F"]$. 

On v\'erifie imm\'ediatemment que $\Delta_1$ et $\Delta_2$ sont bim\'eromorphiquement \'equivalentes (\`a $(X'\vert 0)$ et $(X\vert 0)$ respectivement, et donc aussi \`a $(X"\vert 0))$. 

On en d\'eduit que $h^0(X"\vert\Delta_i, S^N(\Omega^1(X"\vert\Delta_i)))=1,\forall N>0$. 

Par ailleurs, $t:(X"\vert\Delta")\to (X\vert\Delta^+)$ est un morphisme orbifolde \'el\'ementaire si $\Delta^+:=(1-\frac{1}{m}).[F_1+F_2]$. 
D'o\`u l'on d\'eduit que: $h^0(X"\vert\Delta", S^N(\Omega^1(X"\vert\Delta")))= h^0(B,S^N(\Omega^1(B\vert\Delta_B)))\cong N.(1-\frac{1}{m}),$ puisque $p:(X"\vert\Delta")\to (B\vert\Delta_B)$ est un morphisme orbifolde, si $\Delta_B:=(1-\frac{1}{m}).\{b\}$.

En particulier, $(X"\vert\Delta")$ n'est pas bim\'eromorphe \`a $(X"\vert\Delta_i)$, puisque l'\'equivalence bim\'eromorphe pr\'eserve les sections de $S^N(\Omega^q)$. 

Remarquons que l'on peut faire de tels exemples au-dessus d'une courbe $B$ de genre $g\geq 1$ quelconque, et au-dessus de $B=\Bbb P^1$ si on fait des \'eclatements au-dessus d'un nombre suffisant de points $b_i$.
\end{example}

\begin{re} Dans les deux exemples pr\'ec\'edents, nous avons deux orbifoldes lisses $(X\vert \Delta_j),j=1,2$ birationnellement \'equivalentes, mais non birationnellement \'equivalentes \`a $(X\vert\Delta^+)$, si $\Delta^+:=sup\{\Delta_1,\Delta_2\}$. Observer cependant qu'elles sont birationnellement \'equivalentes \`a $(X\vert\Delta^-)$, si $\Delta^-:=inf\{\Delta_1,\Delta_2\}$. Cette observation pourrait peut-\^etre permettre de d\'ecrire l'\'equivalence birationnelle orbifolde \`a l'aide d'une unique relation binaire (\'egale \`a la relation d'\'equivalence qu'elle engendre, en rempla\c cant $X$ par une modification $g:X'\to X$, et $\Delta$ par $\Delta'$ sur $X'$, minimale ad\'equate telle que $(X'\vert \Delta')$ soit birationnelle \`a $(X_i\vert \Delta_i), i=1,2)$. Le principe \'etant de dominer les $X_i$, mais de diminuer les $\Delta_i$. \end{re}

\subsection{Restriction \`a une sous-vari\'et\'e.} \label{catorb}

\begin{definition}\label{fdefrestr} Soit $(Y\vert\Delta)$ une orbifolde g\'eom\'etrique lisse, et $j:V\to Y$  l'inclusion d'une sous-vari\'et\'e \footnote{ie: un sous-ensemble analytique ferm\'e irr\'eductible.}non contenue dans le support de $\Delta$. Soit $j':V'\to V$ une r\'esolution de $V$, telle que $(j')^{-1}(supp(\Delta))$ soit un diviseur de $V'$ \`a croisements normaux. Soit $\Delta'$ un diviseur orbifolde sur $V'$, de support contenu dans $(j')^{-1}(supp(\Delta))$.

L'orbifolde g\'eom\'etrique $(V'\vert\Delta')$ est {\bf une restriction de $\Delta$ \`a $V$} si $j': (V'\vert\Delta')\to (Y\vert\Delta)$ est un morphisme orbifolde, et si $\Delta'$ est le plus petit diviseur orbifolde sur $V'$ tel que $j'$ soit un morphise orbifolde\footnote{On peut formuler cette notion de restriction de fa\c con analogue pour les morphismes divisibles.}.

\end{definition}

\begin{example}\label{rescoutr}Un cas particulier est celui dans lequel $V$ est une courbe lisse, et rencontre transversalement en des points lisses le diviseur $Supp(\Delta)$. Dans ce cas, il existe une restriction $\Delta_V$ de $\Delta$ \`a $V$, obtenue en affectant chacun des points d'intersection de $V$ avec $Supp(\Delta)$ de la $\Delta$-multiplicit\'e de l'unique composante irr\'eductible de $Supp(\Delta)$ \`a laquelle appartient ce point. On a donc simplement: $\Delta_V:=j^*(\Delta)$ dans ce cas. On dit alors que $\Delta$ est {\bf transverse \`a $V$}, et que l'orbifolde $(V\vert \Delta_V)$ est {\bf la restriction transverse de $\Delta$ \`a $V$}. 
\end{example}

\begin{proposition}\label{rescoumin} Soit $(Y\vert\Delta)$ une orbifolde g\'eom\'etrique lisse, et $j:V\to Y$  l'inclusion d'une sous-vari\'et\'e non contenue dans le support de $\Delta$. Soit $(V'\vert\Delta_{V'})$ une restriction de $\Delta$ \`a $V$.

1. Si $u:V"\to V'$ est une modification propre telle que $(j'\circ u)^{-1}(supp(\Delta))$ est \`a croisements normaux, et si $(V"\vert\Delta_{V"})$ est la restriction de $\Delta$ \`a $V$ de support $V"$, alors $u_*(\Delta_{V"})=\Delta_{V'}$ (mais $u$ n'est pas un morphisme orbifolde en g\'en\'eral). 

2. Deux restrictions de $\Delta$ \`a $V$ sont \'egales si $V$ est une courbe. On l'appellera donc la restriction de $\Delta$ \`a $V$.
 \end{proposition}

{\bf D\'emonstration:} Assertion 1. Montrons que  $u_*(\Delta_{V"})=\Delta_{V'}$. Soit donc $F'\in W(V')$ et $E\in W(Y)$ choisi tel que $(j')^*(E)=t.F'+....$, avec $t>0$ et $t.m_{\Delta_{V'}}(F')=m_{\Delta}(E)$, l'\'egalit\'e r\'esultant de la minimalit\'e de $\Delta_{V'}$. Si $F"\in W(V")$ est le transform\'e strict de $F'$ dans $V"$, on a donc: $(j'\circ u)^*(E)=t.F"+....$, et par hypoth\`ese, $t.m_{\Delta_{V"}}(F")\geq m_{\Delta}(E)$. Donc $m_{\Delta_{V"}}(F")\geq m_{\Delta_{V'}}(F')$, et puisque $\Delta_{V"}$ est minimale, on a \'egalit\'e, puisque $m_{\Delta_{V"}}(F"):=\inf\{m_{\Delta}(E)/t\}$, pour $E,t$ comme ci-dessus. Donc $u_*(\Delta_{V"})=\Delta_{V'}$. 

Assertion 2. R\'esulte de l'assertion 1 si $V$ est une courbe $\square$

\

Dans le cas d'une courbe, on peut calculer comme suit la restriction (par application directe des d\'efinitions):

\begin{proposition}\label{rescou} Soit $(X\vert \Delta)$ une orbifolde lisse, $X$ compacte connexe, et $f: C\to X$ un morphisme, birationnel sur son image, d'une courbe lisse et compacte connexe $C$ telle que $V:=f(C)$ ne soit pas contenue dans $Supp(\Delta)$. Alors, pour tout point $a\in C$, la multiplicit\'e de $a$ dans la restriction minimale de $\Delta$ \`a $C$ est $m(a):=\sup_{j\in J(a)}\{\frac{m_j}{t_{j,a}}\}$, si $\Delta:=\sum_j(1-\frac{1}{m_j}).D_j$, $J(a):=\{j\vert f(a)\in D_j\}$, et $f^*(D_j)=t_{j,a}.\{a\}+\dots$, pour $j\in J(a)$\footnote{Dans $Georb^{div}$, on a: $m(a):=ppcm_{j\in J(a)}\{\frac{m_j}{pgcd(m_j,t_{j,a})}\}$.}. 
\end{proposition}

On peut \'egalement calculer la restriction en situation immerg\'ee. Nous traitons le cas d'une courbe.

\begin{theorem}\label{rescouim} Soit $(X\vert \Delta)$ une orbifolde lisse, $X$ compacte. Soit $V\subset X$ une courbe compacte et irr\'eductible non contenue dans $Supp(\Delta)$. Il existe une modification propre $u:X'\to X$ telle que les propri\'et\'es A et B suivantes soient satisfaites:

A. le diviseur $u^{-1}(Supp (\Delta))$ est \`a croisements normaux sur $X'$, qui est lisse.

B.  La transform\'ee stricte $V'$ de $V$ dans $X'$ est lisse, et chacun de ses points d'intersection avec $u^{-1}(Supp (\Delta))$ est un point lisse de ce diviseur (donc contenu dans une unique composante irr\'eductible de $u^{-1}(Supp (\Delta)))$ en lequel l'intersection est transversale. 

On munira alors $X'$ du diviseur orbifolde $\Delta'$ minimum tel que $u:(X'\vert \Delta')\to (X\vert \Delta)$ soit un morphisme orbifolde bim\'eromorphe \'el\'ementaire. On notera que $Supp(\Delta')\subset u^{-1}(Supp(\Delta))$.
Dans cette situation, on dira que la modification $u:(X'\vert \Delta')\to (X\vert \Delta)$ est {\bf transverse \`a $V$}.

Alors: la restriction (transverse) de $\Delta'$ \`a $V'$ est \'egale \`a la restriction de $\Delta$ \`a $V$.
\end{theorem}

{\bf D\'emonstration:}  L'existence de $u$ est une cons\'equence des th\'eor\`emes d'Hironaka. On note: $\Delta:=\sum_j(1-\frac{1}{m_j}).D_j$.
Soit $a\in V'$ un point d'intersection, transversale, avec une unique composante $E$ de $u^{-1}(Supp(\Delta))$. On a donc: $m(E)=max_{j\in J(E)}\{\frac{m_j}{t_{j,E}}\}$, avec: $J(E):=\{j$ tels que: $ u(E)\subset D_j\}$, et $u^*(D_j)=t_{j,E}.E+\dots$ , pour $j\in J(E)$. 

On note:  $J(a):=\{j$ tels que: $ u(a)\in D_j\}$. 

On a donc: $J(a)=J(E)$, puisque $u(a)\in u(E)\subset Supp(\Delta)$.

La multiplicit\'e de $a$ dans la restriction de $\Delta$ est: $m(a):=max_{j\in J(a)}\{\frac{m_j}{s_{a,j}}\}$, tandis que celle de $a$ dans la restriction de $\Delta'$ est: $m'(a):=max_{j\in J(a)}\{\frac{m_j}{t_{j,E}}\}$. 

On note $\nu:V'\to X$ la restriction de $u$ \`a $V'$, qui coincide donc avec l'injection de $V$ dans $X$ compos\'ee avec la normalisation de $V$. On a, pour tout $j$ tel que $\nu(a)\in D_j$, pour les nombres d'intersection locaux en $a$: $$s_{a,j}:=V'.\nu^*(D_j)=V'.u^*(D_j)=V'.t_{j,E}.E)=t_{j,E},$$ par les conditions de transversalit\'e et d'unicit\'e A et B ci-dessus.

L'\'egalit\'e pr\'ec\'edente montre que les deux multiplicit\'es $m(a)$ et $m'(a)$ coincident donc $\square$

\begin{re}\label{restmer} La restriction ne commute pas avec l'\'equivalence bim\'eromorphe. Soit $u:(X'\vert \Delta)\to (X\vert \Delta)$ un morphisme orbifolde bim\'eromorphe \'el\'ementaire, et $V$ une courbe lisse rencontrant transversalement le support de $\Delta$. En g\'en\'eral, les restrictions de $\Delta$ \`a $V$ et de $\Delta'$ \`a $V'$ ne coincident pas. (Consid\'erer par exemple $X=\Bbb P^2, \Delta=0, V=$ une droite, et $u:X'\to \Bbb P^2$ l'\'eclatement d'un point de $V$, $E$ le diviseur exceptionnel, $\Delta':=(1-\frac{1}{m}).E,m>1)$. M\^eme si $V$ est une courbe, et si $v:(X'\vert \Delta')\to (X\vert \Delta)$ est un morphisme orbifolde bim\'eromorphe, avec $\Delta'$ minimale rendant $v$ un morphisme orbifolde, il peut se faire que les restrictions de $\Delta'$ \`a $V'$ et de $\Delta$ \`a $V$ ne coincident pas (si $\Delta'$ n'est pas transverse \`a $V')$. Voir l'exemple \ref{str}.

\subsection{Restriction \`a une courbe g\'en\'erique d'une famille couvrante.}

\begin{definition}\label{coucou} Soit $(C_t)_{t\in T}$ une famille analytique de courbes complexes compactes de la vari\'et\'e $X$. On suppose $T$ compact, normal et irr\'eductible, et $C_t$ irr\'eductible et r\'eduite pour $t\in T$ g\'en\'erique. On note $Z\subset T\times X$ le graphe d'incidence ($T$-propre et $T$-connexe). On note $p:Z\to X$ et $q:Z\to T$, les projections. La famille est dite {\bf couvrante} si $p$ est surjective.

La famille est dite {\bf sans point-base} si, pour tout $A\subset X$, analytique ferm\'e de codimension au moins $2$, $C_t$ ne rencontre pas $A$, pour $t\in T$ g\'en\'erique. 

\end{definition}

Remarquons que si $(C_t)_{t\in T}$ est comme ci-dessus, il existe une modification $u:X'\to X$ telle que la famille $(C'_t)_{t\in T}$ de courbes de $X'$ dont le membre g\'en\'erique est la transform\'ee stricte de $C_t$ pour $t\in T$ g\'en\'erique soit sans point-base sur $X'$. On obtient $X'$ en aplatissant le morphisme propre $p:Z\to X$, et en prenant pour graphe d'incidence de la famille $(C'_t)_{t\in T}$ la composante principale de $Z\times _XX'$. L'absence de points-base est pr\'eserv\'ee pour toute famille obtenue sur une modification de $X'$, que l'on supposera lisse.

Soit alors donn\'e, de plus, un diviseur orbifolde $\Delta$ sur $X$, tel que $(X\vert \Delta)$ soit lisse. On supposera aussi, quitte \`a modifier $X'$, que les propri\'et\'es suivantes sont satisfaites:

A. le diviseur $u^{-1}(Supp (\Delta))$ est \`a croisements normaux.

B. Pour $t\in T$ g\'en\'erique, $C'_t$ est lisse, et chacun de ses points d'intersection avec $u^{-1}(Supp (\Delta))$ est un point lisse de ce diviseur (donc contenu dans une unique composante irr\'eductible de $u^{-1}(Supp (\Delta))$ en lequel l'intersection est transversale. 

On munira alors $X'$ du diviseur orbifolde $\Delta'$ minimum tel que $u:(X'\vert \Delta')\to (X\vert \Delta)$ soit un morphisme orbifolde bim\'eromorphe \'el\'ementaire. 

On notera que $Supp(\Delta')\subset u^{-1}(Supp(\Delta))$, et que la restriction transverse de $\Delta'$ \`a $C'_t$ est bien d\'efinie pour $t\in T$ g\'en\'erique (par la condition B de transversalit\'e). On dira que la modification $u$ est {\bf transverse \`a la famille $T$}.

\begin{theorem}\label{rescougen} Si $(C_t)_{t\in T}$ est une famille analytique de courbes de $X$ comme ci-dessus, si $(X\vert \Delta)$ est lisse, fix\'ee, et si $u:(X'\vert \Delta')\to (X\vert \Delta)$ est une modification choisie comme ci-dessus, alors la restriction de $\Delta'$ \`a $C'_t$ coincide avec la restriction minimale de $\Delta$ \`a $C_t$ (d\'efinie comme en \ref{rescoumin} et \ref{rescou}).
\end{theorem}

{\bf D\'emonstration:} C'est une cons\'equence de \ref{rescouim}, appliqu\'ee au membre g\'en\'erique $C'_t$ de la famille $\square$

\begin{re} Le r\'esultat pr\'ec\'edent et sa d\'emonstration (aux modifications usuelles pr\`es) restent valables dans $Georb^{div}$.
\end{re}

\subsection{Fibres orbifoldes d'une fibration m\'eromorphe.}\label{fibmor}

\begin{definition}\label{fiborb}
Soit $f: (X\vert\Delta)\to Y$ une fibration (propre). On suppose l'orbifolde g\'eom\'etrique $(X\vert\Delta)$ lisse. On d\'efinit alors la {\bf fibre g\'en\'erique orbifolde} $(X\vert\Delta)_y:=(X_y\vert\Delta_y)$ comme \'etant, pour $y\in Y$ g\'en\'erique, la restriction de $\Delta$ \`a $X_y:=f^{-1}(y)$. C'est-\`a-dire que si $\Delta=\sum_{D\in W(X)}(1-\frac{1}{m_{\Delta}(D)}).D$, on a: $\Delta_y:=\sum_{D\in W(X)}(1-\frac{1}{m_{\Delta}(D)}).(D\cap X_y)$. 

Le th\'eor\`eme de Sard montre alors que $(X_y\vert\Delta_y)$ est lisse. 
\end{definition}

\begin{definition}\label{bimeq} On dit que $f$ et $f'$ sont (bim\'eromorphiquement) {\bf \'el\'ementairement \'equivalentes} s'il existe un diagramme commutatif:

\

\centerline{
\xymatrix{ (Y'\vert\Delta' )\ar[r]^{v}\ar[d]_{f'} & (Y\vert\Delta)\ar[d]^{f}\\
X'\ar[r]^u&X\\
}}

\

\

dans lequel $u,v$ sont bim\'eromorphes, et $v$ un morphisme orbifolde tel que $v_*(\Delta')=\Delta$.

Plus g\'en\'eralement, $f$ et $f'$ sont \'equivalentes (on note alors: $f\sim f')$ si elles le sont pour la relation d'\'equivalence engendr\'ee par de tels diagrammes.
\end{definition}

\begin{proposition}\label{fib'} Si $f:(X\vert\Delta)\to Y$ et $f':(X'\vert\Delta')\to Y'$ sont deux fibrations bim\'eromorphiquement \'equivalentes, et si $(X\vert\Delta)$ et $(X'\vert\Delta')$ sont lisses, leurs fibres orbifoldes g\'en\'eriques correspondantes sont bim\'eromorphiquement \'equivalentes. 
\end{proposition}

{\bf D\'emonstration.} Il suffit de d\'emontrer l'assertion lorsque $f$ et $f'$ sont \'el\'ementairement \'equivalentes.  Dans ce cas, il suffit de choisir $x\in X$ tel que les fibres orbifoldes de $f$ et $f'$ au-dessus de $x$ soient lisses. Alors $v$ induit une \'equivalence orbifolde \'el\'ementaire entre ces deux fibres $\square$

 \begin{re}\label{dimfib} Par contre, si $f$ n'est pas presque-holomorphe\footnote{Voir d\'efinition \ref{ph}.}, la fibre orbifolde g\'en\'erique $(X_y\vert\Delta_{X_y})$ de $f$ n'est pas bien d\'efinie, \`a \'equivalence orbifolde pr\`es, et sa dimension canonique d\'epend en g\'en\'eral du mod\`ele bim\' eromorphe $(X\vert \Delta)$ choisi. Consid\'erer, par exemple, $f:X:=\Bbb P^2\dasharrow \Bbb P^1$ la fibration d\'efinie par les coniques passant par $4$ points $a_j\in \bP^2$ en position g\'en\'erale. Les fibres de cette fibration sont donc les coniques de la famille, sans structure orbifolde. Si l'on \'eclate les points $a_j$, et si l'on d\'enote par $A_j$ les diviseurs exceptionnels de cet \'eclatement, alors la fibration $f': X'\to \Bbb P^1$ obtenue est holomorphe. Munissant $X'$ de la structure orbifolde $\Delta_m':= \sum_j(1-\frac{1}{m}). A_j$, les fibres lisses de $f'$ rencontrent transversalement les $A_j$, et sont donc munies de la structure orbifolde induite consistant en quatre points munis de la multiplicit\'e $m$, et leur dimension canonique est donc $-\infty$ (resp. $0$, resp. $1)$ si $m=1$ (resp. $m=2$, resp. $m\geq 3)$. Observer que, pour tout $m$, $(X\vert 0)$ et $(X'\vert\Delta_m')$ sont cependant bim\'eromorphes. 
 \end{re}

Pour les fibres g\'en\'eriques d'une application m\'eromorphe, nous allons cependant d\'efinir une notion de restriction, unique \`a \'equivalence bim\'eromorphe pr\`es, si $(X\vert\Delta)$ est fix\'ee (mais cette restriction d\'epend, \`a \'equivalence bim\'eromorphe orbifolde pr\`es) du choix de $(X\vert\Delta)$, comme le montre \ref{dimfib}.\end{re}

\begin{proposition}\label{ibrest} Soit $(X\vert \Delta)$ une orbifolde g\'eom\'etrique lisse, avec $X$ compacte et connexe. Soit $f:X\dasharrow Y$ une application m\'eromorphe surjective \`a fibres g\'en\'eriques $X_y$ connexes\footnote{La fibre g\'en\'erique $X_y$ est l'image sur $X$ de la fibre correspondante du graphe de $f$. Elle est bien d\'efinie et ind\'ependante de $X,f,Y$, \`a \'equivalences bim\'eromorphes pr\`es.}. On peut alors d\'efinir alors la restriction de $\Delta$ \`a $X_y$ comme la restriction de $\Delta'$ \`a la fibre g\'en\'erique de $f':(X'\vert \Delta')\to Y$ si $g:(X'\vert\Delta')\to (X\vert \Delta)$ est un morphisme orbifolde bim\'eromorphe \'el\'ementaire tel que $f':=f\circ g$ soit holomorphe. Cette restriction g\'en\'erique est alors bien d\'efinie \`a \'equivalence bim\'eromorphe orbifolde pr\`es.
\end{proposition}

{\bf D\'emonstration:} Soit $g:(X'\vert\Delta')\to (X\vert \Delta)$ un morphisme orbifolde bim\'eromorphe \'el\'ementaire (avec $(X'\vert \Delta')$ lisse) tel que $f':=f\circ g:X'\to Y$ soit holomorphe. On note $\Delta_y$ la restriction de $\Delta$ \`a $X'_y$ (au sens usuel de \ref{fiborb}). Alors $g_y:(X'_y\vert\Delta'_y)$ est une restriction de $\Delta$ \`a $X_y$, minimale si $\Delta'$ est minimal rendant $g$ un morphisme orbifolde. Il nous suffit donc de montrer que si $v:Y'\to Y$ est une modification propre, et $h:(X"\vert \Delta")\to (X'\vert \Delta')$ une modification orbifolde bim\'eromorphe \'el\'ementaire telle que $f":=v^{-1}\circ f'\circ h:X"\to Y'$ soit holomorphe, alors $(X_y"\vert \Delta_y")\to (X_y'\to \Delta_y')$ est un morphisme orbifolde bim\'eromorphe (\'el\'ementaire si $\Delta"$ et $\Delta'$ sont minimales rendant $g$ et $h$ des morphismes orbifoldes. Remarquer que si $D"_{min}$ est minimal rendant $g\circ h$ un morphisme orbifolde, alors $(X"\vert \Delta")$ et $(X"\vert \Delta"_{min})$ sont bim\'eromorphiquement \'equivalentes). C'est une v\'erification imm\'ediate, puisque $h^*$ et $r_y$ commutent, si $r_y$ d\'esigne la restriction des diviseurs au-dessus de $X_y$ $\square$

\begin{re}\label{rfo}

1. Lorsque $f$ est presque holomorphe (ie: si le lieu d'indetermination de $f$ ne rencontre pas sa fibre g\'en\'erique $X_y$), alors $(X_y\vert \Delta_y)$ est d\'ej\`a une restriction de $\Delta$ \`a $X_y$ (par le th\'eor\`eme de Sard-Bertini).

2. Les notions relatives aux fibres (orbifoldes) se comportent de mani\`ere ``duale" de celles relatives aux bases (orbifoldes, d\'efinies dans le \S \ref{basfibmor} suivant): les fibrations int\'eressantes sont celles dont les fibres ont une dimension canonique petite ($0,-\infty$), ou sont $\Delta$-rationnelles (voir \S\ref{dr}), ou celles dont la base est de type g\'en\'eral ou du moins de dimension canonique nulle ou positive. Ceci justifie la d\'efinition suivante:
\end{re}

\begin{definition}\label{dfo} Soit $f:(X\vert \Delta)\dasharrow Y$ une fibration m\'eromorphe, avec $(X\vert \Delta)$ lisse et $X$ compacte et connexe. On dit que la fibre orbifolde g\'en\'erique (not\'ee $(X_y\vert \Delta_y)$) poss\`ede la propri\'et\'e $P$ {\bf s'il existe une} fibration holomorphe $f':(X'\vert \Delta')\to Y'$ bim\'eromorphiquement \'equivalente \`a $f$ (avec, bien s\^ur, $(X'\vert \Delta')$ lisse), et dont la fibre orbifolde g\'en\'erique poss\`ede la propri\'et\'e $P$.
\end{definition}

\begin{re}\label{rrfo} 1. Les propri\'et\'es $P$ que nous consid\`ererons sont les suivantes: $\kappa=0$, $\kappa=-\infty$, $ \kappa_+=-\infty$, ou bien: $\Delta$-unir\'egl\'ee ou $\Delta$-rationnellement connexe, ou encore: \^etre ``sp\'eciale". Voir le \S \ref{dr} et le  \S\ref{orbspec} pour ces six derni\`eres notions.

2. Dans cette situation, il n'est pas toujours vrai que pour les propri\'et\'es $P$ pr\'ecedentes, et pour {\bf tout} mod\`ele bim\'eromorphe holomorphe $ f":(X"\vert \Delta")\to Y"$ de $f$ les fibres orbifoldes g\'en\'eriques poss\`edent encore la propri\'et\'e $P$. Voir l'exemple \ref{dimfib}.

3. Cependant, cette propri\'et\'e $P$ sera pr\'eserv\'ee sur {\bf tout} mod\`ele bim\'eromorphe pour les trois fibrations fondamentales que nous construirons et consid\'ererons ici:

La fibration de Moishezon-Iitaka ($P$ est: $\kappa=0)$ (voir le th\'eor\`eme \ref{klo}).

Le quotient $\kappa$-rationnel ($P$ est: $\kappa_+=-\infty)$. En effet cette fibration est presque-holomorphe (si l'orbifolde consid\'er\'ee est finie). Voir \S\ref{krat}.

Le ``coeur": ( $P$ est: \^etre sp\'eciale). En effet: le ``coeur" est toujours presque-holomorphe. Voir le th\'eor\`eme \ref{c}.

En effet, lorsque $\kappa(f\vert \Delta)\geq 0$, et si $(X\vert \Delta)$ est finie, toutes les mod\`eles bim\'eromorphes de $(f\vert \Delta)$ ont un mod\`ele bim\'eromorphe presque-holomorphe, et les fibres orbifoldes g\'en\'eriques ont donc une classe d'\'equivalence bim\'eromorphe bien d\'efinie. La m\^eme conclusion est vraie (sans restriction de finitude si $(f\vert \Delta)$ est de type g\'en\'eral. Voir le \S \ref{gtph}. 
\end{re}

\subsection{R\'esolution d'une orbifolde g\'eom\'etrique.}

On peut d\'efinir comme suit la notion de mod\`ele lisse d'une orbifolde g\'eom\'etrique arbitraire $(X\vert\Delta)$ , d\'efinie sur un espace complexe normal $X$ alg\'ebriquement $\bQ$-factoriel, c'est-\`a-dire tel que tout diviseur de Weil sur $X$ soit $\bQ$-Cartier. On note $Sing(X)$ le lieu singulier de $X$.

   Soit $\Delta=\sum_j (1-\frac{1}{m_j}).D_j$ un diviseur orbifolde de $X$ de support $D,$ o\`u $m_j\in(\bQ\cup \infty),\forall j.$

\begin{definition} \label{res} Une {\bf r\'esolution} $p:(Y\vert\Delta_Y)\to (X\vert\Delta)$ de l'orbifolde g\'eom\'etrique $(X\vert\Delta)$ est une d\'esingularisation $p:Y\to X$ de $X$ telle que $p^{-1}(D\cup Sing(X))$ soit un diviseur \`a croisements normaux sur $Y$, et $(Y\vert\Delta_Y)$ une orbifolde g\'eom\'etrique lisse de $Y$ de support contenu dans $p^{-1}(D\cup Sing(X))$, telle que:

1. $p_*(\Delta_Y)=\Delta$, et:

2. $p:(Y\vert\Delta_Y)\to (X\vert\Delta)$ soit un morphisme orbifolde \footnote{au sens divisible si l'on est dans $Georb^{div}$.}.
\end{definition}

\begin{re} \label{resmin} 

\

1.Si $p:(Y\vert\Delta_Y)\to (X\vert\Delta)$ et $p:(Z\vert\Delta_Z)\to (X\vert\Delta)$ sont deux telles r\'esolutions, et si $u:V\to Y$ et $v:V\to Z$ sont bim\'eromorphes (propres), , telles que les images r\'eciproques de $(D\cup Sing(X))$ par $p\circ u$ et $p\circ v$ soient des diviseurs \`a croisements normaux, il existe alors une plus petite orbifolde g\'eom\'etrique $\Delta_V$ sur $V$ telle que $u:(V\vert\Delta_V)\to (Y\vert\Delta_Y)$ et $v:(V\vert\Delta_V)\to (Z\vert\Delta_Z)$ soient des morphismes orbifoldes bim\'eromorphes, de sorte que $p\circ u:(V\vert\Delta_V)\to (X\vert\Delta)$ et $p\circ v:(V\vert\Delta_V)\to (X\vert\Delta)$ sont aussi des r\'esolutions. En g\'en\'eral, $u$ et $v$ ne sont pas des morphismes bim\'eromorphes {\it \'el\'ementaires}, et  il n'est donc pas \'evident que $(Y\vert \Delta_Y)$ et $(Z\vert \Delta_Z)$ soient bim\'eromorphiquement \'equivalentes (dans la classe des orbifoldes g\'eom\'etriques lisses). Voir cependant la question \ref{qsing} ci-dessous.

2. Lorsque $X$ n'est pas $\bQ$-factorielle, la notion de r\'esolution d'une orbifolde g\'eom\'etrique n'est pas d\'efinie ici (et il ne semble pas imm\'ediat d'en donner une d\'efinition naturelle). La situation est donc, ici encore, analogue \`a celle des r\'esolutions logarithmiques du programme des (Log) mod\`eles minimaux.

3. La notion de r\'esolution d\'efinie ici sera utilis\'ee pour definir les orbifoldes l.c ou klt. 
\end{re}

\begin{question}\label{qsing} Soit $p:(Y\vert\Delta_Y)\to (X\vert\Delta)$ une r\'esolution d'une orbifolde g\'eom\'etrique, avec $X$ $\Bbb Q$-factorielle. 

Pour tous $N>0,q>0$, le faisceau $S^N(\Omega^q(X\vert\Delta)):=p_*(S^N(\Omega^q(Y\vert\Delta_Y))$ est-il ind\'ependant de la r\'esolution $p$?

Si oui, alors:  $p^*(H^0(X,S^N(\Omega^q(X\vert\Delta))))=H^0(Y, (S^N(\Omega^q(Y\vert\Delta_Y))))$ est ind\'ependant de la r\'esolution choisie, pour tous $N,q$.
\end{question}

\subsection{Orbifoldes g\'eom\'etriques log-canoniques et klt}

Ces notions sont d\'efinies comme dans le cadre du PMML (=LMMP en Anglais), avec la condition additionnelle de $\bQ$-factorialit\'e. Il est crucial de pouvoir \'etendre \`a ce cadre \'elargi les r\'esultats du pr\'esent texte.

\begin{definition} 

\

$\bullet$  L'orbifolde g\'eom\'etrique $(X\vert\Delta)$ est dite {\bf log-canonique} si:

1. $X$ est $\bQ$-factorielle.

2. Le faisceau $\omega_X=K_X$ est $\bQ$-Cartier.

3.Pour toute (ie: pour une) r\'esolution $p:(Y\vert\Delta_Y)\to (X\vert\Delta)$, on a: $K_X+\Delta_Y\geq p^*(K_X+\Delta)$. 

$\bullet$ L'orbifolde g\'eom\'etrique $(X\vert\Delta)$ est dite {\bf klt} si, pour une (ie: pour toute) r\'esolution $p:(Y\vert\Delta_Y)\to (X\vert\Delta)$, on a: $(K_X+\Delta_Y)-p^*(K_X+\Delta):=E$ est un diviseur effectif support\'e sur le diviseur exceptionnel de $p$ tout entier.

\end{definition}

\begin{example} Les orbifoldes g\'eom\'etriques lisses sont log-canoniques, et klt si les multiplicit\'es sont finies.

\end{example}

\begin{re}\label{icanlc} Soit $p:(Y\vert\Delta_Y)\to (X\vert\Delta)$ une r\'esolution d'une orbifolde g\'eom\'etrique log-canonique. Alors $p^*(H^0(X,N(K_{(X\vert\Delta)}))=H^0(Y, N.K_{(Y\vert\Delta_Y)})$ est ind\'ependant de la r\'esolution choisie, pour tous $N>0$. Ceci r\'esulte de \ref{invbir'}.
\end{re}

\begin{question}\label{qlc} Supposons que $(X\vert\Delta)$ soit une orbifolde g\'eom\'etrique log-canonique, et que les faisceaux $S^N(\Omega^q(X\vert\Delta)))$ soient biens d\'efinis (ie: ind\'ependants de la r\'esolution  de $(X\vert\Delta)$, comme dans \ref{qsing}). Les faisceaux $S^N(\Omega^q(X\vert\Delta)))$ sont-ils alors r\'eflexifs?

Pour $q=n:=dim(X)$, c'est une cons\'equence de la d\'efinition. Lorsque $q=1,$ ou $q=(n-1)$, et si $N=1$, une r\'eponse positive est fournie dans [G-K-K].
\end{question}

\newpage

\section{BASE ORBIFOLDE}\label{basfibmor}

Nous allons ici d\'efinir la base orbifolde d'une fibration, en \'etablir certaines propri\'et\'es bim\'eromorphes et calculer celle d'une compos\'ee. Les r\'esultats sont \'etablis et \'enonc\'es dans la cat\'egorie $Georb^{Q}$ (voir \ref{defbsorb} ci-dessous), mais restent valables, ainsi que leurs d\'emonstrations (en y rempla\c cant in\'egalit\'e par divisibilit\'e), dans $Georb^{div}$.

\subsection{Base orbifolde d'un morphisme}\label{basmor}

\begin{definition} On appellera {\bf morphisme} toute application  $f:Y\to X$, holomorphe propre et surjective, entre espaces analytiques complexes $Y,X$ normaux. Une {\bf fibration} est un morphisme \`a fibres connexes.

\end{definition}

Soit $\Delta_Y:=\Delta$ un diviseur orbifolde sur $Y$. 

On note $m(\Delta_Y):W(Y)\to \bQ \cup \{\infty\}:=\bar{\Bbb Q}$ la multiplicit\'e de $\Delta_Y$.

Si $D\in W(X)$,  et si $f:Y\to X$ est un morphisme, on \'ecrit: $f^*(D)=\sum_{j\in J(f,D)}m_j.D_j+R$, o\`u $R$ est le plus grand diviseur de $X$ de support contenu dans $f^*(D)$, et $f$-exceptionnel (ie: tel que la codimension de $f_*(R)$ dans $Y$ soit au moins $2)$. Remarquer que les entiers $m_j$ sont bien d\'efinis, m\^eme si $X$ n'est pas $\bQ$-factorielle.

\begin{definition}\label{defbsorb} (Voir [Ca04, 1.29, p. 523]) Soit $(Y\vert\Delta_Y)$ une orbifolde g\'eom\'etrique. La {\bf base orbifolde} $(X\vert\Delta(f,\Delta_Y))$ du morphisme $f:(Y\vert\Delta_Y)\to X$ est l'orbifolde g\'eom\'etrique $(X\vert\Delta(f,\Delta_Y))$ d\'efinie par la multiplicit\'e $m(f,\Delta_Y):W(X)\to \bQ\cup \{\infty\}$ telle que, pour tout $D\in W(X)$: $m(f,\Delta_Y;D):=inf_{j\in J(f,D)}\lbrace m_j.m_{\Delta_Y}(E_j)\rbrace$. 

\

Lorsque $\Delta_Y=0$, $\Delta(f,0)$ est simplement not\'ee $\Delta(f)$.

\

La d\'efinition pr\'ec\'edente est celle des cat\'egories $Georb^Q$ et $Georb^Z$ (voir d\'efinition \ref{morphorb*}). 

Dans $Georb^{div}$, on remplace:

 $m(f,\Delta_Y;D):=inf_{j\in J(f,D)}\lbrace m_j.m_{\Delta_Y}(E_j)\rbrace$ ci-dessus par: 
 
 $m(f,\Delta_Y;D):=pgcd_{j\in J(f,D)}\lbrace m_j.m_{\Delta_Y}(E_j)\rbrace$.
 
On pr\'ecisera si n\'ecessaire la cat\'egorie consid\'er\'ee en notant: $(X\vert\Delta(f,\Delta_Y))^*$, ou: $\Delta(f,\Delta_Y))^*$, avec: $*=Z,Q,div$ selon le cas.
 \end{definition}

\begin{re} 

\

1. M\^eme si $X$ et $(Y\vert\Delta_Y)$ sont lisses, $(X\vert\Delta(f,\Delta_Y))$ n'est pas lisse, en g\'en\'eral.

2. Si $f:(Y\vert\Delta)\to (X\vert\Delta')$ est un morphisme orbifolde, on a: $\Delta'\leq \Delta(f,\Delta)$. En g\'en\'eral, si $X$ est lisse, $f:(Y\vert\Delta)\to (X\vert\Delta(f,\Delta))$ n'est pas un morphisme orbifolde (\`a cause des diviseurs $f$-exceptionnels de $Y$, n\'eglig\'es dans la d\'efinition de $\Delta(f,\Delta_Y))$. Si $f$ est un morphisme fini, c'est cependant le cas (ceci r\'esulte de \ref{comporb}, assertion 2,  ci-dessous). On peut toujours, en augmentant les multiplicit\'es des diviseurs $f$-exceptionnels, faire de $f:(Y\vert\Delta^+_Y)\to(X\vert\Delta(f,\Delta))$ un morphisme orbifolde sans changer $\Delta(f,\Delta)$. 

3. Si $X$ est une courbe (lisse), $f:(Y\vert\Delta)^*\to (X\vert\Delta(f,\Delta_Y)^*)$ est un morphisme dans la cat\'egorie $Georb^*$, avec $*=Z,Q,div$, puisqu'il n'y a alors pas de diviseur $f$-exceptionnel. Dans ce cas, $(X\vert\Delta(f,\Delta_Y)^*$ est aussi la plus grande structure orbifolde sur $X$ rendant $f$ un morphisme dans la cat\'egorie $Georb^*$.

4. Si $f=g\circ h$, avec $h:Y\to X'$ une fibration, et $g:X'\to X$ finie (donc $g\circ h$ est la factorisation de Stein du morphisme $f)$, alors: $\Delta(f,\Delta)=\Delta(g,\Delta(h,\Delta))$ peut \^etre calcul\'ee en $2$ \'etapes. (Ceci r\'esulte aussi de \ref{comporb} ci-dessous). Si $g$ n'est pas finie, on a seulement: $\Delta(f,\Delta)\geq \Delta(g,\Delta(h,\Delta))$. Voir ci-dessous pour les cas d'\'egalit\'e.

5. La d\'efinition \ref{defbsorb} est motiv\'ee, entre autres, par la propri\'et\'e suivante ([Ca04, 1.30]): si $\Delta_Y=\Delta(g)$, pour une fibration $g:Z\to Y$, alors $\Delta(f,\Delta_Y)=\Delta(f,\Delta(g))\geq\Delta(f\circ g)$, la diff\'erence provenant des diviseurs de $Z$ qui sont $g$-exceptionnels, mais non $f\circ g$-exceptionnels. Sur des modifications ad\'equates de $Z,X$ et $Y$, on obtient ([Ca04, 1.31]) des fibrations $g':Z'\to Y', f':Y'\to X'$ telles que $\Delta(f'\circ g')=\Delta(f', \Delta(g'))$.

Nous allons g\'en\'eraliser cette relation au cas o\`u l'on a un diviseur orbifolde initial sur $Z$. Les arguments sont cependant essentiellement les m\^emes. Pour cela, quelques r\'esultats pr\'eliminaires doivent \^etre \'etablis.\end{re}

\begin{lemma}\label{invbir} Soit $v,f,f'=f\circ v$ des applications holomorphes, avec $v:(Y'\vert\Delta')\to (Y\vert\Delta)$ un morphisme orbifolde bim\'eromorphe tel que $\Delta=v_*(\Delta')$:

\centerline{
\xymatrix{ (Y'\vert\Delta' )\ar[r]^{v}\ar[rd]_{f'} & (Y\vert\Delta)\ar[d]^{f}\\
& X\\
}}

Alors: $\Delta(f,\Delta)=\Delta(f',\Delta')$.
\end{lemma}

{\bf D\'emonstration:} Soit $D\in W(X)$, alors:

$f'^*(D)=v^*(f^*(D))=v^*(\sum_jm_j.D_j+R)=\sum_jm_j.(v^*(D_j))+v^*(R)=\sum_j m_j.\overline{D_j}+\sum_{(j,h)}m_j.n_{(j,h)}E_{(j,h)}+v^*(R)=\sum_jm_j.\overline{D_j}+R'$, o\`u $\overline{D_j}$ est le transform\'e strict de $D_j$ par $v$. Bien s\^ur, $v^*(R)$ est $f'$-exceptionnel. 

D'autre part, puisque $v$ est un morphisme orbifolde, on a, pour tous $(j,h):$
$$n_{(j,h)}.m'(E_{(j,h)})\geq m(D_j),$$ et donc: $n_{(j,h)}.m_j.m'(E_{(j,h)})\geq m_j.m(D_j)$. On a not\'e $m:W(Y)\to \bar{\Bbb Q}$, et  $m':W'(Y)\to \bar{\Bbb Q}$ les fonctions de multiplicit\'e associ\'ees \`a $\Delta$ et $\Delta'$ respectivement. 

Puisque $v_*(\Delta')=\Delta$, on a de plus: $m(D_j)=m'(\overline{D_j})$.

Notant $m_f$ et $m'_f$ les multiplicit\'es d\'efinissant $\Delta(f,\Delta)$ et $\Delta(f',\Delta')$ respectivement, on a donc:
$$m_f(D)=inf_j\lbrace m_j.m(D_j)\rbrace\leqno (1)$$
tandis que:
$$m'_f(D)=inf_{(j,h)}\lbrace m_j.m'(\overline{D_j}), m_j.n_{(j,h)}.m'(E_{(j,h)})\rbrace\leqno (2)$$

On d\'eduit alors des relations (1) et (2) pr\'ec\'edentes que $m_f(D)=m'_f(D)$, ce qui est l'assertion du lemme $\square$

\begin{corollary}\label{invdelta'} Dans la situation du lemme \ref{invbir} pr\'ec\'edent, si $h:X\to W$ est une application holomorphe, avec $W$ lisse, alors $\Delta(h\circ f', \Delta')=\Delta(h\circ f,\Delta)$ est ind\'ependant de $\Delta'$, avec $\Delta=v_*(\Delta')$, $v$ \'etant un morphisme d'orbifoldes g\'eom\'etriques.
\end{corollary}

\begin{lemma}\label{compfib} Supposons que, dans  le diagramme commutatif suivant, $w$ {\bf et $g$} soient des morphismes orbifolde, que $w$ et $v$ soient bim\'eromorphes, et que les vari\'et\'es consid\'er\'ees soient $\bQ$-factorielles, connexes et compactes.

\centerline{
\xymatrix{ (Z'\vert\Delta' )\ar[r]^{w}\ar[d]_{g'} & (Z\vert\Delta)\ar[d]^{g}\\
Y'\ar[r]^v& Y\\
}}

Alors: $v:(Y'\vert\Delta(g',\Delta'))\to (Y\vert\Delta(g,\Delta))$ est un morphisme orbifolde. De plus, si $w_*(\Delta')=\Delta$, alors  $v_*(\Delta(g',\Delta'))=\Delta(g,\Delta)$.
\end{lemma}

{\bf D\'emonstration:} Soit $E\in W(Y), E'\in W(Y')$ tels que:

 $v^*(E)=s.E'+\dots$, avec $s\geq 1$. 

On veut montrer que $s.m'_g(E')\geq m_g(E)$, si $m_g,m'_g$ d\'esignent les fonctions de multiplicit\'e associ\'ees \`a $\Delta(g,\Delta)$ et $\Delta(g',\Delta')$ respectivement. De m\^eme, $m,m'$ sont celles de $\Delta$ et $\Delta'$.

Si $g'^*(E')=\sum_jm_j.F'_j+R'$, on a: $m'_g(E'):=\inf_j\lbrace m_j.m'(F'_j)\rbrace$. 

On choisit et fixe un $j$ tel que $g'(F'_j)=E'$ et $m_j.m'(F'_j)=m'_g(E')$. On a donc: $$(v\circ g')^*(E)=(g')^*(s.E'+...)=s.m_j.F'_j+...,\leqno(1)$$

puisque $E'$ est la seule composante irr\'eductible de $v^*(E)$ contenant $g'(F'_j)$. 

Si $g^*(E)=\sum_{k\in K} n_k.F_k$, on a donc: $m_g(E)=\inf_k\lbrace n_k.m(F_k)\}$, puisque $g$ est un morphisme orbifolde, par hypoth\`ese. 

Posant: $w^*(F_{k})=s_{(j,k)}.F'_j+\dots$, on d\'eduit de (1), en calculant la multiplicit\'e de $F'_j$ dans $(g\circ w)^*(E)$, utilisant $g\circ w=v\circ g'$ et le fait que $w$ est un morphisme orbifolde: 

$$s.m_j\geq n_k.s_{(j,k)},\forall k\leqno(2)$$

On d\'eduit de (1) et (2) que: $$s.m'_g(E')=s.m_j.m'(F'_j)\geq n_k.s_{(j,k)}.m'(F'_j)\geq n_k.m(F_k)\geq m_g(E),$$ et donc la conclusion:  $s.m'_g(E')\geq m_g(E)$ par application directe des d\'efinitions; l'avant-derni\`ere in\'egalit\'e r\'esultant du fait que $w$ est un morphisme orbifolde.

La seconde assertion r\'esulte du lemme \ref{invbir} $\square$

\begin{re} La condition ``$g$ morphisme orbifolde" est essentielle. Elle ne peut \^etre affaiblie en: ``$g'$ morphisme orbifolde".
\end{re}

\subsection{Fibrations nettes.}

\begin{definition}\label{nette} Soit $g:(Z\vert\Delta)\to Y$ une fibration, avec $(Z\vert\Delta)$ et $Y$ lisses. On dira que $g$ est {\bf nette}  (relativement \`a $w)$ s'il existe un diagramme:

\centerline{
\xymatrix{ (Z\vert\Delta )\ar[r]^{w}\ar[d]_{g} & (Z'\vert\Delta')\\
Y&\\
}}

dans lequel:

1. $w$ est un morphisme orbifolde, $v$ et $w$ sont bim\'eromorphes,  $Z',Y,Y'$ lisses, et $w_*(\Delta)=\Delta'$.

2. Tout diviseur $g$-exceptionnel de $Z$ est $w$-exceptionnel.

On dira que $g$ est {\bf nette} si elle est nette relativement \`a une fibration $g'$ comme ci-dessus.

 On dira que $g$ est  {\bf nette et haute} si elle est nette, et si $g: (Z\vert\Delta)\to (Y\vert\Delta(g,\Delta))$ est un morphisme orbifolde.

On dira que $g$ est {\bf strictement nette} si elle est nette, et si, de plus, sa base orbifolde est lisse.
\end{definition}

\begin{example}\label{netcourb} Toute fibration sur une courbe est strictement nette et haute.
\end{example}

\begin{proposition}\label{stneth}  Si $g':(Z'\vert\Delta')\to Y'$ est donn\'ee, il existe un diagramme commutatif:

\centerline{
\xymatrix{ (Z\vert\Delta )\ar[r]^{w}\ar[d]_{g} & (Z'\vert\Delta')\ar[d]^{g'}\\
Y\ar[r]^v& Y'\\
}}

 tel que:

1. $g:(Z\vert\Delta)\to Y$ est strictement nette relativement \`a $w$. 

2. $g:(Z\vert\Delta)\to (Y\vert\Delta(g,\Delta))$ est un morphisme orbifolde (ie: $g$ est strictement nette et haute).\end{proposition}

{\bf D\'emonstration:} On construit le diagramme commutatif ci-dessous:

\centerline{
\xymatrix{ Z\ar[r]^{w_1}\ar[rd]_g& Z_1\ar[r]^{w'}\ar[d]_{g_1} & (Z'\vert\Delta')\ar[d]^{g'}\\
&Y\ar[r]^v& Y'\\
}}

dans lequel: $v$ est une modification telle que le morphisme d\'eduit de $g'$ par $v$ soit plat sur la composante pincipale ([R72]), et tel que $Y$ soit lisse, avec $v^{-1}(Supp(\Delta(g',\Delta')))$ soit \`a croisements normaux sur $Y$. La d\'esingularisation de $Z_1$ fournit alors $Z$. On pose: $w:=w'\circ w_1$. Les diviseurs $g$-exceptionnels de $Z$ sont alors $w$-exceptionnels, par platitude de $g_1$. On munit $Z$ du diviseur orbifolde $\Delta$ tel que $w_*(\Delta)=\Delta'$. Pour les diviseurs $F\in W(Z)$ qui sont $w$-exceptionnels, leurs multiplicit\'e peut \^etre choisie arbitrairement assez grande. Ceci \'etablit la premi\`ere assertion.

La seconde assertion r\'esulte de \ref{compfib}, si les multiplicit\'es des diviseurs $g$-exceptionnels de $Z$ (qui sont donc $w$-exceptionnels) sont choisies assez grandes $\square$

Du lemme \ref{compfib}, on d\'eduit:

\begin{corollary}\label{compfib'} Supposons que, dans  le diagramme commutatif suivant, $w$ soit un morphisme orbifolde bim\'eromorphe \'el\'ementaire, et que $g$ soit une fibration haute et strictement nette, et enfin que $v$ soit bim\'eromorphe. 

\
\centerline{
\xymatrix{ (Z'\vert\Delta' )\ar[r]^{w}\ar[d]_{g'} & (Z\vert\Delta)\ar[d]^{g}\\
Y'\ar[r]^v& Y\\
}}

Alors: $v:(Y'\vert\Delta(g',\Delta'))\to (Y\vert\Delta(g,\Delta))$ est un morphisme orbifolde bim\'eromorphe \'el\'ementaire.

En particulier, la classe d'\'equivalence bim\'eromorphe des bases orbifoldes de fibrations \'equivalentes pour la relation d'\'equivalence engendr\'ee par de tels diagrammes est bien d\'efinie (ie: est ind\'ependante du repr\'esentant choisi). 
\end{corollary}

\begin{question}\label{eqbimbasorb}{\bf: \' Equivalence bim\'eromorphe des bases orbifoldes.}

Soit $(Z\vert\Delta)$ et $(Z'\vert\Delta')$ deux orbifoldes lisses bim\'eromorphiquement \'equivalentes, avec $Z,Z'$ compactes, k\" ahler et connexes. Si $g:Z\to Y$ et $g':Z'\to Y'$ sont deux fibrations bim\'eromorphiquement \'equivalentes (ie: avec la m\^eme famille de fibres g\'en\'eriques), avec $Y,Y'$ lisses, et si $g$ et $g'$ sont strictement nettes, leurs bases orbifoldes sont-elles bim\'eromorphiquement \'equivalentes? 
\end{question}

Nous \'etablirons dans le chapitre \S\ref{dimcanfib} suivant une propri\'et\'e plus faible, mais centrale pour les r\'esultats du pr\'esent texte: l' \'egalit\'e de la dimension canonique de deux telles bases orbifoldes ``stables".

\subsection{Composition de fibrations}\label{comfib}

Dans cette section, toutes les orbifoldes et vari\'et\'es consid\'er\'ees sont lisses, compactes et connexes.

On consid\`ere un diagramme commutatif de morphismes orbifolde entre orbifoldes g\'eom\'etriques lisses (compactes et connexes). On suppose que les fl\`eches verticales sont des fibrations.

\centerline{\xymatrix{  (Z\vert\Delta)\ar[d]^{g}\\
 Y\ar[d]^{f}\\X\\}}

\begin{proposition}\label{comporb} Soit $g:Z\to Y$ et $f:Y\to X$ deux fibrations, avec $Y,X$ lisses, et $\Delta$ un diviseur orbifolde sur $Z$. 

On peut donc construire trois diviseurs orbifoldes: $\Delta_Y:=\Delta(g,\Delta)$ sur $Y$, ainsi que $\Delta_{fg}:=\Delta(f\circ g,\Delta)$ et $\Delta_{f/g}:=\Delta(f,\Delta(g,\Delta))$ sur $X$. Alors:

1. $\Delta(f\circ g, \Delta)\leq \Delta(f,\Delta(g,\Delta))$.

2. On a $\Delta(f\circ g, \Delta)= \Delta(f,\Delta(g,\Delta))$ si $m_{\Delta}(F)$ est assez grand, pour tout $F\in W(Z)$ qui est $g$-exceptionnel, mais non pas $f\circ g$-exceptionnel.

3. Si $g:(Z\vert\Delta)\to (Y\vert\Delta(g,\Delta))$ est un morphisme orbifolde, on a: $\Delta(f\circ g, \Delta)= \Delta(f,\Delta(g,\Delta))$.

\end{proposition}

{\bf D\'emonstration:} Remarquons tout d'abord que si $F\in W(Z)$ est $fg$-exceptionnel, ou bien il est $g$-exceptionnel, ou bien $g(F)$ est $f$-exceptionnel. Soit alors $D\in W(X)$. Alors $(g^*(f^*))(D)=(\sum_jg^*(m_jE_j)+g^*(R)=\sum_{j,h}m_j.n_{j,h}F_{j,h}+R'$. Ici $R'$ est $fg$-exceptionnel, mais aucun des $F_{j,h}$ ne l'est. Notons $m_{fg}$ la multiplicit\'e de $\Delta(f\circ g,\Delta)$, et $m_{f/g}$ celle de $\Delta(f,\Delta( g,\Delta)).$

Donc: $m_{fg}(D):=\inf_{j,h}\lbrace m_j.n_{j,h}.m_{\Delta}(F_{j,h})\rbrace$. La somme porte sur les composantes $F_{j,h}$ de $(f\circ g)^*(D)$ qui ne sont pas $f\circ g$-exceptionnelles.

Par d\'efinition, $m_{f/g}(D)=\inf_{j,h}\lbrace m_j.n_{j,h}.m_{\Delta}(F_{j,h})\rbrace$, la somme portant sur les composantes de $(f\circ g)^*(D)$ qui ne sont exceptionnelles ni pour $g$, ni pour $f\circ g$. D'o\`u la premi\`ere assertion, puisque $m_{f/g}(D)\geq m_{fg}(D)$.

Pour \'etablir la seconde assertion, il suffit donc d'observer que $m_{fg}(D)=m_{f/g}(D)$. Or, pour chaque $F$ qui est une composante $g$-exceptionnelle mais non $f\circ g$-exceptionnelle de $(f\circ g)^*(D)$, on peut choisir: $m_{\Delta}(F)\geq m_j.n_{j,h}.m_{\Delta}(F_{j,h})$, pour (au moins) une composante $F_{j,h}=F'$ qui n'est ni $g$-, ni $f\circ g$-exceptionnelle, et telle que, de plus,  $g(F)\subset g(F')=g(F_{j,h})=E_j$, avec: $f(E_j)=D$ (et il existe toujours une telle $F'$).

Pour \'etablir la troisi\`eme assertion, nous allons, pour tout diviseur $F\in W(Z)$ qui est $g$-exceptionnel, mais non $f\circ g$-exceptionnel, montrer l'existence d'un tel diviseur $F'$, si $g$ est un morphisme orbifolde.

Soit donc $D:=f(g(F))\in W(X)$. Soit $E\in W(Y)$ tel que $g^*(E)=t_{E,F}.F+\dots$, avec $t_{E,F}>0$. 

Puisque $g:(Z\vert\Delta)\to (Y\vert\Delta_Y)$ est un morphisme orbifolde, on a:
$$t_{E,F}.m_{\Delta}(F)\geq m_{g}(E):=m_{\Delta_Y}(E).$$

De plus, il existe $F'\in W(Z)$ tel que $g(F')=E$, et $t_{E,F'}. m_{\Delta}(F')=m_{g}(E)$, avec $g^*(E)=t_{E,F'}.F'+\dots$, par d\'efinition de $\Delta_Y=\Delta(g,\Delta)$. 

On a donc aussi: $t_{E,F}.m_{\Delta}(F)\geq t_{E,F'}.m_{\Delta}(F')$.

Donc $D:=f\circ g(F)=f\circ g(F')$. Si $f^*(D)=s.E+\dots$, on a donc: $(f\circ g)^*(D)=s.t_{E,F'}.F'+t^+.F+\dots$, o\`u $t^+\geq s.t_{E,F}$, puisque $f^*(D)$ peut avoir un support dont plusieurs composantes distinctes contiennent $g(F)$. Il en r\'esulte que $t^+.m_{\Delta}(F)\geq s.t_{E,F'}.m_{\Delta}(F')\geq s.m_{g}(E)\geq m_{f/g}(D)$. Puisque ceci est vrai pour tout $F\in W(Z)$ qui est $g$-exceptionnel, mais non $f\circ g$-exceptionnel, on a bien: $m_{f\circ g}(D)\geq m_{f/g}(D),\forall D\in W(X)$. En effet, $m_{f/g}(D)$ est le minimum des quantit\'es $s.t_{E,F'}.m_{\Delta}(F')$, pour $F'$ tel que $g(F')=E'$, avec $f(E')=D$ $\square$

\

De \ref{comporb}, on d\'eduit maintenant:

\begin{proposition}\label{netdelta} Soit $g:Z\to Y$ et $f:Y\to X$ deux fibrations. Soit $\Delta$ un diviseur orbifolde sur $Z$. Si $g:(Z\vert\Delta)\to Y$ est nette et haute (ce que l'on peut r\' ealiser par un morphisme bim\'eromorphe \'el\'ementaire), alors $\Delta(f\circ g,\Delta)=\Delta(f,\Delta(g,\Delta))$.
\end{proposition}

Ce r\'esultat montre que, quitte \`a effectuer sur $(Z\vert\Delta)$ une transformation bim\'eromorphe \'el\'ementaire, on peut calculer la base orbifolde de $f\circ g$ en deux \'etapes.

De \ref{comporb} et \ref{stneth}, on d\'eduit plus g\'en\'eralement, par r\'ecurrence sur $r$:

\begin{proposition}\label{iternet} Soit $g'_j:Z'_j\dasharrow Z'_{j+1}$, pour $j=0,1,\dots,r$ des fibrations m\'eromorphes dominantes connexes, avec $Z'_0$ lisse, compacte et connexe. Soit $\Delta'_0$ un diviseur orbifolde sur $Z'_0$. Il existe alors un morphisme orbifolde bim\'eromorphe \'el\'ementaire $w:(Z_0\vert\Delta_0)\to (Z'_0\vert\Delta'_0)$, des applications bim\'eromorphes $w_j:Z_j\to Z'_j$, et des applications holomorphes $g_j:Z_j\to Z_{j+1}$ telles que, pour $j=0,1,\dots, r$:

1. $g'_j\circ w_j=w_{j+1}\circ g_j$. 

2. $h_j:=g_j\circ g_{j-1}\circ\dots\circ g_1\circ g_0:Z_0\to Z_{j+1}$ est strictement nette et haute.

3. $\Delta(h_j,\Delta_0)=\Delta(g_j,\Delta(h_{j-1},\Delta_0)).$

\end{proposition}

\newpage

\section{DIMENSION CANONIQUE D'UNE FIBRATION.}\label{dimcanfib}

La dimension canonique (``de Kodaira") de la base orbifolde d'une fibration $f:(X\vert \Delta)\to Y$ n'est pas un invariant bim\'eromorphe de $(f\vert \Delta)$. Nous d\'efinissons un invariant bim\'eromorphe de cette fibration (le minimum sur tous les mod\`eles bim\'eromorphes de $(f\vert \Delta)$), qui coincide avec le pr\'ec\'edent lorsque $f$ est ``nette", et montrons comment calculer cet invariant, \`a l'aide d'un faisceau diff\'erentiel de rang $1$, sur tout mod\`ele bim\'eromorphe de $(f\vert \Delta)$.

\subsection{Dimension canonique d'une fibration.}

\begin{definition}\label{kodf} Soit $f:(Y\vert\Delta)\to X$ une fibration avec $(Y\vert \Delta)$ et $X$ lisses. On d\'efinit: $L_f:=f^*(K_X)\subset \Omega^p_Y$, si $p:=dim(X)$.

Rappelons que l'on a aussi d\'efini dans le lemme \ref{kodL} les invariants  suivants:

$\kappa(Y\vert\Delta,L_f):=\kappa(f,\Delta)$, et $p_N(Y\vert\Delta,L_f):=p_N(f,\Delta),\forall N>0$.
\end{definition}

\begin{re} Soit $v:(Y'\vert\Delta')\to (Y\vert\Delta)$ un morphisme orbifolde bim\'eromorphe entre orbifoldes lisses tel que $v_*(\Delta')=\Delta$, et $f':=f\circ v$. Alors il r\'esulte du lemme \ref{invbir} que: $p_N(f',\Delta')=p_N(f,\Delta),\forall N>0$, et $\kappa(f,\Delta)=\kappa(f',\Delta')$.
\end{re}

L'objectif de cette section est de d\'emontrer le:

\begin{theorem}\label{kappa=L}Soit $f:(Y\vert\Delta)\to X$ une fibration, avec $(Y\vert\Delta)$ et $X$ lisses. Alors: 

1. $\kappa(Y\vert\Delta,L_f)\leq\kappa(X\vert\Delta(f,\Delta))$. 

2. Si $f$ est nette, la derni\`ere in\'egalit\'e est une \'egalit\'e.
\end{theorem}

{\bf D\'emonstration:} Elle utilise le lemme \ref{partsup} ci-dessous.

\begin{definition} Soit $F$ un diviseur entier effectif sur $Y$. On dit  que $F$ est {\bf partiellement support\'e} sur les fibres de $f$ s'il existe un diviseur effectif $D$ sur $X$ tel que $F\leq f^*(D)$, et tel que pour toute composante irr\'eductible $D'$ de $D$, il existe une composante irr\'eductible $E'$ de $f^{-1}(D')$non contenue dans $F$, et telle que $f(E')=D'$. \end{definition}

On a:

\begin{lemma}\label{partsup}([Ca04, lemma 1.22]) Soit $F$ un diviseur de $Y$ partiellement support\'e sur les fibres de $f$, et $L\in Pic(X)$. 

Alors l'injection naturelle: $H^0(Y,f^*(L))\subset H^0(Y,f^*(L)+F)$ est surjective.\end{lemma}

Le th\'eor\`eme \ref{kappa=L} r\'esulte alors de la proposition \ref{k=L} suivante:

\begin{proposition}\label{k=L} Soit $N_0=N_0(f,\Delta)$ le plus petit commun multiple des num\'erateurs des multiplicit\'es finies de $\Delta(f,\Delta):=\Delta_f$. Soit $N:=N_0.M$, et $\overline{L_N}:=\overline{L_{N_0.M}}\subset S_{N,p}(Y\vert\Delta), \forall M>0,p:=dim(X)$. Alors:

1. Il existe des diviseurs effectifs (d\'ependants de $M)$ $F,E^+,E^-$ sur $Y$ tels que: $\overline{L_{N}}=f^*(N(K_X+\Delta_f)+F+E^+-E^-$, et tels que: $F$ soit partiellement support\'e sur les fibres de $f$, $E^+$ et $E^-$ sont $f$-exceptionnels et sans composante commune.

2. Si $f$ est nette, il existe un isomorphisme naturel: $$j_M:H^0(Y,\overline{L_{N}})\to f^*(H^0(X,(N(K_X+\Delta_f)))).$$
\end{proposition}

{\bf D\'emonstration (de \ref{k=L}):}

{\bf  Assertion 1.}  Il suffit de voir que, en codimension $1$ au-dessus de $X$, on a: $\overline{L_N}=f^*(N.(K_X+\Delta_f)+F$, et m\^eme seulement que, si $D\in W(X)$, et $E\in W(Y)$ sont tels que $f(E)=D$, alors au-dessus du point g\'en\'erique de $D$, on a \'egalit\'e entre ces deux faisceaux (pour un choix convenable de $F$ d\'etermin\'e ci-dessous). 

On note $\Delta_f:=\Delta(f,\Delta)$, $m:W(Y)\to \overline{\Bbb Q^+}$ la multiplicit\'e d\'efinissant $\Delta$, et $m_f$ celle d\'efinissant $\Delta_f$.

Si $f^*(D)=t.E+\dots$, avec $t>0$ entier, alors $t.s:=t.m(E)\geq m_f(D):=r$, par hypoth\`ese. Dans des coordonn\'ees locales $y=(y_1,\dots,y_n)$ adapt\'ees au voisinage du point $b\in E$, g\'en\'erique dans $E$, $E$ est d\'efini par l'\'equation $y_1=0$, et $f(y)=(x_1,\dots,x_p)=(y_1^t,f_2,\dots,f_p)$, avec $w'(b):=df_2(b)\wedge \dots\wedge df_p(b)\neq 0$. De plus, $D$ est d\'efini par l'\'equation $x_1=0$. Enfin, $\Delta$ (resp. $\Delta_f)$ a en $b$ (resp. en $a:=f(b))$ pour \'equation: $y_1^{1-1/s}=0$ (resp.  $x_1^{1-1/r}=0$ (avec $s,r$ d\'efinis ci-dessus, et $s.t\geq r)$.

Un g\'en\'erateur de $N.(K_X+\Delta_f)$ en $a$ est donc: $w:=(dx_1\wedge\dots\wedge dx_p)/x_1^{N-N/r}$. On en d\'eduit que $f^*(w)=y_1^{((t/r)-s)N-N}.(dy_1\wedge w')^{\otimes N}$, o\`u $w'$ est une $(p-1)$-forme holomorphe non nulle et sans z\'ero sur $Y$ pr\`es de $b$. On a, par hypoth\`ese, $t/r\geq 1/s$, avec \'egalit\'e pour au moins l'une de composantes $E$ de $f^*(D)$. On en d\'eduit bien la premi\`ere assertion.

{\bf Assertion 2.} On suppose que l'on a un diagramme commutatif:

\centerline{
\xymatrix{ (Y'\vert\Delta' )\ar[r]^{v}\ar[d]_{f'} & (Y\vert\Delta)\ar[d]^{f}\\
X'\ar[r]^u&X\\
}}

dans lequel les fl\`eches horizontales sont bim\'eromorphes, les orbifoldes lisses, et que les diviseurs $f'$-exceptionnels de $Y'$ sont $v$-exceptionnels. Donc $f'$ est nette (relativement \`a $f)$, et on veut montrer l'assertion 2 pour $f'$ (et non pas pour $f)$. On note $m',\Delta'_f=\Delta'_{f'}, F',\dots$ les donn\'ees relatives \`a $f'$ qui sont analogues \`a celles introduites pour $f$. 

Soit $A\subset Y$ le lieu au-dessus duquel $v$ n'est pas un isomorphisme: il est analytique de codimension au moins deux dans  $Y$. Observons tout d'abord que l'on a une injection naturelle de faisceaux $v_*(f'^*(N(K_{X'}+\Delta'_f))+F' )\to (f'^*(N(K_X+\Delta_f))+F)$ qui est un isomorphisme au-dessus de $(Y-A)$. Pour tout diviseur $v$-exceptionnel $E'$ (pas n\'ecessairement effectif) de $Y'$, on a donc une bijection naturelle: $H^0(Y', N(K_{X'}+\Delta'_f)+F' +E')\cong H^0(Y', N(K_{X'}+\Delta'_f)+F' )$, puisque les sections ainsi d\'efinies peuvent \^etre vues comme des sections de $f^*(N(K_X\vert\Delta_f))$ sur $(Y-A)$, la conclusion r\'esultant du th\'eor\`eme d'Hartogs. 

On a donc, en particulier, un isomorphisme: $H^0(Y', f'^*(N(K_{X'}+\Delta'_f)+F' )\cong H^0(Y', f'^*(N(K_{X'}+\Delta'_f))+F' +E'^+-E'^-)\cong H^0(Y',\overline{L_{f',N}})$.

Par le lemme \ref{partsup}, on a donc: $H^0(Y', f'^*(N(K_{X'}+\Delta'_f))+F' )\cong H^0(Y', f'^*(N(K_{X'}+\Delta'_f)))\cong H^0(Y',\overline{L_{f',N}})$, et la conclusion $\square$

 \subsection{\'Equivalence bim\'eromorphe de fibrations.}\label{eqbifi}

Soit $f:(Y\vert\Delta)\to X$, et $f:(Y'\vert\Delta')\to X'$ des  fibrations, avec $X,Y,X',Y'$ lisses compactes et connexes. 

$\square$ On supposera dans tout ce \S\ref{eqbifi} que les orbifoldes g\'eom\'etriques $(Y\vert\Delta)$ et $(Y'\vert\Delta')$ sont lisses.

\

 Rappelons (\ref{bimeq}) la:

\begin{definition} On dit que $f$ et $f'$ sont (bim\'eromorphiquement) {\bf \'el\'ementairement \'equivalentes} s'il existe un diagramme commutatif:
\centerline{
\xymatrix{ (Y'\vert\Delta' )\ar[r]^{v}\ar[d]_{f'} & (Y\vert\Delta)\ar[d]^{f}\\
X'\ar[r]^u&X\\
}}

dans lequel $u,v$ sont bim\'eromorphes, $u,v$ holomorphes, et $v$ un morphisme orbifolde tel que $v_*(\Delta')=\Delta$.

Plus g\'en\'eralement, $f$ et $f'$ sont \'equivalentes (on note alors: $f\sim f')$ si elles le sont pour la relation d'\'equivalence engendr\'ee par de tels diagrammes.
\end{definition}

\begin{re} En g\'en\'eral, $u$ n'est pas un morphisme orbifolde sur les bases orbifoldes de $f$ et $f'$, m\^eme s'il est holomorphe.
\end{re}

\begin{proposition}\label{deckap}Dans le diagramme pr\'ec\'edent, on a:

1. $u_*(\Delta(f',\Delta'))=\Delta(f,\Delta)$.

2. $\kappa(X\vert\Delta(f,\Delta))\geq \kappa(X'\vert\Delta(f',\Delta))$

3. On a \'egalit\'e si $\kappa(X)\geq \kappa(X')\geq 0$
\end{proposition}

{\bf D\'emonstration:} La premi\`ere assertion r\'esulte du lemme \ref{compfib}. La seconde s'ensuit imm\'ediatement. La troisi\`eme est \'etablie dans [Ca04,1.14,p. 514] $\square$

\

$\bullet$ La dimension canonique de la base d'une fibration n'est pas, en g\'en\'eral, un invariant bim\'eromorphe. On va maintenant attacher \`a une fibration $f:(Y\vert\Delta)\to X$ comme ci-dessus un invariant bim\'eromorphe (dans la cat\'egorie des orbifoldes lisses) fondamental d'une fibration, gr\^ace \`a la propri\'et\'e de d\'ecroissance ci-dessus.

\begin{definition} Soit $f:(Y\vert\Delta)\to X$ une fibration avec $(Y\vert \Delta)$ et $X$ lisses. On d\'efinit:
  $\kappa(X, f\vert\Delta):=\inf_{f'\sim f}\lbrace\kappa(X'\vert\Delta(f,\Delta)\rbrace$. C'est la {\bf dimension canonique} d'une fibration orbifolde.

\end{definition}

\begin{corollary}\label{kappa=L'}Soit $f:(Y\vert\Delta)\to X$ une fibration, avec $(Y\vert\Delta)$ et $X$ lisses. Alors: 

1. $\kappa((Y\vert\Delta),L_f)=\kappa((Y'\vert\Delta'),L_{f'})$ si $f\sim f'$.

2. $\kappa(f\vert\Delta)=\kappa((Y\vert\Delta),L_f)\leq\kappa(X\vert\Delta(f,\Delta))$. 

3. Si $f$ est nette, la derni\`ere in\'egalit\'e est une \'egalit\'e.\end{corollary}

{\bf D\'emonstration:} La premi\`ere assertion r\'esulte de l'invariance bim\'eromorphe de la dimension canonique de $(L_f,\Delta)$, par le th\'eor\`eme \ref{invbir'}.

Pour l'assertion 2, on choisit un diagramme:

\centerline{
\xymatrix{ (Y'\vert\Delta' )\ar[r]^{v}\ar[d]_{f'} & (Y\vert\Delta)\ar[d]^{f}\\
X'\ar[r]^u&X\\
}}

 comme ci-dessus, avec $f'$ nette relativement \`a $f$, $f$ choisie telle que $\kappa(f\vert\Delta)=\kappa(X\vert\Delta(f,\Delta))=\kappa(X'\vert\Delta(f',\Delta'))$, la premi\`ere \'egalit\'e par le choix de $f$, la seconde par la d\'ecroissance \ref{deckap}.2 et minimalit\'e de $\kappa(X\vert\Delta(f,\Delta))$. Puisque $f'$ est nette, on a: $\kappa(X'\vert\Delta(f',\Delta'))=\kappa((Y'\vert\Delta'),L_{f'})$, par le th\'eor\`eme \ref{kappa=L}.  Puisque $\kappa((Y'\vert\Delta'),L_{f'})=\kappa((Y\vert\Delta),L_f)$, par invariance birationnelle, on a:  $\kappa(f\vert\Delta)=\kappa(X'\vert\Delta(f',\Delta'))=\kappa((Y'\vert\Delta'),L_{f'})=\kappa((Y\vert\Delta),L_f)$. L'in\'egalit\'e r\'esulte de ce que $\kappa(X'\vert\Delta(f',\Delta'))\leq \kappa(X\vert\Delta(f,\Delta))$, par \ref{deckap}

 L'assertion 3 n'est autre que \ref{kappa=L'} $\square$

 \begin{corollary} \label{kappainv}Soit un diagramme commutatif de morphismes orbifoldes:

\centerline{
\xymatrix{ (Y'\vert\Delta' )\ar[r]^{v}\ar[d]_{f'} & (Y\vert\Delta)\ar[d]^{f}\\
X'\ar[r]^u&X\\
}}

dans lequel les fl\`eches horizontales sont g\'en\'eriquement finies, les fl\`eches verticales des fibrations, avec: $(Y'\vert\Delta'),(Y\vert\Delta),X'$ et $X$ lisses. Alors:

1. $\kappa(f\circ v\vert\Delta')\geq \kappa(f\vert\Delta)$

2. Si $v$ est \'etale en codimension $1$, on a \'egalit\'e.
 
 \end{corollary}
 
 {\bf D\'emonstration:} Il suffit, par \ref{kappa=L} pour l'assertion 1, et \ref{etale} pour l'assertion 2, de montrer (par unicit\'e de la saturation) que $v^*(L_f)=L_{f'}\subset \Omega^p_{Y'}, p:=dim(X)$ au point g\'en\'erique de $Y'$.  Or ceci est \'evident, puisqu'en un tel point:  $v^*(L_f)=(f\circ v)^*(K_X)=f'^*(u^*(K_X))=f'^*(K_{X'})$ $\square$
 
 \begin{re} Les corollaires \ref{k=L} et \ref{kappainv} (et leurs d\'emonstrations) subsistent lorsque $X',X$ sont seulement normaux, et $f,f'$ m\'eromorphes, pourvu que l'on d\'efinisse $L_f,L_{f'}$ comme les images dans $\Omega^p_Y,\Omega^p_{Y'}$ de fibr\'es canoniques de mod\`eles lisses arbitraires. \end{re}

  De \ref{kappa=L} on d\'eduit que:

    \begin{corollary}\label{dimcan} Si $p:=dim(X)$, si $f\sim f'$sont comme ci-dessus\footnote{Voir la d\'efinition \ref{bimeq}.}, alors $\kappa(f\vert\Delta)=\kappa(f'\vert\Delta'):=\kappa((Y\vert\Delta), L_f)$ sont bien d\'efinis, o\`u $L_f\subset \Omega_Y^p$ est coh\'erent de rang $1$, et coincide avec $f^*(K_X)$ aux points g\'en\'eriques de $Y$ (en lesquels $f$ est r\'eguli\`ere).
    
    En particulier, si $f:(X\vert\Delta)\dasharrow Y$ est une fibration m\'eromorphe dominante, $\kappa(f\vert\Delta)$ est bien d\'efinie, sur tout mod\`ele holomorphe net de $f$, et est ind\'ependante de ce mod\`ele. On l'appelle la {\bf dimension canonique de $f$}. 
    
    On notera $[Y\vert\Delta(f\vert\Delta)]$ la base orbifolde d'une tel mod\`ele bim\'eromorphe strictement net et haut de $f$. On l'appellera une {\bf base orbifolde stable} de $(f\vert\Delta)$. 
 \end{corollary}

La dimension canonique de $f$ est bien d\'efinie, m\^eme pour les fibrations non presque-holomorphes, bien que la dimension canonique de la fibre orbifolde g\'en\'erique ne soit pas un invariant bim\'eromorphe. La propri\'et\'e suivante justifie partiellement ce fait, puisque les composantes orbifoldes non invariantes bim\'eromorphiquement sont horizontales dans le sens ci-dessous.

  \begin{corollary}\label{dvert} Soit $f:(Y\vert\Delta)\to X$ une fibration holomorphe nette, avec $(Y\vert\Delta)$ une orbifolde g\'eom\'etrique lisse. Soit $\Delta^{vert}$ (resp. $\Delta^{hor})$ le diviseur orbifolde d\'eduit de $\Delta$ par suppression de ses composantes $f$-horizontales (resp. $f$-verticales), c'est-\`a-dire celles qui se projettent surjectivement (resp. ne se projettent pas surjectivement) sur $X$ par $f$. 
  
  Alors: $\kappa(f\vert\Delta)=\kappa(f\vert\Delta^{vert})$.
   \end{corollary}

{\bf D\'emonstration:}  Par construction, $\Delta(f,\Delta)=\Delta(f,\Delta^{vert})$. Puisque $f$ est nette, $\kappa(f\vert\Delta)=\kappa(X\vert\Delta(f,\Delta))=\kappa(X\vert\Delta(f,\Delta^{vert}))=\kappa(f\vert\Delta^{vert}))$$\square$

\subsection{Fibrations de type g\'en\'eral et orbifoldes sp\'eciales: d\'efinition.}

On d\'efinit maintenant ces deux notions, utilis\'ees constamment dans la suite.

   \begin{definition}\label{deftg} Soit $f:(Y\vert\Delta)\dasharrow X$ une fibration m\'eromorphe avec $X,Y$ compacts et ir\'eductibles, et $(Y\vert\Delta)$ une orbifolde g\'eom\'etrique lisse. On dit que $(f\vert\Delta)$ (ou $f$ s'il n' y a pas d'ambiguit\'e sur $\Delta)$est une {\bf fibration de type g\'en\'eral} si $\kappa(f\vert\Delta)=dim(X)>0$.

   Cet ensemble ne d\'epend donc que de la classe d'\'equivalence bim\'eromorphe de $(Y\vert\Delta)$.
 \end{definition}

 \begin{definition}\label{defspec} Une orbifolde g\'eom\'etrique lisse $(Y\vert\Delta)$, avec $Y$ compacte et connexe est dite {\bf sp\'eciale} si:
 
 1. $Y\in \sC$ (ie: $Y$ est bim\'eromorphe \`a une vari\'et\'e K\" ahl\'erienne compacte). 
 
 2. Il n'existe pas de fibration $f:(Y\vert\Delta)\dasharrow X$ de type g\'en\'eral. (De mani\`ere \'equivalente: $(Y\vert\Delta)$ n'a pas d'application m\'eromorphe dominante ``stable" sur une orbifolde g\'eom\'etrique de type g\'en\'eral de dimension strictement positive. Cette notion se formule donc naturellement dans la cat\'egorie bim\'eromorphe des orbifoldes g\'eom\'etriques).
  \end{definition}

\

$\square$ La notion d'orbifolde g\'eom\'etrique sp\'eciale est antith\'etique de celle d'orbifolde de type g\'en\'eral. Nous verrons ci-dessous en quels sens.

\newpage

\section{COURBES $\Delta$-RATIONNELLES.} \label{dr}

Nous introduisons ici les notions de base de la g\'eom\'etrie des courbes rationnelles dans le contexte orbifolde (lisse), et posons un certain nombre de questions, analogues orbifolde de celles d\'ej\`a connues ou conjectur\'ees dans le cadre non-orbifolde. Nous montrons enfin que lorsque $(X\vert\Delta)=X'/G$ est un {\it quotient global} (voir la d\'efinition \ref{dqg}) ces propri\'et\'es peuvent \^etre d\'eduites des r\'esultats connus lorsque $\Delta=0$. On montre en effet que, dans ce cas tr\`es particulier, les courbes rationnelles usuelles de $X'$ et les courbes $\Delta$-rationnelles de $X$ (voir la d\'efinition \ref{d-rat}) se correspondent par images directe et r\'eciproque. On esp\`ere pouvoir \'etendre ces r\'esultats au cas g\'en\'eral (o\`u $(X\vert\Delta)$ est lisse, finie et enti\`ere) en rempla\c cant le rev\^etement orbifolde $X'$ par le champ alg\'ebrique lisse de Deligne-Mumford $\sX\to X$ associ\'e \`a $(X\vert\Delta)$, et en \'etendant au cas des champs alg\'ebriques de Deligne-Mumford les arguments maintenant classiques de la th\'eorie des courbes rationnelles, gr\^ace au champ alg\'ebrique d'Abramovich-Vistoli (construit dans [AV 98]), puisque les courbes rationnelles de $\sX$ et les courbes $\Delta$-rationnelles de $X$ se correspondent alors comme dans le cas des ``quotients globaux". La plupart des r\'esultats attendus (hormis le ``bend-and-break") semblent d'ailleurs pouvoir \^etre obtenus directement, gr\^ace aux techniques usuelles de d\'eformation des courbes rationnelles. Voir la remarque \ref{defrat}.

\subsection{Notion de $\Delta$-courbe rationnelle.}

\

Nous consid\'ererons ici presque exclusivement la version {\it divisible} de courbe $\Delta$-rationnelle. Les versions $\Delta^Z$ et $\Delta^Q$ seront d\'efinies en \ref{d^z} ci-dessous.

\

Soit $(X\vert\Delta)$ une orbifolde {\it lisse} et {\it enti\`ere}, avec: $\Delta:=\sum_j(1-\frac{1}{m_j}).D_j$. 

\

\begin{definition}\label{d-rat}

\

1. Une {\bf courbe rationnelle orbifolde}\footnote{il s'agit ici de la version ``divisible".} est une orbifolde g\'eom\'etrique $(C\vert\Delta')$ telle que $deg(K_{(C\vert\Delta')})<0$. Cette courbe rationnelle orbifolde est dite {\it enti\`ere} (resp. {\it finie}) si ses multiplicit\'es le sont. On a donc: $C\cong \bP^1$.
\end{definition} 

\

\begin{example}\label{dcr}

\

1. Si $C=\bP_1$, et si $\Delta'=\sum_{j} (1-\frac{1}{m_j}). {p_j},$ avec  $m_j>1, \forall j\in J$, les $p_j$ \'etant des points distincts de $X$, alors $(\bP_1\vert\Delta)$ est rationnelle orbifolde si et seulement si $\sum_j(1-\frac{1}{m_j})< 2$. 

2. Lorsque $\Delta'$ est enti\`ere, c'est le cas si et seulement si la suite ordonn\'ee croissante des $m_j,j\in J$ est l'une des suivantes:

$\vert J\vert\leq 2:$ quelconque, $(+\infty,+\infty)$ exclue.

$\vert J\vert=3:$ $(2,2,m),\forall m<+\infty$$; (2,3,3),(2,3,4);(2,3,5).$

3. Lorsque $\Delta'$ est enti\`ere et finie, on a donc: 

$(\bP^1\vert\Delta')$ est rationnelle $\Longleftrightarrow$  $\exists f:\bP^1\to (\bP^1\vert\Delta')$, morphisme orbifolde divisible surjectif.

4. Si $\Delta'$ est enti\`ere, non n\'ecessairement finie, on a l'\'equivalence:

$(\bP^1\vert\Delta')$ est rationnelle $\Longleftrightarrow$  ou bien: $\exists f:\bP^1\to (\bP^1\vert\Delta')$, morphisme orbifolde divisible surjectif, ou bien $\Delta'$ a un support ayant au plus $2$ points, dont un seul au plus a une multiplicit\'e infinie.

\end{example}

\begin{definition} 
 Soit $(X\vert\Delta)$ une orbifolde g\'eom\'etrique lisse et enti\`ere. Une {\bf $\Delta^{div}$-courbe rationnelle} (ou encore: une courbe $\Delta^{div}$-rationnelle) $R\subset (X\vert\Delta)$ est l'image d'un morphisme orbifolde divisible, non-constant et birationnel sur son image $\nu:(\bP^1\vert\Delta')\to (X\vert\Delta)$, dans lequel $(\bP^1\vert\Delta')$ est une courbe orbifolde rationnelle lisse\footnote{On peut d\'efinir de m\^eme une courbe $\Delta$-elliptique comme l'image d'un morphisme orbifolde $\nu:(C\vert \Delta_C)\to (X\vert \Delta)$, birationnel sur son image, tel que $deg(K_C+\Delta_C)=0$, si $C$ est une courbe lisse projective, et $\Delta_C$ minimal rendant $\nu$ un morphisme orbifolde.}.
\end{definition}

De l'exemple \ref{dcr}, on d\'eduit:

\begin{proposition}\label{ratdiv} Soit $(X\vert\Delta)$, $\Delta:=(\sum_j(1-\frac{1}{m_j}).D_j)$ une orbifolde g\'eom\'etrique lisse et enti\`ere, et $R\subset X$ une courbe irr\'eductible non contenue dans le support de $\Delta$, de normalisation $\nu:\hat R\to X$. On a \'equivalence entre les trois propri\'et\'es suivantes:

1. $R\subset (X\vert\Delta)^{div}$ est une $\Delta^{div}$-courbe rationnelle.

2.  $deg(K_{\hat R}+\Delta')=-2+\sum_{p\in \bP^1}(1-\frac{1}{\mu_p})<0$, $\Delta':=\sum_{p\in \bP^1} (1-\frac{1}{\mu_p}).\{p\}$ \'etant le plus petit diviseur orbifolde sur $\hat R$ qui fait de $\nu:(\hat R\vert\Delta')\to (X\vert\Delta)$ un morphisme orbifolde divisible.

On a donc alors: $\hat R\cong \bP^1$, les $\mu_p$ \'etant d\'efinis comme suit.

Pour tous  $j,p\in \hat R$, on  d\' efinit: $t_{j,p}$ par: $\nu^*(D_j)=\sum_{p\in \hat R} t_{j,p}.\{p\}$, et $\mu_p:=ppcm_{\{j\vert  t_{j,p}>0\}}\{ m'_{j,p}:=\frac{m_j}{pgcd (t_{j,p}, m_j)}\}$,
 $ j$ d\'ecrivant l'ensemble des indices tels que: $ t_{j,p}>0$.

 Si $\Delta$ est, de plus, finie, ces conditions sont aussi \'equivalentes \`a:
 
 3. Il existe un morphisme orbifolde divisible $g:\bP^1\to (X\vert \Delta)$ dont l'image est $R$ ($g$ est, en g\'en\'eral fini, mais non birationnel sur son image).

\end{proposition}

\begin{re}\label{rcm} La notion de courbe rationnelle orbifolde d\'efinie ci-dessus coincide donc avec la notion de restriction minimale de $\Delta$ (au sens divisible) telle que d\'efinie en \ref{rescou} et \ref{rescoumin} dans la cat\'egorie $Georb^{div}$. En particulier, le th\'eor\`eme \ref{rescougen} est applicable aux familles couvrantes de $\Delta$-courbes rationnelles (divisibles).
\end{re}

\begin{re}\label{d^z} On d\'efinit les variantes $\Delta^Z$ et $\Delta^Q$ de $\Delta$-courbe rationnelle en rempla\c cant dans les d\'efinitions pr\'ec\'edentes les morphismes orbifoldes par des morphismes orbifoldes quelconques, et en supposant, de plus, enti\`eres les orbifoldes consid\'er\'ees pour la variante $\Delta^Z$. 

Si les multiplicit\'es sont enti\`eres ou $+\infty$, les morphismes orbifoldes $\nu:(\bP^1\vert\Delta')\to (X\vert\Delta)$ consid\'er\'es peuvent \^etre de deux types:

1. Soit {\it divisibles} (ou encore: {\it classiques}). On notera alors $(X\vert\Delta)=(X\vert\Delta)^{div}$ pour pr\'eciser la cat\'egorie de morphismes consid\'er\'es.

2. Soit {\it non-classiques}. On notera alors $(X\vert\Delta)=(X\vert\Delta)^{Z}$ pour pr\'eciser la cat\'egorie de morphismes consid\'er\'es. 

\

$\square$ Si les multiplicit\'es sont rationnelles, les morphismes sont {\it non-classiques}. On le pr\'ecisera n\'eammoins en notant: $(X\vert\Delta)=(X\vert\Delta)^{Q}$.

$\square$ On notera $Georb^*$, $*=Q,Z,div$, l'une des trois cat\'egories dont les objets sont les orbifoldes g\'eom\'etriques lisses, les morphismes \'etant de l'un des trois types: ${Q,Z}$ ou ${div}$.

\

Il y a donc trois types de courbes orbifoldes, selon la cat\'egorie (d'orbifoldes et de morphismes) consid\'er\'ee:

\

1. Les $\Delta^{div}$-courbes rationnelles {\bf divisibles} si $(X\vert\Delta)^{div}$ est enti\`ere et si $\nu$ est un morphisme orbifolde divisible. Leur ensemble est not\'e $Ratl(X\vert\Delta)^{div}$.

2. Les $\Delta^{Z}$-courbes rationnelles {\bf enti\`eres} si $(X\vert\Delta)^{div}$ est enti\`ere et si $\nu$ est un morphisme orbifolde non-classique. Leur ensemble est not\'e $Ratl(X\vert\Delta)^{Z}$. On a donc: $Ratl(X\vert\Delta)^{div}\subset Ratl(X\vert\Delta)^{Z}$: 

3. Les $\Delta^{Q}$-courbes rationnelles, en g\'en\'eral, si $\nu$ est un morphisme orbifolde non-classique. Leur ensemble est not\'e $Ratl(X\vert\Delta)^{Q}$.

Pour les orbifoldes enti\`eres, on a donc: 

$Ratl(X\vert\Delta)^{div}\subset Ratl(X\vert\Delta)^{Z}\subset Ratl(X\vert\Delta)^{Q}$. Ces inclusions sont strictes en g\'en\'eral. Lorsque $(X\vert\Delta)$ est logarithmique (i.e: si $\Delta=\lceil\Delta\rceil$), ces trois notions coincident cependant.

\end{re}

\begin{re} La proposition \ref{ratdiv} est encore valable pour les $\Delta^*$-courbes rationnelles, avec $*=Z$ ou $Q$, \`a condition d'y supprimer le mot ``divisible", et l'hypoth\`ese d'int\'egralit\'e lorsque $*=Q$.

 1. Si l'on veut que $R$ soit une $\Delta^Z$-courbe rationnelle, et si $(X\vert \Delta)$ est enti\`ere, on doit donc choisir: $\mu_p:=max_{\{j\vert  t_{j,p}>0\}}\{ \lceil \frac{m_j}{t_{j,p}}\rceil\}$.

2. Si l'on veut que $R$ soit une $\Delta^{Q}$-courbe rationnelle, on doit choisir: $\mu_p:=max_{\{j\vert  t_{j,p}>0\}}\{ \frac{m_j}{t_{j,p}}\}$.

\end{re}

\begin{example}\label{Q}

\

1. Si $R$, rationnelle, a tous ses ordres de contact avec chacun des $D_j$ d'ordre au moins $m_j$, alors $R$ est $\Delta^Q$-rationnelle, avec: $\Delta'=0$. Lorsque $\Delta$ est enti\`ere, si $R$, rationnelle, a tous ses ordres de contact avec chacun des $D_j$ divisibles par $m_j$, alors $R$ est $\Delta^{div}$-rationnelle, avec: $\Delta'=0$.

2. Si $(X\vert\Delta)$ est une orbifolde g\'eom\'etrique lisse logarithmique ($\Delta=supp(\Delta))$, les $\Delta$-courbes rationnelles sont les courbes rationnelles $R$ de $X$ dont la normalisation rencontre $\Delta$ en, au plus, $1$ point.

Par exemple, si $X=\bP^2$, et si $\Delta=D$, une droite projective affect\'ee de la multiplicit\'e infinie, alors: les droites $L\neq D$ sont $\Delta$-rationnelles, les coniques irr\'eductibles $\Delta$-rationnelles sont celles qui sont tangentes \`a $D$, et une cubique irr\'eductible singuli\`ere est $\Delta$-rationnelle si et seulement si elle est cuspidale, et tangente \`a $D$ en son unique point singulier.

3. Soit $X$ une vari\'et\'e projective lisse et torique. Et $D$ son diviseur anticanonique torique (\`a croisements normaux). On affecte chacune des composantes de $D$ d'une multiplicit\'e finie, et on note $\Delta$ le diviseur orbifolde r\'esultant. Si $R$ est une courbe rationnelle torique (adh\'erence de l'orbite d'un sous-groupe alg\'ebrique \`a un param\`etre du tore agissant sur $X$), non contenue dans $D$, alors $R$ rencontre $D$ en au plus $2$ points en lesquels elle est unibranche. C'est donc une $\Delta^{div}$-courbe rationnelle .

4. Si $X=\bP^2$, et $\Delta$ la r\'eunion de $3$ droites en position g\'en\'erale, affect\'ees de multiplicit\'es enti\`eres arbitraires $a,b,c$, alors $(\bP^2\vert\Delta)$ est Fano. En effet, son fibr\'e canonique est de degr\'e $\delta:= -3+[(1-\frac{1}{a})+(1-\frac{1}{b})+(1-\frac{1}{c})]=-[\frac{1}{a}+\frac{1}{b}+\frac{1}{c}]<0.$ On dispose de trois familles de droites de droites $\Delta^{div}$-rationnelles couvrant $X$, celles passant par l'un des trois points d'intersection de deux des trois droites du support de $\Delta$. Nous verrons en \ref{eqg} et \ref{eqg'} que pour tout sous-ensemble fini $E$ de $\bP^2$, il existe une courbe $\Delta^{div}$-rationnelle irr\'eductible contenant $E$. Cet exemple se g\'en\'eralise imm\'diatemment \`a $\Bbb P^n$, avec $n+1$ droites munies de multiplicit\'es finies.

5. Si $X=\bP^2$, et $\Delta$ la r\'eunion de $4$ droites en position g\'en\'erale, affect\'ees de multiplicit\'es $2,2,a,b$, pour des entiers $4\leq a\leq b$. Alors $(\bP^2\vert\Delta)$ est Fano. En effet, son fibr\'e canonique est de degr\'e $\delta:= -3+[(1-\frac{1}{2})+(1-\frac{1}{2})+(1-\frac{1}{a})+(1-\frac{1}{b})]=-[\frac{1}{a}+\frac{1}{b}]<0$. Dans ce cas, une seule des trois familles de droites pr\'ec\'edentes est form\'ee de courbes $\Delta^{div}$-rationnelles: celles passant par l'intersection des deux droites de multiplicit\'es $a$ et $b$. Si l'on remplace les deux droites de multiplicit\'e $2$ par une conique g\'en\'erale affect\'ee de la multiplicit\'e $2$, on a une seconde famille de droites $\Delta^{div}$-rationnelles: la famille des tangentes \`a $C$.

6. Consid\'erons maintenant $X=\bP^2$, et $\Delta$ la r\'eunion de $4$ droites en position g\'en\'erale, affect\'ees de multiplicit\'es $a,b,c,d$. Alors  son fibr\'e canonique est de degr\'e $\delta:= -3+[(1-\frac{1}{a})+(1-\frac{1}{b})+(1-\frac{1}{c})+(1-\frac{1}{d})]=1-[\frac{1}{a}+\frac{1}{b}+\frac{1}{c}+\frac{1}{d}]<0.$  Donc $(\Bbb P^2\vert \Delta)$ est Fano si et seulement si: $[\frac{1}{a}+\frac{1}{b}+\frac{1}{c}+\frac{1}{d}]>1$. Il n'existe qu'un nombre fini de tels quadruplets. Des exemples sont: $(3,3,5,7)$ et $(2,3,7,41)$ pour lesquels $\delta=-\frac{1}{3.5.7}=- \frac{1}{105}$ et $\delta=-\frac{1}{2.3.7.41}=-\frac{1}{1722}$ respectivement. 

Ce dernier exemple se g\'en\'eralise au cas de $\Bbb P^n$ muni de $n+2$ droites munies de multiplicit\'es $(a_0,a_1,\dots, a_n,a_{n+1}-2)$, avec: $a_0=2, a_1=7,\dots, a_{k+1}=a_0.a_1.\dots .a_k+1$ si $k=1,2,\dots,n+1$. On a alors: $\delta_n=-\frac{1}{a_0.a_1.\dots .a_n.(a_{n+1}-2)}$, probablement maximum parmi les $(X\vert \Delta)$ Fano lisses et enti\`eres de dimension $n$.(On voit facilement que $a_n>(\frac{3}{2})^{2^n}$, et donc que $-\frac{1}{\delta_n}>-(\frac{2}{3})^{2^{n+1}}$.

Si $X=\Bbb P^2$, et si $(a,b,c,d)=(3,3,5,7)$ (resp. $(2,3,7,41))$, il n'existe qu'un nombre fini de droites $\Delta^Z$-rationnelles, celles passant par $2$ des $6$ points d'intersection des $4$ doites du support de $\Delta$, puisqu'une droite rencontrant $Supp(\Delta)$ en au moins $3$ points, y acquiert des multiplicit\'es au moins $3,3,7$ (resp. $(2,3,41))$.  Par contre, $(\bP^2\vert\Delta)$ a dans ces deux exemples, une famille \`a un param\`etre de coniques irr\'eductibles $\Delta^Q$-rationnelles: les coniques tangentes \`a chacune des $4$ droites de $Supp(\Delta)$. Elles acqui\`erent en effet en leurs points de contact des multiplicit\'es rationnelles $a/2,b/2,c/2,d/2$. Puisque $-2+-(1-2/a)+(1-2/b)+(1-2/c)+(1-2/d))=2-2(1/a+1/b+1/c+1/d)<0$, ces coniques sont bien $\Delta^Q$-rationnelles. 

Par contre, ces coniques ne sont pas $\Delta^Z$-rationnelles, car $(\lceil a/2\rceil,\lceil b/2\rceil,\lceil c/2\rceil,\lceil d/2\rceil)=(2,2,3,4)$ (resp. $(1,2,4, 21)$. Il ne semble pas imm\'ediat de fournir pour ces deux exemples une famille couvrante de courbes $\Delta^{div}$-rationnelles. Un d\'ecompte de dimensions indique qu'il semble possible d'avoir des familles \`a nombre positif de param\`etres de $\Delta^{div}$-courbes rationnelles, et m\^eme avec $\Delta'=0$, mais seulement pour de grands degr\'es, divisibles par $105$ et $1722=41.42$ respectivement. 

En effet: les courbes rationnelles planes irr\'eductibles de degr\'e $d=N. 105$ d\'ependent de  $p_N:=3(d+1)-1-3=3d-1$ param\`etres. La condition d'avoir avec une droite des points de contacts d'ordres tous divisibles par $d'$, diviseur de $d$, fournit $\frac{d}{d'}.(d'-1)= d.(1-\frac{1}{d'})$ conditions.

Dans notre cas, ceci fournit donc (avec $d'=3,3,5,7$ successivement), un nombre de conditions total $c_N:=d.[(1-\frac{1}{3})+(1-\frac{1}{3})+(1-\frac{1}{5})+(1-\frac{1}{7})]=N.3.5.7.(3+\frac{1}{105})=N.105.3+3.N=3d+3.N$. On peut donc escompter avoir une famille \`a: $p_N-c_N=3N-1$ param\`etres de $\Delta^{div}$-courbes rationnelles (sans structure orbifolde). Le calcul est enti\`erement analogue dans le cas de multiplicit\'es $(2,3,7,41)$ (ou quelconques telles que $(\Bbb P^2\vert \Delta)$ soit lisse Fano).

Ce calcul montre que le nombre ``attendu" de param\`etres dont d\'epend la courbe $\Delta$-rationnelle $C$ de degr\'e $d$ tel que $\delta.d$ est entier est: $(3-[(1-\frac{1}{3})+(1-\frac{1}{3})+(1-\frac{1}{5})+(1-\frac{1}{7})]).d-1=\delta.d-1=-(K+\Delta).C-1$. 

Remarquons que si les droites consid\'er\'ees sont respectivement d'\'equations homog\`enes $X=0,Y=0,Z=0, uX+vY+wZ=0$ respectivement, et si une telle courbe rationnelle $R(t)$ de degr\'e $105$ existe bien, elle est donn\'ee par un param\'etrage $X(t)=A^3(t), Y(t)=B^3(t), Z(t)=C^5(t)$, les polyn\^omes $A,B,C$ \'etant de degr\'es au plus $35,35,21$ respectivement, et satisfaisant l'\'egalit\'e: $u.A^3+v.B^3+w.C^5=D^7$, $D$ \'etant un polyn\^ome de degr\'e au plus $15$. Si $u,v,w$ sont alg\'ebriques, et les coefficients de $A,B,C,D$ peuvent \^etre alors choisis dans un corps de nombres $k$, et le param\'etrage pr\'ec\'edent fournit une infinit\'e de solutions dans $k$ de l'\'equation: $u.A^3+v.B^3+w.C^5=D^7$. Le probl\`eme de densit\'e potentielle de $(\bP^2\vert \Delta)$ est de savoir si l'on peut trouver une infinit\'e de telles courbes $R$ d\'efinies sur un m\^eme corps $k$. Ces remarques s'appliquent aux multiplicit\'es $(2,3,7,41)$ et \`a l'\'equation: $u.A^2+v.B^3+w.C^7=D^{41}$

\end{example} 

\begin{re}\label{defrat}Des arguments conceptuels de th\'eorie des d\'eformations permettent de g\'en\'eraliser ce calcul aux d\'eformations de courbes $\Delta$-rationnelles sur des orbifoldes g\'eom\'etriques  lisses finies arbitraires.

Soit en effet: $f:(\Bbb P^1\vert \Delta')\to (X\vert \Delta)$ un morphisme orbifolde divisible d\'efinissant une courbe $\Delta$-rationnelle. Par composition avec un morphisme fini $g:\Bbb P^1\to \Bbb P^1$ ramifiant de fa\c con ad\' equate, on obtient un morphisme orbifolde $h:=f\circ g: \Bbb P^1\to (X\vert \Delta)$.

On peut v\'erifier que l'espace tangent en $h$ \`a la d\'eformation de $h$ en tant que morphisme orbifolde (resp. l'espace des obstructions \`a cette d\'eformation) est simplement $H^0(\Bbb P^1, h^*(TX(-log(\Delta))$ (resp. $H^1(\Bbb P^1, h^*(TX(-log(\Delta)))$. Ici $TX(-log(\Delta)$ est le $\Bbb Q$-faisceau des champs de vecteurs holomorphes $(1-\frac{1}{m_j})$-tangents \`a $D_j$, pour chacun des $D_j$. Lorsque $h:\Bbb P^1\to (X\vert \Delta)$ est un morphisme orbifolde divisible, $h^*(TX(-log(\Delta))$ est un faisceau entier sur $\Bbb P^1$.  Et l'espace des d\'eformations \`a $h$ (en tant que morphisme orbifolde) est de dimension au moins $\chi(\Bbb P^1,  h^*(TX(-log(\Delta)))=-(K_X+\Delta).h_*(\Bbb P^1)+ n$.

Cette observation devrait permettre de traiter les questions d'\'evitement des lieux de codimension deux, et le ``glueing lemma" dans le cadre orbifolde, \'evoqu\'es ci-dessous.
Voir aussi l'exemple \ref{ex} ci-dessous.
\end{re}

\

\subsection{Unir\'eglage, Connexit\'e Rationnelle orbifolde}\label{u1}

\begin{definition}\label{u-rc} Soit $(X\vert\Delta)$ une orbifolde g\'eom\'etrique lisse et enti\`ere, avec $X\in \sC$ connexe. On dit que $(X\vert\Delta)$ est:

1. {\bf unir\'egl\'ee} (abr\'eg\'e en: $UR$) si le point g\'en\'erique de $X$ est contenu dans une $\Delta^{div}$-courbe rationnelle.

2. {\bf rationnellement connexe par chaines} (en abr\'eg\'e: $RCC$) si tout couple de points g\'en\'eriques de $X$ est contenu dans une {\it chaine} (ie: une r\'eunion finie {\it connexe}) de $\Delta^{div}$-courbes rationnelles.

2'. {\bf faiblement rationnellement connexe par chaines} (en abr\'eg\'e: $fRCC$) si tout couple de points de $X$ est contenu dans une {\it chaine} de $\Delta^{div}$-courbes rationnelles ou de leurs d\'eg\'en\'erescences dans $Chow(X)$.

3. {\bf rationnellement connexe} (en abr\'eg\'e: $RC$) si tout couple de points g\'en\'eriques de $X$ est contenu dans une $\Delta^{div}$-courbe rationnelle.

5.  {\bf absolument rationnellement connexe} (en abr\'eg\'e: $ARC$) si tout ensemble fini de points g\'en\'eriques de $X$ est contenu dans une $\Delta^{div}$-courbe rationnelle.

On dira aussi que $X$ est $\Delta$-UR (resp. $\Delta$-RCC,RC,ARC), au lieu de: $(X\vert \Delta)$ est UR (resp. RCC,RC,ARC).

\end{definition}
 
\begin{re} Les deux autres notions de $\Delta^*$-courbes rationnelles, avec $*=Z$ ou $Q$ (divisible, enti\`ere ou non) peuvent \^etre aussi employ\'ees ci-dessus, et conduisent d'ailleurs peut-\^etre \`a des notions \'equivalentes dans nombre de cas importants (mais pas toujours: si $X=\Bbb P_1$, par exemple, les notions $UR^Q$ et $UR^Z$ diff\`erent si $\Delta$ est \`a multiplicit\'es rationnelles et si son support a au moins $4$ points). 

On notera en effet, $UR^*$, $RCC^*$, etc..., avec $*=div, Z,Q$ les propri\'et\'es correspondantes, lorsque le contexte n\'ecessite de pr\'eciser la notion consid\'er\'ee. 

On pourra aussi dire, par exemple, de mani\`ere \'equivalente,  que $(X\vert \Delta)^{div}$ est $RC$, ou que $(X\vert \Delta)$ est $RC^{div}$, ou encore: que $X$ est $\Delta^{div}$$-RC$.
\end{re}

\begin{example} 

\

1. Soit $(X\vert\Delta)$ l'exemple torique de \ref{Q}.(3) ci-dessus. Alors $(X\vert \Delta)$ est rationnellement connexe, puisque toutes ses courbes toriques rationnelles non-contenues dans $D$ sont $\Delta$-rationnelles.

2. L'orbifolde $(\bP_2\vert \Delta)$ de l'exemple \ref{Q}.(6) ci-dessus est $RC^Q$ (mais avec des coniques $\Delta$-rationnelles non-enti\`eres). L'argument de comptage de cet l'exemple \ref{Q}.(6) semble indiquer qu'elle devrait \^etre $ARC^{div}$.

\end{example}

\begin{re}

\

 On a les implications \'evidentes: 
 
 $ARC\Longrightarrow RC\Longrightarrow RCC\Longrightarrow fRCC\Longrightarrow UR$.

 Lorsque $\Delta=0$, et si $X$ est projective, on a les implications r\'eciproques:
$RCC\Longrightarrow RC\Longrightarrow ARC$ ([KoMiMo92]).

\end{re}

\begin{question}\label{qqrat}

\

1. Soit $(X\vert \Delta)$ une orbifolde g\'eom\'etrique lisse, enti\`ere et finie, avec $X\in \sC$. Les propri\'et\'es $fRCC^{div}$ et $ARC^{div}$ sont-elles \'equivalentes pour $(X\vert \Delta)$? La condition de finitude ne peut \^etre supprim\'ee, en vertu de l'exemple \ref{ninvblog} ci-dessous.

2. Les propri\'et\'es pr\'ec\'edentes (d\'efinies en \ref{u-rc}) sont-elles stables par d\'eformation (au sens de \ref{defdef})?

3. Nous donnerons ci-dessous au \S\ref{QG} une r\'eponse affirmative \`a la question 1 pr\'ec\'edente dans le cas o\`u $(X\vert \Delta)$ est un ``quotient global", et indiquerons une approche possible, sugg\'er\'ee par cet exemple, \`a la solution de cette question en utilisant la th\'eorie des champs alg\'ebriques de DM lisses.

\end{question}

\begin{question}\label{Qdiv} Supposons $(X\vert \Delta)$ lisse et enti\`ere, avec $X\in \sC$. Les diff\'erentes notions $UR, RC, ARC$ sont-elles \'equivalentes dans $Georb^Q$ et dans $Georb^{div}$?

Explicitement: si $(X\vert \Delta)$ est $UR^{Q}$, est-t'elle $UR^{div}$? La r\'eciroque est \'evidente. M\^eme question pour $RC$ et $ARC$.
\end{question}

\begin{example}\label{ninvblog}

\

1. La propri\'et\'e $RCC$ n'est pas pr\'eserv\'ee par \'equivalence bim\'eromorphe (avec les d\'efinitions ci-dessus) pour les orbifoldes g\'eom\'etriques logarithmiques lisses. Par exemple: $(\bP^2\vert D)$ est (Fano et) $RCC$, si $D$ est la r\'eunion de $2$ droites distinctes concourantes en $1$ point $a$. En effet, les droites passant par $a$ sont $D$-rationnelles et sont, en fait les seules $D$-courbes rationnelles. En effet, par \'eclatement de $a$, on obtient une orbifolde g\'eom\'etrique qui n'est pas $RCC$: elle admet une fibration orbifolde \'evidente sur l' orbifolde g\'eom\'etrique $(\bP^1\vert D')$, o\`u $D'$ est la r\'eunion (r\'eduite) de $2$ points. Donc $\kappa(\bP^1\vert D')=0$, et toutes les courbes $\Delta$-rationnelles de cet \'eclat\'e sont les fibres de cette fibration. Les droites passant par $a$ sont donc bien aussi les seules courbes $D$-rationnelles de $\bP^2$ . 

En particulier, l'orbifolde logarithmique $(\bP_2\vert D)$ est Fano, sans \^etre RE.

2. Le ph\'enom\`ene pr\'ec\'edent provient de ce que les singularit\'es du couple $(\bP^2,D)$ sont log-canoniques $(lc)$, mais pas $klt$. En effet, en consid\'erant l'orbifolde $(\bP^2\vert C)$, o\`u $C$ est maintenant une conique lisse de multiplicit\'e infinie, on voit que cette orbifolde admet pour tout degr\'e $d$ des familles \`a $d$ param\`etres de courbes $C$-rationnelles $R_P$ de degr\'e $d$, donn\'ees explicitement en coordonn\'ees affines $(x,y)$ par le param\'etrage: $x(t):=\frac{t.(P(t)+t^{d-1})}{(P(t)+2t^{d-1})}, y(t):=\frac{t^2.P(t)}{(P(t)+2t^{d-1})}$, si $P(t)$ est un polyn\^ome de degr\'e au plus $(d-2)$ tel que: $P(0)\neq 0$. On suppose ici que $C$ a pour \'equation affine: $y=x^2$. La courbe rationnelle ainsi param\'etr\'ee admet un contact d'ordre $2d$ (donc unique et unibranche) au point $(0,0)$, choisi arbitrairement sur $C$. Cette courbe est donc bien de degr\'e exact $d$. On doit pouvoir v\'erifier que si l'on se donne $d-1$ points g\'en\'eriques de $\bP^2$, un choix ad\'equat des coefficients de $P$ permet de faire passer la courbe $R_P$ par ces $d-1$ points. Lorsque $d=2$, il est \'evident que c'est bien le cas, et $(\bP^2\vert C)$ est $RC$.

\end{example}

La propri\'et\'e suivante est analogue \`a celle du cas $\Delta=0$. On peut la d\'emontrer, \`a l'aide du th\'eorème \ref{rescougen}, en adaptant l'argument de [Ca 91]. Un r\'esultat plus g\'en\'eral sera \'etabli dans \ref{repif}), nous donnons donc une simple indication de preuve.

\begin{proposition}\label{pirc} Soit $(X\vert \Delta)^{div}$ lisse, enti\`ere, finie et $RC^{div}$, avec $X$ projective\footnote{$X\in \sC$ suffit.}. Alors $\pi_1(X\vert \Delta)$ est fini.
\end{proposition}

{\bf D\'emonstration:} On peut adapter l'argument de [Ca 91]. Par le th\'eor\`eme \ref{rescougen}, nous pouvons supposer (quitte \`a modifier $(X\vert \Delta))$ qu'il existe une famille $(C_t)_{t\in T}$ de courbes $\Delta^{div}$-rationnelles  couvrant $X$, et dont le membre g\'en\'erique $C_t$ passe par un point fix\'e $a\notin Supp(\Delta)$, et rencontre $Supp(\Delta)$ transversalement et seulement en des points lisses de $Supp(\Delta)$. On raisonne sur $Z\subset T\times X$ le graphe d'incidence, muni des projections $p:Z\to X,q:Z\to T$. On munit $Z$ de la plus petite structure orbifolde faisant de $p$ un morphisme orbifolde. Les fibres orbifoldes g\'en\'eriques de $q$ sont $\Delta$-rationnelles par \ref{rescougen}, donc de $\pi_1$ fini. On conclut comme dans [Ca 91], \`a l'aide de la suite exacte (voir \ref{sepi1}) des groupes fondamentaux orbifolde associ\'ee \`a $q$, en utilisant la section de $q$ d\'etermin\'ee par le point $a$ $\square$

\subsection{Invariance bim\'eromorphe}

\

Nous allons montrer sur un exemple que, contrairement au cas $\Delta=0$, la notion de $\Delta$-courbe rationnelle n'est pas un invariant bim\'eromorphe en g\'en\'eral.

\begin{example}\label{str} Nous allons montrer que l'inclusion $g_*(Ratl(X_1\vert\Delta_1)^{div})\subset Ratl(X\vert\Delta)^{div}$ peut \^etre stricte, m\^eme si $(X_1\vert\Delta_1)^{div}$ est minimum rendant $g: (X_1\vert\Delta_1)^{div}\to (X\vert\Delta)^{div}$ un morphisme orbifolde. Ceci montre (voir aussi la remarque \ref{restmer}) que la condition de transversalit\'e du th\'eor\`eme \ref{rescouim} n'est pas superflue pour pr\'eserver la restriction des courbes.

Prenons $(X\vert\Delta):=(\Bbb P_2\vert\Delta)$, avec $\Delta:=(1-1/2).(D_1+D_3)+(1-1/3).D_2$, $D_1$ (resp. $D_2$) \'etant la droite d'\'equation affine: $x=0$ (resp. $y=0$), tandis que $D_3$ est une droite g\'en\'erique rencontrant la cubique rationnelle $R$ d'\'equation affine: $y^2=x^3$ en $3$ points distincts $a,b,c$, dont le point \`a l'infini $a=(0,1,0)$ en coordonn\'ees projectives. Alors $R$ est une courbe $\Delta^{div}$-rationnelle de $\Bbb P_2:=X$, puisqu'elle a un contact d'ordre $3$ (resp. $2$) avec $D_2$ (resp. $D_1$), de sorte que le diviseur orbifolde minimum $\Delta_R$ sur $R$ rendant l'inclusion normalis\'ee $(R\vert\Delta_R)\to (\Bbb P_2\vert\Delta)$ un morphisme orbifolde divisible est $\Delta_R:=(1-1/2).(\{a\}+\{b\}+\{c\})$.

Soit alors $g:X_1\to \Bbb P_2=X$ l'\'eclatement du point de coordonn\'ees affines $(0,0)$, et $E$ le diviseur exceptionnel de cet \'eclatement. Soit $R_1$ la transform\'ee stricte de $R$ dans $X_1$. Elle a donc un ordre de contact avec $E$ (resp. la transform\'ee stricte $D'_1$ de $D_1$; resp. $D'_2$ de $D_2$) \'egal \`a $2$ (resp. $0$ (sur $E$); resp. $1$). La structure orbifolde minimum $\Delta_1^{div}$ sur $X_1$ faisant de $g: (X_1\vert\Delta_1)\to (X\vert\Delta)$ un morphisme orbifolde divisible est: 

$\Delta_1^{div}:=(1-1/6). E+(1-1/2). (D'_1+D'_3)+(1-1/3). D'_2$.

Tandis que la structure orbifolde minimum $\Delta_{R_1}^{div}$ sur $R_1$ faisant de l'inclusion $ (R_1\vert\Delta_{R_1})\to (X_1\vert\Delta_1)$ un morphisme orbifolde divisible est: 

$\Delta_{R_1}^{div}:=(1-1/2). (\{d\}+\{a\}+\{b\}+\{c\})$, en identifiant $a,b,c$ avec leurs images r\'eciproques dans $X_1$, tandis que $d:=(E\cap D'_2)$. Il en r\'esulte que $R_1$ n'est pas une $\Delta_1^{div}$-courbe rationnelle (mais elle reste cependant ``elliptique"). 

On peut d'ailleurs v\'erifier que $R_1$ n'est pas non plus une courbe $(X_1\vert\Delta')^Q$-rationnelle si $g:(X_1\vert\Delta')^Q\to (X\vert\Delta)^Q$ est un morphisme orbifolde.
\end{example}

\begin{question}\label{ib} Soit $(X\vert \Delta)$ une orbifolde g\'eom\'etrique lisse, enti\`ere et finie, avec $X\in \sC$. Les propri\'et\'es $UR, RCC, RC, ARC$ sont-elles alors pr\'eserv\'ees par \'equivalence bim\'eromorphe (dans Georb$^*$, avec $*=div, Z$, ou $Q$)?

Nous indiquerons en \ref{comp'} une solution partielle \`a ce probl\`eme lorsque $(X\vert\Delta)$ est un quotient global, et en d\'eduirons une approche possible gr\^ace aux champs alg\'ebriques.
\end{question}
\

Il est clair que si $g:(X'\vert \Delta')\to (X\vert \Delta)$ est un morphisme orbifolde (dans Georb$^*$), avec $g(X')=X$, et si $X'$ poss\'ede l'une des propri\'et\'es $UR, RCC, RC, ARC$, alors $(X\vert \Delta)$ la poss\`ede aussi. Le probl\`eme est donc le rel\`evement des $\Delta^*$-courbes rationnelles de $X$ \`a $(X'\vert \Delta')^{*}$ lorsque $g$ est bim\'eromorphe, avec $g_*(\Delta')=\Delta$. 

\

Nous pouvons, gr\^ace au th\'eor\`eme \ref{rescougen}, donner une r\'eponse partielle affirmative \`a cette question:

\begin{theorem}\label{ibrcp} Soit $(X\vert \Delta)$ une orbifolde lisse, avec $X\in \sC$. Si $u:(X'\vert \Delta')\to (X\vert \Delta')$ est une modification orbifolde \'el\'ementaire minimale, avec $(X'\vert \Delta')$ lisse, et si $(X\vert \Delta)$ poss\`ede l'une des propri\'et\'es $UR,RC,ARC$, alors $(X'\vert\Delta')$ la poss\`ede aussi, si $u$ est transverse \`a une famille ad\'equate de courbes rationnelles de $(X\vert \Delta)$ \footnote{La conclusion ne subsiste pas pour la propri\'et\'e $RCC$, par \ref{ninvblog}.}.
\end{theorem}

{\bf D\'emonstration:} Montrons l'assertion pour la propri\'et\'e $RC$ dans le cadre des morphismes divisibles. Soit $(C_t)_{t\in T}$ une famille couvrante de courbes $\Delta$-rationnelles de $X$. Par le th\'eor\`eme \ref{rescougen}, nous pouvons supposer le transform\'e strict dans $X'$ du membre g\'en\'erique de cette famille est une $\Delta'$-courbe rationnelle rencontrant $Supp(\Delta')$ transversalement et en des points lisses de ce diviseur, puisque l'on a suppos\'e que la structure orbifolde $\Delta'$ est minimale rendant $u$ un morphisme orbifolde, et que $u$ est transverse \`a la famille ``ad\'equate" $T$. 

La d\'emonstration pour les propri\'et\'es $RC$ et $ARC$ est analogue $\square$

\

Le r\'esultat \'evident suivant \ref{evit} montre que l'\'evitement des lieux de codimension $2$ ou plus permet de r\'esoudre affirmativement la question pr\'ec\'edente, lorsque la structure orbifolde $\Delta'$ n'est plus minimale. 

\

Soit $(R_{t})_{t\in T}$ une famille analytique de courbes de $X$, param\'etr\'ee par l'espace analytique irr\'eductible compact $T\subset Chow(X)$. On dira que cette famille est {\bf couvrante} (resp. {\bf bicouvrante}) si son membre g\'en\'erique $R_t$ est irr\'eductible, et si son graphe d'incidence $Z\to T\times X$ est surjectif sur $X$ (resp. si $Z\times_{T}Z\to X\times X$ est surjectif). 

Supposons que le membre g\'en\'erique de cette famille de courbe poss\`ede une propri\'et\'e $\Pi$.Une telle famille est dite {\bf maximale pour 
$\Pi$} si toute famille analytique $(R_s)_{s\in S}$, avec $T\subset S\subset Chow(X)$ dont le membre g\'en\'erique poss\`ede $\Pi$ est telle que $S=T$.

\

\begin{proposition}\label{evit} Soit $(X\vert \Delta)$ lisse, enti\`ere et finie, avec $X\in \sC$. Soit $(R_t)_{t\in T}$ une famille analytique couvrante (resp. bicouvrante) de courbes $\Delta^{div}$-rationnelles de $X$, et maximale pour la propri\'et\'e d'\^etre $\Delta^{div}$-rationnelles. 

1. Si, pour tout sous-ensemble analytique ferm\'e $A\subset X$, de codimension au moins $2$, le membre g\'en\'erique $R_t$ de cette famille ne rencontre pas $A$. Alors, pour tout morphisme orbifolde divisible bim\'eromorphe $g:(X'\vert\Delta')\to (X\vert \Delta)$, avec $(X'\vert \Delta')$ lisse, enti\`ere et finie telle que $g_*(\Delta')=\Delta$, $(X'\vert \Delta')$ est aussi $\Delta^{div}$-$UR$ (resp. $\Delta^{div}$-$RC$).

2. Si la propri\'et\'e pr\'ec\'edente d'\'evitement des lieux de codimension au moins $2$ est satisfaite pour toute telle $(X\vert \Delta)$, les propri\'et\'es $UR^{div}$ et $RC^{div}$ sont alors pr\'eserv\'ees par \'equivalence bim\'eromorphe dans cette classe d'orbifoldes g\'eom\'etriques.
\end{proposition}

{\bf D\'emonstration:} La transform\'ee stricte du membre g\'en\'erique $R_t$ ne rencontre pas le diviseur exceptionnel de $g$, et est donc une $\Delta'$-courbe rationnelle de $X'$. La seconde assertion en r\'esulte imm\'ediatement $\square$

\begin{re} 

1. Lorsque $\Delta=0$, la propri\'et\'e d'\'evitement de $A$ pr\'ec\'edente est satisfaite par toute famille couvrante de courbes rationnelles (puisque le fibr\'e normal du membre g\'en\'erique est semi-positif).

2. La proposition pr\'ec\'edente et sa d\'emonstration, restent valables pour les courbes $\Delta^*$-rationnelles, avec $*=Z,Q$.

3. Nous montrerons dans \ref{comp'} la propri\'et\'e d'\'evitement de \ref{evit} pour les ``quotients globaux", et en d\'eduirons qu'elle est valable pour les $(X\vert \Delta)$ lisses enti\`eres et finies dans $\sC$ si elle l'est pour les champs de DM lisses. 

4. Cette propri\'et\'e d'\'evitement n'est cependant pas satisfaite par les $(X\vert \Delta)$ lisses, enti\`eres et finies g\'en\'erales: si $(X\vert \Delta)=(\bP^2\vert \Delta)$, o\`u $\Delta$ est support\'ee par $3$ droites en position g\'en\'erale et affect\'ees de multiplicit\'es $(a,b,c)$ dont la somme des inverses est au plus $1$ (telles que, par exemple $(2,3,m)$ avec $m\geq 6)$, alors les droites passant par l'un des trois points d'intersection des $3$ droites forment trois familles couvrantes maximales de courbes $\Delta$-rationnelles n'\'evitant pas l'un des trois points. 

5. Rempla\c cant la famille donn\'ee par un ``multiple", obtenu en d\'eformant la compos\'ee $h:=g\circ f:\bP^1\to (X\vert \Delta)$, o\`u $g:\bP^1\to \bP^1$ est finie, de grand degr\'e et telle que $h$ soit un morphisme orbifolde, la propri\'et\'e d'\'evitement doit \^etre satisfaite par $h$, m\^eme si elle ne l'est pas pour l'application ``initiale" $f$. Ceci doit pouvoir \^etre \'etabli par l'observation de la remarque \ref{defrat}.

5. L'invariance bim\'eromorphe des propri\'et\'es $UR^*, RC^*,\dots$ peut cependant \^etre d\'eduite de la propri\'et\'e d'\'evitement $(C)$ plus faible \'enonc\'ee dans le th\'eor\`eme \ref{rerc}. 
\end{re}

\

\subsection{Unir\'eglage et Dimension Canonique: Conjectures}

\

\begin{definition}\label{k-rc} Si $(X\vert\Delta)$ est une orbifolde g\'eom\'etrique lisse, on pose: 

$\kappa_+(X\vert\Delta):=max_f\{\kappa(f\vert\Delta)\}$, $f:(X\vert\Delta)\dasharrow Y$ parcourant l'ensemble des applications m\'eromorphes surjectives (ie: dominantes), avec $dim(Y)>0$. 

C'est un invariant bim\'eromorphe, avec: $\kappa_+\geq \kappa$.

Donc $\kappa_+(X\vert\Delta)=-\infty$ signife: $\kappa(f\vert\Delta)=-\infty,\forall f:(X\vert\Delta)\dasharrow Y$.

\end{definition}

\begin{example} Si $X$ est rationnellement connexe, alors $\kappa_+(X):=\kappa_+(X\vert 0)=-\infty$. 

D\'efinissant: $\kappa_{++}(X\vert \Delta):=max_{L\subset \Omega_X^p,p>0, rg(L)=1}\{\kappa(L\vert \Delta)\},$ on a m\^eme: $\kappa_{++}(X):=\kappa_{++}(X\vert 0)=-\infty$ si $X$ est RC.\end{example}

La d\'emonstration de l'\'enonc\'e suivant est facile, gr\^ace au th\'eor\`eme \ref{rescougen}:

\begin{proposition}\label{k+} Soit $(X\vert\Delta)^Q$ une orbifolde g\'eom\'etrique lisse, avec $X\in \sC$ connexe.

1. Si $(X\vert\Delta)^Q$ est $UR$, alors $\kappa(X\vert\Delta)=-\infty$. 

2. Si $(X\vert\Delta)^Q$ est $RC$, alors $\kappa_+(X\vert\Delta)=-\infty$.

\end{proposition}

{\bf D\'emonstration:} Soit $Z$ un mod\`ele lisse du graphe d'incidence d'une famille alg\'ebrique couvrante (resp. couvrante passant par un point $a\in X-Supp(\Delta)$) de $X$, not\'ee $(C_t)_{t\in T}$ de courbes $\Delta^{Q}$-rationnelles, pour d\'emontrer l'assertion 1 (resp. 2). Les projections naturelles sont not\'ees $p:Z\to X$ et $q:Z\to T$. Pour l'assertion 1, on suppose de plus que $dim Z=dim X$.

Le th\'eor\`eme \ref{rescougen} nous permet, quitte \`a remplacer $(X\vert \Delta)$ par une modification orbifolde \'el\'ementaire minimale transverse \`a la famille $T$, de supposer que le membre g\'en\'erique $C_t$ de cette famille est lisse, et rencontre transversalement $Supp(\Delta)$, et seulement en des points lisses de ce diviseur.

On munit $Z$ de la plus petite structure orbifolde (que l'on peut supposer lisse) pour laquelle $p:(Z\vert \Delta_Z)\to (X\vert \Delta)^Q$ est un morphisme orbifolde. La fibre orbifolde g\'en\'erique $(C'_t\vert \Delta_t)$ de $q$ est alors rationnelle.

Pour l'assertion 1, supposons par l'absurde l'existence d'une section non nulle: $w\in H^0(X, m.(K_X+\Delta)), m>0$. 

Alors $0\neq p^*(w)\in H^0(Z, m.(K_Z+\Delta_Z))$. Contredit le fait que la restriction \`a $C'_t$ de $p^*(w)$, qui est non-nulle, doit s'annuler, puisque $:(C'_t\vert \Delta_t)$ est rationnelle, et que le fibr\'e normal de $C'_t$ dans $X$ est trivial.

Pour l'assertion 2, on note $s:T\to X$ la section telle que: $s(t):=(t,a),\forall t\in T$. 

Soit $0\neq w\in H^0(X, S^m\Omega^r(X\vert \Delta))$, avec $r>0,m>0$. Soit $0\neq p*(w) \in S^N(\Omega(X\vert \Delta))$.

Alors $p^*(w)=q^*(w')$, puisque $:(C'_t\vert \Delta_t)$ est rationnelle de fibr\'e normal trivial, et $w'=0$, puisque $w'=s^*(q^*(w'))=s^*(p^*(w))=0$, puisque $p\circ s:T\to X$ est l'application constante. Donc $w=0$ $\square$

\

La conjecture suivante est la version orbifolde d'une conjecture standard (qui en est le cas o\`u $\Delta=0)$ de la g\'eom\'etrie alg\'ebrique:

\begin{conjecture}\label{-infty} Soit $(X\vert\Delta)^Q$ une orbifolde g\'eom\'etrique lisse, avec $X\in \sC$ connexe. 

Si $\kappa(X\vert\Delta)=-\infty$, alors $(X\vert\Delta)^Q$ est $UR^Q$. Si, de plus, $(X\vert \Delta)$ est \`a multiplicit\'es enti\`eres, elle est $\Delta^{div}$-unir\'egl\'ee.

Si $\kappa_+(X\vert\Delta)=-\infty$, alors $(X\vert\Delta)^Q$ est $RC^Q$. Si, de plus, $(X\vert \Delta)$ est \`a multiplicit\'es enti\`eres, elle est $\Delta^{div}-RC$.
\end{conjecture}

\begin{re}

\

1. La condition $\kappa=-\infty$ \'etant un invariant bim\'eromorphe, la conjecture \ref{-infty} implique l'invariance bim\'eromorphe de l'unir\'eglage, et aussi l'\'equivalence des conditions $UR^Q$ et $UR^{div}$ lorsque $\Delta$ est \`a multiplicit\'es enti\`eres. Ainsi que des conclusions analogues pour la connexit\'e rationnelle et la condition $\kappa_+=-\infty$.

2. Il est facile de montrer (\`a l'aide du ``quotient rationnel" et de [GHS 03]) que si la premi\`ere partie de la conjecture est vraie lorsque $\Delta=0$, la seconde l'est aussi.

3. Le seul cas non-trivial avec $\Delta\neq 0$ dans lequel cette conjecture est connue est celui des surfaces projectives avec diviseur orbifolde logarithmique ($\Delta=Supp(\Delta))$, par [K-M98]. Il s'agit donc dans ce cas de recouvrir $X$ par des courbes rationnelles $R$ rencontrant $\Delta$ en un seul point, en lequel $R$ est unibranche.
\end{re}

\begin{example} Soit $X=\bP^2$, et $\Delta=C$, avec $C$ une conique lisse. Donc $(\bP^2\vert C)$ est Fano. Les courbes rationnelles orbifoldes $R$ sont alors les courbes rationnelles (de degr\'e $d)$, qui coupent $C$ en un unique point en lequel $R$ est unibranche.

Lorsque $d=1$, les tangentes \`a $C$ sont donc de tels exemples. Lorsque $d=2$, $R$ est une conique lisse osculatrice \`a $C$. On a vu en \ref{ninvblog} que l'on peut trouver de telles courbes $C$-rationnelles pour tous les degr\'es $d$. On devrait pouvoir en d\'eduire que $(\bP^2\vert C)$ est $ARC$. \end{example}

Lorsque les multiplicit\'es sont finies, de nombreux nouveaux cas se pr\'esentent, m\^eme lorsque $n=2$, d'orbifoldes g\'eom\'etriques lisses $(\bP^2\vert\Delta)$ qui sont Fano, et pour lesquelles la v\'erification de la conjecture \ref{-infty} n'est pas imm\'ediate.

\begin{example}\label{ex} Soit $X=\bP^2$, $\Delta=2/3(L_3+M_3)+4/5L_5+6/7.L_7$, o\`u  $L_3,M_3,L_5,L_7$ sont quatre droites en position g\'en\'erale. L'orbifolde g\'eom\'etrique lisse $(\bP^2\vert\Delta)$ est Fano, puisque $2/3+1/5+1/7=1+1/105>1$. La conjecture pr\'ec\'edente implique qu'elle devrait \^etre $RC^{div}$. Voir l'exemple \ref{Q} pour une possible v\'erification directe de cette propri\'et\'e. Il serait int\'eressant d'avoir une approche plus conceptuelle (par d\'eformation) de ce probl\`eme. Voir la remarque \ref{defrat}.

\end{example}

\subsection{Quotients globaux: descente et rel\`evement de courbes rationnelles}

Soit $(X\vert\Delta)$  lisse, finie et enti\`ere, avec $X\in \sC$. Nous allons introduire une classe (tr\`es restreinte) de telles orbifoldes g\'eom\'etriques pour lesquelles les probl\`emes concernant les courbes $\Delta^{div}$-rationnelles peuvent \^etre r\'eduits au cas usuel o\`u $\Delta=0$. La consid\'eration des champs alg\'ebriques de DM associ\'es devrait permettre d'\'etendre cette r\'eduction \`a toutes les orbifoldes g\'eom\'etriques $(X\vert\Delta)$  lisse, finie et enti\`ere, avec $X\in \sC$.

\begin{definition}\label{dqg} Soit $(X\vert\Delta)$  lisse, finie et enti\`ere, avec $X\in \sC$, et soit $f:X'\to X$ une application holomorphe propre et finie surjective, avec $X'$ lisse. On dit que $f$ {\bf ramifie au moins} (resp. {\bf ramifie exactement}) au-dessus de $\Delta^{div}$ si $f$ est \'etale au-dessus du compl\'ementaire du support de $\Delta$, et si, pour chaque $j$, et chaque  point $x'$ de $X'$ tel que $f(x')\in D_j$, $f$
ramifie en $x'$ \`a un ordre $m'_j$ multiple de $m_j$ (resp. tel que $m'_j=m_j$).

On dit que $f$ ramifie au moins au-dessus de $\Delta^Z$ si $m'_j\geq m_j$, pour tous $j,x'$ comme ci-dessus.

On dit que $(X\vert\Delta)$ est un {\bf quotient global} s'il existe un $f:X'\to X$ comme ci-dessus, qui ramifie exactement au-dessus de $\Delta^{div}$ (ou $\Delta^Z)$, les notions coincidant alors)\end{definition}

\begin{example}\label{eqg} Soit $X=\bP^2$, $\Delta=\sum_{j=1}^{j=3}(1-\frac{1}{m_j}).D_j$, o\`u les $m_j>1$ sont entiers, et les $D_j$ les droites d'\'equation $T_j=0$, dans les coordonn\'ees homog\`enes $(T_1,T_2,T_3)$. On note $m:=pgcd_j\{m_j\}$. Soit $f:\bP^2\to \bP^2$ le morphisme d\'efini par $f(U_j)=T_j:=U_j^m$, pour $j=1,2,3$. Alors ce morphisme ramifie au moins (resp.exactement) au-dessus de $\Delta^{div}$ pour tous $m_j$ (resp. si $m_j=m, \forall j)$.
\end{example}

\begin{theorem}\label{d-orb'} Soit $(X\vert\Delta)$  lisse, finie et enti\`ere, avec $X\in \sC$. On suppose qu'il existe un morphisme propre et fini $f:X'\to X$ ramifiant au moins au-dessus de $\Delta^{div}$ (resp. $\Delta^Z)$.

1. Soit $R'\subset X'$ une courbe rationnelle, et $R:=f(R')\subset X$. Alors: $R$ est une $\Delta^{div}$-courbe rationnelle (resp. une $\Delta^{Z})$-courbe rationnelle.

2. Si $f$ ramifie exactement au-dessus de $\Delta^{div}$, et si $R\subset X$ est une $\Delta^{div}$-courbe rationnelle, toute composante irr\'eductible $R'\subset X'$ de $f^{-1}(R)$ est une courbe rationnelle de $X'$.
\end{theorem}

Nous donnerons dans les trois sections suivantes des cons\'equences imm\'ediates de ce r\'esultat.

{\bf D\'emonstration:} Assertion 1. Nous ne montrerons que le cas divisible (le cas non-classique est similaire, plus simple). 

Lorsque la restriction $g:R'\to R$ est birationnelle, la conclusion r\'esulte de ce que la compos\'ee de deux morphismes orbifoldes est un morphisme orbifolde (notant que $f:X'\to (X\vert\Delta)$ est un morphisme orbifolde). Dans le cas g\'en\'eral, composant les inclusions de $R,R'$ dans $X,X'$ avec les normalisations, nous pouvons supposer que $R,R'$ sont lisses rationnelles (pour simplifier les notations).

Soit maintenant $b'\in R$ tel que $b:=f(b')\in D_j$. 

Soit $e(b')$ l'ordre de ramification de $g:R'\to R$ en $b'$, et $D'_j:=f^{-1}(D_j)$.

On note aussi $e(b):=pgcd\{e( b'), b'\in f^{-1}(b)\cap R'\}$.

La formule de projection fournit pour les nombres d'intersection {\bf locaux} pr\`es de $b,b'$:

$m'_j.(D'_j.R')_{b'}=(f^*(D_j).R')_{b'}=(D_j.g_*(R))_b=e(b').(D_j.R)_b:=e(b').t_{j,b}$

D'o\`u l'on d\'eduit, par le th\' eor\`eme de Bezout, que, pour tout $j$ tel que $g(b')\in D_j$, $m'_j$, et donc aussi $m_j$, divise $e(b). t_{j,b}$.

Donc: $\frac{m_j}{pgcd(t_{j,b},m_j)}$ divise $e(b)$$, \forall j$. Et $ppcm_j\{\frac{m_j}{pgcd(t_{j,b},m_j)}\}$ divise donc $e(b)$.

Et donc, d\'esignant par $\Delta'$ le plus petit diviseur orbifolde sur $R$ faisant de l'inclusion (normalis\'ee) de $R$ dans $X$ un morphisme $\Delta$-orbifolde {\bf divisible}:
$m_{\Delta'}(b)$ divise $e(b)$. Donc: $K_{R'}\geq f^*(K_R+\Delta')$, ce qui ach\`eve la d\'emonstration de l'assertion 1 du th\'eor\`eme $\square$

\

Remarquons que l'assertion 1 r\'esulte aussi du lemme suivant, qui en donne une version quantitative, puisque $g(R')=g(R)=0$:

\begin{lemma} Soit $f:R'\to R$ un morphisme surjectif de degr\'e $d>0$ entre courbes projectives lisses et connexes. Pour tout $b'\in R'$, on note $e(b')$ l'ordre de ramification de $f$ en $b'$, et pour $b\in R$, par $e(b):=pgcd\{e(b'), b'\in f^{-1}(b)\}$, et enfin par $e'(b):=inf\{e(b'), b'\in f^{-1}(b)\}\leq e(b)$.

Alors: $\frac{2(g(R')-1)}{d}\geq \frac{2(g(R)-1)}{d}+\sum_b(1-\frac{1}{e(b)})\geq \frac{2(g(R)-1)}{d}+\sum_b(1-\frac{1}{e'(b)})$.
\end{lemma}

{\bf D\'emonstration:} Appliquant la formule de Riemann-Hurwitz \`a $f:R'\to R$, on obtient:
$$2g(R')-2=2d.(g(R)-1)+\sum_b(\sum_{\bar b\in f^{-1}(b)}(e(\bar b)-1))=2d.(g(R)-1)+\sum_b(d-\sum_{b'\in f^{-1}(b)}1)$$

Divisant par $d$, nous obtenons: 

$$\frac{2(g(R')-1)}{d}= \frac{2(g(R)-1)}{d}+\sum_b(1-\frac{\sum_{b'}1}{\sum_{b'}e(b')})$$

La conclusion r\'esulte alors de ce que: $\frac{\sum_{b'}1}{\sum_{b'}e(b')}\leq \frac{1}{e(b)}\leq \frac{1}{e'(b)}$ $\square$

\

$\square$ D\'emonstration de l'assertion 2. Consid\'erons le diagramme commutatif:

\centerline{
\xymatrix{ R'\ar[r]^{g'}\ar[d]_{v} & X'\ar[d]^f\\
R\ar[r]^{g}&(X\vert\Delta)\\
}}

dans lequel les fl\`eches horizontales sont les normalisations compos\'ees avec les inclusions, tandis que $v$ est g\'en\'eriquement la restriction de $f$.

Soit $b'\in R', b:=v(b'),$ et $ t_j,t'_j$ respectivement les ordres de contacts en $b$ et $b'$ respectivement de $D_j$ avec $R$, et de $f^{-1}(D_j)$ avec $R'$. Notons enfin $e(b'):=e$ l'ordre de ramification en $b'$ de $v$. La formule de projection (ou la commutativit\'e du diagramme pr\'ec\'edent) montre que: 
$$f_*(R').D_j=e.g^*(D_j)=e.t_j=R'.f^*(D_j)=m_j.t'_j$$

Puisque $X'$ est suppos\'e lisse, $t'_j$ est entier.

Notons $d_j:=pgcd(t_j,m_j), u_j:=t_j/d_j,m'_j:=m_j/d_j$. 

Donc $pgcd(u_j,m'_j)=1$.

On d\'eduit donc de l'\'egalit\'e ci-dessus que, pour tout $j\in J(b)$, ensemble des $j$ tels que $g(b)\in D_j$: 
$$e.u_j=m'_j.t'_j,$$

De sorte que: $u_j$ divise $t'_j$ et $m'_j$ divise $e$ pour tout $j\in J(b)$.

Soit alors: $m':=ppcm\{m'_j, j\in J(b)\}$. 

Alors $m'$ divise $e$ par ce qui pr\'ec\`ede. Posons: $e:=e'.m'$. 

L'assertion du th\'eor\`eme sera \'etablie si l'on montre que $v:R'\to (R\vert\Delta_R)$ est \'etale (ie: si $e=m'$, ce qui \'equivaut ici \`a: $v^*(K_R+\Delta_R)=K_{R'}$, $\Delta_R$ \'etant le diviseur orbifolde sur $R$ attribuant \`a (tout) $b\in R$ la multiplicit\'e $m'$ pr\'ec\'edente, qui est la plus petite faisant de $g:(R\vert\Delta_R)\to (X\vert\Delta)$ un morphisme orbifolde divisible).

Nous allons montrer que $e'=1$, ce qui \'etablira donc le th\'eor\`eme. 

Dans l'\'egalit\'e suivante, les trois termes sont des entiers, et $e'$ divise donc $e":=pgcd(e, t'_j,j\in J(b))$:
$$e'. (\frac{m'}{m'_j})=(\frac{t'_j}{u_j}),\forall j\in J(b).$$

Dans des coordonn\'ees locales adapt\'ees, nous pouvons donc supposer (par lissit\'e et rev\^etement \'etale local de $X'$) que, si $s\in \Bbb D$ est une coordonn\'ee locale sur $R'$, alors:
$$g'(s)=(s^{t'_1}.(1+s.w_1(s)),\dots,s^{t'_n}.(1+s.w_n(s)),$$
 tandis que: $$f\circ g'(s)=p_1(s^e),\dots,p_n(s^e),$$
 
les $w_k,p_k$ et $x_k$ ci-dessous \'etant des fonctions analytiques au voisinage de $0$.
 
 On a donc, sur un voisinage de $0\in \Bbb D$, disque unit\'e de $\Bbb C$:
 $$p_k(s^e)=s^{t'_k.m_k}.(1+s.w_k(s)),\forall k:=1,\dots,n.$$
 
 On en d\'eduit que $e.t_k=t'_k.m_k$ et que $(1+s.w_k(s))=1+s^e.x_k(s^e)$, pour $s$ assez voisin de $0\in \Bbb C$, et $k=1,\dots,n$. Donc:
 $$g'(s)=(s^{t'_1}.(1+s^e.x_1(s^e)),\dots,s^{t'_n}.(1+s^e.x_n(s^e)).$$
 
 Puisque $g'$ est g\'en\'eriquement injective pr\`es de $0$, on a bien:
 $$pgcd(e,t'_k,k=1,\dots,n)=1.$$
 
 Ce qui ach\`eve la preuve $\square$

\begin{re}\label{relv'} Sous les hypoth\`eses pr\'ec\'edentes, la d\'emonstration montre, plus g\'en\'eralement, que si $(R\vert \Delta_R)\to (X\vert\Delta)$ est un morphisme orbifolde, avec $R$ une courbe irr\'eductible, et si $R'\subset f^{-1}(R)\subset X'$ est une composante irr\'eductible, alors $v:R'\to (R\vert \Delta)$ est orbifolde-\'etale au sens suivant: \end{re}

\begin{definition}\label{eta} Soit $f:(B'\vert\Delta_{B'})\to (B\vert\Delta_B)$ un morphisme fini surjectif entre courbes lisses projectives et connexes $B',B$. On dit que $f$ est {\bf \'etale} si $f^*(K_B+\Delta_B)=K_{B'}+\Delta_{B'}$. Ceci \' equivaut \`a: $e(b').m'(b')=m(b)$ pour tout $b'\in B'$, $e(b')$ (resp. $m(b')$, resp. $m(b)$) d\'esignant l'ordre de ramification en $b'$ (resp. la multiplicit\'e de $\Delta_{B'}$ en $b'$, resp. la multiplicit\'e de $\Delta_B$ en $b:=f(b')$).

Si $f:(B'\vert\Delta_{B'})\to (B\vert\Delta_B)$ est \'etale, alors $(B'\vert\Delta_{B'})$ est rationnelle (resp. elliptique) si et seulement si $(B\vert\Delta_{B})$ l'est. (On dit que $(B'\vert\Delta_{B'})$ est elliptique si $deg(K_{(B'\vert\Delta_{B'})})=0$).
\end{definition}

\subsection{Quotients globaux: unir\'eglage et connexit\'e rationnelle}\label{QG}

\

Nous r\'eduisons ici, gr\^ace au th\'eor\`eme \ref{d-orb'}, l'\'etude des courbes $\Delta^{div}$-rationnelles \`a celle des courbes rationnelles usuelles lorsque $(X\vert\Delta)$ est un quotient global au sens de \ref{dqg}.

\begin{theorem}\label{comp} Soit $f:X'\to (X\vert\Delta)$ un quotient global au sens de la d\'efinition \ref{dqg}. Alors:

1. $(X\vert\Delta)^{div}$ est unir\'egl\'ee si et seulement si $X'$ est unir\'egl\'ee.

2. $(X\vert\Delta)^{div}$ est $fRCC$ si et seulement si $X'$ est $fRCC$.

3. $(X\vert\Delta)^{div}$ est $RC$ si et seulement si $X'$ est $RC$.

4. $(X\vert\Delta)^{div}$ est $ARC$ si et seulement si $X'$ est $ARC$.

\

Si $f:X'\to (X\vert\Delta)$ ramifie au moins au-dessus de $(X\vert\Delta)^{div}$, alors:

1. $(X\vert\Delta)^{div}$ est unir\'egl\'ee si $X'$ est unir\'egl\'ee.

2. $(X\vert\Delta)^{div}$ est $ARC$ si $X'$ est $RCC$.
\end{theorem}

{\bf D\'emonstration:} Elle est imm\'ediate, d'apr\`es les d\'efinitions, et le th\'eor\`eme \ref{d-orb'}.

\begin{example}\label{eqg'} Soit $(X\vert\Delta)$ l'example \ref{eqg}. Par \ref{comp} pr\'ec\'edent, $(X\vert\Delta)^{div}$ est $ARC$.
\end{example}

\begin{corollary}\label{comp'} Soit $f:X'\to (X\vert\Delta)$ un quotient global au sens de la d\'efinition \ref{dqg}. Alors les propri\'et\'es suivantes sont \'equivalentes:

1. $(X\vert\Delta)^{div}$ est $fRCC$

2. $(X\vert\Delta)^{div}$ est $RC$

3. $(X\vert\Delta)^{div}$ est $ARC$

\end{corollary}

{\bf D\'emonstration:} Evidente d'apr\`es \ref{comp} et l'\'equivalence entre les propri\'et\'es $ARC$ et $fRCC=RCC$ pour $X'$ (qui r\'esulte de [Ko-Mi-Mo92]).

\

\begin{corollary}\label{fano} Soit $f:X'\to (X\vert\Delta)$ un quotient global au sens de la d\'efinition \ref{dqg}. Si $(X\vert \Delta)$ est Fano (i.e: $-(K_X+\Delta)$ est ample sur $X$), alors:

1. $(X\vert\Delta)^{div}$ est $ARC$.

2. $\pi_1(X\vert \Delta)$ est fini.

\end{corollary}

\begin{re} La d\'emonstration de ce corollaire en supprimant l'hypoth\`ese ``quotient global" semble accessible \`a l'aide des champs de DM lisses. Voir aussi la remarque \ref{defrat}.
\end{re}

Nous en d\'eduisons maintenant une version faible de l'\'equivalence bim\'eromorphe (qui n'est pas aussi \'evidente que lorsque $\Delta=0$):

\begin{corollary}\label{comp"} Soit $f:X'\to (X\vert\Delta)$ un quotient global au sens de la d\'efinition \ref{dqg}, avec $X\in \sC$. Alors:

1. La propri\'et\'e d'\'evitement \ref{evit} est satisfaite par $(X\vert \Delta)$.

\

De plus, si $g:(X_1\vert\Delta')\to (X\vert \Delta)$ est un morphisme orbifolde divisible bim\'eromorphe, avec $g_*(\Delta_1)=\Delta$,  et $(X_1\vert\Delta_1)$ lisse, alors:

2. $(X\vert\Delta)^{div}$ est unir\'egl\'ee si et seulement si $(X_1\vert\Delta_1)^{div}$ est unir\'egl\'ee.

3. $(X\vert\Delta)^{div}$ est $RC$ si et seulement si $(X_1\vert\Delta_1)^{div}$ est $RC$.

\end{corollary}

{\bf D\'emonstration:} Les propri\'et\'es 2. et 3. d\'ecoulent de 1., gr\^ace \`a \ref{evit}. On va \'etablir 1.

Soit $A\subset X$, analytique ferm\'e de codimension au moins $2$, et $A'\subset X'$ son image r\'eciproque dans $X'$, qui y est de codimension au moins $2$ puisque $f$ est finie. Si on a une famille alg\'ebrique couvrante maximale de courbes $\Delta^{div}$-rationnelles de $X$, elle fournit (gr\^ace \`a \ref{d-orb'} par image r\'eciproque par $g$ une famille alg\'ebrique couvrante maximale de courbes rationnelles de $X'$ dont le membre g\'en\'erique \'evite donc $A'$. Leurs images par $g$ \'evitent donc $A$, et sont des courbes rationnelles de $(X\vert\Delta)^{div}$ $\square$

\begin{re} Pour \'etablir l'\'equivalence des conditions $RCC^{div}$ et $ARC^{div}$, ainsi que l'invariance bim\'eromorphe des conditions $UR^{div}$ et $RC^{div}$, il suffit donc de le faire pour les champs de DM (et les conditions $UR$ et $RC$ dans ette cat\'egorie). Voir aussi la remarque \ref{defrat}.
\end{re}

\begin{theorem} Admettons la conjecture \ref{-infty} pr\'ec\'edente lorsque $\Delta=0$. Alors cette conjecture est encore valable pour tout quotient global $f: X'\to (X\vert\Delta)$ au sens de la d\'efinition \ref{dqg} ci-dessous.\end{theorem}

{\bf D\'emonstration:} En effet: $f^*(K_X+\Delta)=K_{X'}$, de telle sorte que $\kappa(X\vert\Delta)=\kappa(X')=-\infty$. Donc $X'$ est unir\'egl\'ee (par \ref{-infty}), et donc aussi $(X\vert\Delta)$, par \ref{d-orb'} $\square$

\begin{corollary} Admettons la conjecture \ref{-infty} lorsque $\Delta=0$. Soit $(X\vert\Delta)$ un quotient global (au sens de \ref{dqg}), avec $X\in \sC$ connexe. Si $\kappa_+(X\vert\Delta)=-\infty$, alors $(X\vert\Delta)$ est $ARC^{div}$. \end{corollary}

{\bf D\'emonstration:} Soit $f:X'\to (X\vert\Delta)$ le morphisme fini faisant de $(X\vert\Delta)$ un quotient global. Alors $\kappa_+(X')=-\infty$, puisque $\kappa_+(X\vert \Delta)=-\infty$. Donc $X'$ est $ARC$, par la conjecture \ref{-infty} pour $\Delta=0$. Donc $(X\vert \Delta)$ est aussi $ARC$, par \ref{comp} $\square$

\subsection{Quotients globaux: le lemme de ``scindage" orbifolde.}

\

\begin{proposition}\label{negk} Soit $(X\vert\Delta)$ une orbifolde lisse avec $X\in \sC$ connexe, et $C_t$ une famille alg\'ebrique couvrante de courbes $\Delta^{Q}$-rationnelles. Alors: $(K_X+\Delta).C_t<0$. 

\end{proposition}

{\bf D\'emonstration:} D'apr\`es le th\'eor\`eme \ref{rescougen}, on peut, quitte \`a modifier $(X\vert\Delta)$ dans $Georb^Q$, supposer que le membre g\'en\'erique de cette famille $(C_t)_{t\in T}$ est lisse et rencontre $Supp(\Delta)$ transversalement en des points lisses. Si la propri\'et\'e est vraie sur l'orbifolde modifi\'ee, elle l'est sur $(X\vert \Delta)$, puisque l'effectivit\'e num\'erique du fibr\'e canonique orbifolde cro\^it par une telle modification. 

On consid\`ere le graphe d'incidence normalis\'e $q:Z\to X$ de la famille, et en utilisant le fait que $q^*(K_X)\leq K_Z+q^*(\Delta)$, puisque la famille est couvrante. On conclut avec la formule d'adjonction: $(K_X+\Delta).C_t=q^*(K_X+\Delta).C'_t\leq (K_Z+q^*(\Delta)). C'_t=deg(K_{C'_t})+q^*(\Delta).C'_t<0$, puisque $C'_t$ est $q^*(\Delta)$-rationnelle, par \ref{rescouim} et \ref{rescoutr}. (On a not\'e $C'_t$ l'image dans $Z$ de $C_t)$ $\square$

\

La question centrale d'existence de la g\'eom\'etrie des $\Delta$-courbes rationnelles est la r\'eciproque:

\begin{question}\label{cass} Soit $(X\vert\Delta)$ une orbifolde lisse avec $X\in \sC$ connexe. Soit $C\subset X$ une courbe projective  irr\'eductible non contenue dans le support de $\Delta$. Soit $a\in C$, $a\notin Supp(\Delta)$. Si $(K_X+\Delta).C<0$, existe-t'il une $\Delta$-courbe rationnelle $R$ passant par $a$?

Si $\Delta$ est enti\`ere, on peut demander \`a $R$ d'\^etre, soit $\Delta^{div}$-rationnelle (version forte), soit seulement $\Delta^Z$-rationnelle, ou m\^eme seulement $\Delta^Q$-rationnelle.
\end{question}

Nous allons d\'eduire du classique ``lemme de scindage" (``Bend and break") de Miyaoka-Mori [Mi-Mo 86] une r\'eponse positive \`a la version divisible dans le cas (tr\`es) particulier suivant.

\begin{theorem}\label{cass'} Soit $(X\vert\Delta)$ une orbifolde g\'eom\'etrique lisse, enti\`ere et \`a multiplicit\'es finies. On suppose $X$ projective. On suppose aussi qu'il existe un morphisme fini $f:X'\to X$, avec $X$ lisse, ramifiant au moins au-dessus de $\Delta^{div}$.

Soit $C\subset X$ une courbe projective irr\'eductible, non contenue dans $Supp(\Delta)$, et telle que $(K_X+\Delta).C<0$. 

Alors: si $a\in C$ n'appartient pas au support de $\Delta$, il existe une $\Delta^{div}$-courbe rationnelle $R$ de $X$ passant par $a$.

Si $H$ est un diviseur ample sur $X$, on peut choisir $R$ telle que: $$H.R\leq \frac{2. H.C}{-(K_X+\Delta).C}$$.
\end{theorem}

{\bf D\'emonstration:} Soit $C'$ une composante irr\'eductible de $f^{-1}(C)$. Alors $K_{X'}.C'=d.f^*(K_X+\Delta).C<0$ si $f_*(C')=d.C$. Soit $H':=f^*(H)$. Si $ a'\in C'$ est au-dessus de $a$, il existe, d'apr\`es [Mi-Mo 86], une courbe rationnelle $R'$ passant par $a'$ et telle que: $$H'.R'\leq 2.\frac{H'.C'}{-(K_{X'}.C')}$$

Posant $R:=f(R')$. Il suffit donc de montrer que $R$ est une $\Delta^{div}$-courbe rationnelle. Le th\'eor\`eme \ref{d-orb'} ci-dessus ach\`eve donc la d\'emonstration $\square$

\begin{corollary}\label{d-u} Soit $(X\vert \Delta)$ et $f:X'\to X$ satisfaisant les hypoth\`eses du th\'eor\`eme pr\'ec\'edent. On suppose que $(K_X+\Delta).C_t<0$ pour une famille $C_t$ de courbes irr\'eductibles de $X$ recouvrant un ouvert (analytique) non vide de $X$. Alors un ouvert de Zariski non-vide de $X$ est recouvert par une famille alg\'ebrique de courbes $\Delta^{div}$-rationnelles $R_s$ (i.e: $(X\vert \Delta)^{div}$ est unir\'egl\'ee).

Si la famille $C_t$ est alg\'ebrique, le $H$-degr\'e des $R_s$ satisfait la borne de \ref{cass'}.
\end{corollary}

\begin{re}\label{cass"} L'hypoth\`ese tr\`es restrictive de l'existence de $f:X'\to X$ dans les deux r\'esultats pr\'ec\'edents devrait pouvoir \^etre  supprim\'ee en rempla\c cant ce rev\^etement  par le champ alg\'ebrique (de Deligne-Mumford) associ\'e \`a $(X\vert \Delta)$, et en \'etablissant un ``lemme de scindage" dans cette cat\'egorie. \end{re}

\subsection{Quotients globaux: semi-positivit\'e g\'en\'erique.}

\

La pr\'esente section est sugg\'er\'ee par une question orale d'Ama\"el Broustet ([Br 08]): ``peut-t-on \'etendre aux orbifoldes g\'eom\'etriques le th\'eor\`eme de semi-positivit\'e g\'eom\'etrique de Miyaoka (voir [Miy 87] et [SP 92])?"

\

On peut pr\'eciser:

\begin{question}\label{sempos} Soit $(X\vert\Delta)$ une orbifolde lisse et finie, avec $X\in \cal C$. S'il existe des entiers strictement positifs $m,p$ tels que $S^m\Omega^p(X\vert\Delta)$ ait un quotient $\cal Q$ tel que $c_1(\cal Q)$ ne soit pas pseudo-effectif, alors $(X\vert\Delta)$ est-elle unir\'egl\'ee (i.e: g\'en\'eriquement recouverte par une famille de $\Delta$-courbes rationnelles)?
\end{question}

On a une r\'eponse positive lorsque $\Delta=0$ si $X$ est projective ([C-Pe 06]), et aussi dans la version divisible dans le cas (tr\`es) particulier suivant:

\begin{theorem}\label{sempos'} Soit $(X\vert\Delta)$ une orbifolde g\'eom\'etrique lisse, enti\`ere et \`a multiplicit\'es finies. On suppose $X$ projective. On suppose aussi qu'il existe un morphisme fini $f:X'\to X$, avec $X'$ {\bf lisse}, ramifiant {\bf exactement} au-dessus de $\Delta'$.

Alors $(X\vert\Delta)$ est unir\'egl\'ee s'il existe des entiers strictement positifs $m,p$ tels que $S^m\Omega^p(X\vert\Delta)$ ait un quotient $\cal Q$ tel que $c_1(\cal Q)$ ne soit pas pseudo-effectif, et si $m$ est divisible par toutes les multiplicit\'es de $\Delta$.
\end{theorem}

\

{\bf D\'emonstration:} C'est une r\'eduction imm\'ediate au cas o\`u $\Delta=0$. En effet: $f^*(S^m\Omega^p(X\vert\Delta))=Sym^m(\Omega^p_{X'})$ a pour quotient $f^*(\cal Q)$, qui n'est pas pseudo-effectif. Donc $X'$ est unir\'egl\'ee, par [C-Pe 06], theorem 0.3, et donc aussi $(X\vert\Delta)^{div}$, par \ref{d-orb'} $\square$

\begin{re} Ici encore, on pourrait s'affranchir de l'hypoth\`ese que $(X\vert \Delta)$ est un quotient global en \'etablissant le r\'esultat de [C-Pe 06] pour les champs alg\'ebriques de DM.
\end{re}

\subsection{Sections orbifoldes.}

\

Cette section tente de formuler le th\'eor\`eme de [GHS 03] dans le cadre des orbifoldes g\'eom\'etriques. Rappelons-le:

\begin{theorem}\label{ghs} ([G-H-S 03]) Soit $g':X'\to B'$ une fibration, avec $X',B'$ projectives, lisses et connexes, et $B'$ une courbe. Si la fibre g\'en\'erique de $g'$ est $RC$, alors $g'$ admet une section.
\end{theorem}

Les bases orbifoldes dans $Georb^{div}$ ne discernent pas l'existence sections locales (voir l'exemple \ref{orbdiv'} ci-dessous). On doit donc ici, se placer dans la cat\'egorie $Georb^Q$. Lorsque $(X\vert \Delta)$ est enti\`ere, une r\'eponse affirmative \`a la question \ref{Qdiv} permettrait de remonter les $\Delta^Q$-courbes rationnelles obtenues en des $\Delta^{div}$-courbes rationnelles.

\begin{definition}\label{dsect}  Soit $g:(X\vert \Delta)\to B$ une fibration sur une courbe projective lisse et connexe, $(X\vert \Delta)$ \'etant lisse, avec $X\in \sC$. Soit $\Delta_B=\Delta^Q_B:=\Delta(g\vert \Delta)^Q$ la base orbifolde de $(g,\Delta)$ dans $Georb^Q$.

Soit $j_C:C\subset X$ une courbe projective irr\'eductible surjective (donc finie) sur $B$, de normalisation $\nu:\hat C\to C$, et non contenue dans $Supp(\Delta)$. On munit $\hat C$ de la plus petite structure orbifolde enti\`ere $\Delta_C=\Delta_C^Q$ qui fait de l'inclusion normalis\'ee $j_C\circ \nu:\hat C\to (X\vert \Delta)^Q$ un morphisme orbifolde. 

Donc, $g_C:=g\circ j_C\circ \nu:(\hat C\vert \Delta_C)\to (B\vert \Delta_B)^Q$ est un morphisme orbifolde.

On dit que $C$ est une {\bf multisection} (resp. une {\bf section}) de $(g\vert \Delta)$ si $g_C$ est orbifolde-\'etale (i.e: si $(g_C)^*(K_B+\Delta_B)=K_{\hat C}+\Delta_C)$ (resp. et si, de plus, $g_C$ est de degr\'e $1)$. 
\end{definition}

\begin{question}\label{qs} Soit $g:(X\vert \Delta)\to B$ une fibration sur une courbe projective lisse et connexe, $(X\vert \Delta)$ \'etant lisse, avec $X\in \sC$. On suppose que la fibre orbifolde g\'en\'erique $(X_b\vert \Delta_b)$ de $g$ est $(RC^{Q})$. (Donc $X$ est Moishezon, et projective si $X$ est K\" ahler).
Alors: 

1. $(g\vert \Delta)$ a-t-elle une multisection $\sigma: B\to X$?

2. A-t-on: $\Delta^Z_B=0$ si $\Delta^{vert}=0$ (i.e: si $Supp(\Delta)$ ne contient aucune composante d'une fibre de $g$)?

3. $(g\vert \Delta)$ a-t'elle une section orbifolde si, de plus, $\Delta^{vert}=0$? 
\end{question}

\begin{re} Si, dans la situation de \ref{qs}, on suppose $(X\vert \Delta)$ enti\`ere, et si on consid\`ere $\Delta^{div}_B:=\Delta(g\vert \Delta)^{div}$, base orbifolde de $(g,\Delta)$ dans $Georb^{div}$, et de mani\`ere analogue, $\Delta_C^{div}$, alors $(g,\Delta)$ n'a en g\'en\'eral pas de multisection orbifolde (voir l'exemple \ref{orbdiv'} ci-dessous).
\end{re}

\begin{example}\label{orbdiv'} Consid\'erons par exemple la seconde projection $f_0: X_0:=\Bbb P_1\times \Bbb P_1\to \Bbb P_1:=B$. On \'eclate d'abord $X_0$ aux trois points $(0,0), (0,1)$ et $(0,\infty)$. On obtient ainsi $u:X_1\to X_0$ et sur $X_1$, six $(-1)-$courbes $F_0,F_1,F_{\infty}$ et $E_0,E_1,E_{\infty}$, o\`u les $F_i$ sont les transform\'es stricts des fibres $f^{-1}(i)$, et $E_i$ le diviseur exceptionnel correspondant, image r\'eciproque de $(0,i)$ pour $i=0,1,\infty$. On note $g:=f_0\circ u:X_1\to \bP_1$.

Soit $3\leq p<q$ deux entiers premiers entre eux. On munit $X_1$ du diviseur orbifolde $\Delta_1:=(1-1/p).(E_0+E_1+E_{\infty})+(1-1/q).(F_0+F_1+F_{\infty})$. Alors $(X_1\vert\Delta_1)^Q$ n'a pas de courbe rationnelle orbifolde $g$-horizontale $(B_1\vert\Delta_{B_1})$. En effet: $\Delta_{g,\Delta_1}^Q=(1-1/p).(\{0\}+\{1\}+\{\infty\})$, et donc: $\kappa((\bP_1\vert\Delta^Q_{g,\Delta_1}))=1$. L'inclusion normalis\'ee de $B_1$ dans $X_1$ compos\'ee avec $g$ fournirait alors un morphisme $g_1:(B_1\vert\Delta_{B_1})^Q\to (\bP_1\vert\Delta^Q_{g,\Delta_1})$ qui induirait pour tout $N$ une inclusion $g_1^*: H^0(S^N(\Omega^1(\bP_1\vert\Delta^Q_{g,\Delta_1}))\to H^0(S^N(\Omega^1(B_1\vert\Delta_{B_1}))$. Ce qui implique que $\kappa(B_1\vert\Delta_{B_1})\geq \kappa((\bP_1\vert\Delta^Q_{g,\Delta_1}))=1$ et contredit la rationalit\'e de $(B_1\vert\Delta_{B_1})^Q$.

Remarquer que, par contre,  $\Delta_{g,\Delta_1}^{div}=0$, puisque $pgcd(p,q)=1$. 
\end{example}

Dans le cas des quotients globaux, on peut encore r\'eduire la solution de la question \ref{qs} au cas o\`u $\Delta=0$, r\'esolu dans [GHS 03]:

\begin{proposition}\label{eta'} Soit $f:X'\to (X\vert\Delta)$ un quotient global au sens de \ref{dqg}. On suppose qu'il existe un diagramme commutatif:

\

\centerline{
\xymatrix{ X'\ar[r]^{f}\ar[d]_{g'} & (X\vert\Delta)\ar[d]^g\\
B'\ar[r]^{h}&B\\
}}

dans lequel les fl\`eches verticales sont des fibrations, et o\`u $B',B$ sont des courbes projectives lisses et connexes. Alors:

1. $K_{B'}\leq h^*(K_B+\Delta(g,\Delta)^{div})$.

2. Si $g'$ a une section $\sigma'$, alors: $K_{B'}= h^*(K_B+\Delta(g,\Delta)^{div})$, et $\Delta(g,\Delta)^{div}=\Delta(g,\Delta)^{Z}$. De plus:

3. Si $g'$ a une section $\sigma'$, et si $B_1:=$ est la normalis\'ee de $f\circ \sigma'(B')$, munie de la plus petite structure orbifolde $\Delta_{B_1}^{div}$ faisant de l'inclusion normalis\'ee $i:(B_1\vert \Delta_{B_1})^{div}\to (X\vert\Delta)^{div}$ un morphisme, alors $g: (B_1\vert \Delta_{B_1})^{div}\to (B\vert \Delta_{B})^{div}$ est \'etale.
\end{proposition}

{\bf D\'emonstration:}

{\bf 1.} Soit $b'\in B'$, et $e(b'):=e$ l'indice de ramification de $h$ en $b'$. Soit $b:=h(b')\in B$. Soit $x'\in X'$ tel que $g'(b')=b'$, et $x:=f(x')$. Soit $D$ une composante irr\'eductible de $g^{-1}(b)$ contenant $x$, et $D'$ une composante irr\'eductible de $f^{-1}(D)$ contenant $x'$.

Alors pr\`es de $x'$ et de $x$, on a: $$t_D.m_{\Delta}(D).D'+...=f^*(g^*(b))=(g\circ f)^*(b)=(h\circ g')^*(b)=g'^*(e.b')=e.d.D'+...,$$

en d\'esignant par $t_D$ (resp. $d$) la multiplicit\'e de $D$ (resp. $D'$) dans la fibre de $g$ au-dessus de $b$ (resp. de $g'$ au-dessus de $b'$).

Puisque, par d\'efinition, $m_{\Delta_B^{div}}(b)=pgcd_{\{D\}} \{t_D.m_{\Delta}(D)\}$, on voit que $e(b')$ divise $m_{\Delta_B^{div}}(b)$, ce qui est l'assertion 1.

{\bf 2.} Si $g'$ admet une section $\sigma'$ telle que $\sigma'(b')=x'$, on a: $d=1$. Donc: $e=t_D.m_{\Delta}(D)\geq m_{\Delta_B^{Z}}(b)$ divise $m_{\Delta_B^{div}}(b)$. Puisque $m_{\Delta_B^{div}}(b)$ divise (toujours) $m_{\Delta_B^{Z}}(b)$, on a l'\'egalit\'e: $e=t_D.m_{\Delta}(D)\geq m_{\Delta_B^{Z}}(b)=m_{\Delta_B^{Z}}(b)$, pour tout $b\in B$, ce qui est l'assertion 2.

Remarquons que $e=t_D$ si $m_{\Delta}(D)=1$.

{\bf 3.} On en d\'eduit donc que $f\circ \sigma':B'\to (B_1\vert \Delta_{B_1})^{div}$ est orbifolde-\'etale au sens de \ref{eta}. En effet, si $u,v$ sont des morphismes divisibles de courbes dont le compos\'e est orbifolde-\'etale, chacun des deux l'est aussi.

Puisque la compos\'ee $h=(g\circ f\circ \sigma'):B'\to (B\vert\Delta_B)^{div}$ est orbifolde-\'etale (par 2. ci-dessus), on en d\'eduit que $g:(B_1\vert \Delta_{B_1})^{div}\to (B\vert \Delta_{B})^{div}$ est aussi orbifolde-\'etale $\square$

\begin{corollary}\label{ghsqg}  Dans la situation du diagramme de \ref{eta'}, si les fibres de $g'$ sont RC, alors $g'$ a une section $\sigma'$ (d'apr\`es \ref{ghs}), et donc (d'apr\`es \ref{eta'}):

1. $K_{B'}= h^*(K_B+\Delta(g,\Delta)^{div})$ (autrement dit: $h$ est orbifolde-\'etale au sens de \ref{eta}).

2. $\Delta(g,\Delta)^{div}=\Delta(g,\Delta)^{Z}$.

3. $g: (B_1\vert \Delta_{B_1})^{div}\to (B\vert \Delta_{B})^{div}$ est \'etale au sens de \ref{eta}.

4. En particulier, si $(B\vert \Delta_{B})^{div}$ est rationnelle, $(X\vert \Delta)^{div}$ et $X'$ sont $RC$.

5. Si $\Delta^{vert}=0$, $h:B'\to B$ est \'etale, et $ \Delta_{B}^{Z}=0$. Si $B\cong \bP^1$, alors $(g,\Delta)$ a une section.
\end{corollary}

{\bf D\'emonstration:} Seule l'assertion 5. ne r\'esulte pas de \ref{eta'}, et est d\'emontr\'ee. Puisque $\Delta^{vert}=0$, $\Delta(g,\Delta)^Z=\Delta(g)^Z$. Les fibres de $g'$ \'etant $RC$, celles de $g$ le sont aussi, et $\Delta(g)^Z=0$, d'apr\`es [GHS 03]. Les autres assertions en sont des cons\'equences directes $\square$

\begin{re} L'\'enonc\'e de [GHS 03] dans la cat\'egorie des champs de DM lisses fournirait donc, avec les arguments ci-dessus, une r\'eponse affirmative \`a a question \ref{qs}.
\end{re}

Nous allons maintenant renforcer la question \ref{qs} pos\'ee ci-dessus:

\begin{conjecture}\label{ht'} Soit $g:(X\vert \Delta)\to B$ une fibration, avec $X,B$ projectives, lisses et connexes, $(X\vert \Delta)$ enti\`ere, et $B$ une courbe. On suppose que la fibre orbifolde g\'en\'erique $(X_b\vert\Delta_b)^{div}$ de $g$ est lisse et $RC$, et que le support de $\Delta$ n'a pas de composante irr\'eductible $g$-verticale (ie: contenue dans une fibre de $g$). Soit $E\subset X$ est un sous-ensemble fini en chaque point duquel la fibre de $g$ est lisse et r\'eduite.

Alors $(g,\Delta)$ admet une section orbifolde $\sigma$ dont l'image contient $E$.

De mani\`ere \'equivalente, pour toute composante irr\'eductible $D_j$, de multiplicit\'e $m_j$, du support de $\Delta$, et tout $b$ de $B$, l'ordre de contact en $\sigma(b)$ de $D_j$ avec $\sigma(B)$ est divisible par $m_j$, si $\sigma(b)\in D_j$.
\end{conjecture}

\begin{re} Il est clair que la conjecture \ref{ht'} fournit une section orbifolde de $(g,\Delta)$ (au sens de \ref{ghs}) dans la situation particuli\`ere o\`u $\Delta$ n'a pas de composante $g$-verticale. 
\end{re}

La conjecture \ref{ht'} est la version orbifolde g\'eom\'etrique de la suivante:

\begin{conjecture}\label{ht}([H-T 04]) Soit $g:X\to B$ une fibration, avec $X,B$ projectives, lisses et connexes, et $B$ une courbe. Si la fibre g\'en\'erique de $g$ est $RC$, et si $E\subset X$ est un sous-ensemble fini en chaque point duquel la fibre de $g$ est lisse et r\'eduite, alors $g$ admet une section dont l'image contient $E$.\end{conjecture}

\begin{proposition}\label{ghsorb} Soit $g:(X\vert \Delta)\to B$ une fibration, avec $X,B$ projectives, lisses et connexes, $(X\vert \Delta)$ enti\`ere, et $B$ une courbe. On suppose que la fibre orbifolde g\'en\'erique $(X_b\vert\Delta_b)$ de $g$ est $RC^{div}$. 

Admettons la conjecture \ref{ht'}. 

1. Alors $(g,\Delta)$ a une multisection orbifolde au sens de \ref{ghs}. 

2. En particulier, si $(X\vert \Delta)$ est lisse et si $(B\vert \Delta_{g,\Delta})^Q$ est rationnelle, alors $(X\vert\Delta)$ est $RCC^{div}$.
\end{proposition}

{\bf D\'emonstration:}

Soit $\Delta_B:=\Delta_{g,\Delta}^Q=\sum_{b\in B}(1-\frac{1}{m(b)}).\{b\}=\sum_k(1-\frac{1}{m_k}).\{b_k\}$, avec: $m_k:=m(b_k):=inf_{j}\{t_j.m_{\Delta}(F_j)\}$, si $g^*(b_k)=\sum_jt_j.F_j$. Pour chacun des $b_k$, soit $j_k$ un indice tel que $t_{j_k}.m_{\Delta}(F_{j_k})=m_k$. Posons: $m_k(g):=t_{j_k}$, et $\Delta_g:=\sum_k(1-\frac{1}{m_k(g)}).\{b_k\}$. 

Il existe alors un rev\^etement orbifolde \'etale: $u:B'\to (B\vert \Delta_g)$\footnote{Sauf si $B=\bP_1$, et si le support de $\Delta_g$ a $1$ ou $2$ \'el\'ements. Dans ce cas, on choisit simplement $B'=B$. La courbe $B_1$ r\'esultant de la construction qui suit ne sera pas une section orbifolde, mais sera une courbe $\Delta^{div}$-rationnelle, ce qui suffit pour nos applications.}. On consid\`ere $g': X':=\widehat{X\times_BB_1}\to B_1$, d\'eduit de $g$ par changement de base et normalisation. Alors l'image r\'eciproque $F'_{j_k}$ de chacun des $F_{j_k}$ est une composante r\'eduite de la fibre de $g'$ au-dessus de $b'_k:=u^{-1}(b_k)$. Notons $f: X'\to X$ la projection naturelle. Choisissant pour chaque $k$ un point lisse $E_k$ de $F'_{j_k}$, et notant $E$ la r\'eunion des $E_k$, on peut appliquer la conjecture \ref{ht'} \`a $(X'\vert f^*(\Delta^{hor}))$, o\`u $\Delta^{hor}$ est la partie $g$-horizontale de $\Delta$. On obtient ainsi une section $\sigma'$ de $g$ contenant $E$, et ayant les ordres de tangence requis le long de $f^*(\Delta)^{hor})$. 

On v\'erifie imm\'ediatement que $B_1:=f(\sigma'(B'))$ est une (multi)section orbifolde de $g: (X\vert\Delta)\to B$. La derni\`ere assertion est alors \'evidente $\square$

\subsection{Quotients rationnels orbifoldes}

\begin{theorem}\label{qrat}Soit $(X\vert\Delta)$ une orbifolde g\'eom\'etrique lisse et enti\`ere, avec $X\in \sC$ connexe. Il existe une unique fibration m\'eromorphe dominante $R:=R_{(X\vert\Delta)}:X\dasharrow R(X\vert\Delta)$, qui est presque-holomorphe, telle que\footnote{Ici encore, le r\'esultat et sa d\'emonstration sont valables pour les trois notions de $\Delta$-courbe rationnelle. On pr\'ecisera $R^*, *=div, Z, Q$ si besoin est.}:

1. Les fibres orbifoldes de $R$ sont $fRCC$.

2. Si $y\in R(X\vert\Delta)$ est g\'en\'eral, la fibre $X_y$ de $R$ au-dessus de $y$ ne rencontre aucune courbe rationnelle orbifolde de $(X\vert\Delta)$.

\end{theorem}

\begin{re} L'exemple \ref{ninvblog} montre que cette fibration n'est pas toujours un invariant bim\'eromorphe de l'orbifolde $(X\vert \Delta)$, lorsque celle-ci a une singularit\'e non klt.
\end{re}

{\bf D\'emonstration:} Nous utilisons ici (bri\`evement) l'appendice technique \ref{qm} ci-dessous. 

Si $(X\vert\Delta)$ n'est pas unir\'egl\'ee, $R=id_X$. Sinon, on d\'esigne par $A\subset \sC(X)$ l'ensemble des points param\'etrant une courbe rationnelle orbifolde (r\'eduite) de $(X\vert\Delta)$. On va montrer que $A$ est Z-r\'egulier (au sens de \ref{zreg}). La conclusion r\'esultera alors de \ref{ared}. Soit donc $B\subset \sC(X)$ un sous-ensemble analytique ferm\'e irr\'eductible (donc compact, par [Lieb78]) tel que $Z_a$ soit une courbe rationnelle orbifolde r\'eduite de $(X\vert\Delta)$ pour $a\in A':=A\cap B$, $A'$ non contenu dans une r\'eunion d\'enombrable de sous-ensembles analytiques ferm\'es stricts de $B$. Quitte \`a consid\'erer une modification $u:(X'\vert \Delta')\to (X\vert \Delta)$, bim\'eromorphe orbifolde minimale transverse \`a la famille $B$, le membre g\'en\'erique $Z'_b$ de la famille, param\'etr\'ee par $B$, de $X'$ obtenue par transform\'ee stricte par $u$ est donc une courbe rationnelle de $X'$, par le th\'eor\`eme \ref{rescougen}. Soit $Z'\subset B\times X'$ le graphe d'incidence (normalis\'e) de la famille $Z'_b,b\in B$, muni de ses projections $p':Z'\to X'$ et $q':Z'\to B$. Au-dessus d'un ouvert de Zariski dense $B^*$ de $B$, $q'$ est une submersion de fibre $\bP_1$. Soit $\Delta_{Z'}$ la plus petite structure orbifolde sur $Z'$ rendant $p':(Z'\vert \Delta_{Z'})\to (X'\vert \Delta')$ un morphisme orbifolde. Par transversalit\'e, la restriction transverse de $\Delta_{Z'}$ \`a $Z'_b$ est aussi la restriction de $\Delta$ \`a $Z_b$ (voir \ref{rescougen}). Cette restriction est donc une courbe $\Delta$-rationnelle pour tout $a\in A"$. Puisque $A"$ n'est pas contenu dans un sous-ensemble analytique ferm\'e strict de $B$, la restriction $(Z'_b\vert \Delta'_{Z',b})$ de $\Delta_{Z'}$ \`a $Z'_b$ est donc une courbe rationnelle orbifolde pour $b\in B$ g\'en\'erique. En fait, il existe une structure orbifolde fixe $\Delta_P$ sur $\bP_1$ (\`a automorphisme pr\`es de $\bP_1)$ telle que $(Z'_b\vert \Delta'_{Z',b})$ soit une courbe rationnelle orbifolde isomorphe \`a $(\Bbb P^1\vert \Delta_P)$, pour $b\in B$ g\'en\'erique. Il en est donc de m\^eme pour la restriction de $\Delta$ \`a $Z_b$, $b\in B$ g\'en\'erique, puisque ces deux restrictions coincident $\square$

\begin{question}\label{rir} La fibration $R$ pr\'ec\'edente est-elle un invariant bim\'eromorphe de l'orbifolde $(X\vert \Delta)$ si celle-ci est finie? Ceci semble pouvoir \^etre d\'emontr\'e en observant que les fibres d'un \'eclatement centr\'e en l'intersection de $m\geq 2$ composantes du support de $\Delta$ sont des projectifs munis d'une structure orbifolde finie support\'ee par les hyperplans de coordonn\'ees, et donc $ARC$, par (la g\'en\'eralisation imm\'ediate de) l'exemple \ref{eqg'}.
\end{question}

\begin{definition}\label{dre} Soit $(X\vert \Delta)$ lisse, avec $X\in \sC$. On dit que $(X\vert \Delta)$ est 

1. {\bf fortement rationnellement connexe} si, toute orbifolde lisse bim\'eromorphiquement \'equivalente \`a $(X\vert \Delta)$ est $fRCC$.

2. {\bf rationnellement engendr\'ee} si, pour toute fibration $f: (X\vert \Delta)\dasharrow Y$, toute base orbifolde stable de $(f\vert\Delta)$ est $UR$ \footnote{Ici encore, on peut se placer dans l'une des trois cat\'egories $Georb^*$.}. 

On note alors: $(X\vert \Delta)^*$ est $FRCC$ (resp. $RE)$, ou $(X\vert \Delta)$ est $FRCC^*$ (resp. $RE^{*})$, si l'on veut pr\'eciser la cat\'egorie $Georb^*$ dans laquelle on s'est plac\'e.
\end{definition}

\begin{re} 

\

1. Si $(X\vert \Delta)$ est $RE^*$ avec $dim(X)>0$, elle est $UR^*$.

2. Les propri\'et\'es $RE^{*}$ et $FRCC^*$ sont, par d\'efinition, pr\'eserv\'ees par \'equivalence bim\'eromorphe. On a l'implication: $FRCC^{*}\Longrightarrow RE^{*}$. Soit, en effet, $g:(X'\vert \Delta')\to Y$ une fibration holomorphe nette, avec $(X'\vert \Delta')$ bim\'eromorphe \`a $(X\vert \Delta)$. Puisque $(X'\vert \Delta')$ est $fRCC$, elle est recouverte par une famille de $\Delta'$-courbes rationnelles non contenues dans les fibres de $g$, et dont les images par $g$ fournissent donc un unir\'eglage de sa base orbifolde.

3. Lorsque $\Delta=0$, et $X$ projective (ou K\"ahler compacte), on a \'equivalence entre les conditions $RE$ et $RC$. On \'etablit cette \'equivalence en consid\'erant le ``quotient" rationnel de $X$, et en appliquant [GHS 03]. Voir \ref{rerc} pour le cas g\'en\'eral.
\end{re}

\begin{question}\label{rir'} Les propri\'et\'es $FRCC$ et $RE$ sont-elles \'equivalentes?
\end{question}

\begin{corollary}\label{barr} Soit $(X\vert\Delta)$ une orbifolde g\'eom\'etrique lisse et enti\`ere, avec $X\in \sC$ connexe. Il existe une unique fibration m\'eromorphe dominante $\bar R:=\bar R_{(X\vert\Delta)}:X\dasharrow R(X\vert\Delta)$ telle que\footnote{Ici encore, le r\'esultat et sa d\'emonstration sont valables pour les trois notions de $\Delta$-courbe rationnelle. On pr\'ecisera $\bar R^*, *=div, Z, Q$ si besoin est.}:

1. Les fibres orbifoldes de $R$ sont $FRCC$.

2. Si $y\in R(X\vert\Delta)$ est g\'en\'eral, la fibre $X_y$ de $R$ au-dessus de $y$ ne rencontre aucune courbe rationnelle orbifolde de $(X\vert\Delta)$. 
\end{corollary}

{\bf D\'emonstration:} Si $u:(X'\vert \Delta')\to (X\vert \Delta)$ est une modification bim\'eromorphe orbifolde, l'image par $u$ de la fibre g\'en\'erique de $R_{(X'\vert \Delta')}$ est contenue dans une fibre de $R_{(X\vert \Delta)}$, puisque sont image est $fRCC$. On choisit pour $\bar R$ la fibration dont la fibre g\'en\'erique est l'image dans $X$ d'une fibration $R_{(X'\vert \Delta')}$ dont la dimension des fibres est minimum. Cette fibration est un invariant bim\'eromorphe, car si $(X\vert \Delta)$ et $(X_1\vert \Delta_1)$ sont (\'el\'ementairement) domin\'ees par $(X^+\vert \Delta^+)$, les fibrations $\bar R$ d\'efinies ci-dessus coincident pour $(X\vert \Delta)$ et $(X^+\vert \Delta^+)$, ainsi que pour $(X_1\vert \Delta_1)$ et $(X^+\vert \Delta^+)$, donc aussi pour $(X\vert \Delta)$ et $(X_1\vert \Delta_1)$. Les fibres orbifoldes de $\bar R$ sont $FRCC$: on peut supposer que $\bar R=R$ (ie: que la dimension des fibres de $R$ est minimum). Si la fibre orbifolde $(X_y\vert \Delta_y)$ g\'en\'erique de $R$ n'est pas $FRCC$, il existe un mod\`ele bim\'eromorphe $(X'_y\vert \Delta'_y)$ de cette fibre telle que la fibration $R$ associ\'ee (qui est presque-holomorphe)  $R'_y:(X'_y\vert \Delta'_y)\to Z'_y$, avec $dim(Z'_y)>0$. Les arguments de \ref{gtgen} fournissent alors une factorisation: $R=h\circ g$, avec $g:(X\vert \Delta)\to Z$ et $h:Z\to Y$ , avec $dim(Z)>dim(Y)$ telle que les fibres g\'en\'eriques de $g$ soient les images dans $X$ de celles des $R'_y$. On en d\'eduit, par transform\'ee stricte de $\Delta$, l'existence d'une orbifolde lisse $(X'\vert \Delta')$ bim\'eromorphe \`a $(X\vert \Delta)$ et pour laquelle la fibration $R$ coincide avec $g$, contredisant la minimalit\'e de la dimension des fibres de $R$. La propri\'et\'e 2 est \'evidente, puisque $\bar R=R$ $\square$

\begin{proposition}\label{re}Soit $f:(X\vert\Delta)\dasharrow Y$ une fibration m\'eromorphe dominante, avec $X\in \sC$ connexe, et $(X\vert\Delta)$ lisse. Si la fibre orbifolde (g\'en\'erique) $(X_y\vert\Delta_y)^*$ et la base orbifolde stable de $f$, sont $RE^{*}$, alors $(X\vert\Delta)$ est $RE^{*}$.
\end{proposition}

{\bf D\'emonstration:} Le r\'esultat est trivial si $Y=X$. On suppose donc que $dim(Y)<dim(X)$. Les fibres de $f$ sont $RE$, donc unir\'egl\'ees. Soit $g:(X\vert\Delta)\dasharrow Z$ une fibration m\'eromorphe dominante. On peut supposer, par le th\'eor\`eme \ref{rescougen}, que $g$ est holomorphe nette. Si $g$ ne se factorise pas par $f$, la base orbifolde stable de $g$ est unir\'egl\'ee, puisque recouverte par les images par $g$ des fibres orbifoldes de $f$, qui sont $RE$, par hypoth\`ese, et dont les images (munies de leur structures orbifoldes restreintes) par $g$ sont  donc unir\'egl\'ees. Sinon, $g=h\circ f$, avec $f:X\to Y$ et $h:Y\to Z$ que l'on peut supposer holomorphes et nettes, par \ref{rescougen}. La base orbifolde de $h:(Y\vert \Delta_{f,\Delta})\to Z$ est aussi celle de $f$ (par la propri\'et\'e \ref{comporb}). Donc la base orbifolde de $f$ est $RE$, puisque celle de $g$ est RE, par hypoth\`ese $\square$

\begin{corollary}\label{qre} Soit $(X\vert\Delta)$ une orbifolde g\'eom\'etrique lisse, avec $X\in \sC$ connexe. Il existe\footnote{Le r\'esultat est valable pour les trois notions de $\Delta$-courbe rationnelle.} une unique fibration m\'eromorphe dominante 

$R_*:=R_{*(X\vert\Delta)}:X\dasharrow R_*(X\vert\Delta)$ telle que:

1. Les fibres orbifoldes de $R_*$ sont $RE$.

2. L'une des bases orbifoldes stables de $R_*$ n'est pas unir\'egl\'ee.
\end{corollary}

{\bf D\'emonstration:} Si $(X\vert \Delta)$ n'est pas $UR$, on prend $R_*:=id_X$. Sinon, on consid\`ere la fibration $\bar R$ donn\'ee par \ref{barr}. Ses fibres sont $FRCC$, donc $RE$. 

Si l'une des bases orbifoldes stables de $R_0$ n'est pas unir\'egl\'ee, on d\'efinit: $R=R_0$. Sinon, on d\'efinit, par r\'ecurrence sur $dim(X)$: $R_*:=(R_B)_*\circ R_0$, si $B$ est une base orbifolde stable de $R_0$, et $(R_B)_*$ la fibration $R_*$ correspondante. Les fibres orbifoldes de la fibration compos\'ee sont bien $RE$, par \ref{re} $\square$

\begin{re} Si l'implication $RE\Longrightarrow fRCC$ est vraie lorsque les multiplicit\'es sont finies, alors $R=R_*$. Ce serait le cas si les questions \ref{qqrat} avaient une r\'eponse positive.
\end{re}

On peut du moins \'etendre au cas $RE$ certaines des propri\'et\'es valables dans le cas $RC$.

\begin{proposition}\label{kre} Soit $(X\vert \Delta)$ lisse, avec $X\in \sC$. Alors:

1. Si $(X\vert \Delta)^Q$ est $RE$, on a: $\kappa_+(X\vert \Delta)=-\infty$.

2. Si $(X\vert \Delta)^{div}$ est $RE$ et finie, on a: $\pi_1(X\vert \Delta)$ est fini.
\end{proposition}

{\bf D\'emonstration:} L'assertion 2 sera \'etablie en \ref{repif}. Pour d\'emontrer l'assertion 1, il suffit, par la construction de $R_*$, de montrer que si $r:(X\vert \Delta)\to R$ est une fibration dont les fibres orbifoldes g\'en\'eriques $(F\vert \Delta_F)$ sont $UR^Q$, et dont la base orbifolde stable $(R\vert \Delta_R)$ satisfait: $\kappa_+(R\vert \Delta_R)=-\infty$, alors: $\kappa_+(X\vert \Delta)=-\infty$. Mais ceci r\'esulte des arguments fournis en \ref{k+} $\square$

\begin{corollary} Admettons la conjecture \ref{-infty}. Soit $(X\vert\Delta)$ une orbifolde g\'eom\'etrique lisse, avec $X\in \sC$ connexe. Si $\kappa_+(X\vert\Delta)=-\infty$, alors $(X\vert\Delta)$ est $RE$. \end{corollary}

{\bf D\'emonstration:} Soit $R_*:=R_{*(X\vert\Delta)}:X\dasharrow R_*(X\vert\Delta)$ la fibration construite en \ref{qre}. Sa base est de dimension nulle, puisque $\kappa_+(X\vert\Delta)=-\infty$. Donc $(X\vert\Delta)$ est RE $\square$

\begin{re} Admettant une autre conjecture ($C_{n,m}^{orb}$, voir \S\ref{cnmorb}), nous construirons en \ref{k-qr} une variante de la fibration $R_*$: ses fibres orbifoldes g\'en\'erales ont $\kappa_+=-\infty$, et sa base orbifolde stable est de dimension canonique $\kappa\geq 0$ (elle n'est donc pas unir\'egl\'ee). Si \ref{-infty} est vraie, ces deux fibrations coincident.\end{re}

\subsection{Quotients globaux: l'implication $RE\Longrightarrow RC.$}

\

Lorsque $\Delta=0$, et $X$ projective (ou K\"ahler compacte), on a \'equivalence entre les conditions $RE$ et $RC$. On \'etablit cette \'equivalence en consid\'erant le ``quotient" rationnel de $X$, et en appliquant [GHS 03].

\begin{question}\label{rerc} Soit $(X\vert \Delta)$ lisse, finie et enti\`ere, avec $X\in \sC$. Si $(X\vert \Delta)$ est $RE^{div}$, est-t-elle $RC^{div}$?
\end{question}

La question n\'ecessite, pour \^etre r\'esolue affirmativement, de pouvoir:

1. montrer que la condition $RC^{div}$ est bim\'eromorphiquement pr\'eserv\'ee. Ce qui r\'esulterait de la condition d'\'evitement des lieux de codimension au moins $2$ (\'enonc\'ee dans \ref{evit}).

2. montrer que si $g:(X\vert\Delta)\to B$ est une fibration sur une courbe, avec des fibres orbifoldes g\'en\'eriques $RC^{div}$, alors il existe une courbe $\Delta^{div}$-rationnelle surjective sur $B$ si $(B\vert\Delta^{Z}_{g,\Delta})$ est rationnelle. Ce qui est le cas si $(g,\Delta)$ a une multisection (au sens de \ref{dsect}).

Nous allons donner, conditionnellement en la condition d'\'evitement \ref{evit}, une r\'eponse affirmative pour les quotients globaux.

\begin{theorem}\label{rercc} Soit $X'\to (X\vert\Delta)$ un quotient global (voir \ref{dqg}). Si $(X\vert\Delta)$ est $RE^{div}$, il est $(ARC)^{div}$ si la condition (C) suivante est satisfaite:

(C) Toute $(Z\vert \Delta_Z)$ lisse, enti\`ere et finie qui est $UR^{div}$ contient, pour tout $A\subset Z$, Zariski-ferm\'e de codimension au moins $2$, une famille alg\'ebrique couvrante de courbes $\Delta_Z^{div}$-rationnelles dont le membre g\'en\'erique ne rencontre pas $A$.
\end{theorem}

{\bf D\'emonstration:} Il nous suffit de montrer que $X'$ est $RC$, par \ref{comp}. Nous pouvons supposer que le rev\^etement $f: X'\to X$ est Galoisien de groupe $G$, fini (voir \ref{dqg}). Soit donc $g':X'\dasharrow Y'_0$ le quotient rationnel de $X'$ ([Ca 92], [Ko-Mi-Mo 92]): ses fibres sont $RC$ et sa base n'est pas unir\'egl\'ee  si ce n'est pas un point. Supposons donc $dim(Y')>0$. Nous allons montrer que $Y'$ est unir\'egl\'e, ce qui \'etablira le r\'esultat. La fibration $g'$ est pr\'eserv\'ee par le groupe $G$. Soit $h_0:Y'_0\to Y_0:=Y/G$ le quotient naturel. Il existe donc une application $g:X\dasharrow Y_0:=Y/G$ telle que $g\circ f=h\circ g'$.

Il existe donc des modifications $u:X_1\to Y$ et $u':X'_1\to X'$, ainsi que des applications holomorphes: $f_1:X'_1\to X_1$ et $h:Y'\to Y$ telles que:

1. $h\circ g'_1=g_1\circ f_1:X'_1\to Y$, les applications $u'$ et $f_1$ \'etant $G$-\'equivariantes.

2. Pour toute structure orbifolde $\Delta_1$ sur $X_1$ telle que $u_*(\Delta_1)=\Delta$ et telle que $u:(X_1\vert\Delta_1)^{div}\to (X\vert\Delta)^{div}$ soit un morphisme orbifolde, alors: (puisque $(X\vert\Delta)$ est $RE^{div}$ par hypoth\`ese), la base orbifolde $(Y\vert\Delta_Y)$ est lisse et unir\'egl\'ee, avec: $\Delta_Y:=(\Delta_{g_1,\Delta_1})^{Z}$.

3. $h:Y'\to Y$ est finie (mais $Y'$ est seulement normal), et est un rev\^etement Galoisien de groupe $H$, quotient de $G$.

Il existe donc sur $Y$ une unique structure orbifolde $\Delta'_Y$ telle que $h:Y'\to Y$ ramifie {\bf exactement} au-dessus de $\Delta'_Y$. On peut, de plus, puisque $(X\vert \Delta)$ est, par hypoth\`ese, $RE^{div}$, supposer que:

4. $\Delta'_Y$ divise  $\Delta_Y$ (i.e: $m_{\Delta'_Y}(E)$ divise $m_{\Delta_Y}(E),\forall E\in W(Y))$.

En effet: soit $E\in W(Y)$ telle que $m'(E)$ ne divise pas $m(E)$, notant $m'(E)$ (resp. $m(E))$ la multiplicit\'e de $E$ relative \`a $\Delta'_Y$ (resp. $\Delta_Y$). 

Il existe donc $D\in W(X_1)$ tel que $m(E)=t_E(D).m_{\Delta_1}(D)$, avec $g_1^*(E)=t_E(D).D+...$, et $t_E(D)>0$, i.e: $g_1(D)=E$. 

On a alors $2$ cas: ou bien $D$ est $u$-exceptionnel, ou bien $u(D)=D_j\in W(X)$. 

Dans le premier cas, on peut augmenter $m_{\Delta_1}(D)$ de telle sorte que $t_E(D).m_{\Delta_1}(D)$ soit un multiple de $m'(E)$: $(Y\vert\Delta_Y)$ reste unir\'egl\'ee, puisque, par hypoth\`ese, toute base orbifolde stable de $(g,\Delta)$ est unir\'egl\'ee. 

Dans le second cas, si $D'\in W(X')$ et si $g'_1(D')=E'$, et $h(E')=E$, on a: $(h\circ g'_{1})$$^*(E)=g'_{1}$$^*(m'(E).E')+\dots=(t'_{E'}(D').m'(E)).D'+...=(g_1\circ f_{1})$$^*(E)=f_{1}$$^*(t_E(D).D+...)=(t_E(D).m_{\Delta}(D_0).D+...$, o\`u $D_0:=u(D)$, puisque l'ordre de ramification de $f_1$ au-dessus de $D$ coincide avec celui de $f$ au-dessus de $D_0$.

On a donc (puisque $f_1$ et $h$ sont Galoisiens, et que seuls des diviseurs non $g_1\circ f_1$-exceptionnels interviennent): 

$t'_{E'}(D').m'(E)=t_E(D).m_{\Delta}(D_0)=t_E(D).m_{\Delta_1}(D)=:m(E)$. 

Donc: $m'(E)$ divise $m(E)$.

En r\'ep\'etant un nombre fini de fois cette construction (au plus autant de fois que le nombre de diviseurs $D$ qui sont $u$-exceptionnels dans les images r\'eciproques par $g_1$ des composantes du support de $\Delta_1)$, on obtient la divisibilit\'e annonc\'ee.

Soit alors $B\subset Y$ une courbe $\Delta^{div}_Y$-rationnelle, membre g\'en\'erique d'une famille couvrante de telles courbes. C'est, a fortiori, une courbe $(\Delta'_Y)^{div}$-rationnelle, et on a donc: $(K_Y+\Delta'_Y).B<0$. Soit $Z'\subset Y'$ le lieu singulier de $Y'$, et $Z:=h(Z')\subset Y$. Donc, $Z$ est de codimension au moins $2$ dans $Y$. Puisque $(Y\vert\Delta_Y)$ satisfait la condition (C) par hypoth\`ese, nous pouvons supposer que $B$ \'evite $Z$. Soit $B'$ une composante irr\'eductible de $h^{-1}(B)$, membre d'une famille couvrante de $Y'$, avec $h_*(B')=d.B$, $d>0$: on a donc: $K_{Y'}.B'=d.(K_Y+\Delta'_Y).B<0$. Par [Mi-Mo 86], $Y'$ est unir\'egl\'e, et le th\'eor\`eme d\'emontr\'e $\square$

\subsection{Hyperbolicit\'e alg\'ebrique.}\label{sha}

Cette section est inspir\'ee par une discussion avec Aaron Levin.

\begin{question}\label{qha} Soit $(X\vert\Delta)$ une orbifolde g\'eom\'etrique lisse et enti\`ere, $X$ projective. On suppose que $\kappa(K_X+\Delta)=n:=dim(X)>0$. Alors existe-t-il un sous-ensemble alg\'ebrique $A\subsetneq X$ contenant toutes les courbes $\Delta^{div}$-rationnelles et $\Delta^{div}$-elliptiques (au sens divisible, donc\footnote{ On peut se poser la question, plus difficile, pour le cas non-divisible \'egalement.}) de $X$? Et plus g\'en\'eralement, existe-t'il $A'\subsetneq X$ alg\'ebrique, contenant toutes les sous-vari\'et\'es $V\subset X$ dont les restrictions minimales $(V'\vert \Delta')$ de $\Delta$  (au sens de \ref{fdefrestr}) ne sont pas toutes de type g\'en\'eral?
\end{question}

\begin{re}\label{eha} Un cas connu est celui des orbifoldes logarithmiques $(\Bbb P^n\vert H_d\vert D)$, si $H_d\subset \Bbb P^n$ est une hypersurface lisse g\'en\'erale de dimension $n$ et de degr\'e $d\geq 2n+1$. (Voir [P-R 06], et [B 95, thm 1.5.2] pour une approche diff\'erente dans le cas des surfaces). Le cas propre (o\`u $X=H_d$, $\Delta=0)$ a \'et\'e trait\'e auparavant par H. Clemens et C. Voisin). \end{re}

 \subsection{Appendice: Quotients m\'eromorphes.}\label{qm}

  Nous rappelons ici les r\'esultats de [Ca04, appendice] (auquel nous renvoyons pour les d\'emonstrations et plus de d\'etails).  Ils seront utilis\'es dans le pr\'esent texte dans les \S\ref{qrat} et \ref{gamred}.

  Soit $X\in \sC$, normal et connexe. On note $\sC(X)$ ou $Chow(X)$ la vari\'et\'e de Chow (ou espace des cycles) de $X$ construit dans [Ba75]. Pour $a\in \sC(X)$, on note $Z_a\subset X$ le support du cycle param\'etr\'e par $a$.
  
  On note $A\subset \sC(X)$ un sous-ensemble (ensembliste).  On supposera toujours que la famille $A$ est {\bf couvrante}, c'est-\`a-dire que la r\'eunion des $Z_a,a\in A$ est $X$. Soit $R_A\subset X\times X$ la relation d'\'equivalence pour laquelle deux points de $X$ sont \'equivalents si et seulement s'ils sont contenus dans $A$-chaine, ie: une r\'eunion finie {\it connexe} de $Z_a, a\in A$. 
  
  Si $V\subset X$ est analytique ferm\'e irr\'eductible, on dira que $V\in A$ s'il existe $a\in A$ tel que $Z_a=V$.

  \begin{theorem}\label{aq} On suppose $A$ analytique ferm\'e dans $\sC(X)$. Il existe alors une unique fibration $q_A:X\dasharrow Q_A$, qui est presque-holomorphe (voir \ref{ph}), et telle que pour $b\in Q_A$ g\'en\'eral (g\'en\'erique si $A$ a un nombre fini de composantes), la fibre $X_b=q_A^{-1}(b)$ est la classe de $R_A$-\'equivalence de chacun de ses points. La fibration $q_A$ est appel\'ee le $A$-quotient de $X$.
  \end{theorem}

\begin{definition} \label{zreg} On dit que $A\subset \sC(X)$ est Z-r\'egulier si, pour tout $B\subset \sC(X)$, analytique ferm\'e et irr\'eductible, $A\cap B$ soit contient une intersection d\'enombrable d'ouverts de Zariski denses de $B$, soit est contenu dans une r\'eunion d\'enombrable de sous-ensembles analytiques ferm\'es stricts de $B$. 

Si $A$ est Z-r\'egulier, il existe une unique r\'eunion finie ou d\'enombrable de sous-ensembles analytiques ferm\'es irr\'eductibles sans inclusions $B_n$ de $\sC(X)$ (appel\'es les composantes de $A)$, tels que, pour chaque $n$, $A\cap B_n$ contienne une intersection d\'enombrable d'ouverts de Zariski denses de $B_n$, et tels que $A$ soit contenu dans la r\'eunion des $B_n$. On note $B$ la r\'eunion des $B_n$.
\end{definition}

\begin{theorem}\label{ared} Soit $A\subset \sC(X)$, Z-r\'egulier et couvrant. Soit $q_A:X\dasharrow Q_A$ le $B$-quotient, $B$ \'etant la r\'eunion des composantes de $A$ dans $\sC(X)$. Si $b\in Q_A$ est g\'en\'eral, et si $Z_a,a\in A$, rencontre $X_b$, alors $Z_a\subset X_b$. 
\end{theorem}

\begin{definition}\label{dstab} Soit $A\subset \sC(X),X\in \sC$. On dit que $A$ est {\bf stable} si, pour tout $V\subset X$ analytique ferm\'e irr\'eductible muni d'une fibration m\'eromorphe dominante $g:V\dasharrow W$, alors $V\in A$ si:

1. Les fibres g\'en\'erales de $g$ sont dans $A$.

2. Il existe $Z\subset V$, $Z\in a$ tel que $g(Z)=W$.\end{definition}

\begin{theorem}\label{stab} Soit $A\subset \sC(X)$, Z-r\'egulier et couvrant. Soit $q_A:X\dasharrow Q_A$ le $B$-quotient, $B$ \'etant la r\'eunion des composantes de $A$ dans $\sC(X)$. Si $A$ est stable, et si $b\in Q_A$ est g\'en\'eral, alors $X_b\in A$. 
\end{theorem}

\newpage

\section{ADDITIVIT\' E ORBIFOLDE}

On va rappeler ici certaines conjectures et r\'esultats de [Ca04, \S 4, pp. 564-574] auquel nous renvoyons pour les d\'emonstrations et d\'etails. Ces r\'esultats sont techniquement essentiels pour \'etablir les r\'esultats principaux du pr\'esent texte (le ``coeur" et sa d\'ecompostion).

\subsection{La conjecture $C_{n,m}^{orb}$.}

\begin{conjecture} (Conjecture $C_{n,m}^{orb}$) \label{cnmorb}Soit $f:(Y\vert\Delta)\to S$ une fibration holomorphe, l'orbifolde g\'eom\'etrique $(Y\vert\Delta)$ \'etant lisse avec $Y$ compacte et connexe dans la classe $\sC$. 

Alors $\kappa(Y\vert\Delta)\geq \kappa(Y_s\vert\Delta_s)+\kappa(f\vert\Delta)$

On a not\'e $(Y_s\vert\Delta_s)$ la fibre orbifolde de $f$ au-dessus du point g\'en\'eral $s\in S$. (Remarquons que cette orbifolde g\'eom\'etrique est lisse, par le th\'eor\`eme de Sard, et $ \kappa(Y_s\vert\Delta_s)$ est ind\'ependant de $s\in S$, g\'en\'eral, par le th\'eor\`eme de coh\'erence des images directes de Grauert, appliqu\'e aux $f_*(m(K_Y+\Delta))$, $m>0$ assez divisible). \end{conjecture}

 \begin{re} \label{fonctm}

 1. Ici comme partout ailleurs\footnote{\`A l'exception des d\'efinitions du groupe fondamental, du rev\^etement universel, de la pseudom\'etrique de Kobayashi, et des points entiers au sens des orbifoldes g\'eom\'etriques.}, les coefficients des composantes de $\Delta$ sont rationnels dans [0,1], et pas n\'ecessairement de la forme ``standard" $(1-\frac{1}{m})$.
 
 2. Lorsque $f$ est nette (au sens de \ref{nette}), on a aussi: $\kappa(f\vert\Delta)=\kappa(S\vert\Delta(f,\Delta))$.
 
 3. Cette conjecture est \'evidemment la version orbifolde de la conjecture $C_{n,m}$ d'Iitaka, qui affirme que $\kappa(Y)\geq \kappa(Y_s)+\kappa(S)$ si $Y$ est projective.
 
 4. Cette conjecture a de nombreuses cons\'equences, dont certaines seront developp\'ees ci-dessous au \S\ref{deccoeur}. L'une d'entre elles est l'in\'egalit\'e: $\kappa(X\vert \Delta)\geq \kappa(Y\vert \Delta_Y)$ si l'on a un morphisme orbifolde surjectif $f:(X\vert \Delta)\to (Y\vert \Delta_Y)$ entre orbifoldes lisses, compactes et connexes, et si $\kappa(X\vert \Delta)\geq 0$. Une cons\'equence de cette in\'egalit\'e est la fonctorialit\'e des applications de Moishezon-Iitaka si, de plus, $\kappa(Y\vert \Delta_Y)\geq 0$: il existe $M_f:M(X\vert \Delta)\to M(Y\vert \Delta_Y)$ telle que $M_Y\circ f=M_f\circ M_X$. Voir le \S\ref{moiit} pour ces notions.
  \end{re}

\subsection{Le cas des fibrations de type g\'en\'eral.}

Le r\'esultat principal est ici le:

 \begin{theorem}\label{addorb} Lorsque $f:(Y\vert\Delta)\to S$  est, de plus, une fibration de type g\'en\'eral, la conjecture pr\'ec\'edente $C_{n,m}^{orb}$ est vraie, et on a alors: 
 
 $\kappa(Y\vert\Delta)= \kappa(Y_s\vert\Delta_s)+dim(S)$.  \end{theorem}

 Ce th\'eor\`eme est ais\'ement d\'eduit du suivant, adaptation au cadre orbifolde de r\'esultats de E. Viehweg, initi\'es par T. Fujita et Y. Kawamata (Voir [Fuj78], [Kaw80], [Vie83]; un r\'esultat similaire pour les fibrations telles que $\kappa(XY_x\vert \Delta_x)=0$ est d\^u \`a Y. Kawamata, dans le contexte ``num\'erique" [Kaw98]):

  \begin{theorem}\label{fporb}[Ca04, 4.13, p. 568] Soit $f:Y\to S$ une fibration, avec $Y,S$ lisses, $Y$ dans la classe $\sC$, et $S$ projective. Soit $D=\sum m_j.D_j$ un diviseur entier et effectif sur $Y$ dont le support est \`a croisements normaux au-dessus du point g\'en\'erique de $S$. Soit $m>0$ un entier tel que $m\geq m_j$, pour tout $j$ tel que $D_j$ soit une composante $f$-horizontale de $D$ (ie: tel que $f(D_j)=S)$. Alors: $f_*(mK_{Y/S}+D)$ est un faisceau coh\'erent faiblement positif \footnote{Voir [Vi83] pour cette notion, due \`a E. Viehweg. Des rappels se trouvent aussi dans [Ca04, \S 4]} sur $S$.
  \end{theorem}

 \begin{re} 
 
 \
 
 1. Bien qu'obtenue par les m\^emes m\'ethodes que celles de [Vi 83], cette g\'en\'eralisation en \'etend consid\'erablement le champ d'application.
 
 2. Il r\'esulte de \ref{addorb} que si $f:(X\vert\Delta)\dasharrow Y$ est une fibration de type g\'en\'eral (ie: la base orbifolde d'un mod\`ele holomorphe net est de type g\'en\'eral), alors les fibres g\'en\'eriques orbifoldes de deux mod\`eles nets holomorphes de $f$ ont la m\^eme dimension canonique, \'egale \`a la diff\'erence entre $\kappa(X\vert\Delta)$ et $dim(Y)$. Nous verrons en fait (en \ref{phftg}) un r\'esultat plus pr\'ecis: une telle fibration est presque-holomorphe. 
 
 3. Dans le Lemme 4.10, p.567 de [Ca04], l'hypoth\`ese (par exemple) que $g_*(E)$ est localement libre a \'et\'e omise, comme me l'a signal\'e O. Debarre. Cette hypoth\`ese est difficilement v\'erifiable en pratique, mais on a la version suivante, qui couvre les applications pr\'esentes: 
  \end{re}

  \begin{proposition} Soit $f:Y\to S$ une fibration avec $Y$ compact, normal et connexe. Soit $A$ un $\bQ$-diviseur ample sur $S$ et $D$ un fibr\'e en droites sur $Y$ tel que $f_*(D)$ soit faiblement positif. Alors $\kappa(Y,D+f^*(A)+E)=\kappa(Y_s,D_s)+dim(S)$, pour $E$ un diviseur effectif $f$-exceptionnel ad\'equat sur $Y$. S'il existe un morphisme birationnel $v:Y\to Y'$ contractant tous les diviseurs $f$-exceptionnels de $Y$, et un fibr\'e en droites $D'$ sur $Y'$ tel que $D=v^*(D')$, on peut prendre $E=0$.
  
   \end{proposition}

  Lorsque $f$ est nette (c'est la situation pr\'esente, et aussi celle consid\'er\'ee dans [Ca04]), l'hypoth\`ese de la seconde assertion de la proposition est satisfaite (et le lemme 4.10 peut donc bien \^etre appliqu\'e tel quel).

  {\bf D\'emonstration:} (C'est, sous une forme simplifi\'ee, celle de \ref{k=L} et de \ref{kappa=L} ci-dessus). Il suffit de montrer que $H^0(Y,m(L+E+f^*(A))\neq 0$ pour $E$ et $m>0$ ad\'equats. Par hypoth\`ese, $H^0(S, S^m({F})\otimes mA)\neq 0$, si $S^m({F})$ est le bidual de $Sym^m(f_*(D))$. Une section non nulle de ce faisceau fournit donc une section de $m(D+f^*(A))$ ayant des p\^oles sur un diviseur $f$-exceptionnel de $Y$ (on suppose $mA$ entier). D'o\`u la premi\`ere assertion. La seconde assertion r\'esulte de ce que, sous l'hypoth\`ese de contractibilit\'e de l'\'enonc\'e, toute section de $m(D+f^*(A))$ d\'efinie sur le compl\'ementaire du lieu exceptionnel de $f$ se prolonge \`a $Y$ tout entier $\square$

   \subsection{Premi\`ere application: $\kappa=0$ et orbifoldes Fano.}

   Notre premi\`ere application du th\'eor\`eme \ref{addorb} est l'exemple fondamental suivant d'orbifolde g\'eom\'etrique lisse ``sp\'eciale" (voir aussi la d\'efinition \ref{defspec} pour ce terme):

\begin{theorem}\label{k=ospec} Soit $(Y\vert\Delta)$ une orbifolde g\'eom\'etrique lisse, avec $Y$ compacte et connexe dans la classe $\sC$. Si $\kappa(Y\vert\Delta)=0$, alors $(Y\vert\Delta)$ est ``sp\'eciale" (ie: il n'existe pas de fibration de type g\'en\'eral  $f:(Y\vert\Delta)\dasharrow X)$.
\end{theorem}

   {\bf D\'emonstration:} Sinon, soit $f:(Y\vert\Delta)\dasharrow S$ une fibration de type g\'en\'eral. On a donc $dim(S)>0$, et (c'est imm\'ediat): $\kappa(Y_s\vert\Delta_s)\geq 0$. Donc: $0=\kappa(Y\vert\Delta)=\kappa(Y_s\vert\Delta_s)+dim(S)\geq dim(S)>0$, par \ref{addorb}. Contradiction $\square$

  \
  
  Le m\^eme argument fournit un r\'esultat un peu plus g\'en\'eral:

  \begin{theorem} \label{essk>o}Soit $f:(Y\vert\Delta)\to S$ une fibration holomorphe de type g\'en\'eral, l'orbifolde g\'eom\'etrique $(Y\vert\Delta)$ \'etant lisse avec $Y$ compacte et connexe dans la classe $\sC$. 
  
  Si $\kappa(Y\vert\Delta)\geq 0$, alors $dim(S)\leq \kappa(Y\vert\Delta)$, et on a \'egalit\'e si et seulement si $\kappa(Y_s\vert\Delta_s)=0$ (auquel cas $f$ est la fibration de Moishezon-Iitaka de $(X\vert \Delta))$\end{theorem}

\begin{corollary}\label{fanorb} Soit $(X\vert\Delta)$ lisse et Fano (ie: $-(K_X+\Delta)$ est ample sur $X)$. 
Alors $(X\vert\Delta)$ est sp\'eciale.
\end{corollary}

{\bf D\'emonstration:} Soit $H$ une section lisse de $-m(K_X+\Delta)$ intersectant transversalement $\Delta$. Soit $\Delta':=\Delta+(1/m).H$. Alors: $(X\vert\Delta')$ est lisse, et $K_X+\Delta'$ est de $\bQ$-torsion, donc $\kappa(X\vert\Delta')=0$. Donc $(X\vert\Delta')$ est sp\'eciale, par \ref{k=ospec}. Donc aussi $(X\vert\Delta)$, puisque $\Delta\leq \Delta'$$\square$

\

Le m\^eme argument fournit l'exemple suivant (sugg\'er\'e par une question de M. Musta\c ta):

\begin{example}\label{torspec} Soit $X$ une vari\'et\'e torique de fibr\'e anticanonique $D$, et $\Delta\leq D$. Alors $(X\vert\Delta)$ est sp\'eciale. En effet: le diviseur anticanonique, compl\'ementaire de l'orbite ouverte, est \`a croisements normaux. 
\end{example}

Un cas particulier utilis\'e dans la suite est le suivant:

\begin{example}\label{projspec} Soit $(Y\vert\Delta):=(\bP^r\vert D_r)$ l'orbifolde g\'eom\'etrique (logarithmique) lisse de l'exemple \ref{exproj'}, obtenue de $\bP^r$ en munissant les $(r+1)$ hyperplans de coordonn\'ees de la multiplicit\'e $+\infty$. Cette orbifolde g\'eom\'etrique (torique) est donc sp\'eciale.\end{example}

\subsection{Compos\'ees de fibrations de type g\'en\'eral}\label{comptg}

Le r\'esultat suivant est essentiel pour \'etablir les propri\'et\'es du ``coeur", au \S9 ci-dessous.

\begin{theorem}\label{cftg} Soient $f:(Z\vert \Delta_Z)\to Y$ et $g:Y\to X$ des fibrations, avec $(Z\vert \Delta_Z)$ lisse, et $Z\in \sC$. Soit $g_x:(Z_x\vert \Delta_{Z_x})\to Y_x$ la restriction de $g$ au-dessus de $x\in X$, g\'en\'eral. On suppose $f\circ g:(Z\vert \Delta_Z)\to X$ de type g\'en\'eral. Alors: $\kappa(g\vert \Delta_Z)=\kappa(g_x\vert \Delta_{Z_x})+dim(X)$. 

En particulier: si $(g_x:(Z_x\vert \Delta_{Z_x})\to Y_x$ est de type g\'en\'eral, alors $g$ est de type g\'en\'eral.
\end{theorem}

{\bf D\'emonstration:} Les invariants en jeu sont bim\'eromorphes. On peut donc supposer que $g, f,$ et $gf$ sont nettes et hautes. Les dimensions canoniques de ces trois fibrations sont donc celles de leurs bases orbifoldes, et $\Delta(gf,\Delta_Z)=\Delta(g,\Delta(f,\Delta_Z))$. On pose: $\Delta_Y:=\Delta(f\vert \Delta_Z)$. Le th\'eor\`eme \ref{addorb} montre alors la seconde des \'egalit\'es suivantes: $\kappa(g\vert \Delta_Z)=\kappa(Y\vert \Delta_Y)= \kappa(Y_x\vert \Delta_{Y_x})+dim(X)=\kappa(g_x\vert\Delta_{Z_x})+dim(X)$ $\square$

\subsection{Le quotient $\kappa$-rationnel (conditionnel)}\label{krat}

{\bf On suppose dans tout ce \S\ref{krat} que  $C_{n,m}^{orb}$ est vraie. (Voir \ref{cnmorb}).}

\begin{lemma}

Soit $(X\vert\Delta)$ une orbifolde g\'eom\'etrique lisse, avec $X\in \sC$. 

Soit $f:(X\vert\Delta)\dasharrow Y$ et $g:(X\vert\Delta)\dasharrow Z$ deux fibrations avec $dim(Y)>0$ et $dim(Z)>0$, telles que $\kappa(f\vert\Delta)\geq 0$ et $\kappa(g\vert\Delta)\geq 0$.

Il existe alors une fibration $h:(X\vert\Delta)\dasharrow V$ telle que $\kappa(h\vert\Delta)\geq 0$ qui domine $f$ et $g$ (ie: il existe $u:V\dasharrow Y$ et $v:V\dasharrow Z$ telles que $v\circ h=g$ et $u\circ h=f)$.
\end{lemma}

{\bf D\'emonstration:} On peut supposer $f$ et $g$ holomorphes, nettes, avec bases orbifoldes lisses. Soit $W\subset Y\times Z$ l'image du morphisme produit $k:f\times g:X\to Y\times Z$ d\'efini par $k(x)=(f(x),g(x))$.  On note $u':W\to Y$ et $v':W\to Z$ les projections naturelles, telles que $u'\circ k=f$, et $v'\circ k=g$. 

Soit $h:X\to V$ la factorisation  de Stein de $k:X\to W$. Observons que les projections naturelles $u:V\to Y$ et $v:V\to Z$ sont bien des fibrations, puisque $V_y=h(X_y)$ et $V_z=h(X_z)$ sont connexes, pour tous $y,z\in Y,Z$.

On peut encore, quitte \`a modifier encore $X,Y$ et $Z$, supposer la fibration $h:(X\vert\Delta)\to V$ nette et \`a base orbifolde lisse, et supposer aussi (par \ref{compfib}) que: $\Delta(f,\Delta)=\Delta(u,\Delta(h,\Delta))$. 

La famille $g(X_y)=(Z_y)_{y\in Y}$ forme une famille couvrante de sous-vari\'et\'es de $Z$. Notons $\Delta_Z:=\Delta(g,\Delta)$. L'orbifolde g\'eom\'etrique stable  $[(Z\vert\Delta_Z)_{Z_y}]$ induite par restriction est donc bien d\'efinie, et est telle que $\kappa([(Z\vert\Delta_Z)_{Z_y}])\geq 0,$ puisque $\kappa(g\vert\Delta)=\kappa(Z\vert\Delta(g,\Delta))=\kappa(Z\vert\Delta_Z)\geq 0$, par hypoth\`ese.

La fibration $v:(V\vert\Delta(h,\Delta))\to Y$ a pour fibre orbifolde g\'en\'erique $(V\vert\Delta(h,\Delta))_{V_y}$, qui est g\'en\'eriqement finie sur $[(Z/(\Delta_Z)_{V_y}]$. Donc $\kappa((V\vert\Delta(h,\Delta))_{V_y}\geq 0$. Par hypoth\`ese $\kappa(Z\vert\Delta(v,\Delta(h,\Delta)))=\kappa(Z\vert\Delta(v\circ h,\Delta))=\kappa(Z\vert\Delta_Z)\geq 0$.

Appliquant $C_{n,m}^{orb}$ (suppos\'ee vraie) \`a la fibration $v:(V\vert\Delta_V)\to Y$, on obtient donc: 
$\kappa(V\vert\Delta_V)\geq 0$ $\square$

\begin{corollary}\label{k-qr} {\bf On suppose que  $C_{n,m}^{orb}$ est vraie. (Voir \ref{cnmorb}).}

Soit $(X\vert\Delta_X)$ une orbifolde g\'eom\'etrique lisse avec $X\in \sC$, et $dim(X)>0$. 

Il existe une unique fibration $r^+_{X\vert\Delta}:(X\vert\Delta)\dasharrow R^+(X\vert\Delta)$ telle que:

1. $\kappa([R^+(X\vert\Delta)\vert\Delta(r^+_{X\vert\Delta},\Delta)])\geq 0$.

2. $\kappa_+(X\vert\Delta)_r=-\infty$, si $[(X\vert\Delta)_r]$ est la fibre orbifolde stable g\'en\'erique de $r^+_{X\vert\Delta_X}$.

Cette fibration, bien d\'efinie \`a \'equivalence bim\'eromorphe pr\`es sur $(X\vert\Delta)$, est appel\'ee le {\bf $\kappa$-quotient rationnel (conditionnel)} de $(X\vert\Delta)$. 

\end{corollary}

{\bf D\'emonstration:} L'unicit\'e est claire, puisque $\kappa(h\vert\Delta)=-\infty$, pour toute fibration $h:(X\vert\Delta)\dasharrow T$ telle que $dim(h(X_r))>0$ pour $r\in R^+(X\vert\Delta)$ g\'en\'erique.

Etablissons l'existence. Si $\kappa_+(X\vert\Delta)=-\infty$, on choisit pour $r^+_{X\vert\Delta}$ l'application constante sur un point. Sinon, on choisit une fibration $f:(X\vert\Delta)\dasharrow Z$ avec $dim(Z)>0$ maximum, telle que $\kappa(f\vert\Delta)\geq 0$. Le lemme pr\'ec\'edent montre que $f$ domine toute fibration $g:(X\vert\Delta)\dasharrow T$ telle que $\kappa(g\vert\Delta)\geq 0$.

Il reste \`a montrer que $\kappa_+([(X\vert\Delta)_z])=-\infty$, pour $z\in Z$ g\'en\'erique.

Sinon, par r\'ecurrence sur $dim(X)>0$, toute fibre g\'en\'erale $[(X\vert\Delta)_z]$ de $f$ admet un (unique) $\kappa$-quotient rationnel (conditionnel) $r^+_z:(X\vert\Delta)_z\dasharrow R^+_z$, holomorphe sur un mod\`ele bim\'eromorphe ad\'equat de $f$, et dont la famille des fibres forme donc une composante irr\'eductible de $\sC$$how(X_z)$. Il r\'esulte alors (par des arguments exactement similaires) de la d\'emonstration de \ref{gtgen} qu'il existe des  fibrations $r^+:(X\vert\Delta)\dasharrow V$ et $s:V\dasharrow Z$ telles que:

1. $s\circ r=f:(X\vert\Delta)\dasharrow Z$.

2. $r^+_{\vert X_z}=r^+_z$, pour $z\in Z$ g\'en\'eral.

En particulier, \`a la fois la base orbifolde stable et la fibre orbifolde stable de $s:(V\vert\Delta(r^+,\Delta))\dasharrow Z$ ont une dimension canonique positive ou nulle. Il r\'esulte alors de $C_{n,m}^{orb}$ que $\kappa(V\vert\Delta(r^+,\Delta))\geq 0$. Ce qui contredit la maximalit\'e de $dim(Z)$, puisque $dim(V)>dim(Z)$, par construction. Les fibres orbifoldes g\'en\'eriques de $f$ ont donc bien $\kappa_+=-\infty$ $\square$

\begin{re}\label{r=M=Y} 

\

1. On a donc: $R^+(X\vert\Delta)=(X\vert\Delta)$ si et seulement si $\kappa(X\vert\Delta)\geq 0$, et dans ce cas, $M(X\vert\Delta)=(X\vert\Delta)$ si et seulement si $\kappa(X\vert\Delta)=dim(X)\geq 0$. 

2. Il peut se faire que $r^+_{(X\vert\Delta_X)}$ ne soit pas presque-holomorphe lorsque $\Delta\neq 0$. Par exemple si $X=\bP^2$, et si $\Delta$ est l'orbifolde g\'eom\'etrique logarithmique (multiplicit\'es $+\infty)$ dont le support consiste en deux droites (projectives, distinctes). Ceci est cependant peut-\^etre particulier aux orbifoldes g\'eom\'etriques logarithmiques.
\end{re}

Ici encore, $r^+_{(X\vert\Delta_X)}$ jouit de la propri\'et\'e de fonctorialit\'e suivante:

\begin{lemma} \label{fonctr} Soit $f:(X\vert\Delta_X)\dasharrow (Y\vert\Delta_Y)$ est un morphisme dans la cat\'egorie m\'eromorphe des orbifoldes g\'eom\'etriques lisses. On suppose que $X\in \sC$.

Notons (pour simplifier les notations) $[R^+_X\vert\Delta_{R^+_X}]$ et $[R^+_Y\vert\Delta_{R^+_Y}]$ les bases orbifoldes stables de $r^+_{(X\vert\Delta_X)}$ et de $r^+_{(Y\vert\Delta_Y)}$ respectivement.

Il existe alors un (unique) morphisme $r^+_f:[R^+_X\vert\Delta_{R^+_X}]\dasharrow [R^+_Y\vert\Delta_{R^+_Y}]$ tel que:
$r^+_f\circ r^+_{(X\vert\Delta_X)}=r^+_{(Y\vert\Delta_Y)}\circ f$.
\end{lemma}

{\bf D\'emonstration:} Notons $F_X$ la fibre orbifolde g\'en\'erique de $r^+_{(X\vert\Delta_X)}$, $R^+_Y$ la base orbifolde stable de $r^+_{(Y\vert\Delta_Y)}$, et $r^+_Y:=r^+_{(Y\vert\Delta_Y)}$. Alors $r^+_Y:F_X\dasharrow R^+_Y$ d\'efinit une famille couvrante de $R^+_Y$. Il s'agit de montrer que $dim(F_X)=0$. Supposons le contraire. L' orbifolde g\'eom\'etrique obtenue par restriction de $[\Delta(r^+_Y,\Delta_X)]$ \`a $r^+_Y(F_X)$ a donc $\kappa=-\infty$, puisque quotient de $F_X$ telle que $\kappa_+(F_X)=-\infty$. Puisque cette famille est couvrante, ceci contredit $\kappa(R^+_Y\vert\Delta(r^+_Y,\Delta_Y))\geq 0$ $\square$

\newpage

  \section{ORBIFOLDES SP\'ECIALES I}\label{orbspec}

\subsection{Fibre et base orbifolde stables d'une fibration}

{\bf Rappels.} Dans cette section, on consid\`erera uniquement des orbifoldes g\'eom\'etriques $(X\vert\Delta)$ lisses, avec $X\in \sC$, lisse et connexe.

Si $f:(X\vert\Delta)\dasharrow Y$ est m\'eromorphe surjective, on d\'efinira sa {\bf base orbifolde stable}, not\'ee $[Y\vert\Delta(f,\Delta)]$ ou $[f\vert\Delta]$, obtenue comme base orbifolde d'un repr\'esentant holomorphe {\it net} arbitraire de $f$. Nous {\it ne savons pas} si la classe d'\'equivalence bim\'eromorphe (au sens orbifolde) de $[Y\vert\Delta(f,\Delta)]= [f\vert\Delta]$ est ind\'ependante du mod\`ele choisi. Mais, par \ref{dimcan}, sa dimension canonique $\kappa([f\vert\Delta]):=\kappa(f\vert \Delta)$ est bien d\'efinie. C'est donc sur cet unique invariant que sont bas\'ees toutes les consid\'erations qui suivent.

On dit (definition \ref{deftg}) que $f$ est de {\bf type g\'en\'eral} si sa base orbifolde stable est de type g\'en\'eral et de dimension strictement positive.

Si $f:(X\vert\Delta)\dasharrow Y$ est une fibration, on d\'esignera par $(X\vert\Delta)_y$ sa {\bf fibre orbifolde g\'en\'erique stable} , d\'efinie, par modification de $X$, sur un mod\`ele holomorphe de $f$ (\'egalement not\'e $f)$, dont la fibre orbifolde g\'en\'erique $(X_y\vert\Delta_y)$ est lisse (par le th\'eor\`eme de Sard). Voir \ref{ibrest}. Remarquons que la classe d'\'equivalence bim\'eromorphe de $(X\vert\Delta)_y$ g\'en\'erique ne d\'epend pas du mod\`ele bim\'eromorphe de $f$ choisi, rendant $f$ holomorphe. Mais elle peut d\'ependre du repr\'esentant bim\'eromorphe de $(X\vert \Delta)$. Voir cependant \ref{ph} et \ref{phftg}  ci-dessous: si $f$ est presque-holomorphe (d\'efinition \ref{ph}), alors $(X\vert\Delta)_y$ ne d\'epend pas du mod\`ele choisi, et $f$ est presque-holomorphe si elle est de type g\'en\'eral ou plus g\'en\'eralement, si $(X\vert\Delta)$ est finie, et si $\kappa(f\vert \Delta)\geq 0$. Si $f$ est la fibration de Moishezon-Iitaka, la dimension de Kodaira de $(X\vert \Delta)_y$ ne d\'epend pas du mod\`ele bim\'eromorphe de $(X\vert \Delta)$.

Lorsque $f$ est holomorphe, on n'a pas toujours: $(X\vert\Delta)_y=(X_y\vert\Delta_y)$, $\Delta_y$ \'etant la restriction de $\Delta$ \`a $X_y$. Ce sera cependant le cas pour les fibrations que nous consid\`ererons ici (voir la remarque \ref{rrfo}).

$\bullet$ Nous montrerons en, \ref{specdomtg}, si $f:(X\vert\Delta)\dasharrow Y$ est une fibration de type g\'en\'eral, et si $g:(X\vert\Delta)\dasharrow T$ est une fibration dont la fibre g\'en\'erique (orbifolde stable) $(X\/\vert\Delta)_t$ est sp\'eciale, il existe une (unique) factorisation $h:T\dasharrow Y$ de $f$ telle que $f=h\circ g$.

\begin{definition} L'orbifolde g\'eom\'etrique lisse $(X\vert\Delta)$ est dite {\bf sp\'eciale} si:

1. $X\in \sC$.

2. Aucune fibration $f:(X\vert\Delta)\dasharrow Y$ n'est de type g\'en\'eral.

(En particulier, $X\in \sC$ est dite sp\'eciale si $(X\vert 0)$ est sp\'eciale).

Cette propri\'et\'e est pr\'eserv\'ee par \'equivalence bim\'eromorphe.
\end{definition}

\subsection{Premiers exemples}

Les exemples et  contre-exemples les plus simples sont les suivants:

\begin{example}\label{exspec}

\

{\bf 1.} Si $(X\vert\Delta)$ est lisse de type g\'en\'eral, de dimension non nulle, alors $(X\vert\Delta)$ n'est pas sp\'eciale.

{\bf 2.} Si $X$ est une courbe lisse projective de genre $g$, et si $\Delta=\sum_{j} (1-\frac{1}{m_j}). {p_j},$ avec  $m_j>1, \forall j\in J$, les $p_j$ \'etant des points distincts de $X$, alors $(X\vert\Delta)$ est sp\'eciale si et seulement si $2(g-1)+\sum_j(1-\frac{1}{m_j})\leq 0$. 

C'est le cas si et seulement si $g=1$ et $\Delta=0$, ou $g=0$, et si la suite ordonn\'ee croissante des $m_j,j\in J$ est l'une des suivantes:

$\vert J\vert\leq 2:$ quelconque.

$\vert J\vert=3:$ $(2,2,m),\forall m\leq +\infty$$; (2,3,4);(2,3,5);(2,3,6); (2,4,4).$

$\vert J\vert=4:$  $(2,2,2,2)$.

{\bf 3.} Si $f:(X\vert\Delta)\dasharrow (X'\vert\Delta')$ est m\'eromorphe surjective (ie: telle que $f$ ait un mod\`ele bim\'eromorphe lisse qui est un morphisme orbifolde sujectif), et si $(X\vert\Delta)$ est sp\'eciale, alors $(X'\vert\Delta')$ est sp\'eciale, par \ref{composgt}.

{\bf 4.} Si $(X\vert\Delta)$ est une orbifolde g\'eom\'etrique lisse, si $X\in \sC$, et si $\kappa(X\vert\Delta)=0$, alors $(X\vert\Delta)$ est sp\'eciale. (Ceci r\'esulte du th\'eor\`eme \ref{k=ospec}).

{\bf 5.} Si $(X\vert\Delta)$ est lisse et Fano, alors $(X\vert\Delta)$ est sp\'eciale (par le corollaire \ref{fanorb}). Plus g\'en\'eralement (du moins conjecturalement):

{\bf 6.} Si $\kappa_+(X\vert\Delta)=-\infty$ (voir d\'efinition \ref{k-rc}), et en particulier si $(X\vert\Delta)$ est $RE$ (voir d\'efinition \ref{u-rc}), alors $(X\vert\Delta)$ est sp\'eciale. C'est \'evident, par d\'efinition.

{\bf 7.} Si $(X\vert\Delta)$ est $\cal S$-connexe (voir d\'efinition en \ref{s-conn}) , alors $(X\vert\Delta)$ est sp\'eciale (par le corollaire \ref{s-conn} du chapitre suivant).

{\bf 8.} Pour tout $n\geq 0$ et tout $k\in\{-\infty,0,1,\dots,n-1\}$, il existe des orbifoldes g\'eom\'etriques sp\'eciales de dimension $n$ avec $\kappa=k$ (Voir [Ca04] pour des exemples simples).

{\bf 9.} Si $X$ est une vari\'et\'e quasi-projective lisse, et $\overline{X}$ une compactification lisse de $X$ telle que $\overline{X}-X:=D$ soit un diviseur \`a croisements normaux de $\overline{X}$, on dira que $X$ est sp\'eciale si et seulement si $(\overline{X}\vert D)$ est sp\'eciale. Cette condition ne d\'epend pas de la compactification choisie.

{\bf 10.} $X$ est sp\'eciale s'il existe une application m\'eromorphe non-d\'eg\'en\'er\'ee $\varphi:\bC^n\dasharrow X$ (voir [Ca 04], qui traite d'une situation plus g\'en\'erale). De mani\`ere similaire, s'il existe un morphisme orbifolde $\varphi:\bC^n\to (X\vert \Delta)$, alors $(X\vert \Delta)$ n'est pas de type g\'en\'eral ([Sak 74]. Je remercie E. Rousseau pour m'avoir indiqu\'e cette r\'ef\'erence).

\end{example}

\subsection{Composition de fibrations sp\'eciales}

Le r\'esultat simple suivant est crucial:

\begin{theorem}\label{compspec} Soit $f:(X\vert\Delta)\dasharrow Y$ une fibration m\'eromorphe surjective telle que:

1. $X\in \sC$

2. La base orbifolde stable $[Y\vert\Delta(f,\Delta)]$ de $f$ est sp\'eciale.

3. La fibre g\'en\'erale orbifolde stable $(X\vert\Delta)_y$ de $f$ est sp\'eciale.

Alors: $(X\vert\Delta)$ est sp\'eciale.
\end{theorem}

{\bf D\'emonstration:} On supposer que $f$ est nette, \`a base orbifolde lisse sp\'eciale, et \`a fibres orbifoldes g\'en\'erales sp\'eciales. Soit $g:(X\vert\Delta)\dasharrow T$ une fibration de type g\'en\'eral. Il r\'esulte de \ref{specdomtg} que $g=h\circ f$ pour une fibration $h:Y\dasharrow T$. Quitte \`a changer de mod\`eles bim\'eromorphes pour $f,g,h$, on peut supposer par \ref{compfib} que $\Delta_T(g,\Delta)=\Delta_T(h,\Delta_Y(f,\Delta))$ est la base orbifolde stable de $h:(Y\vert\Delta(f,\Delta))\to T$. Puisque $g$ est de type g\'en\'eral, il en est donc de m\^eme pour $h$. Mais ceci contredit le fait que $[(Y\vert\Delta(f,\Delta)]$ est sp\'eciale. Donc $(X\vert\Delta)$ est sp\'eciale $\square$

\begin{re} Ce r\'esultat est \'evidemment faux (surfaces elliptiques de base $\bP^1$ ayant au moins $5$ fibres multiples, par exemple) sans l'introduction de structures orbifoldes g\'eom\'etriques. Ce r\'esultat simple est l'un de ceux qui justifient la n\'ecessit\'e de travailler dans la cat\'egorie des orbifoldes g\'eom\'etriques, plut\^ot que dans celle des vari\'et\'es. Un autre r\'esultat fournissant une conclusion similaire est celui portant sur les groupes fondamentaux des fibrations nettes: le $\pi_1$ orbifolde de l'espace total est extension de celui de la base orbifolde par (un quotient de) celui de la fibre (orbifolde).
\end{re}

Par it\'eration, on obtient, \`a l'aide de \ref{iternet}:

\begin{corollary}\label{specdec} Soit $f_j:(X_j\vert\Delta_j)\dasharrow X_{j+1}, j=0,\dots,k-1,$ une suite de fibrations m\'eromorphes telle que:

1. $X_0\in \sC$.

2. $X_k$ est un point.

3. Pour tout $j=0,\dots,k-1$, $[(X_{j+1} \vert\Delta_{j+1})]=[(X_{j+1}\vert\Delta(f_j,\Delta_j)]$.

4. Pour tout $j=0,\dots, k-1$, la fibre {\bf orbifolde} stable g\'en\'erique $F_j$ de $f_j$ est telle que: ou bien $\kappa(F_j)=0$, ou bien $\kappa_+(F_j)=-\infty$.

Alors: $(X_0\vert\Delta_0)$ est sp\'eciale.
\end{corollary}

Nous verrons au \S\ref{deccoeur} que, r\'eciproquement, sous r\'eserve de la validit\'e de la conjecture $C_{n,m}^{orb}$, on peut canoniquement d\'ecomposer les orbifoldes g\'eom\'etriques sp\'eciales en tours de fibrations \`a fibres orbifoldes g\'en\'eriques ayant soit $\kappa=0$, soit $\kappa_+=-\infty$.

\subsection{Orbifoldes sp\'eciales ``divisibles".}

Nous avons consid\'er\'e dans les sections pr\'ec\'edentes la notion d'orbifolde sp\'eciale dans la cat\'egorie $Georb^Q$. Lorsque l'orbifolde consid\'er\'ee est enti\`ere, on peut donner les m\^emes d\'efinitions, mais dans $Georb^{div}$. On dira alors que $(X\vert \Delta)$ est $div$-sp\'eciale, ou sp\'eciale$^{div}$ si elle est lisse, $X\in \sC$, et n'admet pas de fibration de tpe g\'en\'eral dans $Georb^{div}$. 

On a l'implication \'evidente: $(X\vert \Delta)$ sp\'eciale $\Longrightarrow (X\vert \Delta)$ sp\'eciale$^{div}$, puisqu'une fibration de type g\'en\'eral dans $Georb^{div}$ l'est aussi dans $Georb^Q$.

\begin{question}\label{qspec} Est-il vrai que: $(X\vert \Delta)$ sp\'eciale $^{div}\Longrightarrow (X\vert \Delta)$ sp\'eciale?
\end{question}

Cette implication r\'eciproque r\'esulterait d'une r\'eponse affirmative \`a la question \ref{qmult}. (Appliquer le ``coeur" (voir \ref{c}) \`a $(X\vert \Delta)$, suppos\'ee \^etre sp\'eciale$^{div}$: la base orbifolde est de type g\'en\'eral, mais coincide avec la base orbifode dans $Georb^{div}$ si la r\'eponse \`a \ref{qmult} est affirmative. Cette base orbifolde est donc un point, par hypoth\`ese).

Le d\'evissage conditionnel de \ref{sdeviss} reste {\it a fortiori} valable dans la cat\'egorie $Georb^{div}$ (puisque la sur-additivit\'e des dimensions de Kodaira y est-{\it a fortiori }-valable (puisque les bases orbifoldes divisibles divisent les bases orbifoldes calcul\'ees dans $Georb^Q$.  Pour \'etablir l'implication r\'eciproque, il suffit donc (sous r\'eserve de $C_{n,m}^{orb})$, de montrer (puisque les orbifodes avec $\kappa=0$ sont sp\'eciales et en utilisant \ref{compspec}) de montrer qu'une orbifolde telle que $\kappa_+^{div}=-\infty$   est sp\'eciale (dans $Georb^Q)$. On d\'efinit \'evidemment la condition $\kappa_+^{div}=-\infty$ en imposant \`a toute base orbifolde (calcul\'ee dans $Georb^{div})$ de toute fibration d'avoir $\kappa=-\infty$.

\newpage

 \section{FIBRATIONS DE TYPE G\'EN\'ERAL.} \label{fgt}

  Nous \'etablissons ici un certain nombre de propri\'et\'es de ces fibrations, utilis\'ees dans la construction du ``coeur", au \S\ref{coeur}

  \begin{notation}\label{notat} D\'esormais\footnote{C'est-\`a-dire: dans le pr\'esent \S\ref{fgt}.}, $(Y\vert\Delta)$ d\'esignera une orbifolde g\'eom\'etrique lisse avec $Y$ compacte et connexe.

On dira que $Y\in \cal C$ si $Y$, compact, normal et connexe, est bim\'eromorphe \`a une vari\'et\'e K\" ahl\'erienne compacte.

On notera $f:(Y\vert\Delta)\dasharrow X$ une application m\'eromorphe surjective.\end{notation}

   \subsection{Fibrations de type g\'en\'eral: compos\'ees} \label{fgt1}

  \
  
   \begin{definition} Soit $f:(Y\vert\Delta)\dasharrow X$ une fibration m\'eromorphe avec $X,Y$ compacts et ir\'eductibles, et $(Y\vert\Delta)$ une orbifolde g\'eom\'etrique lisse. On dit que $(f\vert\Delta)$ (ou $f$ s'il n' y a pas d'ambiguit\'e sur $\Delta)$est une {\bf fibration de type g\'en\'eral} si $\kappa(f\vert\Delta)=dim(Y)>0$.

   Si $(Y\vert\Delta)$ est une orbifolde g\'eom\'etrique lisse, avec $Y$ compacte et connexe. On note $TG(Y\vert\Delta)$ l'ensemble des classes d'\'equivalence bim\'eromorphe de fibrations de type g\'en\'eral sur $(Y\vert\Delta)$. Si $f$ est une telle fibration on notera $[f]$ sa classe d'\'equivalence dans $TG(Y\vert\Delta)$.

   Cet ensemble ne d\'epend donc que de la classe d'\'equivalence bim\'eromorphe de $(Y\vert\Delta)$.
 \end{definition}

\begin{question}\label{fintg} $TG(Y\vert \Delta)$ est-il fini si $Y\in \sC$?
\end{question}

\begin{proposition} \label{composgt}Soit $g:(Z\vert\Delta_Z)\to (Y\vert\Delta_Y)$ un morphisme surjectif  d'orbifoldes g\'eom\'etriques lisses ($Y$ et $Z$ \'etant compactes et connexes) , et $f:(Y\vert\Delta_Y)\to X$ holomorphe surjective. Si $f$ est de type g\'en\'eral, la compos\'ee $f\circ g$ l'est aussi.
La propri\'et\'e subsiste lorsque $f$ et $g$ sont m\'eromorphes.
  \end{proposition}

  {\bf D\'emonstration:} Si $f$ est de type g\'en\'eral, $\kappa(Y,L_f)=p=dim(X)>0$. Or $L_{f\circ g}=g^*(L_f)$. Donc $\kappa(Z, L_{f\circ g})=\kappa(Z, g^*(L_f))=\kappa(Y,L_f)$. D'o\`u la premi\`ere assertion. La seconde assertion se d\'emontre de la m\^eme mani\`ere, puisque la dimension canonique des faisceaux $L_f$ est un invariant bim\'eromorphe des orbifoldes g\'eom\'etriques lisses $\square$
  
\begin{re} Le th\'eor\`eme \ref{cftg} est une r\'eciproque partielle de cette proposition.
\end{re}

Rappelons la notion d'orbifolde (g\'eom\'etrique lisse) sp\'eciale:

  \begin{definition} Une orbifolde g\'eom\'etrique lisse $(Y\vert\Delta)$, avec $Y$ compacte et connexe est dite {\bf sp\'eciale} si:
 
 1. $Y\in \sC$ (ie: $Y$ est bim\'eromorphe \`a une vari\'et\'e K\" ahl\'erienne compacte). 
 
 2. Il n'existe pas de fibration $f:(Y\vert\Delta)\dasharrow X$ de type g\'en\'eral.  \end{definition}

 \begin{corollary} Si $u:(X\vert\Delta)\dasharrow (X'\vert\Delta')$ est m\'eromorphe surjective entre orbifoldes g\'eom\'etriques lisses, et si $(X\vert\Delta)$ est sp\'eciale, alors $(X'\vert\Delta')$ est aussi sp\'eciale.
  \end{corollary}

\subsection{Faisceaux de Bogomolov}\label{boggt}

\begin{definition} Soit $(Y\vert\Delta)$ une orbifolde g\'eom\'etrique lisse, et $L\subset \Omega_X^p,$ avec $p>0$, un sous-faiscau coh\'erent de rang $1$. Rappelons que, d'apr\`es \ref{bogorb}, $\kappa(Y\vert\Delta, L)\leq p$.

On dit que $L$ est un {\bf Faisceau de Bogomolov} si $\kappa(Y\vert\Delta, L)=p$.

Deux tels faisceaux sont {\bf \'equivalents} s'ils coincident sur un ouvert non vide de $X$. On note $Bog(Y\vert\Delta)$ l'ensemble (\'eventuellement vide) des classes d'\'equivalence de tels faisceaux. On note $[L]$ la classe d'\'equivalence de $L$.
\end{definition}

Si $f:(Y\vert\Delta)\dasharrow X$ est une fibration m\'eromorphe, de faisceau $L_f:=f^*(K_X)$ au point g\'en\'erique de $Y$, alors le corollaire \ref{kappa=L} montre que $L_f$ est un faisceau de Bogomolov si et seulement si $f$ est de type g\'en\'eral. D'o\`u une application naturelle (clairement injective): $F:=F_{(Y\vert\Delta)}:TG(Y\vert\Delta)\to Bog(Y\vert\Delta)$.

Il r\'esulte du corollaire \ref{bogorb} que cette application est surjective lorsque $Y$ est bim\'eromorphe \`a une vari\'et\'e K\" ahl\'erienne compacte (ie: dans la classe $\sC$). D'o\`u le:

\begin{theorem}\label{bijbog} Si $(Y\vert\Delta)$ est une orbifolde g\'eom\'etrique lisse, avec $Y$ dans la classe $\sC$,  l'application $F_{(Y\vert\Delta)}:TG(Y\vert\Delta)\to Bog(Y\vert\Delta)$ pr\'ec\'edente est bijective. 

En particulier, $(Y\vert \Delta)$ est sp\'eciale si et seulement si $Bog(Y\vert\Delta)$ est vide.
\end{theorem}

\subsection{Restriction \`a une sous-vari\'et\'e g\'en\'erique.}\label{restrssv}

On \'etend au cas des orbifoldes g\'eom\'etriques g\'eom\'etriques certains des r\'esultats de [Ca 04].

Pour les notions de restriction au sens orbifolde, voir le \S\ref{catorb}, ainsi que \ref{rfo} et \ref{rrfo}.

\begin{proposition}\label{proprestr} Soit $f:(Y\vert\Delta)\dasharrow X$ une application m\'eromorphe de type g\'en\'eral, avec $(Y\vert\Delta)$ lisse, $Y$ compacte et connexe. Soit $V\subset Y$ une sous-vari\'et\'e non contenue dans le support de $\Delta$, et non contenue dans le lieu d'ind\'etermination de $f$. Si $f(V)=X$, alors $f_V:(V\vert\Delta_V)\dasharrow X$ est aussi de type g\'en\'eral, pour tout choix de la restriction $f_V$.
\end{proposition}

{\bf D\'emonstration:} Soit $p:=dim(X)>0$. Soit $(V'\vert\Delta_{V'})$ une restriction arbitraire de $\Delta$ \`a $V$. Voir \S\ref{catorb} pour cette notion. On obtient $L_{f_V}\subset \Omega^p_{V'}$ en composant l'inclusion $L_f:=f^*(K_X)\subset \Omega^p_Y$ avec la restriction naturelle $ \Omega^p_Y\to  \Omega^p_{V'}$. 

On en d\'eduit l'in\'egalit\'e: $p:=\kappa(f\vert\Delta)\leq \kappa(f_V\vert\Delta_V)$, et l'assertion $\square$

 \

Rappelons la:
 
 \begin{definition}\label{defgen} Un point $t\in T$ est dit {\bf g\'en\'eral} s'il appartient \`a un sous-ensemble de $T$ contenant une intersection d\'enombrable d'ouverts de Zariski non vides. Une fibre $Y_t$ de $h$ est dite g\'en\'erale si $t\in T$ est g\'en\'eral. 
 \end{definition}

  \begin{theorem}\label{thmrestr} Soit $(Y\vert\Delta)$ une orbifolde g\'eom\'etrique lisse, $f:(Y\vert\Delta)\dasharrow X$ une fibration m\'eromorphe, et $j':(V'\vert \Delta_{V'})\to (X\vert \Delta)$ la restriction (au sens de \ref{fdefrestr}) de $\Delta$ \`a une sous-vari\'et\'e $V$ de $X$ passant par un point g\'en\'eral de $X$. Si $f$ est de type g\'en\'eral, alors $f\circ j':(V'\vert\Delta_{V'})\to W$ l'est aussi, si  $dim(W)>0$, avec $W:=f(V)$.
  
  En particulier, $dim(W)=0$ si $(V'\vert \Delta_{V'})$ est sp\'eciale.
    \end{theorem}

    {\bf D\'emonstration:} On peut, par le th\'eor\`eme \ref{rescouim}, remplacer $(Y\vert\Delta)$ par une modification \'el\'ementaire $(Y'\vert\Delta')$, et supposer de plus que:
    
    1.  $f$ est nette et haute,
    
    2. la restriction $(V\vert \Delta_V)$  de $\Delta$ \`a $V$ est bien d\'efinie (au sens de \ref{catorb}, de sorte que $(V\vert \Delta_V)$ est lisse),
    
    3. $W:=f(V)\subset X$ est lisse (et est donc le membre g\'en\'erique d'une famille couvrante de $X)$. 
    
    4. La restriction $g:(V\vert \Delta_V)\to W$ de $f$ \`a $V$ est nette et haute.
    
    5. L'application $\varphi: X\to Z$ d\'efinie par $L_f$ est bim\'eromorphe, et que sa rstriction \`a $W$ est aussi bim\'eromorphe sur son image.

   On a alors: $K_{X\vert W}\leq K_W$. Notons $g: V\to W$ a restriction de $f$, et $j:V\to Y$ l'inclusion. Alors: $dim(W)=\kappa(V\Delta_V, j^*(L_f))\leq \kappa(V\vert \Delta, L_g)=\kappa(g\vert \Delta_V)\leq dim(W)$. Donc $(g\vert \Delta_V)$ est bien de type g\'en\'eral $\square$

 De \ref{thmrestr} on d\'eduit le:

 \begin{corollary}\label{specdomtg} Soit $f:(Y\vert\Delta)\dasharrow X$ une application m\'eromorphe surjective de type g\'en\'eral, l'orbifolde g\'eom\'etrique $(Y\vert\Delta)$ \'etant lisse, $Y$ compacte et connexe. 
 
 Supposons aussi qu'il existe une fibration m\'eromorphe $h:(Y\vert\Delta)\dasharrow T$ telle que les fibres orbifoldes g\'en\'erales de $h$ soient sp\'eciales\footnote{Voir \ref{dfo} pour cette notion.}. Alors  $h$ se factorise par $f$ (ie: il existe $g: T\dasharrow X$ telle que $f=g\circ h)$.
 \end{corollary}

  {\bf D\'emonstration:} Sinon, d'apr\`es le th\'eor\`eme \ref{thmrestr} pr\'ec\'edent, la restriction de $f$ \`a $(Y_t\vert\Delta_{Y_t})$ est de type g\'en\'eral, ce qui est impossible, puisque les fibres orbifoldes g\'en\'erales sont supos\'ees sp\'eciales $\square$

\begin{example}\label{projspec'} Soit $(Y\vert\Delta):=(\bP^r\vert D_r)$ l'orbifolde g\'eom\'etrique (logarithmique) lisse de l'exemple \ref{exproj'}. Alors tous les faisceaux $S^N_q(Y\vert\Delta)$ sont triviaux. Il n'existe donc pas de fibration de type g\'en\'eral $f:(Y\vert\Delta)\dasharrow X$ sur cette orbifolde g\'eom\'etrique, et $TG(Y\vert\Delta)=0$ dans ce cas. Cette orbifolde g\'eom\'etrique est donc sp\'eciale. Une seconde d\'emonstration de cette propri\'et\'e est donn\'ee en \ref{torspec}.
\end{example}

  \`A l'aide de l'exemple \ref{projspec'}, on obtient le r\'esultat suivant, utilis\'e dans le \S\ref{gtph} ci-dessous:

\begin{corollary} Soit $f:(Y\vert\Delta)\dasharrow X$ une application m\'eromorphe surjective de type g\'en\'eral, l'orbifolde g\'eom\'etrique $(Y\vert\Delta)$ \'etant lisse, $Y$ compacte et connexe. Supposons aussi qu'il existe une fibration holomorphe $h:Y\to T$ telle que les fibres g\'en\'eriques $Y_t\cong \bP^r$ de $h$ soient des sous-vari\'et\'es lisses connexes de $Y$ rencontrant transversalement le support de $\Delta$ en $s\leq (r+1)$ hyperplans projectifs. Alors la famille des fibres de $h$ se factorise par $f$.
\end{corollary}

  \subsection{Crit\`eres de presque-holomorphie}\label{gtph}

  \begin{definition}\label{ph} Soit $f:Y\dasharrow X$ une application m\'eromorphe propre et surjective (d\'efinition par passage au graphe $Y_f\subset X\times Y$ de $f$) entre espaces analytiques normaux connexes. Soit $I_f\subset Y$ son lieu d'ind\'etermination (lieu au-dessus duquel le graphe de $f$ a des fibres de dimension strictement positive), et $f(I_f)\subset X$ l'image de $I_f$ par $f$ (image par $f':Y_f\to X$ de l'image r\'eciproque de $I_f$ dans $Y_f)$. 
  
  On dit que $f$ est {\bf presque-holomorphe} si $f(I_f)\subsetneq X$.
  \end{definition}

  L'exemple le plus simple de fibration m\'eromorphe non presque-holomorphe est $f:\bP^2\dasharrow \bP^1$ dont les fibres sont les droites passant par un point $a\in \bP^2$. La famille des fibres de cette fibration  (vue dans $\sC(\bP^2)=Chow(\bP^2))$ se d\'eforme (en fonction de $a)$.
  
  Par contraste, les applications presque-holomorphes jouissent de propri\'et\'es similaires \`a celles des applications holomorphes. Par exemple: si $f:Y\dasharrow X$ est une fibration presque-holomorphe, la famille des ses fibres (famille vue dans la vari\'et\'e de Chow $\sC(Y))$ forme une composante irr\'eductible de $\sC(X)$, par le classique ``lemme de rigidit\'e".
  
  En particulier, l'ensemble des fibrations presque-holomorphes (\`a \'equivalence bim\'eromorphe pr\`es) sur une vari\'et\'e complexe fix\'ee $Y$ (d\'enombrable \`a l'infini) est fini ou d\'enombrable. Remarquer que cette propri\'et\'e tombe en d\'efaut sur l'exemple non presque-holomorphe tr\`es simple donn\'e ci-dessus.
  
  De plus, si $f$ est presque-holomorphe, la restriction $(\Delta_Y)_x$ de l'orbifolde g\'eom\'etrique $\Delta_Y$ \`a sa fibre g\'en\'erique $Y_x$ est bien d\'efinie, et $(Y_x\vert(\Delta_Y)_x)$ est lisse si $(Y\vert\Delta_Y)$ est lisse. La classe d'\'equivalence bim\'eromorphe de $(Y_x\vert(\Delta_Y)_x)$ ne d\'epend que de celle de $(Y\vert\Delta_Y)$.

   \begin{theorem}\label{phftg} Si $(Y\vert\Delta)$ est une orbifolde g\'eom\'etrique lisse, avec $Y$ compacte et connexe, et si $f:(Y\vert\Delta)\dasharrow X$ est une fibration m\'eromorphe de type g\'en\'eral, alors $f$ est presque-holomorphe.
   \end{theorem}

   {\bf D\'emonstration:} Soit $v:Y'\to Y$ une compos\'ee d'\'eclatements de centres lisses telle que $f':=f\circ v:Y'\to X$ soit holomorphe. Si $f$ n'est pas presque holomorphe, l'un des diviseurs exceptionnels $E$ de l'un de ces \'eclatements est tel que:
   
   1.  $f'(E)=X$,
   
   2. $v(E):=T$ est de codimension $r+1\geq 2$ dans $Y$.
   
   3. Les fibres g\'en\'eriques $E_t$ de $E$ sur $T$ ont des images $X_t$ dans $X$ de dimension $d>0$. (Autrement dit: $f'_E:E\to X$ ne se factorise pas par $v)$.
   
   4. On peut supposer que $E$ n'est pas contenu dans le support de $\Delta'$ (ceci par le corollaire \ref{dvert}, puisque $f$, et donc $f'$ sont de type g\'en\'eral).
   
   Puisque $(Y\vert\Delta)$ est lisse, $T$ est contenu dans $s\leq (r+1)$ des composantes (lisses, d'intersections transversales) du support de $\Delta$. Pour $t\in T$ g\'en\'erique, l'orbifolde g\'eom\'etrique $(E_t\vert\Delta_{E_t})$ est donc bim\'eromorphe \`a $(\bP^r\vert\Delta_r)$, o\`u $\Delta_r$ est support\'ee par $s$ hyperplans projectifs en position g\'en\'erale. Donc $(E_t\vert\Delta_{E_t})$ est sp\'eciale. Mais d'apr\`es  \ref{proprestr}, $f_E$ est de type g\'en\'eral, donc d'apr\`es \ref{thmrestr}, $f_{E_t}$ aussi. Mais ceci contredit \ref{projspec'} $\square$

  \begin{re} 
  
  \
  
  1. La condition de lissit\'e de $(Y\vert\Delta)$ est essentielle. Soit $f:(\bP^2\vert\Delta)\dasharrow \bP^1$ la fibration dont les fibres sont les droites passant par un point donn\'e $a\in \bP^2$, et soit  $\Delta=(1-\frac{1}{m}).(D_1+D_2+D_3)$, les $D_j$ \'etant $3$ droites concourantes en $a \in \bP^2$. Cette fibration est de type g\'en\'eral si $m\geq 4$ (puisque sa base orbifolde g\'eom\'etrique l'est: consid\'erer $u^+(\Delta)$ sur l'\'eclat\'e de $\bP^2$ en $a$). Elle n'est \'evidemment pas presque-holomorphe. 
  
  2.  Lorsque $\Delta=0$, la condition de lissit\'e de $Y$ est \'egalement essentielle (Exemple 2.23, p. 534 de [Ca04] du c\^one $Y$ sur une vari\'et\'e de type g\'en\'eral). 
   \end{re}
 
 \
 
 La d\'emonstration du th\'eor\`eme \ref{phftg} pr\'ec\'edent montre aussi (puisque $(\bP^r\vert (1-\frac{1}{m}).\Delta_r)$ est unir\'egl\'ee, pour tout $0<m<+\infty$):
 
    \begin{theorem}\label{phftg'} Si $(Y\vert\Delta)$ est une orbifolde g\'eom\'etrique lisse et finie, avec $Y$ compacte et connexe, et si $f:(Y\vert\Delta)\dasharrow X$ est une fibration m\'eromorphe telle que: $\kappa(f\vert \Delta)\geq 0$, alors $f$ est presque-holomorphe.
   \end{theorem}

\subsection{R\'eduction de type g\'en\'eral simultan\'ee.}

 \begin{definition}\label{tgmax} Soit $g:(Y\vert\Delta)\dasharrow X$ une fibration de type g\'en\'eral, avec $(Y\vert\Delta)$ lisse, $Y$ compacte et connexe. On dit que $g$ est {\bf maximum} si toute autre fibration de type g\'en\'eral  $h:(Y\vert\Delta)\dasharrow V$ est domin\'ee par $g$ (ie: telle qu'existe $k:X\dasharrow V$ avec $h=k\circ g)$.
  \end{definition}

  Une telle fibration de type g\'en\'eral maximum est unique (\`a \'equivalence bim\'eromorphe pr\`es) si elle existe. L'existence (qui est l'un des r\'esultats principaux du pr\'esent texte) sera \'etablie pour les vari\'et\'es de la classe $\sC$ dans la section \ref{coeur}, \`a l'aide du th\'eor\`eme d'additivit\'e orbifolde pour les dimensions canoniques. On utilisera en particulier le:

 \begin{theorem}\label{gtgen} Soit $f:(Y\vert\Delta)\to S$ une application holomorphe surjective, avec $(Y\vert\Delta)$ une orbifolde g\'eom\'etrique lisse, $Y$ \'etant compacte connexe et dans la classe $\sC$. 
 On suppose qu'il existe $T\subset S$, un sous-ensemble (arbitraire\footnote{ie: non n\'ecessairement analytique.}) non contenu dans une r\'eunion d\'enombrable de ferm\'es de Zariski {\bf stricts}\footnote{ie: diff\'erents de $S$. Un tel ensemble $T$ est dit ne pas \^etre {\it analytiquement maigre} dans $S.$} de $S$, tel que pour tout $t\in T$, l'orbifolde g\'eom\'etrique (que l'on peut supposer lisse) $(Y_t\vert\Delta_t):=(Y_t\vert\Delta_{Y_t})$ admette une fibration de type g\'en\'eral $\bar{g_t}:(Y_t\vert\Delta_t)\dasharrow X_t$ maximum.
 
 Il existe alors une unique fibration $g:(Y\vert\Delta)\dasharrow X$ au-dessus de $f$ (ie: telle qu'existe $h:X\to S$ avec $f=h\circ g)$ telle que pour $x\in U$ g\'en\'eral dans $S$, $g_s:(Y_s\vert\Delta_s)\dasharrow X_s$ soit la fibration de type g\'en\'eral maximum de $(Y_s\vert\Delta_s)$. De plus, pour tout $s=t\in T\cap U$, $g_s=\bar{g_t}$
  \end{theorem}

{\bf D\'emonstration:} Le th\'eor\`eme 2.35, p. 538 de [Ca04] montre l'existence d'une fibration $g:(Y\vert\Delta)\dasharrow X$ au-dessus de $S$ telle que $g_t=\bar{g_t}$ pour tout $t\in T'\subset T$, $T'$ n'\'etant pas analytiquement maigre dans $S$. 

Les arguments de la d\'emonstration de [Ca04,thm 2.35] s'appliquent sans changement \`a la situation consid\'er\'ee ici, et plus g\'en\'eralement  au cas de fibrations presque-holomorphes poss\'edant une propri\'et\'e d'unicit\'e sur les fibres $(Y_t)_{t\in T}$ de $f$
$\square$

Le th\'eor\`eme r\'esulte alors de la proposition suivante:

\begin{proposition}\label{gengen} Soit $f:(Y\vert\Delta)\to S$ une application holomorphe surjective, avec $(Y\vert\Delta)$ une orbifolde g\'eom\'etrique lisse, $Y$ \'etant compacte connexe. Soit  $f=h\circ g$ une factorisation de $f$ par une application m\'eromorphe $g:Y\dasharrow X$ et une fibration $h:X\to S$.

On suppose que la restriction $g_t:(Y_t\vert\Delta_t)\dasharrow X_t$  est de type g\'en\'eral pour tout $t\in T$, $T$ non analytiquement maigre dans $S$. Alors, $g_s$ est de type g\'en\'eral pour $s\in S$ g\'en\'eral. \end{proposition}

{\bf D\'emonstration:} Pour tout entier $m>0$, consid\'erons le faisceau analytique coh\'erent  $L:=g^*(K_{X/S})\subset \Omega_{Y/S}^p$,  avec $p:=dim(X/S)=dim(X)-dim(S)$, et la saturation $\overline{L_m}\subset S_{m,p}((Y\vert\Delta)/S)$ de $L^{\otimes m}$ dans le faisceau $S_{m,p}((Y\vert\Delta)/S)$ quotient de $S_{m,p}(Y\vert\Delta)$ par le sous-faisceau des (polyn\^omes en les) formes nulles sur les $p$-vecteurs tangents \`a $Y$ qui sont $f$-verticaux.

Par le th\'eor\`eme d'image directe de Grauert, les faisceaux $F_m:=f_*(\overline{L_m})$ sont coh\'erents, et, pour $s\in U$ g\'en\'eral dans $S$, la fibre $F_{m,s}$ de $F_m$ en $s$ coincide avec $H^0(Y_s, (\overline{L_{f_s}})$$_m)$ pour tous les $m>0$. On a not\'e $L_{g_s}:=(g_s)^*(K_{X_s})\subset \Omega^p_{Y_s}$, et $\overline{L_{g_s}}$$_{,m}\subset S_{m,p}(Y_s\vert\Delta_s)$ la saturation de $L_{g_s}^{\otimes m}$ dans $ S_{m,p}(Y_s\vert\Delta_s)$.

L'assertion r\'esulte alors de ce que si, pour un $s\in U$, le syst\`eme lin\'eaire $\vert\overline{ L_{g_s}}$$_{,m}\vert$ d\'efinit une application $\varphi_s:Y_s\dasharrow Z_s$ de rang $p$, il en est de m\^eme pour tout $s'\in S$ g\'en\'erique, et qu'alors $Z_s=X_s$ (quitte \`a remplacer $m$ par un multiple ad\'equat) $\square$

\newpage

 \section{LE ``COEUR" D'UNE ORBIFOLDE}\label{coeur}

\subsection{Construction du Coeur.}

Le r\'esultat central du pr\'esent texte est le:

 \begin{theorem}\label{c} Soit $(Y\vert\Delta)$ une orbifolde g\'eom\'etrique lisse, avec $Y$ compacte connexe et dans la classe $\sC.$\footnote{ie: bim\'eromorphe \`a $Y'$, K\" ahl\'erienne compacte.} Il existe alors une unique\footnote{\`A \'equivalence bim\'eromorphe pr\`es} fibration $c_{(Y\vert\Delta)}: (Y\vert\Delta)\dasharrow C(Y\vert\Delta)$ telle que:
 
 1. $\kappa(c_{(Y\vert\Delta)}\vert\Delta)=dim(C({Y\vert\Delta}))$ (et $c$ est donc presque-holomorphe).
 
 2. Les fibres g\'en\'erales $(Y_c\vert\Delta_c)$ de $c_{(Y\vert\Delta)}$ sont sp\'eciales.
 
 Cette fibration est presque-holomorphe.
 \end{theorem}

\begin{definition} La fibration $c_{(Y\vert\Delta)}$ est appel\'ee le {\bf coeur}. 
 
 On notera $[C(Y\vert\Delta)]$ une base orbifolde stable (voir \ref{dimcan}) de $c_{(Y\vert\Delta)}$\footnote{Nous ne savons pas si c'est un invariant de $c_{(Y\vert\Delta)}$ dans la cat\'egorie bim\'eromorphe des orbifoldes g\'eom\'etriques lisses. }.

 La dimension de $C(X\vert\Delta)$ est appel\'ee {\bf dimension essentielle} de $(X\vert\Delta)$, et est  not\'ee $ess(X\vert\Delta)$.

 On notera aussi $K(C(Y\vert\Delta)):=A(Y\vert\Delta)$ l'alg\`ebre canonique (voir d\'efinition en \ref{kodL}) de la base orbifolde stable de $c_{(Y\vert\Delta)}$. On l'appelle {\bf l'alg\`ebre essentielle de $(Y\vert\Delta)$.}
  \end{definition}

 \begin{re} 
 
 \
 
 1. La fibration  $c_{(Y\vert\Delta)}$ est donc de type g\'en\'eral si $dim(C{(Y\vert\Delta))}>0$. Donc $dim(C{(Y\vert\Delta))}=0$ si et seulement si $(Y\vert\Delta)$ est sp\'eciale. 
 
 2. A l'autre extr\^eme, $(Y\vert\Delta)$ est de type g\'en\'eral si et seulement si $dim(C{(Y\vert\Delta))}=dim(Y)$, c'est-\`a-dire si et seulement si $c_{(Y\vert\Delta)}=id_Y$.
 
 3. D'apr\`es \ref{specdomtg}, le coeur est donc l'unique fibration de type g\'en\'eral maximum (au sens de \ref{tgmax}) sur $(Y\vert\Delta)$ si  $(Y\vert\Delta)$ n'est pas sp\'eciale.

  4. Si $C_{n,m}^{orb}$ est vraie, alors le coeur d\'efini en \ref{c} coincide avec celui d\'efini en \ref{c'} ci-dessous comme aboutissement de l'it\'eration $(M\circ r)^n, n:=dim(X)$ du $\kappa$-quotient rationnel $r$ et de la fibration de Moishezon-Iitaka $M$.
    \end{re}

{\bf D\'emonstration:} L'unicit\'e est une cons\'equence imm\'ediate du th\'eor\`eme \ref{specdomtg}. On montre donc maintenant l'existence. On proc\`ede par r\'ecurrence sur $n:=dim(Y)$. Si $n=0$, la conclusion est (\'evidemment) vraie, en consid\'erant les points comme des orbifoldes g\'eom\'etriques sp\'eciales. Supposant la conclusion vraie pour $dim(Y')<n,$ on choisit pour  $c_{(Y\vert\Delta)}$ une application constante si $(Y\vert\Delta)$ est sp\'eciale. Sinon, on choisit une fibration $f:(Y\vert\Delta)\dasharrow S$ de type g\'en\' eral, avec $p:=m(X)>0$ maximum. Il suffit de montrer que les fibres g\'en\'erales (orbifoldes) de $f$ sont sp\'eciales: on aura alors $f= c_{(Y\vert\Delta)}$. Si les fibres orbifoldes g\'en\'erales de $f$ ne sont pas sp\'eciales, il existe un sous-ensemble $T\subset S$ qui n'est pas analytiquement maigre dans $S$, et tel que pour $s\in T$, $(Y_s\vert\Delta_s)$ n'est pas sp\'eciale, et  $c_{(Y_s\vert\Delta_s)}$ est donc, par l'hypoth\`ese de r\'ecurrence, puique $dim(S)>0$, par hypoth\`ese, l'unique fibration de type g\'en\'eral maximum sur $(Y_s\vert\Delta_s)$. On d\'eduit alors de \ref{gtgen} et \ref{gengen}  l'existence d'une factorisation $f=h\circ g$ de $f$, avec $g:Y\dasharrow X$ et $h:X\dasharrow S$ des fibrations, telles que $dim(X)>dim(S)$, et telle que, pour $s\in S$ g\'en\'eral, la restriction $g_s:(Y_s\vert\Delta_s)\dasharrow X_s$ de $g$ \`a $Y_s$ soit de type g\'en\'eral.

Puisque $f=h\circ g$ est de type g\'en\'eral, on d\'eduit du th\'eor\`eme \ref{cftg} que $g:(Y\vert\Delta)\dasharrow X$ est une fibration de type g\'en\'eral. Ceci contredit le fait que $dim(S)$ soit maximum parmi les bases de fibrations de type g\'en\'eral d\'efinies sur $(Y\vert\Delta)$, et montre donc que les fibres orbifoldes de $f$ sont sp\'eciales $\square$

\begin{example} Soit $(Y\vert\Delta)$ une orbifolde g\'eom\'etrique lisse avec $Y$ compacte et connexe dans la classe $\sC$. Si $\kappa(Y\vert\Delta)\geq 0$, alors $dim(C(Y\vert\Delta))\leq \kappa(Y\vert\Delta)$, et on a \'egalit\'e si et seulement si $c_{(Y\vert\Delta)}=M_{(Y\vert\Delta)}$ est la fibration de Moishezon-Iitaka de l'orbifolde g\'eom\'etrique $(Y\vert\Delta)$(autrement dit: si soit $\kappa(Y\vert\Delta)=0$, soit si $M_{(Y\vert\Delta)}$ est de type g\'en\'eral. Ceci r\'esulte im\'ediatemment du th\'eor\`eme \ref{essk>o}.
 \end{example}

\begin{proposition} Soit $(Y\vert\Delta)$ une orbifolde g\'eom\'etrique lisse, avec $Y\in \sC$. Soit $f:(Y\vert\Delta)\dasharrow S$ une fibration m\'eromorphe \`a fibres orbifoldes g\'en\'erales sp\'eciales, et $g:(Y\vert\Delta)\dasharrow X$ une fibration m\'eromorphe telle que $\kappa(g\vert\Delta)=dim(X)$. On pose $c:=c_{(Y\vert\Delta)}$.

Il existe alors une unique couple de fibrations $\sigma:S\to C(Y\vert\Delta)$ et $\gamma:C(Y\vert\Delta)\to X$ telles que: $c\circ \sigma=f$ et $\gamma\circ c=g$.

\

\centerline{
\xymatrix{ (Y\vert\Delta)\ar[r]^{f}\ar[rd]^c\ar[d]_{g} & S\ar[d]^{\sigma}\\
X&C(Y\vert\Delta)\ar[l]^{\gamma}\\
}}

\end{proposition}

{\bf D\'emonstration:} C'est une cons\'equence imm\'ediate de \ref{specdomtg}: pour construire $\sigma$ (resp. $\gamma)$, on utilise le fait que $c$ est \`a base orbifolde stable de type g\'en\'eral  (resp. \`a fibre orbifolde g\'en\'erale sp\'eciale) $\square$

  \subsection{Fonctorialit\'e.}

  \begin{proposition} \label{c_f}Soit $f:(Y\vert\Delta)\dasharrow (X\vert\Delta')$ une application m\'eromorphe surjective d'orbifoldes g\'eom\'etriques lisses, avec $Y\in \sC$. Alors $f$ induit une (unique) application m\'eromorphe d'orbifoldes g\'eom\'etriques lisses $c_f:C(Y\vert\Delta)\dasharrow C(X\vert\Delta')$ telle que $c_{(X\vert\Delta')}\circ f=c_f\circ c_{(Y\vert\Delta)}$.
 
  De plus, $f$ induit une application m\'eromorphe orbifolde $c_f:[C(Y\vert\Delta)]\to [C(X\vert\Delta')]$ entre bases orbifoldes stables.\end{proposition}

    {\bf D\'emonstration:} La compos\'ee $c_{(X\vert\Delta')}\circ f:(Y\vert\Delta)\dasharrow C(X\vert\Delta')$ est de type g\'en\'eral, par la proposition \ref{composgt}. Sa restriction \`a la fibre g\'en\'erale orbifolde de $c_{(Y\vert\Delta)}$ est donc, par le th\'eor\`eme \ref{proprestr}, soit constante, soit de type g\'en\'eral. Comme cette fibre orbifolde est sp\'eciale, cette application est constante, et $c_{(X\vert\Delta')}\circ f$ se factorise donc par $c_{(Y\vert\Delta)}$. La seconde assertion r\'esulte de \ref{compfib}, en consid\'erant des repr\'esentants strictement nets et hauts $\square$

  \subsection{Connexit\'e par cha\^ines sp\'eciales.}

   \begin{proposition}\label{specrestr} Soit $(Y\vert\Delta)$ une orbifolde g\'eom\'etrique lisse, avec $Y\in \sC$, et $c_{(Y\vert\Delta)}:(Y\vert\Delta)\dasharrow C(Y\vert\Delta))$ son coeur. Si $(V_t)_{t\in T}$ est une famille couvrante de sous-vari\'et\'es de $Y$, et si $(V_t\vert\Delta_{V_t})$ (d\'efinie par restriction de $\Delta$ \`a $V_t$ comme en \ref{fdefrestr}) est sp\'eciale pour $t\in T$ g\'en\'eral, alors $V_t$ est contenue dans la fibre de $c_{(Y\vert\Delta)}$, pour tout $t\in T$.
  \end{proposition}

{\bf D\'emonstration:} Sinon, la restriction de $c_{(Y\vert\Delta)}$ est de type g\'en\'eral, par le th\'eor\`eme \ref{thmrestr}, ce qui contredit le fait que $(V_t\vert\Delta_{V_t})$ est sp\'eciale $\square$

  \begin{corollary}\label{s-conn} Soit $(Y\vert\Delta)$ une orbifolde g\'eom\'etrique lisse, $Y\in C$ compacte et connexe. Supposons que deux points g\'en\'eraux de $Y$ peuvent \^etre joints par une chaine connexe de sous-vari\'et\'es (au sens de \ref{fdefrestr}) de $Y$ telles que les orbifoldes g\'eom\'etriques de $(Y\vert\Delta)$ obtenues par restriction \`a ces sous-vari\'et\'es (au sens de \ref{fdefrestr}) soient sp\'eciales (On dira alors que $(Y\vert\Delta)$ est $\cal S$-connexe). 
  
  Alors $(Y\vert\Delta)$ est sp\'eciale. (Autrement dit: $(Y\vert\Delta)$ est sp\'eciale si elle est lisse et $\cal S$-connexe).
  \end{corollary}

{\bf D\'emonstration:} Soit $y$ l'un des deux points, choisi tel que la fibre $F$ de $c:=c_{(Y\vert\Delta)}$ passant par $y$ ne passe pas par le lieu d'ind\'etermination de $c_{(Y\vert\Delta)}$. Ceci est possible, puisque $c_{(Y\vert\Delta)}$ est presque-holomorphe, en vertu du th\'eor\`eme \ref{gtph}. Soit $V\subset Y$ une sous-vari\'et\'e telle que $(V\vert\Delta_V)$ soit sp\'eciale et rencontre $F$ sans \^etre contenue dans $F$ (on peut supposer que la restriction de $\Delta$ \`a $V$ est d\'efinie, par le th\'eor\`eme \ref{thmrestr}. Alors $c(V)\neq C(Y\vert\Delta)$, par \ref{proprestr}. Par r\'ecurrence sur $dim(Y)$, on en d\'eduit (en consid\'erant les vari\'et\'es de la forme $W:=c^{-1}(V))$ que $(Y\vert\Delta)$ admet une famille couvrante d'orbifoldes g\'eom\'etriques $(W\vert\Delta_W)$ obtenues par restriction de $\Delta$, et contenant strictement les fibres de $c$. Contradiction avec \ref{specrestr} $\square$

\begin{re} La lissit\'e de $(Y\vert\Delta)$ est essentielle (consid\'erer le c\^one sur une vari\'et\'e de type g\'en\'eral). En sens inverse, il existe des $Y$ (avec $\Delta=0$) sp\'eciales n'ayant aucune famille couvrante non-triviale de sous-vari\'et\'es sp\'eciales: par exemple, les vari\'et\'es ab\'eliennes simples.
  \end{re}

\begin{corollary}\label{rccspec} Si $(Y\vert\Delta)$ est lisse, avec $Y\in \sC$, et si $(Y\vert\Delta)$ est rationnellement connexe par chaines d'orbifoldes obtenues par restriction de $\Delta$, alors $(Y\vert\Delta)$ est sp\'eciale. 
\end{corollary}

  \subsection{Invariance par rev\^etement \'etale.}

   \begin{proposition} \label{invetcoeur}Soit $v:(Y'\vert\Delta')\to (Y\vert\Delta)$ un morphisme d'orbifoldes g\'eom\'etriques lisses et connexes dans la classe $\sC$. Si $v$ est \'etale en codimension $1$ (au sens orbifolde, voir \ref{defet} et \ref{etale}), alors $c_v:C(Y'\vert\Delta')\to C(Y\vert\Delta)$ est g\'en\'eriquement fini. En particulier, $ess(Y'\vert \Delta')=ess(Y\vert \Delta)$, et $(Y'\vert \Delta')$ est sp\'eciale si et seulement si $(Y\vert \Delta)$ l'est.
  \end{proposition}

    {\bf D\'emonstration:} On peut supposer que le rev\^etement $v$ est Galoisien de groupe fini $G$ en codimension $1$ (ie: sur le compl\'ementaire de $A$, Zariski ferm\'e de codimension au moins $2$ dans $Y)$.     
    Le groupe $G$ agit naturellement sur $C(Y'\vert\Delta')$, par unicit\'e du coeur. Soit alors $q:C(Y'\vert\Delta')\to C'(Y\vert\Delta)$ le quotient par $G$ de $C(Y'\vert\Delta')$. On a donc une application naturelle $c':Y\dasharrow C'(Y\vert\Delta)$ telle que $c'\circ v=q\circ c_{(Y'\vert\Delta')}$. Les fibres orbifoldes de cette application sont sp\'eciales, puisque images par $v$ de celles de $c_{(Y'\vert\Delta')}$. De plus, $c'$ est une fibration de type g\'en\'eral, par \ref{kappainv}. Donc $C'(Y\vert\Delta)=C(Y\vert\Delta)$ et $c'=c_{(Y\vert\Delta)}$. D'o\`u le r\'esultat $\square$

\subsection{Le coeur ``divisible".}\label{cd}

Si $(X\vert \Delta)$ est lisse et enti\`ere, avec $X\in \sC$, on peut d\'efinir (avec la m\^eme d\'emonstration) le coeur $c^{div}:(X\vert \Delta)\to C^*$ de $(X\vert \Delta)$ dans $Georb^{div}$: ses fibres g\'en\'erales sont sp\'eciales $^{div}$ (voir \ref{qspec}) et sa base orbifolde stable dans $Georb^{div}$ est de type g\'e,n\'eral (ou un point).

Si la question \ref{qmult} a une r\'eponse affirmative, alors $c^{div}=c$.

Puisque l'on ignore la r\'eponse \`a \ref{qmult}, la  consid\'eration de $c^*$ est essentielle pour l'\'etude de $\pi_1(X\vert \Delta)$ qui est (par \ref{sepi1}) extension du $\pi_1$ de la base orbifolde divisible $B^*$ de $c^*$ par celui de la fibre orbifolde g\'en\'erique $F^*$ de $c^*$. De plus, l'orbifolde lisse $F^*$ \'etant sp\'eciale, $\pi_1(F)$ est, par la conjecture \ref{cpi1}, presque ab\'elien.

\newpage

  \section{D\'ECOMPOSITION DU COEUR.}\label{deccoeur}

La conjecture $C_{n,m}^{orb}$ implique que le ``coeur" d'une orbifolde g\'eom\'etrique lisse $(Y\vert\Delta)$ de la classe $\cal C$ n'est autre que la compos\'ee des it\'er\'es, dans le cadre orbifolde g\'eom\'etrique (c'est essentiel), des deux fibrations de base de la classification: la fibration $M$ de Moishezon-Iitaka, et le quotient $\kappa$-rationnel $r$. En bref: $c=(M\circ r)^n, n=dim(Y)$.

\subsection{La d\'ecomposition conditionnelle du coeur.}

On rappelle les propri\'et\'es de $M$ et $r$.

\begin{lemma}  On admet l'existence du $\kappa$-quotient rationnel \footnote{Cons\'equence de $C_{n,m}^{orb}$.} pour les orbifoldes g\'eom\'etriques lisses dans $\sC$. 

Soit $(Y\vert\Delta)$ une orbifolde g\'eom\'etrique lisse, avec $Y\in \sC$. 

Soit $r:Y\dasharrow R^+=R^+(Y\vert\Delta)$ son $\kappa$-quotient rationnel, de base orbifolde stable $[R^+\vert\Delta_R^+]:=(R^+\vert\Delta(r,\Delta))$. (Voir le \S\ref{krat}).

Soit  $M:(R^+\vert\Delta_R^+)\dasharrow M(R^+)$ la fibration de Moishezon-Iitaka de $(R^+\vert\Delta_R^+)$(qui est bien d\'efinie, puisque $\kappa(R^+\vert\Delta_R^+)\geq 0)$.

On notera $s_{(Y\vert\Delta)}:=M\circ r:(Y\vert\Delta)\dasharrow S(Y\vert\Delta)$ la compos\'ee, appel\'ee {\bf r\'eduction sp\'eciale \'el\'ementaire} de $(Y\vert\Delta)$. 

Alors:

1. Les fibres orbifoldes stables g\'en\'erales de $s_{(Y\vert\Delta)}$ sont sp\'eciales.

2. $S(Y\vert\Delta)=(Y\vert\Delta)$ si et seulement si $(Y\vert\Delta)$ est de type g\'en\'eral (ie: si et seulement si $\kappa(Y\vert\Delta)=dim(Y)\geq 0)$. 
\end{lemma}

{\bf D\'emonstration:} La premi\`ere assertion r\'esulte de \ref{compspec}, la seconde de la remarque \ref{r=M=Y} pr\'ec\'edente $\square$

\begin{theorem}\label{c'}

{\bf (On admet l'existence du $\kappa$-quotient rationnel pour les orbifoldes g\'eom\'etriques lisses dans $\sC$, cons\'equence de $C_{n,m}^{orb})$.} 

Soit $(Y\vert\Delta)$ une orbifolde g\'eom\'etrique lisse, avec $Y\in \sC$. 

On d\'efinit une suite de fibrations $s^k:(Y\vert\Delta)\dasharrow S^k(Y\vert\Delta):=S^k$, de base orbifolde stable not\'ee $(S^k\vert\Delta_{S^k})$, par r\'ecurrence sur $k\geq 0$ par les conditions:

1. $s^0=id_Y$.

2. $s^{k+1}:=s_{(S^k\vert\Delta_{S^k})}:(S^k\vert\Delta_{S^k})\dasharrow S^{k+1}$

Alors:

1. Il existe un plus petit entier $k\leq n=dim(Y)$ tel que $S^{k+1}=S^{k}$. 

Cet entier sera appel\'e {\bf la longueur} de $(Y\vert\Delta)$, not\'e $\nu(Y\vert\Delta)$.

2.  Si $k\geq \nu(Y\vert\Delta)$, $s^k=s^n:(Y\vert\Delta)\dasharrow S^k=S^n$ est l'unique fibration de type g\'en\'eral \`a fibres orbifoldes sp\'eciales d\'efinie sur l'orbifolde g\'eom\'etrique $(Y\vert\Delta)$.

La fibration $S^n$ (conditionnelle en $C_{n,m}^{orb})$ est donc le coeur de $(Y\vert\Delta)$. Elle a \'et\'e construite inconditionnellement dans le \S \ref{coeur}.
\end{theorem}

{\bf D\'emonstration:} La premi\`ere assertion r\'esulte de ce que $dim(S^{k+1})<dim(S^k)$ si $S^{k+1}\neq S^k$, et de ce que $S^{k+j}=S^k$ pour tout $j>0$ si cette \'egalit\'e a lieu pour $j=1$. La seconde assertion r\'esulte du lemme pr\'ec\'edent, l'unicit\'e de \ref{compspec} $\square$

Reformulons le r\'esultat pr\'ec\'edent:

\begin{theorem}\label{c=jr^n} {\bf (On admet l'existence du $\kappa$-quotient rationnel (conditionnel) pour les orbifoldes g\'eom\'etriques lisses dans $\sC)$.} Soit $(Y\vert\Delta)$ une orbifolde g\'eom\'etrique lisse, avec $Y\in \sC$ et $n=dim(Y)$. Alors: $c=(M\circ r)^n$, avec les notations pr\'ec\'edentes.
\end{theorem}

Pour les orbifoldes g\'eom\'etriques sp\'eciales, on a donc:

\begin{corollary} {\bf (On admet l'existence du $\kappa$-quotient rationnel (conditionnel) pour les orbifoldes g\'eom\'etriques lisses dans $\sC)$}

Soit $(Y\vert\Delta)$ une orbifolde g\'eom\'etrique lisse, avec $Y\in \sC$, et $n=dim(Y)$.

Alors $(Y\vert\Delta)$ est sp\'eciale si et seulement si $dim(S^{n}(Y\vert\Delta))=0$. 

De mani\`ere \'equivalente: $(Y\vert\Delta)$ est sp\'eciale  si et seulement si $(Y\vert\Delta)$ est (comme en \ref{specdec}) une tour de fibrations \`a fibres orbifoldes g\'en\'erales stables $F$ ayant soit $\kappa=0$, soit $\kappa_+=-\infty$.
\end{corollary}

\begin{re} La suite de fibrations $s^k$ fournit une d\'ecomposition intrins\`eque et fonctorielle dans la cat\'egorie bim\'eromorphe des orbifoldes g\'eom\'etriques lisses de $\sC$ en composantes des trois ``g\'eom\'etries pures" \'el\'ementaires:

1. $\kappa_+=-\infty$ (conjecturalement rationnellement connexes dans un sens orbifolde ad\'equat).

2. $\kappa=0$.

3. $\kappa=dim$ (``type g\'en\'eral").

Les analogues ``num\'eriques" de ces trois g\'eom\'etries bim\'eromorphes sont:

1. $K<0$ (``g\'eom\'etrie de Fano)

2. $K\equiv 0$. (ou $c_1=0)$.

3. $K>0$ (fibr\'e canonique ample).

L'objectif du programme des mod\`eles minimaux (en version ``logarithmique") est pr\'ecis\'ement de r\'eduire la version bim\'eromorphe \`a la version num\'erique. Ce ``programme" pourrait \^etre (en l'absence de m\'ethodes bim\'eromorphes directes) l'outil permettant de d\'emontrer certaines des conjectures \'enonc\'ees ci-dessous concernant les orbifoldes g\'eom\'etriques.

La ``g\'eom\'etrie sp\'eciale" combine donc (conditionnellement, et seulement dans la cat\'egorie orbifolde bim\'eromorphe) les deux premi\`eres ``g\'eom\'etries pures" $(\kappa_+=-\infty$ et $\kappa=0)$, tandis que le ``coeur" d\'ecompose canoniquement et fonctoriellement les orbifoldes g\'eom\'etriques dans $\sC$ en leurs composantes ``sp\'eciale" (les fibres), et de ``type g\'en\'eral" (la base orbifolde).
\end{re}

\subsection{Analogie avec les alg\`ebres de Lie.}

Dans le dictionnaire ci-dessous entre notions de g\'eom\'etrie complexe, et notions de la th\'eorie des alg\`ebres de Lie complexes, le ``coeur" est l'analogue de la fibration de Levi-Mal\v cev, tandis que la d\'ecomposition du coeur $c=(M\circ r)^n$ est l'analogue de la suite des quotients par la s\'erie d\'eriv\'ee.

$$\begin{tabular}{| c || c | }
 \hline
G\'eom\'etrie orbifolde & Alg\`ebres de Lie \\
 
\hline
\hline
Fibration& Extension\\
 
 \hline

Fibre g\'en\'erale orbifolde& Noyau\\
 \hline
 
 Base orbifolde stable&Quotient \\
 \hline
 \hline
 Sp\'eciale&R\'esoluble\\
 \hline
 Type g\'en\'eral&Semi-simple\\
 \hline
 $\kappa=0$ ou $\kappa_+=-\infty$&Ab\'elienne\\
 \hline
 \hline
 Coeur&Levi-Mal\v cev\\
 \hline

\end{tabular}$$

\subsection{Vari\'et\'es de ``type" et dimension donn\'es.}

\begin{definition}\label{typ} Admettons la conjecture $C_{n,m}^{orb}$. Pour toute orbifolde g\'eom\'etrique $(X\vert\Delta)$ de dimension $n$, avec $X\in \cal C$, et $k=0,1,\dots,n$, soit $d_k=d_k(X\vert\Delta):=dim(S^k(X\vert\Delta))$, et pour $k=1,\dots,n$, soit: $d'_k:=dim(R(S^{k-1}(X\vert\Delta)))$.

Les $(2n+1)$ entiers: $d_0:=dim(X)\geq d'_1\geq d_1\geq d'_2\geq d_2\geq \dots\geq d_n=dim((X\vert\Delta)$ forment une suite d\'ecroissante $t(X\vert\Delta)=(d_0,d'_1,\dots,d_n)$ appel\'ee le {\bf type} de $(X\vert\Delta)$.

Soit $\nu=\nu(X\vert\Delta)$ le plus petit entier $k\leq n$ tel que $d_k=d_{k+1}=d_n$ ou $d'_{k-1}=d'_{k}=d_n$. Cet entier $\nu$ est la {\bf longueur} de $(X\vert\Delta)$.

Le type de $(X\vert\Delta)$ est soumis aux conditions: $d_k>d'_{k+1}$ si $k<\nu(X\vert\Delta)$ (puisque la compos\'ee de deux fibrations dont les fibres orbifoldes g\'en\'erales ont $\kappa_+=-\infty$ a aussi cette propri\'et\'e).
\end{definition}

Nous allons voir que ces conditions sont les seules, en g\'en\'eral. Il nous suffira pour cela, par r\'ecurrence sur $\nu$, d'\'etablir le lemme \ref{extyp} suivant:

\begin{lemma}\label{extyp} Soit $X=(X\vert 0)$ une vari\'et\'e projective lisse de dimension $m$, de type $t=(m=\delta_0,\delta'_1,\delta_1,\dots ,\delta_{\nu-1}=\dots=\delta_{m})$ et de longueur $(\nu-1)$. Alors, pour tout entier $r>0$:

1. Si $\kappa(X)\geq 0$, $((X\times \bP^r)\vert 0)$ est de longueur $\nu$ et de type $(d_0=m+r, d'_1=m=\delta_0=\delta_1,d'j=\delta_{j-1}, d_j=\delta_{j-1}$ pour $j=1,\dots,m)$.

2. Il existe une vari\'et\'e projective lisse $Y$ telle que $\kappa(Y)=dim(X)=m$, $dim(Y)=dim(X)+r$, et une fibration $J:Y\to X$ qui est la fibration de Moishezon-Iitaka de $Y$. Le type de $Y$ est alors: $(m+r=d_0=d_1,d'_j=\delta'_{j-1},d_j=\delta_{j-1}$ pour $j=1,\dots,m)$.
\end{lemma}

{\bf D\'emonstration:} L'assertion 1 est \'evidente, puisque la projection de $Y:=(X\times \bP^r)$ sur son premier facteur est le quotient rationnel de $Y$, et coincide donc avec $r$ pour l'orbifolde g\'eom\'etrique $(Y\vert 0)$. 

Pour l'assertion 2, on choisit sur $X$ un fibr\'e r\`es ample $H$, et pour $Y$ un membre g\'en\'erique du syst\`eme lin\'eaire $\vert d.H\times \cal O$$_{\bP^{r+1}}(r+1)\vert$, sur $(X\times \bP^{r+1})$, pour $d>0$ entier assez grand. La projection de $Y$ sur $X$ induite par celle du produit sur son premier facteur n'a pas de fibre multiple en codimension $1$, et satisfait les conditions \'enonc\'ees. (Nous ne le v\'erifions pas ici) $\square$

\begin{example} Pour tout $n\geq 0$, on note $c(n)$ (tout comme en [Ca 04, 6.22]), le nombre de types en dimension $n$.  Par \'etude directe on voit ([Ca 04, 6.22]) que $c(0)=1,c(1)=3,c(2)=8,c(3)=21$.

En dimension $1$, les trois types sont: $(1,0,0)$, $(1,1,0)$ et $(1,1,1)$ correspondant respecivement \`a $\kappa=-\infty,0,1$ pour les courbes $(X\vert 0)$.

En dimension 2, les $8$ types sont: $(2,0,0,0,0)$, $(2, 1,1,0,0)$, $(2,1,1,1,1)$, $(2,2,0,0,0)$, $(2,2,1,0,0)$,$(2,2,1,1,0)$, $(2,2,1,1,1)$, $(2,2,2,2)$. Ils correspondent, pour les surfaces $(X\vert 0)$ aux invariants $(\kappa, \tilde{q})$ suivants respectivement, o\`u $\tilde{q}$ est le maximum des irr\'egularit\'es des  rev\^etements \'etales finis de $X$: $(-\infty, 0)$, $(-\infty, 1)$, $(-\infty,+\infty)$ (ie: $q\geq 2$), $\kappa=0$,$(\kappa=1,0)$,$(\kappa=1,1)$,$(\kappa=1,+\infty)$, $\kappa=2$. 

On trouvera dans [Ca 04,3.38 ] une liste partielle des invariants birationnels des vari\'et\'es $(X\vert 0)$ correspondants aux diff\'erents types non-sp\'eciaux possibles en dimension $3$.

Plus g\'en\'eralement, on va montrer ([Ca 04, 6.22])\footnote{L'expression donn\'ee de $c(n)$ est exacte, bien que la formule de r\'ecurrence y soit incorrecte!} que: 
$c(n)=\frac{a^{n+1}-b^{n+1}}{\root\of 5}$, avec: $a=\frac{3+\root\of 5}{2}$ et $b=\frac{3-\root\of 5}{2}$, car $c(n+1)=3c(n)-c(n-1)$. 

En effet: on a $c(n+1)=(n+2)+\sum_{0}^{n}(n+1-d_1).c(d_1)$. Le premier terme correspond aux types tels que $d'_1=d_1\leq (n+1)$, le deuxi\`eme aux types tels que: $(d'_1>d_1)$. On en d\'eduit: $c(n+1-c(n)=1+\sum_{0}^{n}c(d_1)$, puis l'\'egalit\'e voulue en faisant la diff\'erence avec $c(n)-c(n-1)$. 
\end{example}

\begin{re} Sous r\'eserve de $C_{n,m}^{orb}$, nous avons donc d\'efini une s\'erie de nouveaux invariants bim\'eromorphes des orbifoldes g\'eom\'etriques lisses $(X\vert \Delta)$. Ces invariants sont donc d\'efinis en particulier dans les cas extr\^emes $\Delta=0$ (cas compact) et $\Delta=Supp(\Delta)$ (cas ouvert, qui s'applique \`a toute compactification \`a croisements normaux de toute vari\'et\'e quasi-projective). Nous conjecturons que ces invariants entiers sont stables par d\'eformation (par exemple lorsque $\Delta=0$).

Ces invariants g\'en\'eralisent et pr\'ecisent la dimension canonique (dite ``de Kodaira"). Ils fournissent, d\'ej\`a lorsque $\Delta=0$, (voir, par exemple, le cas de la dimension $2$ ci-dessus) une description beaucoup plus pr\'ecise de la structure de $X$ que la classique dimension ``de Kodaira".
\end{re}

\subsection{Fonctorialit\'e (conditionnelle)}

 Des propri\'et\'es de fonctorialit\'e \ref{fonctm}.4 et \ref{fonctr} pour $r$ et $M$, on d\'eduit imm\'ediatemment celles de $s$ (r\'eduction sp\'eciale \'el\'ementaire) et du coeur:

 \begin{proposition} \label{foncts} Soit $f:(X\vert\Delta_X)\dasharrow (Y\vert\Delta_Y)$ est un morphisme dans la cat\'egorie m\'eromorphe des orbifoldes g\'eom\'etriques lisses. On suppose que $X\in \sC$. Soit $k\geq 0$ un entier. Notons (pour simplifier les notations) $[S^k_X\vert\Delta_{S^k_X}]$ et $[S^k_Y\vert\Delta_{S^k_Y}]$ respectivement les bases orbifoldes stables de $s^k_{(X\vert\Delta_X)}$ et de $s^k_{(Y\vert\Delta_Y)}$, suppos\'ees exister.

Il existe alors un (unique) morphisme $s^k_f:[S^k_X\vert\Delta_{S^k_X}]\dasharrow [S^k_Y\vert\Delta_{S^k_Y}]$ tel que: $s^k_f\circ s^k_{(X\vert\Delta_X)}=s^k_{(Y\vert\Delta_Y)}\circ f$.

On a: $\nu(X\vert\Delta_X)\geq \nu(Y\vert\Delta_Y)$.

Si $\nu=\nu(X\vert\Delta_X)$, $s^{\nu}_f$ induit donc un morphisme orbifolde (dans la cat\'egorie bim\'eromorphe) entre les bases orbifoldes des coeurs de $(X\vert\Delta_X)$ et de $(Y\vert\Delta_Y)$
\end{proposition}

\begin{question} Si $u: (X\vert\Delta)\to (X'\vert\Delta')$ est un morphisme orbifolde \'etale en codimension $1$, entre orbifoldes g\'eom\'etriques lisses de $\sC$, les morphismes naturels induits par $u$ ($r^+_u,M_u,s_u, s^k_u)$ sont-ils \'egalement (sur des repr\'esentants ad\'equats) \'etales en codimension $1$? En particulier, la dimension des orbifoldes g\'eom\'etriques images doit \^etre alors la m\^eme. Nous verrons (en \ref{invetcoeur}) que cette derni\`ere propri\'et\'e est v\'erifi\'ee pour le coeur (voir \ref{invetcoeur}), et l'est donc bien, en particulier, pour chaque \'etape de la d\'ecomposition.\end{question}

\subsection{Rel\`evement de propri\'et\'es par d\'evissage.}\label{sdeviss}

Soit $(X\vert\Delta)$ une orbifolde g\'eom\'etrique lisse sp\'eciale. Admettons la conjecture $C_{n,m}^{orb}$. Soit alors $s^k:(X\vert\Delta)\dasharrow S^k(X\vert\Delta)$ la suite de fibrations d\'efinies ci-dessus, pour $k\leq \nu(X\vert\Delta)$, avec $s^k=(M\circ r)^k$. Cette suite est appel\'ee le {\bf d\'evissage canonique} de $(X\vert\Delta)$.

Soit maintenant (P) une propri\'et\'e susceptible d'\^etre ou non satisfaite par une orbifolde g\'eom\'etrique lisse $(X\vert\Delta)$, avec $X\in \sC$. De mani\`ere \'equivalente, notons $(P)$ la classe des orbifoldes g\'eom\'etriques lisses $(X\vert\Delta)$, avec $X\in \sC$ poss\'edant la propri\'et\'e (P).

On a alors le lemme de ``rel\`evement" \'evident suivant, qui permet de relever aux orbifoldes g\'eom\'etriques sp\'eciales des propri\'et\'es satisfaites pour les classes $\kappa=0$ et $\kappa_+=-\infty$:

\begin{lemma}\label{deviss} Admettons $C_{n,m}^{orb}$.

{\bf A.} Supposons que la classe (P) poss\`ede les propri\'et\'es de stabilit\'e suivantes:

1. Invariance bim\'eromorphe.

2.$(X\vert\Delta)\in (P)$ si $\kappa(X\vert\Delta)=0$.

3.$(X\vert\Delta)\in (P)$ si $\kappa_+(X\vert\Delta)=-\infty$. 

4. Stabilit\'e par extension: si $f:(X\vert\Delta)\dasharrow Y$ est une fibration m\'eromorphe de base orbifolde stable $[Y\vert\Delta(f,\Delta)]\in (P)$, et si $(X\vert\Delta)_y\in (P)$, pour $y\in Y$ g\'en\'eral, alors $(X\vert\Delta)\in (P)$. 

Alors: $(X\vert\Delta)\in (P)$ si $(X\vert\Delta)$ est lisse et sp\'eciale.

{\bf B.} Si, de plus, la classe (P) est telle que:

5. $(X\vert\Delta)\notin (P)$ si $(X\vert\Delta)$ est lisse, de type g\'en\'eral, avec: $dim(X)>0$.

6. Si $f:(X\vert\Delta)\dasharrow Y$ est une fibration, avec $(X\vert\Delta)\in (P)$, alors $[Y\vert\Delta(f,\Delta)]\in (P)$.

Alors: $(P)$ est {\bf exactement} la classe des orbifoldes g\'eom\'etriques sp\'eciales.
\end{lemma}

(Pour d\'emontrer le point {\bf B}, utiliser le ``coeur", en raisonnant par l'absurde).

Ce lemme devrait permettre de r\'eduire la preuve d'un certain nombre de propri\'et\'es conjecturales (voir S.\ref{conj}) des orbifoldes g\'eom\'etriques lisses et sp\'eciales aux cas ``\'el\'ementaires" cruciaux $\kappa=0$ et $\kappa_+=-\infty$. La stabilit\'e par extensions est elle-m\^eme conjecturale dans la plupart des situations. Sa signification est la suivante: si les obstructions \`a \'etendre $(P)$ s'annulent localement sur la base (orbifolde stable) de $f$, elles s'annulent globalement sur cette base. Des exemples de propri\'et\'es auxquelles appliquer ce principe de rel\`evement sont fournies dans le \S{\ref{conj}.

\newpage

  \section{GROUPE FONDAMENTAL}\label{gf}

  \subsection{Groupe fondamental d'une orbifolde g\'eom\'etrique lisse}

  Les orbifoldes g\'eom\'etriques consid\'er\'ees dans cette section sont enti\`eres, {\bf lisses} et connexes. Les morphismes orbifoldes le sont au sens {\bf divisible} (voir d\'efinition \ref{morphorb}).

  On montre ici que le groupe fondamental orbifolde se comporte, sous l'hypoth\`ese usuelle de lissit\'e, comme dans le cas des vari\'et\'es sans structure orbifolde (et m\^eme mieux, puisque les fibrations induisent toujours des suites exactes, voir \ref{sepi1}. Ce fait est, avec \ref{compspec},  une seconde indication du fait que la cat\'egorie orbifolde g\'eom\'etrique puisse \^etre le cadre naturel de la g\'eom\'etrie k\" ahl\'erienne ou alg\'ebrique)\footnote{Remarquons que les r\'esultats de ce chapitre s'\'etendent aux champs de Deligne-Mumford lisses dont ``l'espace des modules" est projectif (ou dans la classe $\sC)$.}.

  \begin{re} \label{pdiv} Remarquons enfin que le groupe fondamental d'une orbifolde lisse et enti\`ere $(X\vert \Delta)$, avec $X\in \sC$ peut \^etre \'etudi\'e en consid\'erant son coeur ``divisible" $c^*$ (voir le \S\ref{cd}), qui exprime $\pi_1(X\vert \Delta)$ comme extension de $\pi_1(B^*)$, la base orbifolde de $c^*$, par $\pi_1(F^*)$, la fibre orbifolde g\'en\'erique de $c^*$. Ce dernier groupe est, conjecturalement (par \ref{cpi1}), presque ab\'elien, puisque $F^*$ est sp\'eciale$^{div}$.
 \end{re} 
  
  La d\'efinition suivante est classique:

  \begin{definition}\label{dpi1} Soit $(X\vert\Delta)$ une orbifolde g\'eom\'etrique lisse, avec $X$ connexe, et $\Delta:=\sum_{j\in J} (1-\frac{1}{m_j}). D_j$.  On note $X^*:=(X-supp(\Delta))$ le compl\'ementaire dans $X$ du support de $\Delta$, et $a\in X^*$. On note $\pi_1(X\vert\Delta,a)$ le quotient de $\pi_1(X^*,a)$ par le sous-groupe {\bf normal} engendr\'e par les lacets $g_j^{m_j}$, d\'esignant par $g_j,\forall j\in j$ le lacet bas\'e en $a$, tournant une fois dans le sens direct autour du diviseur $D_j$ (lacet dit ``\'el\'ementaire" pour $D_j$). (Dans les probl\`emes consid\'er\'es ici, on omettra la mention des points-base, qui n'y jouent aucun r\^ole). 
  
  On notera enfin $X^0$ tout ouvert de Zariski de $X$, compl\'ementaire d'un sous-ensemble analytique ferm\'e $A$ de codimension $2$ au moins, $A$ contenant le lieu singulier $Sing(\lceil\Delta\rceil)$. On note $\Delta^0:=\Delta\cap X^0$.
  \end{definition}

  \begin{re} Si $X$ est compacte, $\pi_1(X\vert\Delta)$ est donc de pr\'esentation finie, puisque quotient de $\pi_1(X^*)$, qui est de type fini, par le sous-groupe normal engendr\'e par un nombre fini d'\'el\'ements.
  \end{re}

  \begin{example}\label{pp} Si $X=\Bbb P^2$, et si $D:=Supp(\Delta)$ est \`a croisements normaux, et r\'eunion de $r\geq 0$ courbes irr\'eductibles de degr\'es $d_1,\dots,d_r$, alors $G^*:=\pi_1(\Bbb P^2-D)$ est le quotient de $\Bbb Z^{\oplus r}$ par les sous-groupe engendr\'e par $(d_1,\dots,d_r)$ (par [De 79]). Donc toute orbifolde enti\`ere $(X\vert \Delta)$ avec $Supp(\Delta)=D$ a un $\pi_1$ ab\'elien fini, quotient de $G^*$. La situation est tr\`es diff\'erente d\`es que l'on a des points triples: par exemple si $D$ est la r\'eunion de trois droites concourantes. Voir aussi [B-H-H 87, 3.1.D, p.108, et 2.2.F, pp. 65-66]. 
  \end{example}

 \begin{proposition}\label{fonctpi1} Soit $f:(X\vert\Delta_X)\to (Y\vert\Delta_Y)$ un morphisme orbifolde {\bf divisible}\footnote{Voir d\'efinition en \ref{morphorb}.} entre orbifoldes g\'eom\'etriques lisses, avec $X$ connexe.

1. Il induit un morphisme fonctoriel naturel de groupes $f_*:\pi_1(X\vert\Delta_X)\to (Y\vert\Delta_Y)$.

 2. Si $f$ est une fibration, $f_*$ est surjectif.

 3. L'injection $j^0:(X^0\vert\Delta^0)\to (X\vert\Delta)$ induit un isomorphisme de groupes: $j^0_*:\pi_1(X^0\vert\Delta^0)\to \pi_1(X\vert\Delta)$

 4. Si $f$ est \'etale en codimension $1$, $f_*$ est injectif. L'indice de $\pi_1(X\vert\Delta_X)$ dans $\pi_1(Y\vert\Delta_Y)$ est \'egal au degr\'e g\'eom\'etrique de $f$.

 5. Si $f$ est propre et surjective, son image est d'indice fini (au plus \'egal au nombre de composantes connexes d'une fibre g\'en\'erique) dans $\pi_1(Y\vert\Delta_Y)$. 
 \end{proposition}

 {\bf D\'emonstration:} 1. L'image par $f$ du support $D_X$ de $\Delta_X$ est contenue dans le support $D_Y$ de $\Delta_Y$. Si $D$ est une composante de $\Delta_X$ de multiplicit\'e $m$, et si $E$ est une composante quelconque de $\Delta_Y$ contenant $f(D)$, de multiplicit\'e $m'$, on a donc, par hypoth\`ese: $t.m=k.m'$, pour un entier $k>0$. Soit donc $g_D$ un lacet de $X-D_X$ \'el\'ementaire pour $D$. Soit $g_E$ un lacet \'el\'ementaire relatif \`a $E$ dans $(Y-D_Y)$. Puisque $f^*(E)=t.D+\dots$, $f_*(g_D)=(g_E)^t$ (\`a conjugaison pr\`es dans $\pi_1(Y-D_Y))$. Donc $f_*(g_D^{m})=(g_E)^{t.m}=(g_E)^{m'.k}=((g_E)^{m'})^{k}$, est donc dans le noyau du quotient $\pi_1(Y-D_Y)\to \pi_1(Y\vert\Delta_Y)$, puisque ceci est vrai pour toute composante $E$ de $D_Y$. La fonctorialit\'e de $f_*$ r\'esulte imm\'ediatemment de sa d\'efinition.

 2. Soit $X^{**}:=(X-f^{-1}(D_Y))\cap X^*\subset X^*\subset X$. Si $f$ est une fibration, le morphisme de groupes $\pi_1(X^{**})\to \pi_1(Y^*)$ est surjectif. En effet: $X$, et donc $X^{**}$ est lisse et connexe, les fibres g\'en\'eriques de $f_{\vert X^{**}}:X^{**}\to Y^*$ sont connexes, et $f_{\vert X^{**}}:X^{**}\to Y^*$ est surjective, par hypoth\`ese. Donc $f_*$ est \'egalement surjectif, puisque $f_*$ est d\'eduit du morphisme de groupes pr\'ec\'edent par passage aux quotients. Plus pr\'ecis\'ement, le diagramme commutatif ci-dessous \'etablit la surjectivit\'e de $f_*$.

\centerline{
\xymatrix{ \pi_1(X^{**})\ar[r]\ar[d] & \pi_1(X^{*})\ar[r]&\pi_1(X\vert\Delta)\ar[ld]^{f_*}\\
\pi_1(Y^*)\ar[r]& \pi_1((Y\vert\Delta_Y))&\\
}}

  3. Les groupes $\pi_1(X^*)$ et $\pi_1(X^0)^*$ sont naturellement isomorphes par $j^0$, et les noyaux des morphismes d\'efinissant $\pi_1(X^0\vert\Delta^0)$ et $\pi_1(X\vert\Delta)$ sont engendr\'es par les (classes des) lacets \'el\'ementaires correspondants.

  4. On est donc r\'eduit, par 3, au cas o\`u $f$ est \'etale au sens orbifolde.

  On peut donc supposer \^etre dans la situation de \ref{etale}. Si $D$ est une composante de $D_X$ de multiplicit\'e $m$, alors $f(D)=E$ est une composante du support de $D_Y$ qui est de multiplicit\'e $m'$, avec $r.m=m'$, si $r$ est l'indice de ramification de $f$ le long de $D$, tel donc que $f^*(E)=r.D+\dots$. Donc $f_*(g_D)^m=((g_E)^r)^{m'}$, avec les notations utilis\'ees dans la preuve de 1. Puisque le morphisme naturel $\pi_1(X^*)\to \pi_1(Y^*)$ induit par la restriction de $f$ est injectif, et que les sous-groupes normaux de $\pi1(Y^*)$ engendr\'es par les lacets \'el\'ementaires autour des composantes de $D_X$ et de $D_Y$ coincident, $\pi_1(X\vert\Delta_X)$ est bien un sous-groupe de $\pi_1(Y\vert\Delta_Y)$ dont l'indice coincide avec celui de $\pi_1(X^*)$ dans $\pi_1(Y^*)$, \'egal au degr\'e (fini ou pas) de $f$.

  5. Le diagramme utilis\'e dans la d\'emonstration de l'assertion 2 ci-dessus montre que l'on peut remplacer $Y$ par l'un quelconque de ses ouverts non vides de Zariski $Y'$, et $X$ par $X':=f^{-1}(Y')$. Soit $f=h\circ s$, $s:X\to Z$ une fibration, et $h:Z\to Y$ finie, la factorisation de Stein de $f$. Quitte \`a remplacer $Y$ par $Y'$ ad\'equat, on supposera que $(Z\vert\Delta(s,\Delta_X))$ est lisse, que $s_{\Delta}:(X\vert\Delta_X)\to(Z\vert\Delta(s,\Delta_X))$ est un morphisme orbifolde (n\'ecessairement divisible), et que $h:Z\to Y$ est \'etale et induit un morphisme orbifolde divisible $h_{\Delta}:(Z\vert\Delta_Z):=(Z\vert\Delta(s,\Delta_X))\to (Y\vert\Delta_Y)$. La seule assertion non imm\'ediate est la divisibilit\'e de $s_{\Delta}$ et $h_{\Delta}$. Puisque $f$ est divisible, $\Delta_Y$ divise $\Delta(f,\Delta_X)$. Par finitude de $h$, $\Delta(f,\Delta_X)=\Delta(h,\Delta(s,\Delta_X)):=\Delta(h,\Delta_Z)$. Ce sont les assertions de divisibilit\'e annonc\'ees. On est donc r\'eduit au cas o\`u $f=s$ est \'etale finie. L'assertion de divisibilit\'e signifie alors que $h^*(\Delta_Y)$ divise $\Delta_Z$. Donc $\pi_1(Z/h^*(\Delta_Y))$ est un quotient de $\pi_1(Z\vert\Delta_Z)$. Puisque $\pi_1(Y\vert\Delta_Y)$ est un quotient de $\pi_1(Z/h^*(\Delta_Y))$ de degr\'e \'egal au nombre de composantes connexes de la fibre g\'en\'erique de $f$, le r\'esultat est \'etabli.

  \begin{re}\label{pi1div}
  
  1. La condition de divisibilit\'e du morphisme orbifolde $f$ est en g\'en\'eral n\'ecessaire, pour l'existence de $f_*$ : consid\'erer par exemple $X=\bP^1$ avec les trois diviseurs orbifoldes $\Delta_k=\{0\}+(1-\frac{1}{a_k}). \{\infty\}$, pour $k=1,2,3$ et: $1<a_1<a_2<a_3=+\infty$, $a_1$ et $a_2$ premiers entre eux. Ils induisent deux morphismes de groupes: $\bZ\to \bZ_{a_1}$ et $\bZ\to \bZ_{a_2}$, mais aucun morphisme $\bZ_{a_2}\to \bZ_{a_1}$ n'est compatible par composition avec les pr\'ec\'edents. 
  
  2. Si $f:(X\vert\Delta)\to (Y\vert\Delta)$ est l'\'eclatement d'un point lisse du support de $\Delta_Y$, situ\'e sur une composante $D$ de $D_Y$ affect\'ee dans $\Delta_Y$ de la multiplicit\'e $m>1$, et si le diviseur exceptionel $D$ de $f$ est affect\'e dans $\Delta_X$ tel que $f_*(\Delta_X)=\Delta_Y$ d'une multiplicit\'e $m'$, premi\`ere avec $m$, et assez grande pour que $f:(X\vert\Delta)\to (Y\vert\Delta)$ soit un morphisme orbifolde (non divisible), alors $\pi_1(X\vert\Delta_X)\cong \pi_1(Y\vert\Delta'_Y)$, o\`u $\Delta'_Y$ est obtenu de $\Delta_Y$ par suppression de la composante $(1-\frac{1}{m}).D$. Ce dernier groupe est en g\'en\'eral un quotient strict de $\pi_1(Y\vert\Delta_Y)$.
  
  En particulier, le groupe fondamental {\bf n'est pas} un invariant bim\'eromorphe dans la cat\'egorie des orbifoldes lisses (munie des morphismes {\bf non} divisibles).
  
  3. Un cas particulier de $f:(X\vert\Delta_X)\to (Y\vert\Delta_Y)$ auquel l'assertion 2 peut \^etre appliqu\'ee est celui dans lequel $X=Y$ et $\Delta_Y$ divise $\Delta_X$ (ie: $m_{\Delta_Y}(D)$ divise $m_{\Delta_X}(D),\forall D\in W(X)=W(Y))$: $\pi_1(Y\vert\Delta_Y)$ est alors un quotient de $\pi_1(X\vert\Delta_X)$.
  \end{re}

  \begin{proposition} $f:(X\vert\Delta_X)\to (Y\vert\Delta_Y)$ un morphisme orbifolde {\bf divisible}, avec $X$ connexe. On suppose que $f$ est bim\'eromorphe au sens orbifolde (ie: bim\'eromorphe de $X$ sur $Y$, et tel que $f_*(\Delta_X)=\Delta_Y)$. Alors $f_*$ est un isomorphisme de groupes.

  En particulier, deux orbifoldes g\'eom\'etriques bim\'eromorphes {\bf au sens divisible}\footnote{ie: s'il existe une chaine d'\'equivalences bim\'eromorphes orbifoldes bim\'eromorphes {\bf divisibles} les reliant.} ont des groupes fondamentaux isomorphes.
  \end{proposition}

  {\bf D\'emonstration:} Par l'observation faite ci-dessus, on ne change pas $\pi_1(Y\vert\Delta_Y)$ en restreignant $Y$ au compl\'ementaire de $A\subset Y$, Zariski ferm\'e de codimension au moins $2$. On choisit $A$ contenant le lieu d'ind\'etermination de $f^{-1}$. Alors $f':X^*\to Y^*$ est un isomorphisme. Et il existe donc un morphisme quotient  $h: \pi_1(Y^*)\to \pi_1(X\vert\Delta_X)$ tel que $f_*\circ h=\delta_Y:\pi_1(Y^*)\to \pi_1(Y\vert\Delta_Y)$ est le quotient naturel. Puisque le noyau de $\delta$ (engendr\'e par les lacets $\Delta_Y$-\'el\'ementaires) est contenu dans celui de $h$ (engendr\'e par les lacets $\Delta_X$-\'el\'ementaires), $f_*$ est un isomorphisme $\square$

   \begin{definition}\label{dreg} Soit $\Delta'$ et $\Delta$ deux diviseurs orbifoldes sur $X$, connexe. Alors $id_X:(X\vert\Delta_X)\to (X\vert\Delta')$ est un morphisme orbifolde {\bf divisible} si et seulement si $\Delta'$ {\bf divise } $\Delta$ (ie: que pour chaque $D\in W(X)$, $m_{\Delta'}(D)$ divise $m_{\Delta}(D))$. Ce morphisme induit un morphisme de groupes $(id_X)_*:\pi_1(X\vert\Delta_X)\to \pi_1(X\vert\Delta')$ surjectif.
   
   On dit que $(X\vert\Delta)$ est {\bf r\'eguli\`ere} si le seul diviseur orbifolde $\Delta'$ divisant $\Delta$, tel que le $(id_X)_*$ soit un isomorphisme est $\Delta'=\Delta$ (autrement dit: si on ne peut pas ``r\'eduire" strictement $\Delta$ sans r\'eduire strictement son $\pi_1)$. 
   \end{definition}

  \subsection{Suite exacte associ\'ee \`a une fibration orbifolde nette}

    \begin{proposition}\label{sepi1} $f:(X\vert\Delta_X)\to (Y\vert\Delta_Y)$ un morphisme orbifolde {\bf divisible}, avec $X$ connexe. Soit $(X_y\vert\Delta_y)$ une fibre orbifolde lisse g\'en\'erique de $f$, notant $\Delta:=\Delta_X$, et $j:(X_y\vert\Delta_y)\to (X\vert\Delta)$ l'inclusion naturelle. Elle induit une morphisme de groupes naturel $j_*:\pi_1(X_y\vert\Delta_y)\to \pi_1(X\vert\Delta)$.

    Si $f$ est nette au sens divisible\footnote{ie: satisfait la condition de \ref{nette}, les morphismes orbifoldes l'\'etant alors tous au sens divisible.}, la suite de groupes suivante induite par $j$ et $f$ est exacte:

\centerline{
\xymatrix{ \pi_1(X_y\vert\Delta_y)\ar[r]^{j_*} & \pi_1(X\vert\Delta)\ar[r]^{f_*}&\pi_1(Y\vert\Delta_Y)\ar[r]&\{1\}\\
}}

    \end{proposition}

  {\bf D\'emonstration:} Seule l'exactitude en $ \pi_1(X\vert\Delta)$ reste \`a \'etablir. Puisque $f_*\circ j_*=(f\circ j)_*$, il faut montrer que le noyau de $f_*$ est contenu dans l'image de $j_*$. Rappelons (voir \ref{nette}) que $f$ ``nette" (au sens divisible) signifie qu'il existe un diagramme commutatif:

\centerline{
\xymatrix{ (X\vert\Delta )\ar[r]^{w}\ar[d]_{f} & (X'\vert\Delta')\ar[d]^{f'}\\
Y\ar[r]^v& Y'\\
}}

dans lequel:

1. $w$ est un morphisme orbifolde divisible, $v$ et $w$ \'etant  bim\'eromorphes,  $Y,Y'$ lisses, et $w_*(\Delta)=\Delta'$.

2. Tout diviseur $g$-exceptionnel de $X$ est $w$-exceptionnel.

  Si $A\subset Y$ est analytique ferm\'e de codimension $2$ ou plus, $f^{-1}(A)=B\cup E\subset X$ est analytique ferm\'e r\'eunion de $B$, et $E'$, analytique ferm\'es dans $X$, $B$ de codimension $2$ ou plus, et $E$ un diviseur $f$-exceptionnel, donc $w$-exceptionnel, puisque $f$ est nette relativement \`a $f'$. Il en r\'esulte que $\pi_1(X"\vert\Delta")=\pi_1(X\vert\Delta)$, si $X"=f^{-1}(Y-A)$, et $\Delta":=\Delta\cap X'$.

  Rempla\c cant $Y$ par un ouvert ad\'equat $Y"$ et $X$ par $X":=f^{-1}(Y")$, on supposera donc d\'esormais que $f$, ainsi que sa restriction au support de $\Delta_X$ sont \`a fibres \'equidimensionnelles, le support de $\Delta_Y$ \'etant lisse.

  On va maintenant montrer que l'on peut supposer que, de plus, $\Delta=\Delta^{vert}$, partie $f$-verticale de $\Delta$, ce dernier diviseur orbifolde \'etant simplement la r\'eunion des composantes du support de $\Delta$ qui ne sont pas envoy\'ees surjectivement sur $Y$ par $f$, les multiplicit\'es des composantes ($f$-verticales) restantes \'etant conserv\'ees.

    Alors: $\Delta_Y:=\Delta(f,\Delta)=\Delta(f,\Delta^{vert})$. Cette derni\`ere \'egalit\'e r\'esulte imm\'ediatemment des d\'efinitions.

  On a, de plus, une surjection naturelle $\pi_1(X\vert\Delta)\to \pi_1(X\vert\Delta^{vert})$, compatible avec les restrictions \`a $X_y$.

  La propri\'et\'e cruciale pour cette seconde r\'eduction est le fait que le noyau du morphisme naturel $\pi_1(X\vert\Delta)\to \pi_1(X\vert\Delta^{vert})$ est engendr\'e par les lacets \'el\'ementaires autour des composantes $f$-horizontales de $\Delta$, et est donc contenu dans $\pi_1(X_y\vert\Delta_{X_y})=\pi_1(X_y),$ puisque $\Delta^{vert}_{\vert X_y}=0$.

  On supposera donc que $f$ et sa restriction \`a $\Delta=\Delta^{vert}$ sont \'equidimensionnels, et que $\Delta_Y=\Delta(f,\Delta)$ est \`a support lisse.

 Soit donc $g \in \pi_1(X^*)$ tel que son image dans $\pi_1(X\vert\Delta)$ soit dans le noyau de $f_*$ ci-dessus. On peut donc \'ecrire: $f(g)=\Pi_{r}h_{r}^{p_r.\mu_r}$, o\`u les $h_{r}^{\mu_r}$ sont des lacets \'el\'ementaires autour des composantes du support de $\Delta_Y$. Puisque $\Delta_Y$ est la base orbifolde de $(f\vert\Delta)$, pour chaque $r$, on peut \'ecrire, par le th\'eor\`eme de Bezout: $h_r^{\mu_r}=\Pi_sf(g_{r,s}^{t_{r,s}.m_{r,s}.q_{r,s}})$, puisque $\mu_r=pgcd_s(t_{r,s}.m_{r,s})$, par la d\'efinition de la base orbifolde (divisible). Les $g_{r,s}$ sont des lacets \'el\'ementaires autour de composantes de $\Delta^{vert}$ dont l'image par $f$ est la composante d'indice $r$, de multiplicit\'e $\mu_r$, du support de $\Delta_Y$.

 Donc: $g=\Pi_{r,s}f(g_{r,s}^{p_r.t_{r,s}.m_{r,s}.q_{r,s}}).g'$, o\`u $f(g')=1$ dans $\pi_1(Y^*)$. Mais la suite exacte d'homotopie:

\centerline{
\xymatrix{ \pi_1(X_y)\ar[r]^{j_*} & \pi_1(X^*)\ar[r]^{f_*}&\pi_1(Y^*)\ar[r]&\{1\}\\
}}

 montre que $g'\in \pi_1(X_y)=\pi_1(X_y\vert\Delta_{X_y})$. Ce qui ach\`eve a d\'emonstration.

  Passant \`a un mod\`ele holomorphe net, on en d\'eduit la version bim\'eromorphe suivante:

  \begin{corollary}\label{sepi} Soit $(X\vert\Delta)$ lisse, avec $X$ connexe, et $f:X\dasharrow Y$ m\'eromorphe dominante connexe et propre, de fibre orbifolde g\'en\'erique $(X\vert\Delta)_y$ (sur un mod\`ele bim\'eromorphe rendant $f$ holomorphe), et de base orbifolde stable $[Y\vert\Delta_Y]$. La suite de groupes suivante est alors exacte:

\centerline{\xymatrix{ \pi_1(X_y\vert\Delta_y)\ar[r]^{j_*} & \pi_1(X\vert\Delta)\ar[r]^{f_*}&\pi_1([Y\vert\Delta_Y])\ar[r]&\{1\}\\}}

  \end{corollary}

   Explicitons le lien entre le groupe fondamental de $X\in \sC$, lisse et connexe, et celui de la base orbifolde (stable) de sa $\Gamma$-r\'eduction.

  \begin{corollary}\label{gammared}Soit $X\in \sC$, lisse et connexe. Soit $\gamma_X:X\to Y:=\Gamma(X)$ sa $\Gamma$-r\'eduction (au sens de [ca94], ou sa r\'eduction de Shafarevich), et $X_g$ sa fibre g\'en\'erique (lisse). Soit $(\Gamma(X)\vert\Delta(\gamma_X))=[Y\vert\Delta_Y]$ sa base orbifolde (stable). Alors on a une suite exacte de groupes: 
  
  \centerline{\xymatrix{ \pi_1(X_g)\ar[r]^{j_*} & \pi_1(X)\ar [r]^{\gamma_*}  &   \pi_1([Y\vert\Delta_Y])\ar[r]&\{1\}\\}}

  \end{corollary}

  \begin{re} Cette suite exacte montre que l'\'etude de $\pi_1(X)$ n\'ecessite celle de l'orbifolde g\'eom\'etrique $(Y\vert\Delta_Y)$.
  \end{re}

  \subsection{Rev\^etement universel d'une orbifolde lisse}\label{revun}

  Soit $(X\vert\Delta)$ une orbifolde g\'eom\'etrique lisse avec $X$ connexe. On note $\Delta=\sum_{j\in J} (1-\frac{1}{m_j}). D_j$. On note aussi $X^*:=(X-\lceil \Delta\rceil)$, et $X^0:=X-A$, avec $A\supset Sing(\lceil\Delta\rceil)$ le lieu singulier du support de $\Delta$, si $A$ est de codimension $2$ au moins dans $X$.

  Nous adoptons ici la terminologie de [N87], et en rappellons certains r\'esultats.

  \begin{definition}\label{drevram} Un {\bf rev\^etement ramifi\'e en au plus $\Delta$} de $X$ est une application holomorphe $g:\bar X\to X$ telle que:
  
  1. $\bar X$ est normal et connexe, et les fibres de $g$ sont discr\`etes.
  
  2. La restriction $g^*:\bar X^*:=g^{-1}(X^*)\to X^*$ de $g$ \`a $\bar X^*$ est un rev\^etement non ramifi\'e (de groupe $G\subset \pi_1(X^*))$. 
  
  3. Si $\bar D_j^0$ est une composante irr\'eductible de $g^{-1}(D_j)\subset \bar X^0:=g^{-1}(X^0)$, alors $g$ ramifie \`a l'ordre $n_j$, diviseur de $m_j$ le long de $\bar D_j^0$. 
  
  On affecte l'adh\'erence de $\bar D_j^0$  de la multiplicit\'e $\bar m_j:=\frac{m_j}{n_j}$. On obtient ainsi un diviseur orbifolde $\bar \Delta$ sur $\bar X$ dont le support est contenu dans la r\'eunion des $g^{-1}(D_j)$.
  
  4. Tout $a\in X$ a un voisinage ouvert connexe $W_a$ tel que, pour toute composante connexe $U_a$ de $g^{-1}(W)$, $U_a\cap g^{-1}(a)$ soit r\'eduit \`a un point $\bar a$, et $g_{\vert U_a}:U_a\to W$ est propre (et finie).
  
  On dira que $\bar X$ est le compl\'et\'e de $\bar X^*$ au-dessus du support de $\Delta$. 
  
  On dit que {\bf $g$  ramifie en $\Delta$} si $n_j=m_j,\forall j\in J$.
  
  Le morphisme $g:(\bar X\vert \bar \Delta)\to (X\vert \Delta)$ est alors un {\bf morphisme orbifolde \'etale}.

  La notion d'isomorphisme de rev\^etements ramifi\'es est la notion usuelle. Les rev\^etements ramifi\'es sont consid\'er\'es \`a isomorphisme pr\`es dans la suite.
   \end{definition}

   \begin{theorem}\label{trevram} [N87, 1.3.8] Dans la situation pr\'ec\'edente \footnote {Le r\'esultat de [N87] couvre, plus g\'en\'eralement, les diviseurs orbifoldes \`a groupes fondamentaux locaux finis. On l'applique ici au-dessus de $X-\Delta_{\infty}$, si $\Delta_{\infty}$ est la r\'eunion des composantes irr\'eductibles du support de $\Delta$ qui sont affect\'ees de la multiplicit\'e $+\infty$.}, il existe une bijection naturelle entre les sous-groupes de $\pi_1(X^*)$ contenant le sous-groupe normal $K$ engendr\'e par les $g_j^{m_j},j\in J$, et les rev\^etements ramifi\'es en au plus $\Delta$ de $X$ (\`a isomorphisme de rev\^etement pr\`es). 
   
   Cette bijection associe au rev\^etement $g:\bar X\to X$ le sous-groupe de $\pi_1(X^*)$ consitu\'e des automorphismes du rev\^etement $g^*:\bar X^*\to X^*$ obtenu par restriction de $g$ au-dessus de $X^*$. R\'eciproquement, \`a un sous-groupe, le th\'eor\`eme associe la compl\'etion au-dessus du support de $\Delta$ du rev\^etement de $X^*$ d\'efini par ce sous-groupe.\end{theorem}

  \begin{corollary}\label{crevram} Il existe un unique rev\^etement $u_{X\vert\Delta}:\overline{(X\vert\Delta)}:=\bar X\to X$ ramifi\'e en au plus $\Delta$, avec $\bar X^0$ (et donc a fortiori $\bar X)$ simplement connexe. On l'appelle le {\bf rev\^etement universel} de $(X\vert\Delta)$. Alors:

 1. Pour tout $j\in J$, soit $n_j$ l'ordre auquel ramifie $u_{X\vert\Delta}$ le long de $\bar D_j^0$: c'est un diviseur de $m_j$. Soit $\bar{\Delta}:=\sum_{j\in J} (1-1/n_j). D_j$: ce diviseur orbifolde divise $\Delta$, et on l'appelle {\bf la r\'egularisation de $\Delta$}. 
  
  2. L'identit\'e de $X$ induit un morphisme orbifolde (divisible) $id_X:(X\vert\Delta)\to (X\vert\bar{\Delta})$ qui est un isomorphisme au niveau des $\pi_1$. De plus:
  
  3. $u_{X\vert\Delta}:\bar X\to X$ est aussi le rev\^etement universel de $(X\vert\bar{\Delta})$, et $u_{X\vert\Delta}$ ramifie en $\bar{\Delta}$. On a: $\overline{(\bar\Delta)}=\bar \Delta$. En particulier: $\bar \Delta=\Delta$ si et seulement si $(X\vert\Delta)$ est r\'eguli\`ere. 
  
  4. Le sous groupe normal de $\pi_1(X^*)$ engendr\'e par les $g_j^{m_j},j\in J$ est aussi le groupe normal engendr\'e par les $g_j^{n_j}$ (notations de \ref{drevram}).

   5. Si $f:(X'\vert\Delta')\to (X\vert\Delta)$ est un morphisme orbifolde divisible bim\'eromorphe, alors le rev\^etement universel de $(X'\vert\Delta')$ est d\'eduit de celui de $(X\vert\Delta)$ par le changement de base $f:X'\to X$ et normalisation. 
   \end{corollary}

{\bf D\'emonstration:} Le rev\^etement universel est donc obtenu en consid\'erant le sous groupe $K$. L'assertion 1 est \'evidente, d'apr\`es les d\'efinitions. L'isomorphisme au niveau des $\pi_1$ r\'esulte de ce que les groupes $K$ coincident pour $\Delta$ et sa r\'egularisation. L'assertion 3 est \'evidente. L'assertion 4 r\'esulte de la d\'efinition du groupe $K$ associ\'e \`a la r\'egularisation de $\Delta$, et de ce qu'il coincide avec celui associ\'e \`a $\Delta$. La propri\'et\'e 5. r\'esulte de ce que les fibres des deux espaces coincident au-dessus de $X^0$, de ce que l'on a un morphisme naturel du changement de base normalis\'e $\overline{X'\times_X\bar X}\to \bar X'$ vers $\bar X'$ au-dessus de $X$, et de ce que les fibres des deux espaces au-dessus de $X$ sont discr\`etes $\square$

 \begin{re} 1. L'orbifolde $(X\vert \Delta)$ et sa r\'egularis\'ee ont donc le m\^eme $\pi_1$, mais pas n\'ecessairement le m\^eme rev\^etement universel si $(X\vert \Delta)$ n'est pas finie. Par exemple: si $D$ est un point de $\Bbb P^1$, alors $(\Bbb P^1\vert D)$ (resp. sa r\'egularis\'ee $(\Bbb P^1\vert 0))$ a pour rev\^etement universel: $(\Bbb P^1-D\vert 0)$ (resp. $(\Bbb P^1\vert 0))$.
 
 2. En g\'en\'eral, $\bar X$ est singulier. Par exemple, si $\Delta:=(1-\frac{1}{m}).(D_1+D_2)$ est le diviseur orbifolde sur $\bP^2$, avec $D_1,D_2$ deux droites projectives distinctes et $2\leq m<+\infty$ entier, alors $\pi_1(\bP^2\vert\Delta)\cong \bZ_m$. De plus, $\bar X$, le rev\^etement universel orbifolde de $(\bP^2\vert\Delta)$ est le c\^one sur la courbe rationnelle normale de degr\'e $m$ dans $\bP^m$ (c'est la conique plane si $m=2)$. Remarquons que le groupe fondamental local de la singularit\'e du rev\^etement universel est ici trivial.\end{re}

  Concernant le rev\^etement universel de $X\in \cal C$, de \ref{gammared}, on d\'eduit:

  \begin{proposition}\label{grev} Soit $X\in \sC$, lisse et connexe. Soit $\gamma_X:X\to Y:=\Gamma(X)$ sa $\Gamma$-r\'eduction (au sens de [Ca94]), ou r\'eduction de Shafarevich. Soit $(\Gamma(X)\vert\Delta(\gamma_X))=[Y\vert\Delta_Y]$ une base orbifolde stable (au sens divisible. Voir \ref{}), et $u_Y:\bar Y\to Y$ son rev\^etement universel orbifolde. On suppose que l'on a modifi\'e $X$ et $Y$ de telle sorte que $\gamma$ soit holomorphe et nette. Donc $(Y\vert\Delta (\gamma_X))=[Y\vert\Delta_Y]$.
  
  D'apr\`es \ref{gammared}, on a une suite exacte, o\`u $\pi_1(X_y)_X$ d\'esigne l'image dans $\pi_1(X)$ de $\pi_1(X_y)$:

    \centerline{\xymatrix{ \pi_1(X_y)_X\ar[r]^{j_*} & \pi_1(X)\ar[r]^{(\gamma_{X})_*}&\pi_1(Y\vert\Delta_Y)\ar[r]&\{1\}\\}}

  Soit $u:\bar X\to X$ le rev\^etement Galoisien de groupe $\pi_1(Y\vert\Delta_Y)$, quotient de $\pi_1(X)$ par le groupe normal fini $\pi_1(X_y)_X$.

  Soit $X'=X\times_Y \bar Y$ le produit fibr\'e normalis\'e. Alors: le morphisme naturel $X'\to X$ se rel\`eve en un morphisme propre et bim\'eromorphe $\beta: X'\to \bar X$.
    \end{proposition}

  {\bf D\'emonstration:} Soit $Y_0\subset Y$ l'ouvert au-dessus duquel $\gamma_X$ est \'equidimensionnel,   et $\Delta(\gamma_X))$ lisse. Soit $\bar Y_0$ son image r\'eciproque dans $\bar Y$. Soit $X_0:=\gamma_X^{-1}(Y_0)$. Puisque $\gamma_X$ est suppos\'ee nette, $\pi_1(X_0)\cong \pi_1(X)$. Il nous suffit de voir que, au-dessus de $Y_0$, l'application naturelle $X'\to X$ est \'etale.   
 Soit $(Y\vert\Delta'_Y)$ l'orbifolde g\'eom\'etrique r\'egularis\'ee de $(Y\vert\Delta_Y)$. Donc $\Delta'_Y$ divise $\Delta_Y$, et $(Y\vert \Delta'_Y)$ a le m\^eme $\pi_1$ et le m\^eme rev\^etement universel que $(Y\vert \Delta_Y)$. Puisque $\Delta'_Y$ divise $\Delta_Y$, un calcul local direct montre alors que la normalisation de $X\times _Y \bar Y$ est \'etale au-dessus de $X_0$. (En effet, chaque composante de chaque fibre de $\gamma_X$ au-dessus de $Y_0$ a une multiplicit\'e multiple de la multiplicit\'e de $\Delta'_Y$ au point image par $\gamma_X$ de cette fibre) $\square$

  \subsection{Groupe fondamental d'une sous-orbifolde.}
  
  Soit $(X\vert\Delta)$ une orbifolde g\'eom\'etrique enti\`ere et lisse, $X$ compacte et connexe. Soit $V\subset X$ une sous-vari\'et\'e complexe ferm\'ee et irr\'eductible , non-contenue dans le support de $\Delta$. Soit $X^*:=X-Supp(\Delta)$, et $V^*:=V-(Supp(\Delta)\cap V)$. Soit $\nu:\hat V\to V$ la normalisation, et $j:\hat V\to X$ la compos\'ee de $\nu$ avec l'inclusion de $V$ dans $X$. On note encore $ j^*:  \hat V^*:=\nu^{-1}(V^*)\to X^*$ la restriction de $j$. On a donc un morphisme naturel de groupes $\hat j_*:\pi_1(\hat V^*)\to \pi_1(X\vert\Delta)$,  compos\'e de $j^*_*:\pi_1(\hat V^*)\to\pi_1(X^*)$  avec la projection: $\pi_1(X^*)\to \pi_1(X\vert\Delta)$. On notera $\pi_1(V)_{(X\vert\Delta)}\subset \pi_1(X\vert\Delta)$ l'image de $\hat j_*$.
  
  \begin{re}\label{rpi} On peut choisir pour $X^*$ tout ouvert de Zariski contenant le compl\'ementaire de $Supp(\Delta)$ sans changer la d\'efinition de $\pi_1(V)_{(X\vert \Delta)}$.
  \end{re}
  
  Soit $u:=u_{(X\vert\Delta)}:\bar X\to X$ le rev\^etement universel, et $\bar V\subset \bar X$ une composante irr\'eductible de $u^{-1}(V)$. Soit $G_{\bar V}\subset \pi_1(X\vert\Delta)$ le sous-groupe stabilisant globalement $\bar V$: ce groupe est d\'efini \`a conjugaison pr\`es dans $\pi_1(X\vert\Delta)$, et on v\'erifie ais\'ement qu'il coincide (\`a conjugaison pr\`es) avec $\pi_1(V)_{(X\vert\Delta)}$.
  
  On a donc la:
  
  \begin{proposition}\label{piv} Si $(X\vert \Delta)$ est finie, $\pi_1(V)_{(X\vert \Delta)}$ est un groupe fini si et seulement si toute composante irr\'eductible $\bar V$ de $u^{-1}(V)$ est compacte. 
  
  Dans le cas g\'en\'eral: $\pi_1(V)_{(X\vert \Delta)}$ est un groupe fini si et seulement si toute composante irr\'eductible $\bar V$ de $u^{-1}(V)$ est $u$-propre sur $V-\Delta_{\infty}$, en notant $\Delta_{\infty}$ la r\'eunion des composantes irr\'eductible de $Supp(\Delta)$  affect\'ees de la multiplicit\'e $+\infty$ \end{proposition}

  On notera $j':(V'\vert \Delta_{V'}) \to (X\vert\Delta)$ une restriction de $\Delta$ \`a $V$, au sens de \ref{fdefrestr}, dans la cat\'egorie $Georb^{div}$. 
  
  Soit $j'_*:\pi_1(V'\vert \Delta_{V'})\to \pi_1(X\vert \Delta)$ le morphisme naturel, et $\pi_1(V')_{(X\vert\Delta)}\subset \pi_1(X\vert\Delta)$ son image. Cette image est d\'efinie seulement \`a conjugaison pr\`es dans $\pi_1(X\vert\Delta)$, et de plus, d\'epend {\it a priori} du choix de $(V'\vert\Delta_{V'})$. Cependant:
  
  \begin{proposition}\label{pirestr} $\pi_1(V')_{(X\vert\Delta)}=\pi_1(V)_{(X\vert\Delta)}$ (\`a conjugaison pr\`es dans $\pi_1(X\vert \Delta))$
  \end{proposition}  
  
  {\bf D\'emonstration:} On a une application naturelle, commutant avec  $\hat j^*$ et $j'^*$, holomorphe propre et \`a fibres connexes $v:V'^*:=j'^{-1}(V^*)\to \hat V^*$. Donc $v_*$ est une surjection au niveau des groupes fondamentaux. Les images dans $\pi_1(X^*)$, et {\it a fortiori} dans $\pi_1(X\vert \Delta)$, de $\pi_1(V'^*)$ et $\pi_1(\hat V^*)$ coincident donc $\square$

  \begin{theorem}\label{ibpir} Soit $g:(Y\vert \Delta_Y)\to (X\vert \Delta)$ un morphisme orbifolde bim\'eromorphe divisible entre orbifoldes lisses, et enti\`eres . Soit $V\subset X$ une sous-vari\'et\'e irr\'eductible, et $W\subset Y$ sa transform\'ee stricte par $g$. On suppose que $V$ (resp. $W)$ n'est pas contenue dans $Supp(\Delta)$ (resp. $Supp(\Delta_Y))$. Alors $g_*(\pi_1(W_{(Y\vert \Delta_Y)})\cong(\pi_1(V_{(X\vert \Delta)})$.
  \end{theorem}
  
  {\bf D\'emonstration:} L'assertion r\'esulte du diagramme commutatif suivant, d\'eduit des applications naturelles dans ce contexte:

  \centerline{
\xymatrix{ \pi_1(\hat W^*)\ar[r]^{\hat j^*}\ar[d]_{(g_W)_*} & \pi_1(Y^*)\ar[d]^{g_*}\ar[r]\ar[d]&\pi_1(Y\vert\Delta_Y)\ar[d]^{\cong}\\
\pi_1(\hat V^*)\ar[r]^{\hat j^*} & \pi_1(X^*)\ar[r]&\pi_1(X\vert\Delta)\\
}}

  Observer que les fl\`eches verticales sont surjectives. On a (\'eventuellement) restreint $X^*$ de telle sorte que $g^{-1}(X^*)$ contienne $Y^*$ (cf. remarque \ref{rpi}) $\square$

    \begin{example}\label{pir} Soit $R\subset (X\vert\Delta)$ une courbe rationnelle orbifolde. Alors  $\pi_1(R)_{(X\vert\Delta)}$ est un groupe fini. En effet, si $R$ admet une restriction divisible $(R'\vert \Delta_{R'})$ de $\Delta$ qui est une courbe rationnelle orbifolde, c'est imm\'ediat,  puisque l'on a un morphisme orbifolde $(R'\vert \Delta_{R'})\to (X\vert\Delta)$ qui induit un morphisme au niveau des groupes fondamentaux, et que le rev\^etement universel de $(R'\vert \Delta_{R'})$ est isomorphe \`a $(\bP^1\vert\Delta')$, avec $\Delta'$ vide ou support\'ee par un seul point. 
    
    De plus, si  $g:(Y\vert \Delta_Y)\to (X\vert \Delta)$ est comme dans \ref{ibpir}, et si $S\subset Y$ est la transform\'ee stricte de $R$, alors $\pi_1(S)_{(Y\vert \Delta_Y)}$ est fini. Remarquer que, en g\'en\'eral, $S$ n'est pas une courbe rationnelle orbifolde, et que si $(S'\vert \Delta_{S'})$ est une restriction de $\Delta_Y$ \`a $S$, alors $\pi_1(S'\vert \Delta_{S'})$ peut-\^etre infini (commensurable \`a un groupe de surface de genre $g\geq 2$ arbitrairement grand), de sorte que l'argument pr\'ec\'edent ne peut \^etre appliqu\'e. 
  \end{example}

  \subsection{$\Gamma$-r\'eduction (ou r\'eduction de Shafarevich) orbifolde.}\label{gamred}

  Nous \'etendons partiellement au cadre orbifolde g\'eom\'etrique certains r\'esultats de [Ca94]: la $\Gamma$-r\'eduction de [Ca94] (voir aussi [Ko93] dans le cas projectif).

 \begin{theorem}\label{gqmer} [Cla08]Soit $(X\vert\Delta)$ une orbifolde g\'eom\'etrique lisse, enti\`ere et {\bf finie}, avec $X\in \sC$. Il existe une unique fibration m\'eromorphe presque holomorphe $g:X\dasharrow Y$ telle que:
 
 1. Pour $y\in Y$ g\'en\'erique, $\pi_1(X_y)_{(X\vert\Delta)}$ est fini.
 
 2. Pour $y\in Y$ g\'en\'eral, si $V\subset X$ est Zariski-ferm\'e irr\'eductible tel que $V$ rencontre $X_y$, mais n'est pas contenu dans $Supp(\Delta)$, et si $\pi_1(V)_{(X\vert\Delta)}$ est fini, alors $V\subset X_y$.
 
 La fibration $g$ est appel\'ee la {\bf $\Gamma$-r\'eduction} (ou {r\'eduction de Shafarevich}) de $(X\vert\Delta)$. On la note: $\gamma_{(X\vert\Delta)}:(X\vert\Delta)\to \Gamma(X\vert\Delta)$. La dimension $dim(Y)$ de la base de cette r\'eduction est appel\'ee la $\Gamma$-dimension de $(X\vert\Delta)$, not\'ee $\gamma d(X\vert\Delta)$. 
 
 On dit que $(X\vert\Delta)$ est de {\bf $\pi_1$-type g\'en\'eral} si $dim(Y)=dim(X)$ ci-dessus.
  \end{theorem}

  {\bf D\'emonstration:} Elle est analogue \`a celle de [94]. On construit l'application de mani\`ere $\pi_1(X\vert \Delta)$-\'equivariante au niveau de $\bar X$, en consid\'erant le quotient m\'eromorphe de $\bar X$ (qui est normal) par la relation d'\'equivalence engendr\'ee par les familles couvrantes de sous-espaces analytiques compacts de $\bar X$. Le seul point nouveau est la necessit\'e, pour \'evaluer le volume des cycles analytiques compacts de $\bar X$, de construire une forme de K\" ahler \'equivariante sur $\bar X$, relev\'ee d'une forme de K\" ahler orbifolde (\`a p\^oles ad\'equats) sur $(X\vert \Delta)$. Voir [Cla08] pour les d\'etails. On obtient alors $\gamma_{(X\vert \Delta)}$ par passage au quotient par $\pi_1(X\vert \Delta)$ $\square$

   \begin{re}\label{pil} Lorsque l'orbifolde $(X\vert \Delta)$ ci-dessus n'est plus finie (par exemple: ``logarithmique"), la conclusion pr\'ec\'edente tombe en d\'efaut: l'application $\gamma_{(X\vert \Delta)}$ n'est pas presque-holomorphe en g\'en\'eral, et ses fibres ne sont pas toujours les images des cycles compacts maximaux ``g\'en\'eraux" de $\bar X$. Ces ph\'enom\`enes sont illustr\'es par l'exemple d\'ej\`a consid\'er\'e de $(\Bbb P^2\vert D)$, $D$ \'etant le diviseur r\'eduit constitu\'e de deux droites distinctes d'intersection $a\in \Bbb P^2$. Dans ce cas, on a bien une (unique) application $\gamma:(\Bbb P^2\vert D)\to \Bbb P^1$ jouissant des propri\'et\'es ci-dessus, dont les fibres sont les droites passant par $a$. Ce ph\'enom\`ene est \'evidemment d\^u au fait que cette orbifolde a, au point $a$, une singularit\'e log-canonique mais non klt. On peut tr\`es certainement construire aussi (un peu diff\'eremment) une $\Gamma$-r\'eduction pour les orbifoldes lisses et enti\`eres non-finies si $X\in \sC$ en consid\'erant les familles couvrantes de cycles $V$ de $X$ dont les composantes irr\'eductibles de l'image inverse dans $\bar X$ sont propres sur $(V-(Supp(\Delta)\cap V))$ (voir proposition \ref{piv}).
   
   \end{re}
 
 \begin{corollary}\label{repif} Soit $(X\vert\Delta)^{div}$ une orbifolde g\'eom\'etrique lisse rationnellement engendr\'ee (ie: $RE^{div})$. Alors $\pi_1(X\vert\Delta)$ est fini.
 \end{corollary}

 {\bf D\'emonstration:} Puisque $(X\vert\Delta)$ est une tour de fibrations m\'eromorphes\footnote{Les bases orbifoldes sont donc consid\'er\'ees ici dans Georb$^{div}$.} \`a fibres orbifoldes RCC, il suffit , par \ref{qre} et \ref{re}, de montrer le r\'esultat lorsque $(X\vert\Delta)$ est RCC. Puisque deux points g\'en\'eriques de $X$ sont alors joints par une chaine connexe de courbes orbifoldes dont les groupes fondamentaux sont finis, l'assertion r\'esulte de \ref{gqmer} de l'exemple \ref{pir}, et du crucial \ref{ibpir} \footnote{On pourrait invoquer aussi \ref{rescougen}.} (qui nous dispense de savoir que la connexit\'e rationnelle est pr\'eserv\'ee par \'equivalence bim\'eromorphe orbifolde) $\square$

  \

  Rappelons ([Ca98]):

  \begin{definition} Une classe $\sG$ de groupes est dite {\bf stable} si:
  
  1. Tout groupe isomorphe \`a $G\in \sG$ est dans $\sG$.
  
  2. Tout quotient, tout sous-groupe d'indice fini de $G\in \sG$ est dans $\sG$. 
  
  3. Tout extension de deux groupes de $\sG$ est dans $\sG$. 
  \end{definition}

  Les exemples les plus simples sont les classes des groupes finis, virtuellement (ou presque) ab\'eliens, virtuellement polycycliques, virtuellement r\'esolubles. La classe des groupes virtuellement nilpotents n'est pas stable.

   \begin{corollary} Soit $(X\vert\Delta)$ une orbifolde g\'eom\'etrique lisse, enti\`ere et finie, avec $X\in \sC$ connexe. Soit $\sG$ une classe stable de groupes. Si deux points g\'en\'eriques de $X$ peuvent \^etre joints par une chaine connexe de sous-vari\'et\'es (analytiques ferm\'ees irr\'eductibles) $V\subset X$, non contenues dans $Supp(\Delta)$, et telles que $\pi_1(V)_{(X\vert\Delta)}\in \sG$, alors $\pi_1(X\vert\Delta)\in \sG$. 
   \end{corollary}

   {\bf D\'emonstration:} Lorsque $\sG$ est la classe des groupes finis, c'est une cons\'equence imm\'ediate de \ref{gqmer}. En g\'en\'eral, on consid\`ere le sous-ensemble $A\in \sC(X)$ constitu\'e des $a$ tels que $\pi_1(Z_a\vert\Delta_{Z_a})_{(X\vert\Delta)}\in \sG$. Il s'agit de montrer que cet ensemble est Z-r\'egulier et stable. Les arguments sont les m\^emes que ceux de la d\'emonstration de \ref{gqmer} \footnote{Il suffit, dans la premi\`ere (resp. derni\`ere) ligne, de remplacer: ``$\pi_1(V)_{(X\vert\Delta)}$ fini" par: ``$\pi_1(V)_{(X\vert\Delta)}\in \sG$" (resp. ``$N,Q$ finis" par: ``$N,Q\in \sG")$}, et sont expos\'es dans [Ca98] et [Ca04] auxquels nous renvoyons $\square$

  \begin{re} Un r\'esultat similaire est certainement vrai aussi lorsque $(X\vert \Delta)$ n'est pas finie. L'\'enonc\'e et la d\'emonstration n\'ecessitent cependant des modifications.
  \end{re}

    \subsection{Finitude r\'esiduelle et crit\`ere d'ab\'elianit\'e.}

  De \ref{crevram} on d\'eduit:

  \begin{corollary}\label{galfin} 1. Il existe un rev\^etement ramifi\'e (Galoisien) fini $g:\bar X\to X$ qui ramifie en $\Delta$ (exactement) si et seulement s'il existe un sous-groupe (normal) $G'$ d'indice fini $G$ de $\pi_1(X\vert\Delta)$ ne contenant, pour chaque $j\in J$, aucun des \'el\'ements $g_j^k$, pour $0<k<m_j$. Cette condition est satisfaite, en particulier, si $\pi_1(X\vert\Delta)$ est r\'esiduellement fini.
  
  2. S'il existe un rev\^etement ramifi\'e (Galoisien) fini $g:\bar X\to X$ qui ramifie en $\Delta$ (exactement), alors la compl\'etion au-dessus de $A\subset X$ du rev\^etement universel (au sens usuel) de $\bar X^0$ est le rev\^etement universel de $(X\vert\Delta)$. 
   \end{corollary}

{\bf D\'emonstration:} 1. Soit $G'<G$ un tel groupe, et $G'^*<\pi_1(X^*)$ son image r\'eciproque. Elle d\'efinit un rev\^etement \'etale fini $X'^*$ de $X^*$. Le groupe $G'^*$ contient $K$, et d\'efinit donc un rev\^etement ramifiant en au plus $\Delta$. Il ramifie en $\Delta$ exactement, puisque, par hypoth\`ese, les lacets $g^{n_j}$ ne sont pas dans $G'$ si $n_j$ est un diviseur strict de $m_j$. 

Si $\pi_1(X\vert\Delta)$ est r\'esiduellement fini, un tel sous-groupe existe (par d\'efinition).

2. Le rev\^etement universel de $\bar X_0$ satisfait alors les conditions de \ref{crevram} caract\'erisant le rev\^etement universel de l'orbifolde g\'eom\'etrique $(X\vert\Delta)$ $\square$

\begin{re}   

1. Dans la situation de \ref{galfin},  si $\sigma:\tilde{X}\to \bar X$ est une r\'esolution des singularit\'es (quotient) de $\bar X$, comme ci-dessus, alors l'inclusion de $\bar X^0\subset \tilde{X}$ induit un morphisme: $j_*: \pi_1(\bar X^0)\to \pi_1(\tilde{X})$ surjectif, mais de noyau infini, en  g\'en\'eral, comme le montre l'exemple (classique) suivant: $E$ une courbe elliptique, $(-1)$ l'involution usuelle sur $E$, $(X\vert\Delta)=((\bP_1\times \bP_1)\vert\Delta)$ le quotient de $E\times E$ par le sous-groupe \`a quatre \'el\'ements engendr\'e par les deux involution agissant sur chacun des facteurs. On choisit pour $\bar X$ la surface de Kummer quotient de $E\times E$ par l'action diagonale des involutions. N\'eammoins, lorsque les groupes fondamentaux locaux des singularit\'es de $\bar X$ sont triviaux, ce morphisme $j_*$ est bijectif. Nous exploiterons ce fait dans \ref{galfinlisse} ci-dessous.
\end{re}

  \begin{theorem}\label{galfinlisse} Soit $(X\vert\Delta)$ lisse, avec $X$ complexe, compacte et connexe. Si $\pi_1(X\vert\Delta)$ est r\'esiduellement fini, il existe un rev\^etement fini Galoisien ramifi\'e en $\Delta$  (exactement) $g:\bar X\to (X\vert\Delta)$  tel que $\pi_1(\bar X^0)\to \pi_1(\bar X)$ soit bijectif, et $\pi_1(\tilde{ X})$ est un sous-goupe d'indice fini de $\pi_1(X\vert\Delta)$, si $\sigma:\tilde{X}\to \bar X$ est une r\'esolution des singularit\'es (quotient) de $\bar X$.
\end{theorem}

{\bf D\'emonstration:} On peut naturellement stratifier $Supp(\Delta)$ comme r\'eunion disjointe d'un nombre fini de sous-vari\'et\'es lisses $V_j$, localement ferm\'ees. A chacune des $V_j$ est attach\'e un sous-groupe fini ab\'elien $G_j$ de $\pi_1(X\vert\Delta)$, bien d\'efini \`a conjugaison pr\`es seulement, et \'egal au groupe fondamental local de $(X\vert\Delta)$ en un point quelconque de $V_j$. Le groupe $\pi_1(X\vert\Delta)$ \'etant suppos\'e r\'esiduellement fini, il admet un sous-groupe normal d'indice fini qui ne rencontre qu'en $\{1\}$ chacun des $G_j$. Le rev\^etement associ\'e \`a un tel sous-groupe satisfait les conditions \'enonc\'ees, puisque les groupes fondamentaux locaux deviennent tous triviaux sur un tel rev\^etement, donc $\pi_1(\bar{X_0})=\pi_1(\bar X)$ $\square$

  \begin{theorem}\label{pifano} Soit $(X\vert \Delta)$ lisse, enti\`ere, finie, et Fano, avec $X$ projective. Si $\pi_1(X\vert \Delta)$ est r\'esiduellement fini, alors il est fini \footnote{On conjecture en \ref{cpi1} que l'hypoth\`ese de finitude r\'esiduelle est superflue.}
  \end{theorem}

{\bf D\'emonstration:} Soit $\bar \Delta$ la r\'egularisation de $\Delta$ au sens de \ref{crevram}. Il existe donc un rev\^etement orbifolde \'etale fini $g:X'\to (X\vert \bar \Delta)$, puisque $\pi_1(X\vert \Delta)=\pi_1(X\vert\bar\Delta)$ est r\'esiduellement fini.  Puisque $\bar \Delta$ divise $\Delta$, et que la paire $(X\vert \Delta)$ est klt, la paire $(X'\vert 0)$ est aussi klt (en fait, $X'$ n'a que des singularit\'es quotient). De plus, $X'$ est ``faiblement Fano", puisque $K_X'=g^*(K_X+\Delta')= g^*(K_X+\Delta)-g^*(\Delta")$, avec $\Delta":=\Delta-\bar\Delta$. Donc $-(K_X+\bar\Delta)=-(K_X+\Delta)+\Delta"$ est somme d'un $\Bbb Q$-ample, et d'un $\Delta"$ tel que la paire $(X\vert \Delta")$ soit aussi klt. Le r\'esultat de [Tak00] s'applique donc, et montre que $X'$ est simplement connexe. Mais le rev\^etement universel de $X'$ (\'egal \`a $X'$, donc) est aussi celui de $(X\vert \Delta)$ et de $(X\vert \bar \Delta)$. Donc $\pi_1(X\vert \Delta)$ est fini $\square$
  
 \begin{re} Le r\'esultat pr\'ec\'edent peut \^etre affin\'e au cas d'orbifoldes Fano klt. Il semble que la th\'eorie $L^2$ puisse \^etre adapt\'ee au cas orbifolde pour \'eliminer l'hypth\`ese de finitude r\'esiduelle
 \end{re}
 
  On a enfin l'extension suivante au cadre orbifolde g\'eom\'etrique de [Delz06] sur les quotients r\'esolubles des groupes de K\" ahler\footnote{Si $\pi_1(X\vert \Delta)$ est r\'esiduellement fini, on peut d\'eduire directement \ref{presolb}.2-3 de [De06].}:
  
  \begin{theorem}\label{presolb}[Ca09] Soit $(X\vert\Delta)$ une orbifolde g\'eom\'etrique lisse, enti\`ere et finie, avec $X\in \sC$, connexe. Alors:
  
  1. $\pi_1(X\vert\Delta)$ a un groupe quotient r\'esoluble, mais non virtuellement nilpotent, si et seulement si l' un des rev\^etements orbifolde \'etales finis de $(X\vert\Delta)$ admet une fibration sur une courbe de genre $g\geq 2$ .
  
  2. Si $\pi_1(X\vert\Delta)$ est r\'esoluble, alors $\pi_1(X\vert\Delta)$ est presque nilpotent. 
  
  3. Si $(X\vert\Delta)$ est sp\'eciale, tout quotient r\'esoluble de $\pi_1(X\vert\Delta)$ est presque ab\'elien; si $\pi_1(X\vert\Delta)$ est lin\'eaire (ie: plongeable dans un $Gl(N,\bC))$, alors $\pi_1(X\vert\Delta)$ est presque ab\'elien.
   \end{theorem}

  \begin{theorem}\footnote{Cet \'enonc\'e est inspir\'e par une discussion avec K. Yamanoi.}\label{ab} Soit $X\in \cal C$ telle que, pour tout rev\^etement fini \'etale $X'$ de $X$, l'application d'Albanese de $X'$ soit surjective. Alors tout groupe quotient $G$ de $\pi_1(X)$ qui est r\'esoluble est presque ab\'elien. 
  
  En particulier, si $X$ est sp\'eciale, tout quotient r\'esoluble de $\pi_1(X)$ est presque ab\'elien. 
  \end{theorem}

  {\bf D\'emonstration:} Il r\'esulte de [Delz06] que $G$ est presque nilpotent (sinon un rev\^etement \'etale fini de $X$ fibre sur une courbe de genre $2$ ou plus, contredisant l'hypoth\`ese sur l'application d'Albanese). Il r\'esulte alors de [Ca95] que $G$ est presque ab\'elien $\square$

  \newpage
  
  \section{CONJECTURES}\label{conj}

  Les conjectures qui suivent sont motiv\'ees par le lemme de d\'evissage \ref{deviss} dans le cas des orbifoldes g\'eom\'etriques sp\'eciales, et par des conjectures standard dans le cas des vari\'et\'es (lisses) avec $\kappa=0$ ou $\kappa_+=-\infty$ (conjecturalement rationnellement connexes). D'autres sont des versions orbifoldes de celles de S. Lang (voir [La86]) concernant l'arithm\'etique et l'hyperbolicit\'e des vari\'et\'es de type g\'en\'eral. Ces conjectures orbifoldes semblent \^etre des interm\'ediaires incontournables pour atteindre les propri\'et\'es conjecturales des vari\'et\'es sp\'eciales (sans structure orbifolde), puis des vari\'et\'es arbitraires, \`a l'aide du ``coeur".
  
  Nous ne cherchons pas \`a \'etablir une liste (pourtant limit\'ee) des cas connus.

  \subsection{Stabilit\'e par d\'eformation et sp\'ecialisation}

  \begin{definition}\label{defdef} Soit $0\in S$ un espace analytique connexe point\'e, et $(X\vert\Delta)$ une orbifolde g\'eom\'etrique lisse, avec $X\in \sC$.  Une {\bf d\'eformation } de $(X\vert\Delta)$ param\'etr\'ee par $S$ est une orbifolde g\'eom\'etrique lisse $(\sX\vert\cal D)$ munie d'une submersion surjective propre et connexe $f:\sX\to S$ dont toutes les fibres sont dans $\cal C$, et telle que:
  
  1. Si $\cal D$$=\sum_j(1-\frac{1}{m_j}).\cal D$$_j$, alors chacune des $\cal D$$_j$ est lisse, et la restriction de $f$ \`a $\cal D$$_j$ est submersive.
  
  2.  $(\sX_0/\cal D$$_0)=(X\vert\Delta)$, posant: $\Delta_s:=\sum_j(1-\frac{1}{m_j}).(\cal D$$_j)_s,\forall s\in S.$\end{definition}

  \begin{re} Nous imposons \`a chacune des $\cal D$$_j$ d'\^etre submersive sur $S$, et non seulement lisse, \`a cause de l'exemple suivant: $(X_s\vert\Delta_s)=(\Bbb P^2\vert C_s),$  o\`u $C_s$ est une conique r\'eduite (affect\'ee de la multiplicit\'e $+\infty)$, $C_s$ \'etant lisse si $s\neq 0$, et r\'eunion de deux droites distinctes si $s=0$. Dans ce cas, le quotient $\kappa$-rationnel de $(\Bbb P^2\vert C_s)$ est un point si $s\neq 0$, et $\Bbb P^1$ si $s=0$. Les conjectures \ref{conjdef'}.(2+4) ci-dessous ne seraient donc pas satisfaites.\end{re}
  
  \begin{conjecture}\label{conjdef}  La dimension essentielle de $(X\vert\Delta)$ est invariante par d\'eformation (ie: $\forall s\in S, ess(X_s\vert\Delta_s)=ess(X_0\vert\Delta_0)$ dans la situation pr\'ec\'edente). En particulier, si $(X_0\vert\Delta_0)$ est sp\'eciale, $\forall s\in S, (X_s\vert\Delta_s)$ est sp\'eciale.
  
  Plus pr\'ecis\'ement: il existe une fibration m\'eromorphe (unique) $c:\sX\to \sC$ au-dessus de $S$ qui induit, pour tout $s\in S,$  une fibration m\'eromorphe $c_s:X_s\to C_s$ qui est le coeur de $(X_s\vert\Delta_s)$. Sur un mod\`ele bim\'eromorphe ad\'equat, la base orbifolde stable de $(\sC/\cal D$$(c,\Delta))$ est une d\'eformation de celle de $(X_0\vert\Delta_0)$. 
  
  L'alg\`ebre essentielle $A(X\vert\Delta)$ est aussi un invariant de d\'eformation. Elle est de type fini. (Voir \ref{c} et \ref{kodL} pour sa d\'efinition).
  \end{conjecture}

  Les deux derni\`eres conjectures sont analogues \`a l'invariance par d\'eformation de l'alg\`ebre canonique, et \`a la finitude de ses g\'en\'erateurs ([Siu01], [Pa05], [BCHM06]).

  Nous allons maintenant voir que ces conjectures peuvent \^etre r\'eduites, dans une certaine mesure, par le d\'evissage de \S.\ref{sdeviss}, au cas de l'invariance par d\'eformation des plurigenres en version orbifolde, et au cas $\kappa_+=-\infty$. Ce d\'evissage sugg\'erant l'invariance par d\'eformation de nouveaux invariants interm\'ediaires d\'efinissant le type (voir d\'efinition \ref{typ}), et en particulier de la longueur $\nu(X\vert\Delta)$. Le lemme \ref{deviss} montre en effet imm\'ediatement que les conjectures \ref{conjdef} r\'esultent des suivantes (et de $C_{n,m}^{orb}$, ou de sa cons\'equence: l'existence du $\kappa$-quotient rationnel):

    \begin{conjecture}\label{conjdef'} 1. $\kappa(X\vert\Delta)$ est invariante par d\'eformation (ie: $\forall s\in S, $   $\kappa(X_s\vert\Delta_s)=\kappa(X_0\vert\Delta_0)$ dans la situation pr\'ec\'edente). L'alg\`ebre canonique $K(X\vert\Delta)$ est de type fini, et invariante par d\'eformation. (Voir \ref{kodL} pour la d\'efinition).
    
   2. La dimension du $\kappa$-quotient rationnel $r^+_{(X\vert\Delta)}$ de $(X\vert\Delta)$ est un invariant de d\'eformation. 
   
  3.  Si $\kappa(X\vert\Delta)\geq 0$, il existe une fibration m\'eromorphe (unique) $\mu:\sX\to \cal M$ au-dessus de $S$ qui induit, pour tout $s\in S,$  une fibration m\'eromorphe $\mu_s:X_s\to M_s$ qui est la fibration de Moishezon-Iitaka de $(X_s\vert\Delta_s)$. Sur un mod\`ele bim\'eromorphe ad\'equat, la base orbifolde stable de $(\cal M/$$\cal D$$(\mu,\Delta))$ est une d\'eformation de celle de $(X_0\vert\Delta_0)$.

  4.  il existe une fibration m\'eromorphe (unique) $r:\sX\to \cal R$ au-dessus de $S$ qui induit, pour tout $s\in S,$  une fibration m\'eromorphe $r^+_s:X_s\to R^+_s$ qui est le $\kappa$-quotient rationnel de $(X_s\vert\Delta_s)$. Sur un mod\`ele bim\'eromorphe ad\'equat, la base orbifolde stable de $(\cal R/$$\cal D$$(r,\Delta))$ est une d\'eformation de celle de $(X_0\vert\Delta_0)$. 
  \end{conjecture}

  \begin{re} Les r\'esultats de [Cla06] et de [BCHM06] \'etablissent respectivement la stabilit\'e par d\'eformation de $K(X\vert\Delta)$ et le fait qu'il soit de type fini (du moins lorsque $\Delta=0)$. La conjecture \ref{-infty} implique facilement l'invariance par d\'eformation du $\kappa$-quotient rationnel.  \end{re}

 \begin{example} (Sugg\'er\'e par une remarque de Y. Tschinkel). Soit $X'_0$ le c\^one sur une hypersurface lisse de degr\'e $(n+2)$ (donc de type g\'en\'eral) de $\bP_n, n\geq 2$. Il est sp\'ecialisation d'une famille d'hypersurfaces lisses de degr\'e $(n+2)$ (donc sp\'eciales) de $\bP_{n+1}$. Soit $X_0$ le transform\'e total de $X'_0$ dans l'\'eclat\'e de $\bP_{n+1}$ en le sommet de $X'_0$. Donc $X_0$ a deux composantes: l'\'eclat\'e de $X'_0$ en son sommet, qui n'est pas sp\'ecial, et une seconde composante isomorphe \`a $\bP_n$. Et $X_0$ est encore sp\'ecialisation de la famille pr\'ec\'edente de vari\'et\'es sp\'eciales. Observons que $X_0$ est $\sS$-connexe (avec la d\'efinition \ref{sconn'} ci-dessous), bien que l'une de ses composantes ne le soit pas.
 \end{example}

\begin{definition}\label{sconn'} Une orbifolde g\'eom\'etrique $(X\vert\Delta)$ est dite $\cal S$-connexe si $X\in \sC$ est de dimension pure, si ses composantes irr\'eductibles sont lisses et se coupent transversalement, si chaque composante $(X_k\vert\Delta_k)$ de $(X\vert\Delta)$ est lisse, et si deux points g\'en\'eriques de $X$ peuvent \^etre joints par une chaine connexe de sous-orbifoldes g\'eom\'etriques sp\'eciales.
\end{definition}

\begin{conjecture}\label{conjspec} Toute sp\'ecialisation d'orbifoldes g\'eom\'etriques lisses sp\'eciales est $\cal S$-connexe.
\end{conjecture}

Cette conjecture ne semble pas pouvoir \^etre simplement d\'eduite par d\'evissage des cas $\kappa=0$ et $\kappa_+=-\infty$.
  
  \
  
 Dans une autre direction:
 
  \begin{conjecture} $\gamma d(X\vert\Delta)$ est invariante par d\'eformation.
  \end{conjecture}
  
  Le cas des familles projectives de vari\'et\'es de dimension $3$ (sans structure orbifolde g\'eom\'etrique) a \'et\'e partiellement r\'esolu par B. Claudon [Cla07].

    \subsection{Groupe fondamental et rev\^etement universel}

    \begin{conjecture}\label{cpi1} Soit $(X\vert\Delta)$ une orbifolde g\'eom\'etrique lisse enti\`ere, avec $X\in \sC$ connexe.
    
    1. Si $(X\vert\Delta)$ est sp\'eciale, alors $\pi_1(X\vert\Delta)$ est presque ab\'elien (not\'e: $\pi_1(X\vert\Delta)\in \widetilde{Ab}$, ceci signifie que ce groupe a un sous-groupe d'indice fini ab\'elien). En particulier:
    
    2. Si $\kappa(X\vert\Delta)=0$, alors $\pi_1(X\vert\Delta)\in \widetilde{Ab}$. 
    
    3. Si $\kappa_+(X\vert\Delta)=-\infty$, ou si $(X\vert \Delta)$ est Fano, alors $\pi_1(X\vert\Delta)$ est presque ab\'elien, et fini si $(X\vert\Delta)$ est finie (voir l'exemple \ref{pil}).   
    
    4. Plus g\'en\'eralement, si $(X\vert \Delta)$ est sp\'eciale $^{div}$ (voir \ref{qspec}), alors $\pi_1(X\vert \Delta)$ est presque ab\'elien.
    
     \end{conjecture}

    Dans ce cas particulier, le lemme de ``d\'evissage" \ref{deviss} permet de r\'eduire la conjecture \ref{cpi1}.(1) \`a ses cas particuliers \ref{cpi1}.(2) et (3).

    \begin{proposition} Supposons la conjecture $C_{n,m}^{orb}$ vraie, ainsi que les \'enonc\'es (2) et (3) de la conjecture \ref{cpi1}. Alors l'\'enonc\'e (1) de la conjecture \ref{cpi1} est \'egalement vrai.
    \end{proposition}

    {\bf D\'emonstration:} Il suffit de montrer que si $f:(X\vert\Delta)\dasharrow Y$ est une fibration dont la base orbifolde stable et la fibre orbifolde g\'en\'erique ont un $\pi_1$ presque ab\'elien, alors $G:=\pi_1(X\vert\Delta)\in \widetilde{Ab}$ aussi. Or $G$ est extension d'un groupe presque ab\'elien de type fini par un autre groupe du m\^eme type, et il est donc polycyclique, et en particulier r\'esoluble lin\'eaire, et r\'esiduellement fini. Il r\'esulte donc de \ref{presolb} que $G$ est presque presque-ab\'elien $\square$

    \begin{re}

    0. On a montr\'e en \ref{presolb} que la conjecture \ref{cpi1}(1) est vraie si $\pi_1(X\vert\Delta)$ est suppos\'e lin\'eaire.
    
    1. Le cas particulier de \ref{cpi1}(3) o\`u $(X\vert\Delta)$ est Fano (ie: $-K_{(X\vert\Delta)}$ ample) est peut-\^etre accessible par les m\'ethodes $L^2$, ou la construction de m\'etriques orbifoldes \`a courbure de Ricci positive lorsque les multiplicit\'es sont finies. (Dans le cas logarithmique, la condition $(X\vert \Delta)$ Fano n'implique en effet pas que: $\kappa_+(X\vert\Delta)=-\infty$).

    2. Le cas particulier o\`u $K_{(X\vert\Delta)}\equiv 0$ est peut-\^etre accessible par l'usage de m\'etriques de K\" ahler-Einstein. (Voir [Ca02] pour le cas tr\`es particulier des vari\'et\'es \`a singularit\'es quotient).

    3. Lorsque $(X\vert\Delta)$ est $RE$, la conjecture \ref{cpi1}.(3) est vraie, par \ref{repif}. La conjecture \ref{cpi1}.(3) est donc une cons\'equence de \ref{repif}, et de la conjecture \ref{-infty}.
    \end{re}

    \begin{example} Soit $(\Bbb P^2\vert\Delta)$, o\`u $\Delta:=\frac{2}{3}.C+\frac{1}{2}.Q$, $C$ (resp. $Q$) \'etant une cubique (resp. conique) lisse, $C$ et $Q$ se coupant transversalement. 
    
    Alors $K_{\Bbb P^2}+\Delta\cong \sO_{\Bbb P^2}$. Si $u:P:=\Bbb P^1\times \Bbb P^1\to \Bbb P^2$ est le rev\^etement double ramifi\'e le long de $Q$, alors $C':=u^{-1}(C)$ est lisse de bidegr\'e $(3,3)$. Le rev\^etement cyclique $S$ de degr\'e $3$ de $P$ ramifi\'e le long de $C'$ est une surface $K3$ (car $h^1(S,\sO_S)=0$, par utilisation de suites spectrales). Donc $\pi_1(\Bbb P^2\vert\Delta)\cong \Bbb Z_6$ (puisque l'on peut aussi d'abord prendre le rev\^etement triple de $\Bbb P^2$ ramifi\'e le long de $C)$. On retrouve ainsi la conclusion de \ref{pp} dans ce cas particulier.
    \end{example}
    
    \begin{question}\label{qrev} Soit $(X\vert\Delta)$ une orbifolde g\'eom\'etrique lisse, enti\`ere et finie, avec $X\in \sC$ connexe. Existe-t-il $Y\in \sC$, normal et \`a singularit\'es terminales, tel que:
    
    1.  $\pi_1(Y)$ et $\pi_1(X\vert\Delta)$ soient commensurables (ie: admettent des sous-groupes d'indices finis isomorphes)?
    
    2. les rev\^etements universels de $Z$ et de $(X\vert\Delta)$ soient analytiquement isomorphes?
    \end{question}

    Lorsque $(X\vert\Delta)$ admet un rev\^etement \'etale fini $\bar X\to (X\vert\Delta)$ comme en \ref{galfin}, et en particulier lorsque $\pi_1(X\vert\Delta)$ est r\'esiduellement fini, la r\'eponse \`a ces deux questions est  ``oui", avec $Z=Y=\bar X$.

     \subsection{Pseudom\'etrique de Kobayashi}\label{pskob}

   On rappelle la notion de pseudom\'etrique de Kobayashi d'une orbifolde g\'eom\'etrique $(X\vert\Delta)$: c'est la plus grande des pseudom\'etriques $d:X\times X\to [0,+\infty[$ telles que $d\leq h^*(d_{\Bbb D})$, pour tout morphisme orbifolde $h:\Bbb D\to (X\vert\Delta)$, $d_{\Bbb D}$ \'etant la m\'etrique de Poincar\'e sur le disque unit\'e $\Bbb D$. On la note $d_{(X\vert\Delta)}$. Consid\'erant seulement les morphismes orbifoldes divisibles $:\Bbb D\to (X\vert\Delta)$, on obtient la pseudom\'etrique {\it classique} (ou {\it divisible)} $d^*_{(X\vert\Delta)}$. On a: $d_{(X\vert\Delta)}\leq d^*_{(X\vert\Delta)}$.

   \begin{conjecture}\label{cd} Soit $(X\vert\Delta)$ une orbifolde g\'eom\'etrique lisse, avec $X\in \sC$, connexe.
   
   1. $d_{(X\vert\Delta)}$ est nulle si et seulement si $(X\vert\Delta)$ est sp\'eciale.
   
   2. Si $(X\vert\Delta)$ est de type g\'en\'eral, il existe un ouvert de Zariski dense $U$ de $X$ tel que la restriction de $d_{(X\vert\Delta)}$ \`a $U\times U$ soit une m\'etrique.
   
   3. Si $c:(X\vert\Delta)\to C(X\vert\Delta)$ est le coeur, alors $d_{(X\vert\Delta)}=c^*(d_{(C(X\vert\Delta)\vert\Delta(c,\Delta))})$.
   \end{conjecture}

     \begin{re} 
     
     \
     
     1. La conjecture (1) pr\'ec\'edente peut \^etre r\'eduite aux cas $\kappa=0$ et $\kappa^+=-\infty$ par d\'evissage si l'on peut montrer qu'une orbifolde g\'eom\'etrique munie d'une fibration a une pseudom\'etrique de Kobayashi nulle s'il en est de m\^ eme pour sa base et ses fibres orbifoldes.
     
     2. L'assertion (2) est une version orbifolde de la conjecture hyperbolique de S. Lang.
     
     3. L'assertion (3) signifie qu'il n'y a pas d'obstruction globale au rel\`evement \`a $(X\vert\Delta)$ des morphismes orbifoldes $h:\Bbb D\to (C(X\vert\Delta)\vert\Delta(c,\Delta))$, base orbifolde du coeur.
     
     4. La conjecture \ref{cd}(1) est \'etablie pour les courbes dans [C-W05], o\`u l'on montre aussi que le lemme de Brody reste valable en version orbifolde g\'eom\'etrique.Dans [C-P 05], une version orbifolde des th\'eor\`emes d'hyperbolicit\'e de Bogomolov et M$^c$Quillan est \'etablie pour certaines surfaces orbifoldes avec $(c_1^2-c_2)>0$. Voir [R 08] pour de nombreux compl\'ements int\' eressants et am\'eliorations sur ce sujet. 
     \end{re}

     \
     
     On peut tenter de relier l'annulation de la pseudom\'etrique de Kobayashi orbifolde \`a l'existence de courbes enti\`eres orbifoldes.

     \begin{conjecture} Soit $(X\vert\Delta)$ une orbifolde lisse, avec $X\in \cal C$. On a \'equivalence entre les propri\'et\'es suivantes:
     
    1. $(X\vert\Delta)$ est sp\'eciale.
    
    2. $d_{(X\vert\Delta)}\equiv 0$.
    
    3. Il existe une $\Delta^{div}$-courbe enti\`ere  (ie: un morphisme $h:\bC\to (X\vert\Delta)^{div})$ dont l'image est Zariski-dense dans $X$.
    
    4. Il existe une $\Delta^{div}$-courbe enti\`ere dont l'image est dense (pour la topologie analytique) dans $X$.
    
    5.  Tout sous-ensemble fini de $X$ est contenu dans une $\Delta^{div}$-courbe enti\`ere dont l'image est dense dans $X$. 
    \end{conjecture}

    \begin{example} 
    \
    
    1. Soit $X$ une vari\'et\'e quasi-projective lisse, et $X=(\overline{X}-D)$ une compactification de $X$ lisse telle que $D$ soit un diviseur \`a croisements normaux de $\overline{X}$. On dit que $X$ est sp\'eciale si $(\overline{X}\vert D)$ est sp\'eciale. Cette d\'efinition ne d\'epend pas de la compactification choisie. La conjecture pr\'ecedente affirme donc en particulier que $X$ est sp\'eciale si et seulement si elle contient une courbe enti\`ere (au sens usuel) Zariski-dense. 
    
    2. Un r\'ecent (Novembre 2008) preprint de J. Winkelmann caract\'erise les surfaces quasi-projectives $X$ dont l'application quasi-Albanese a une image de dimension $2$, et qui contiennent une courbe enti\`ere Zariski-dense. Sa caract\'erisation semble coincider dans ce cas avec le fait d'\^etre sp\'eciale, et donc fournir pour cette classe de surfaces une solution (positive) de la conjecture pr\'ec\'edente.
    
    3. Si $(X\vert \Delta)$ lisse, enti\`ere, avec $X\in \cal C$ contient une $\Delta^{div}$-courbe enti\`ere Zariski-dense, elle est d'apr\`es la conjecture pr\'ec\'edente, sp\'eciale, et son groupe fondamental est donc presque-ab\'elien. En particulier, une vari\'et\'e quasi-projective contenant une courbe enti\`ere (au sens usuel) Zariski-dense a un groupe fondamental presque-ab\'elien. 
    \end{example}

       \subsection{Points rationnels: corps de fonctions}\label{fonct}

      Soit $B$ une courbe projective complexe lisse et connexe, et $k_B:=\bC(B)$ le corps de ses fonctions m\'eromorphes. Soit $f:X\to B$ une application holomorphe surjective et connexe, $X$ \'etant une vari\'et\'e projective complexe lisse (et connexe, donc). Pour $b\in B$, on note $X_b$ la fibre sch\'ematique de $f$ au-dessus de $b$. On suppose $X$ muni d'une structure d'orbifolde g\'eom\'etrique lisse $(X\vert\Delta)$ qui induit, pour $b\in B$ g\'en\'erique, une structure d'orbifolde g\'eom\'etrique lisse $(X_b\vert\Delta_b)=(X\vert\Delta)_b$ sur $X_b$. On dira simplement que $(X\vert\Delta)$ est une orbifolde g\'eom\'etrique lisse d\'efinie sur $k_B$. On dira que $(X\vert\Delta)$ est sp\'eciale (resp. de type g\'en\'eral, etc...) s'il en est de m\^eme pour $(X_b\vert\Delta_b)$, pour $b\in B$ g\'en\'eral. Voir [Ca 01] pour plus de d\'etails.

      Un point $k_B$-rationnel $s$ de $X$ sur $k_B$ est une section $s:B\to X$ de $f$. On note $X(k_B)$ l'ensemble de ces points. (Cet ensemble est essentiellement un invariant birationnel de $(X,f)$: si $(X',f')$ est un second mod\`ele birationnel de $f$, $X(k_B)$ et $X'(k_B)$ coincident sur un ouvert de Zariski non vide commun).
      
      Soit $S\subset B$ un sous-ensemble fini au dessus du compl\'ementaire duquel $(X\vert\Delta)$ a {\it bonne r\'eduction}, c'est-\`a-dire est tel que $(X_b\vert\Delta_b)$ est lisse si $b\notin S$. 
      
      Si $\Delta=\sum_{j\in J}(1-\frac{1}{m_j})D_j$, si $s\in X(k_B)$, si $b\notin S$, et si $j\in J$, on note $(s.D_j)_{b}$ l'ordre de contact en $s(b)$ de $s(B)$ avec $D_j$. C'est un entier positif (ou nul).
      
      On d\'efinit alors $(X\vert\Delta)(k_B,S)$ comme le sous-ensemble des $s\in X(k_B)$ tels que $\forall j,b\notin S$, on ait: $(s.D_j)_{b}\geq m_j$ si $(s.D_j)_{b}\geq 1$. Autrement dit: l'ordre de contact doit \^etre {\it au moins} \'egal \`a $m_j$ si $s(b)\in D_j$, ceci pour tous $j\in J, b\notin S$. 
      
      On note $(X\vert\Delta)(k_B)$ la r\'eunion des $(X\vert\Delta)(k_B,S')$, lorsque $S'\subset B$ est finie, de compl\'ementaire de bonne r\'eduction au sens pr\'ec\'edent. Cet ensemble est essentiellement un invariant birationnel de $((X\vert\Delta), f)$. 
      
      On peut d\'efinir ces notions en version {\it classique}: $(X\vert\Delta)^*(k_B,S)$ est alors le sous-ensemble des $s\in X(k_B)$ tels que $\forall j,b\notin S$, on ait: $(s.D_j)_{b}$ est {\bf  divisible par} $m_j$. On a, bien s\^ur: $(X\vert\Delta)^*(k_B,S)\subset (X\vert\Delta)(k_B,S)$.

      Les extensions finies de corps $k'/k$ correspondent bijectivement aux morphismes finis $B'\to B$, $B'$ courbe projective lisse et connexe, posant: $k'=k_{B'}$. Un morphisme $f:(X\vert\Delta)\to B$ comme ci-dessus d\'efinit alors par changement de base (et d\'esingularisation) un morphisme $f':(X'\vert\Delta')\to B'$, et des inclusions $(X\vert\Delta)(k_B,S)\subset (X\vert\Delta)(k_{B'},S')$, avec $S'$ image inverse de $S$ dans $B'$.

      Si $g:(X\vert\Delta)\to (Y\vert\Delta_Y)$ est un morphisme orbifolde au-dessus de $B$, il induit naturellement une application $f:(X\vert\Delta)(k_B,S)\to (Y\vert\Delta_Y)(k_B,S), \forall S$. Si $f$ est un morphisme orbifolde {\it divisible (ou classique)}, il induit de m\^eme une application $f:(X\vert\Delta)^*(k_B,S)\to (Y\vert\Delta_Y)^*(k_B,S)$.

      On dit que $f:(Y\vert\Delta_Y)\to B$ est {\bf isotrivial} s'il existe $B'\to B$ fini, et $(F\vert\Delta_F)$ une orbifolde g\'eom\'etrique lisse telle que $(Y\vert\Delta_Y)\times _BB'$ soit birationnel au-dessus de $B'$ \`a $(F\vert\Delta_F)\times B'$. 
      
      On dit que $(X\vert\Delta)$ n'a pas de quotient isotrivial (sur $k_B)$ s'il n'existe pas de morphisme orbifolde m\'eromorphe $g:(X\vert\Delta)\dasharrow (Y\vert\Delta_Y)$dominant au-dessus de $B$ tel que $(Y\vert\Delta_Y)$ soit isotrivial.

\begin{conjecture}\label{conjcorfonct} Soit $(X\vert\Delta)$ une orbifolde g\'eom\'etrique lisse d\'efinie sur $k_B$, et sans quotient isotrivial.

1. $(X\vert\Delta)$ est sp\'eciale si et seulement s'il existe une extension finie $k_{B'}/k_B$ telle que $(X\vert\Delta)(k_{B'})$ soit Zariski dense (ie: tel que la r\'eunion des $s(B), s\in (X\vert\Delta)(k_{B'})$ soit Zariski dense dans $X_{B'})$.

2. $(X\vert\Delta)$ est de type g\'en\'eral si et seulement s'il existe un ouvert de Zariski dense $U$ de $X$ tel que, pour toute extension finie $k_{B'}/k_B$, l'ensemble des $s\in(X\vert\Delta)(k_{B'})$ tels que $s(B')$ rencontre $U$ soit fini.
 \end{conjecture}

\begin{re}

\

1. La conjecture \ref{conjcorfonct}.2. est la version orbifolde de la conjecture corps de fonctions de Bombieri-Lang. 

2. Le coeur montre que si $(X\vert\Delta)$ n'est pas sp\'eciale, et si  \ref{conjcorfonct}.2. est vraie, alors $(X\vert\Delta)(k_{B'})$ n'est Zariski dense pour aucune extension finie $k_{B'}/k_B$, ce qui \'etablit  \ref{conjcorfonct}.1. dans ce cas. Plus pr\'ecis\'ement, si  \ref{conjcorfonct}.2. est vraie, alors $(c_{(X\vert\Delta)}(X\vert\Delta)(k'))\cap U$ est fini, pour toute extension finie $k_{B'}/k_B$, et pour un ouvert de Zariski non vide $U$ de $C(X\vert\Delta)$. 

3. La densit\'e potentielle des orbifoldes g\'eom\'etriques lisses sp\'eciales peut \^etre conjecturalement r\'eduite aux cas $\kappa=0$ et $\kappa_+=-\infty$ par le m\^eme ``d\'evissage" et sous les m\^emes hypoth\`eses que dans les cas pr\'ec\'edents.

4. Voir [Ca 05] pour un cas particulier en dimension relative $2$. 
\end{re}

      \subsection{Points rationnels: arithm\'etique}\label{arithm}

  Ce cas est analogue au pr\'ec\'edent (aux questions usuelles d'isotrivialit\'e pr\`es). On renvoie \`a [Abr06] en particulier pour une pr\'esentation et une discussion d\'etaill\'ee des notions pr\'esent\'ees ici.

  Si $X,\Delta=\sum_{j\in J}(1-\frac{1}{m_j}).D_j$ sont d\'efinis sur un corps de nombres $k$, $X$ lisse et $Supp(\Delta)$ \`a croisements normaux, et si $\cal X$, $\cal D$$=\sum_{j\in J}(1-\frac{1}{m_j}).\cal D$$_j$ sont des mod\`eles de $X, D_j,\forall j$ d\'efinis et de bonne r\'eduction (ie: si la r\'eduction de $\Delta$ reste \`a croisements normaux pour toute place $v\notin S$) sur l'anneau des entiers $\cal O$$_{k,S}$, on note alors $(X\vert\Delta)(\cal O$$_{k,S})$ l'ensemble des points $\cal O$$_{k,S}$-int\'egraux  $x\in X(\cal O$$_{k,S})$ tels que pour toute place $v\notin S$ de $\cal O$$_k$, et tout $j\in J$, le nombre d'intersection arithm\'etique $(x.\cal D$$_{j})_{v}$ de $x$ avec $\cal D$$_j$ en $v$ est: soit nul, soit sup\'erieur ou \'egal \`a $m_j$. 
  
  Pour tout $k,S$ fix\'es, les points int\'egraux de $(X/ \Delta)$ sont essentiellement ind\'ependants des mod\`eles choisis: deux mod\`eles \'etant choisis, ces points int\'egraux coincident sur un ouvert de Zariski non vide (``commun") de $X$.

  On peut introduire la notion plus restrictive de points int\'egraux {\bf classiques} en imposant la condition que les nombres d'intersection arithm\'etiques  $(x.\cal D$$_{j})_{v}$ soient divisibles par $m_j$. Voir [DG98] dans le cas des courbes pour les multiplicit\'es ``classiques", et [Ca05] pour les diff\'erences et motivations. On notera $(X\vert\Delta)^*(\cal O$$_{k,S})$ l'ensemble de ces points int\'egraux ``classiques". On a \'evidemment (pour un mod\`ele fix\'e): $(X\vert\Delta)^*(\cal O$$_{k,S})\subset (X\vert\Delta)(\cal O$$_{k,S})$.

\

Si $f: (X\vert\Delta)\to (Y\vert\Delta')$ est un morphisme orbifolde d\'efini sur $k$, il induit une application naturelle $f:(X\vert\Delta)(\cal O$$_{k,S})\to (Y\vert\Delta')(\cal O$$_{k,S})$, sur des mod\`eles sur lesquels $f$ est d\'efini. Si ce morphisme est un morphisme orbifolde {\it divisible}, il induit aussi une application $f:(X\vert\Delta)^*(\cal O$$_{k,S})\to (Y\vert\Delta')^*(\cal O$$_{k,S})$.

\begin{conjecture}\label{conjarith} Soit $(X\vert\Delta)$ une orbifolde g\'eom\'etrique lisse d\'efinie sur $k$, un corps de nombres \footnote{On pourrait formuler cette conjecture, plus g\'en\'eralement, pour $k$ de type fini sur $\bQ$.}. On suppose fix\'e un mod\`ele d\'efini sur $\sO$$_{k,S}$ comme ci-dessus.

1. $(X\vert\Delta)$ est sp\'eciale si et seulement s'il existe une extension finie $k'/k$ telle que $(X\vert\Delta)(k')$ soit Zariski dense pour un, et donc tout mod\`ele. (``Densit\'e potentielle"). 

2. $(X\vert\Delta)$ est de type g\'en\'eral si et seulement s'il existe un ouvert de Zariski dense $U\subset X$ tel que, pour toute extension finie $k'/k$, $(X\vert\Delta)(k')\cap U$ soit fini pour un, et donc tout mod\`ele.

On peut formuler ces deux conjectures \`a la fois pour les points int\'egraux ``classiques" et ``non-classiques". La premi\`ere (resp. la seconde) est plus forte en version ``classique" (resp. ``non-classique").
\end{conjecture}

\begin{re}

\

1. La conjecture \ref{conjarith}.2. n'est qu'une version orbifolde g\'eom\'etrique de la conjecture arithm\'etique de Bombieri-Lang. Remarquons que cette conjecture est ouverte m\^eme pour les courbes, pour lesquelles elle r\'esulte cependant de la conjecture $abc$ (voir [Ca05]). La version corps de fonctions complexes est cependant \'etablie dans [Ca05]. La version ``classique" peut \^etre cependant d\'eduite du th\'eor\`eme de Faltings (voir [DG98]).

2. Le coeur montre que si $(X\vert\Delta)$ n'est pas sp\'eciale, et si  \ref{conjarith}.2. est vraie, alors $(X\vert\Delta)(k')$ n'est Zariski dense pour aucune extension finie $k'/k$, ce qui \'etablit  \ref{conjarith}.1. dans ce cas. Plus pr\'ecis\'ement, si  \ref{conjarith}.2. est vraie, alors $(c_{(X\vert\Delta)}(X\vert\Delta)(k'))\cap U$ est fini, pour toute extension finie $k'/k$, et pour un ouvert de Zariski non vide $U$ de $C(X\vert\Delta)$. 

3. La densit\'e potentielle des orbifoldes g\'eom\'etriques lisses sp\'eciales peut \^etre conjecturalement r\'eduite aux cas $\kappa=0$ et $\kappa_+=-\infty$ par le m\^eme ``d\'evissage" que dans les cas pr\'ec\'edents.

\end{re}

\subsection{Multiplicit\'es ``classiques" et ``non-classiques".}

La consid\'eration des multiplicit\'es ``non-classiques" (bas\'ee sur $inf$ et non $pgcd$) est justifi\'ee par plusieurs raisons (voir aussi divers aspects dans [Ca05]):

1. La compatibilit\'e exacte avec les faisceaux de diff\'erentielles sym\'etriques (prop. \ref{modiff}).

2. La bijection entre faisceaux de Bogomolov satur\'es et fibrations de type g\'en\'eral (th\'eor\`eme \ref{bijbog}), plus g\'en\'eralement, l'\'egalit\'e entre $\kappa(X\vert \Delta)$ et $\kappa(f\vert \Delta)$ (voir th\'eor\`eme \ref{kappa=L}).

3. L'existence de multisections orbifoldes locales d'une fibration est assur\'ee si la base orbifolde est d\'efinie avec les multiplicit\'es $inf$, mais non avec les multiplicit\'es $pgcd$. Ce qui justifie la pr\'eservation conjecturale de certaines propri\'et\'es par ``extension" orbifolde (i.e: si fibres et base orbifolde d'une fibration la poss\`edent, alors l'orbifolde ambiante aussi).

Une propri\'et\'e importante des multiplicit\'es $pgcd$ qui est perdue par les multiplicit\'es $inf$ est la compatibilit\'e avec le groupe fondamental (voir \ref{pi1div}).

En sens inverse, on peut attendre, parfois, un comportement similaire pour les orbifoldes d\'efinies par ces deux types de multiplicit\'es. On a par exemple (question \ref{qspec}):

\begin{question}\label{qspec'} Est-il vrai que: $(X\vert \Delta)$ sp\'eciale $^{div}\Longleftrightarrow (X\vert \Delta)$ sp\'eciale?
\end{question}

On a vu qu'une r\'eponse affirmative \`a \ref{qspec'} r\'esulterait d'une telle r\'eponse \`a la question \ref{qmult} suivante.

\begin{question}\label{qmult}

1. Soit $(X\vert \Delta)$ lisse, enti\`ere et finie, $X$ dans $\sC$. Soit $f:(X\vert \Delta)\to Y$ une fibration. On suppose que les fibres orbifoldes de $f$ sont sp\'eciales. A-t-on alors: $\Delta(f\vert \Delta)=\Delta^*(f\vert \Delta)$? (On suppose que $\Delta^{vert}=0)$\end{question}

Cette propri\'et\'e est satisfaite lorsque $\Delta=0$ si $X_y$ est rationnellement connexe ($\Delta(f)=0)$, ou si $X_y$ est un tore complexe. Le premier cas non-trivial est celui dans lequel $X_y$ est une surface $K3$ et $Y$ une courbe. Remarquons que cette propri\'et\'e est en d\'efaut d\'ej\`a lorsque $X_y$ est une courbe de type g\'en\'eral. Remarquons enfin qu'il suffit de v\'erifier la propri\'et\'e \ref{qmult} lorsque $Y$ est une courbe (projective).

Si cette propri\'et\'e de \ref{qmult} est satisfaite, et si $c:(X\vert \Delta)\to C$ est le ``coeur" d'une orbifolde lisse $(X\vert \Delta)$, avec $X\in \sC$, alors $\pi_1(X\vert\Delta)$ est extension de $\pi_1(F)$ (conjecturalement presque ab\'elien) par $\pi_1(B)$, si $F$ (resp. $B)$ est la fibre orbifolde g\'en\'erique (resp la base orbifolde) de $c$.

Une question apparent\'ee \`a \ref{qmult} est la conjecture \ref{conjarith}, qui affirme l'\'equivalence quantitative des ensembles de points entiers``classiques" et ``non-classiques". De mani\`ere analogue, la conjecture \ref{-infty}  et la question \ref{qha} postulent l'\'equivalence quantitative des ensembles de courbes $\Delta$-rationnelles ``divisibles" et ``non-divisibles".

 \subsection{Formes diff\'erentielles}

  Soit $(X\vert\Delta)$ une orbifolde g\'eom\'etrique lisse avec $X\in \cal C$. 
  
  Soit $q\geq 0$ un entier, et $\Omega^q(X\vert\Delta):=\oplus_{N\geq 0}H^0(X,S^N_q)$ l'alg\`ebre des $q$-formes diff\'erentielles sym\'etriques sur $(X\vert\Delta)$. 
  
  Pour $x\in X, x\notin Supp(\Delta)$, on a une application naturelle d'\'evaluation en $x$: $ev^q_x: \Omega^q(X\vert\Delta)\to \oplus_{N\geq 0}Sym^N(\Omega^q_{X,x})$.

 \begin{conjecture}\label{diff}
 \
 
 1. Si $\kappa_+(X\vert\Delta)=-\infty$, alors $\Omega^q(X\vert\Delta)=\bC,\forall q\geq 0$. 
 
 2. Si $\kappa(X\vert\Delta)=0$, alors $ev^q_x$ est injective pour $x\in X$ g\'en\'erique.
 \end{conjecture}

    \begin{re}\label{diff'}
    \
    
    0. La conjecture pr\'ec\'edente entraine en particulier que: si $L\subset \Omega^p_X,p>0$, on a: $\kappa((X\vert\Delta),L)\leq 0$ si $\kappa(X\vert\Delta)=0$, et: $\kappa((X\vert\Delta),L)=-\infty$ si $\kappa_+(X\vert\Delta)=-\infty.$
    
1. La conjecture \ref{diff}.1 est vraie pour les vari\'et\'es rationnellement connexes (avec $\Delta=0)$. Elle l'est aussi pour les orbifoldes rationnellement connexes (au sens orbifolde). Elle doit pouvoir \^etre \'etablie au moins dans le cas des multiplicit\'es enti\`eres et finies \`a l'aide de m\'etriques orbifoldes \`a courbure de Ricci positive.

2. La conjecture \ref{diff}.2 est vraie pour les vari\'et\'es K\" ahl\'eriennes compactes avec $c_1(X)=0$ (et $\Delta=0)$, par la solution de la conjecture de Calabi sur l'existence de m\'etriques de K\" ahler Ricci-plates et le parall\'elisme des $q$-formes sym\'etriques qui en r\'esulte. La conjecture de K. Ueno ($h^0(X,\Omega^q_X)\leq \frac{n!}{q!(n-q)!}$ si $\kappa(X)=0$, avec $n:=dim(X))$ est un cas particulier de \ref{diff}. Elle doit pouvoir \^etre \'etablie au moins dans le cas des orbifoldes g\'eom\'etrique \`a multiplicit\'es enti\`eres et finies \`a l'aide de m\'etriques orbifoldes \`a courbure de Ricci nulle.

3. Une conjecture sur  l'alg\`ebre $\Omega^q(X\vert\Delta)$ lorsque $(X\vert\Delta)$ est sp\'eciale semble plus difficile \`a formuler, le d\'evissage expos\'e au \S\ref{sdeviss} ne fournissant pas  de structure simple apparente, m\^eme en admettant \ref{diff}. 

\end{re}

\subsection{Familles de vari\'et\'es canoniquement polaris\'ees.}

 De fa\c con vague, il est conjectur\'e que ``l'espace des modules fin" des vari\'et\'es de type g\'en\'eral a des composantes irr\'eductibles qui sont elles-m\^emes de type g\'en\'eral.

 Dans le cas des vari\'et\'es \`a fibr\'e canonique ample (dont les courbes de genre $g\geq 2$ fournissent l'exemple classique), on va pr\'esenter des r\'esultats et une conjecture plus pr\'ecise, dont le cadre naturel semble \^etre justement les notions d'orbifoldes g\'eom\'etriques logarithmiques sp\'eciales et de ``coeur" d\'evelopp\'ees dans le pr\'esent texte.

Soit $g:V\to B$ un morphisme projectif, submersif et \`a fibres connexes de $V$, lisse, sur une base $B$ quasi-projective connexe. On supposera que $B=\bar{B}-D$, o\`u $\bar{B}$ est projective lisse, et $D$ un diviseur \`a croisements normaux de $\bar{B}$. Donc $B$ n'est autre que l'orbifolde g\'eom\'etrique logarithmique $(\bar{B}/D)$. On notera:  $\overline{\kappa}(B):=\kappa(\bar{B}/D)$.

On suppose que $g$ est une {\bf famille de vari\'et\'es canoniquement polaris\'ees}, c'est-\`a-dire satisfait les conditions pr\'ec\'edentes, et que, de plus, le faisceau canonique relatif $K_{V/B}$ est ample sur toutes les fibres $V_b,b\in B$ de $g$.

On note alors:  $Var(g)$ le rang de l'application de Kodaira-Spencer $ks_g(b):TB_b\to H^1(V_b,TV_b)$ au point g\'en\'erique de $B$. Donc $Var(g)=0$ si et seulement si $g$ est isotriviale (ie: ses fibres sont deux-\`a-deux isomorphes).

Notre conjecture \footnote{On pourrait la formuler, plus g\'en\'eralement, lorsque les fibres $V_b$ ont un fibr\'e canonique  semi-ample, ou m\^eme nef.} est la suivante:

\begin{conjecture}\label{conjmod} Si $B=(\bar{B}/D))$ est sp\'eciale, la famille $g$ est isotriviale.

Donc, $Var(g)\leq dim(C(\bar{B}/D))$, la dimension du ``coeur" de la base orbifolde g\'eom\'etrique logarithmique de la famille consid\'er\'ee, dans le cas g\'en\'eral.

En effet, pour $B$ arbitraire, la restriction de $g$ au-dessus des fibres du ``coeur" $c_B:\bar{B}\to C(\bar{B}/D)$ de $B$ serait isotriviale, et la ``variation" de $g$ se factoriserait par $C(\bar{B}/D)$ (ie: au point g\'en\'erique $b$ de $B$, l'application de Kodaira-Spencer $ks_g(b)$ s'annule sur l'espace tangent \`a la fibre de $c_B$ en $b)$.
  \end{conjecture}

  Cette conjecture g\'en\'eralise et renforce en les pr\'ecisant consid\'erablement les $3$ conjectures ant\'erieures $A,B,C$ suivantes:
  
  A. (Viehweg, voir [V-Z02]): si $Var(g)=dim(B)$, alors $B$ est de type g\'en\'eral (ie: $\overline{\kappa}(B):=dim(B))$. En effet, cette \'egalit\'e entrainerait alors: $dim(B)=Var(g)\leq dim(C(\bar{B}/D)\leq dim(B)$, or $dim(B)=dim(C(\bar{B}/D)$ si et seulement si $B$ est de type log-g\'en\'eral.
  
   B. ([Ke-Kov06] Si $\overline{\kappa}(B)=0$, alors $Var(g)=0$ (c'est-\`a-dire que $g$ est isotriviale). En effet, la condition $\overline{\kappa}(B)=0$ entraine que $B$ est Log-sp\'eciale. 
  
    C. ([Ke-Kov06]  Si $\overline{\kappa}(B)=-\infty$, alors $Var(g)\leq (dim(B)-1)$. En effet, si $\overline{\kappa}(B)=-\infty$, $B$ n'est pas de type Log-g\'en\'eral. Et $dim(C(\bar{B}/D)\leq (n-1).$

    La d\'ecomposition conditionnelle du coeur montre que les deux cas cruciaux dans lesquels la conjecture devrait \^etre d'abord \'etablie sont les cas o\`u $\overline{\kappa_+}(B)=-\infty$ (et aussi, en particulier, lorsque $(\bar{B}/D)$ est Fano), et lorsque $\overline{\kappa}(B)=0$ (en particulier lorsque $c_1(\bar{B}/D)=0).$

    Cependant, la conjecture \ref{conjmod} traite aussi les nombreux cas dans lesquels $(\bar{B}/D)$ ``fibre" de mani\`ere it\'er\'ee (au sens orbifolde) avec des fibres de l'un des deux types pr\'ec\'edents. Le premier cas non trait\'e par les conjectures A,B,C ci-dessus \'etant celui dans lequel $dim(B)=2$ et soit $\overline{\kappa}(B)=-\infty$, soit $\overline{\kappa}(B)=1$. Il est r\'esolu ci-dessous, \`a titre d'exemple (\`a une exception pr\`es, si $\overline{\kappa}(B)=1)$.

  \

 La  conjecture \ref{conjmod} est d\'emontr\'ee lorsque $dim(B)=1$ ( [Kov96], [Kov00], g\'en\'eralisant le cas classique dans lequel les $V_b$ sont des courbes de genre $g\geq 2)$.

   Lorsque $dim(B)=2$, les conjectures A, B et C ci-dessus sont d\'emontr\'ees dans [Ke-Kov06]\footnote{Une d\'emonstration de la conjecture \ref{conjmod} lorsque $dim B\leq 3$ vient d'\^etre annonc\'ee dans [J-K 09].}. 
   
 \
 
 Cependant, lorsque $dim(B)=2$, il existe de nombreux cas ($B$ sp\'eciale avec $\overline{\kappa}(B)=1$ ou $-\infty)$ dans lesquels la conjecture \ref{conjmod} renforce les conjectures A,B,C (la conclusion \'etant alors: ``$g$ isotriviale", et non: ``$Var(g)\leq 1")$.

\begin{theorem}\label{B2} Soit $g:V\to B$ une famille de vari\'et\'es canoniquement polaris\'ees (dans le sens ci-dessus). On suppose que $dim(B)=2$, et que $B$ est sp\'eciale. Alors $g$ est isotriviale (sauf, peut-\^etre, si $\overline{\kappa}(B)=1$ et si la fibration de Moishezon-Iitaka de $B$ est isotriviale).
\end{theorem}

  {\bf D\'emonstration:} Tout comme dans [Ke-Ko06], la d\'emonstration r\'esulte essentiellement du r\'esultat suivant de [V-Z]:

  \begin{theorem}\label{vz} Soit $g:V\to B$ une famille de vari\'et\'es canoniquement polaris\'ees. Il existe alors un entier $N>0$ et un sous-fibr\'e $L$ de rang $1$ de $S^N_{1}(B):=S_{N,1}(\bar{B}/D)=Sym^N(\Omega^1_{\bar{B}}(log(D))$ tel que $\kappa(\bar{B}, L)\geq Var(g)$.
  \end{theorem}

    Le th\'eor\`eme \ref{B2} r\'esulte alors du:

    \begin{lemma} \label{sym} Soit $B=(\bar B/D)$ une orbifolde g\'eom\'etrique lisse, projective et logarithmique sp\'eciale de dimension $2$. Pour tout $N>0$, et pour tout sous-fibr\'e $L$ de rang $1$ de  $S^N_{1}(B)$, on a alors: $\kappa(\bar B,L)\leq 0$ (sauf, peut-\^etre, si $\overline{\kappa}(B)=1$ et si la fibration de Moishezon-Iitaka de $B$ est isotriviale).
    \end{lemma}

    {\bf D\'emonstration:} Nous traitons successivement les cas $\overline{\kappa}(B)=0$, $\overline{\kappa}(B)=-\infty$, et $\overline{\kappa}(B)=1$. Les deux premiers cas \'etant d\'ej\`a essentiellement connus.

     Lorsque $\overline{\kappa}(B)=0$, le r\'esultat est, en fait, \'etabli dans [Ke-Kov06]. Il nous reste \`a traiter les cas $\overline{\kappa}(B)=-\infty$ et $\overline{\kappa}(B)=1$.

    Lorsque $\overline{\kappa}(B)=-\infty$, il r\'esulte de [K-McK99] que $\bar B$ est recouverte par des courbes rationnelles rencontrant $D$ en au plus un point (unibranche). Si le point g\'en\'erique de $\bar B$ est contenu dans $2$ telles courbes distinctes (au moins), alors $(\bar B/D)$ est sp\'eciale, et $g$ isotriviale (par [Kov00], par exemple). 
    
    Sinon, il existe (apr\`es \'eventuelle modification) une fibration $f:(\bar B/D)\to C$ dont la fibre g\'en\'erique est l'une des courbes rationnelles pr\'ec\'edentes. Et $D=D^h+D^v$ est alors r\'eunion de sa partie $f$-verticale $D^h$, qui est soit vide, soit une section de $f$ (puisque les courbes $\Delta$-rationnelles de $\bar B$ rencontrent $\Delta$ en au plus un point unibranche), tandis que la partie $f$-verticale $D^v$ est effective, contenue dans une r\'eunion finie de fibres de $f$. Soit $F_c\cong \bP^1$ une fibre g\'en\'erique lisse de $f$. La restriction de $\Omega^1_{(\bar B/D)}$ \`a $F=F_c,c\in C$ est une extension de $\cal O$$_F(-d)$ par $\cal O$$_F=f^*(T^*_{C,c})$, avec $d=1$ si $D^h\neq 0$, et $d=2$ sinon. Donc, pour tout $N>0$, les sections de la restriction de $S^N_{1}(B)$ \`a  $F$ forment un espace vectoriel complexe de dimension $1$ engendr\'e par $f^*((T^*_{C,c})^{\otimes N})$. 
    
    Les sections de $S^N_{1}(B)$ sont donc de la forme: $f^*(N.(K_C+\Delta_C))$, pour $\Delta_C$ un diviseur effectif sur $C$. Par la proposition \ref{k=L}, et l'exemple \ref{netcourb}, on a l'inclusion: $\Delta_C\leq \Delta(f,D)$. 
    
    Si $S^N_{1}(B)$ admet deux sections non nulles $s,t$ telles que $t=u.s$, pour $u$ m\'eromorphe sur $\bar B$, c'est donc que $1=\kappa(C\vert\Delta_C)\leq\kappa(C\vert\Delta(f,D))\leq1$, et $(\bar B/D)=B$ n'est donc pas sp\'eciale.

    Supposons donc d\'esormais que $\overline{\kappa}(B)=1$. Soit $f:(\bar B/D)\to C$ la fibration de Moishezon-Iitaka. Nous la supposons non isotriviale.  Il suffit alors de d\'emontrer, comme ci-dessus, que $L\subset f^*(N.K_{C})$ au-dessus du point g\'en\'erique $c$ de $C$: puisque l'application de Kodaira-Spencer $ks_f(c)$ n'est pas nulle en $c$, les classes successives des extensions d\'eduites de la filtration naturelle de quotients $(TF^*)^{\otimes j}\otimes (f^*(TC_c))^{\otimes (N-j)}$ de $Sym^N(\Omega^1_{\bar B}){\vert F}$ comme extension de $TF^*$ par $f^*(TC_c)$ sont donc aussi non-nulles. Les sections de $Sym^N(\Omega^1_{\bar B}){\vert F}$ se r\'eduisent donc \`a celles de $(f^*(TC_c))^{\otimes N}=f^*(K_{C,c}^{\otimes N})$ $\square$

    \begin{re} Il est int\'eressant de comprendre si la conclusion du lemme \ref{sym} subsiste dans le cas o\`u $f:(\bar B/D)\to C$ n'est pas isotriviale. Cette question est \'etroitement li\'ee a la la remarque \ref{diff'}.3 ci-dessus. Cependant, pour l'\'etude de la conjecture \ref{conjmod}, l'approche de [J-K 09], qui pr\'ecise le th\'eor\` eme \ref{vz} en montrant que $L$ ``provient" de la ``base" de la variation de $g$, est plus naturelle (elle fournit imm\'ediatemment \ref{B2} sans utiliser [K-McK 99] lorsque $\overline{\kappa}(B)=-\infty$, par exemple).
    \end{re}

\subsection{Questions}

On rappelle les questions suivantes:

\subsubsection{Orbifoldes log-canoniques.} 

Voir \ref{qsing} et \ref{qlc}

    \subsubsection{\' Equivalence bim\'eromorphe des bases orbifoldes.}

    Voir \ref{eqbimbasorb}. La question est cruciale pour les trois fibrations fondamentales consid\'er\'ees ici: Moishezon-Iitaka, quotient $\kappa$-rationnel et coeur.

    \subsubsection{Courbes $\Delta$-rationnelles.} 
    
    Voir les questions et conjectures: \ref{qqrat}, \ref{ib}, \ref{-infty}, \ref{cass}, \ref{sempos}, \ref{qs}, \ref{ht'}, \ref{rercc}.
    
    \subsubsection{Additivit\'e orbifolde.}

    Voir \ref{cnmorb}.

    \subsubsection{Finitude des fibrations de type g\'en\'eral.} 
    
    Voir \ref{fintg}.
    
    \subsubsection{Hyperbolicit\'e alg\'ebrique}

     Voir question \ref{qha}

\section{BIBLIOGRAPHIE}


[Abr 07]D. Abramovich. Birational geometry for number theorists. math. AG/0701105v2

[AV 98] D. Abramovich-A. Vistoli.Compactifying the space of stable maps. math. AG/9811059

[Ba75] D. Barlet. Espace analytique r\'eduit des cycles analytiques complexes compacts d'un espace analytique de dimension finie. LNM 482 (1975), 1-158.

[B-H-H 87] G. Barthel, F. Hirzebruch, T. H\" ofer. Geraden-Konfigurationen und Algebraische Fl\" achen. Aspekte der Mathematik. Band D4. Friedr. Vieweg u. Sohn. Braunschweig-Wiesbaden (1987). (Eine Ver\" offentlichung des Max-Planck Instituts f\" ur Mathematik, Bonn).

[B87] F. Beukers. Ternary forms equations. J.Number Theory 54 (1995),113-133 .

[BCHM06]C.Birkar-P.Cascini-C.Hacon-J. McKernan. Existence of minimal models for varieties of log general type. math. AG/0610203

[Bo79] F. Bogomolov. Holomorphic tensors and vector bundles on projective varieties. Math. Ussr Izv. 13 (1979), 499-555.

[Br 08] A. Broustet. Communication orale (Novembre 2008).

[Ca 92] F. Campana. Connexit\'e rationnelle des vari\'et\'es de Fano. Ann. Sc. ENS. 25 (1992), 539-545.

[Ca93] F.Campana.Remarques sur les groupes de k\" ahler nilpotents. Ann. Sc. ENS 28(1993), 307-316.

[Ca94] F.Campana.Remarques sur le rev\^etement universel des vari\'et\'es k\" ahl\'eriennes compactes. Bull. SMF 122(1994), 255-284.

[Ca98] F.Campana. $\sG$-connectedness of compact k\" ahler manifolds. Cont. Math. 241(1999), 85-97.

[Ca01]F.Campana. Special varieties and classification theory. math.AG/0110151.

[Ca04]F.Campana. Orbifolds, special varieties and classification theory. Ann.Inst. Fourier. 54 (2004), 499-665.

[Ca05] F. Campana. Fibres multiples sur les surfaces. Man. Math. 117(2005), 429-461.(arXiv: math/0410469).

[Ca 09] F. Campana. Quotients r\'esolubles ou nilpotents des groupes de K\" ahler orbifoldes. arXiv: 0903.0560

[C-P05] F.Campana-M.P\u aun. Vari\'et\'es faiblement sp\'eciales \`a courbes enti\`eres d\'eg\'en\'er\'ees. Comp. Math.143 (2007), 95-111. (arXiv:math.AG/0512124).

[C-Pe 06] F. Campana-T. Peternell. Geometric stability of the cotangent bundle and the universal cover of a projective manifold. (arXiv: math/0405093).

[C-W05]F.Campana-J. Winkelmann. A Brody theorem for orbifolds. A paraitre \`a Man. Math. (arXiv:math/0604571).

[Cla06]B.Claudon. Invariance for multiples of the twisted canonical bundle. Ann. Inst. Fourier 57 (2007), 289-300. (arXiv: Math. AG/0511736). 

[Cla08]B.Claudon. Gamma-reduction for smooth orbifolds. Manuscripta Math. 127 (2008), 521-532.(disponible sur arXiv: 0801.2894).


[DG98]H. Darmon-A. Granville. On the equations $z^m=F(x,y)$ and $Ax^p+By^q=Cz^r$. Bull. London Math. Soc. 27(19995), 513-543.

[De74] P. Deligne. Th\'eorie de Hodge II. Publ. IHES 40 (1972), 5-57.

[De 79] P. Deligne. Le groupe fondamental du compl\'ementaire d'une courbe plane n'ayant que des points doubles ordinaires est ab\'elien. Expos\'e Bourbaki 543 (Novembre 1979). Tome 22. LNM 642, 1-10. (disponible sur numdam)

[Delz06]T. Delzant. L'invariant de Bieri-Neuman-Strebel des groupes fondamentaux des vari\'et\'es K\" ahl\'eriennes. arXiv:math.DG/0603038

[Fuj78] T.Fujita. On K\" ahler fibre spaces over curves. J. Math. Soc. Jap. 30 (1978), 779-794.

[GHS03] T.Graber-J.Harris-J.Starr. Families of rationally connected varieties. J. Amer. Math. Soc. 16(2003), 57-67.

[G-K-K08] D. Greb, S. Kebekus, S. Kov\`acs. Extension theorems for differential forms and Bogomolov-Sommese vanishing on Log-canonical varieties. arXiv: 0808.3647.

[J-K 09] K. Jabbush, S. Kebekus. Families over special base manifolds and a conjecture of Campana. arXiv:0905.1746

[Kaw81] Y.Kawamata. Characterisation of Abelian Varieties. Comp. Math. (1981), 253-276.

[Kaw98] Y.Kawamata. Subadjunction of log-canonical divisors II. Amer. J. Math. 120 (1998), 893-899.

[Ke-Kov06] S. Kebekus-S. Kov\`acs. Families of canonically polarized varieties over surfaces. (arXiv:math/0511378). To appear in Inv. Math. 

[K-McK99] S. Keel-J. McKernan. Rational curves on quasi-projective surfaces. Memoirs of the AMS 669 (1999).

[Ko93]J. Koll\' ar. Shafarevitch maps and plurigenera of algebraic varieties. Inv. Math. 113 (1993), 177-215.

[KoMiMo92]J.Koll\' ar-Y.Miyaoka--S.Mori. Rationally connected varieties. J. Alg. Geom. 1 (1992), 429-448.

[Kov96] S.Kov\`acs. Smooth families over rational and elliptic curves. JAG 5 (1996), , 369-385.

[Kov00] S.Kov\`acs. Algebraic hyperbolicity of fine moduli spaces. JAG 9 (2000), 169-174.

[La86] S. Lang. Hyperbolic and Diophantine Analysis. Bull. AMS 14(1986), 159-205.

[Lieb78]D.Lieberman. Compactness of the Chow Scheme. LNM 670 (1975), 140-186.

[Miy 87] Y. Miyaoka. Deformation of a morphism along a foliation. 
Proc. Symp. Pure Math. vol. 46, 245-268 (1987)

[Mi-Mo 86] Y. Miyaoka-S. Mori. A numerical criterion for uniruledness. Ann. Math. 124 (1986), 65-69.

[N87]M.Namba. Branched coverings and algebraic functions. Pitman research Notes in Mathematics series 161. Longman Scientific and Technical (1987).

[Pa05]M.P\u aun.Siu's  invariance of plurigenera: a one-tower proof . A paraitre au J. Diff.Geom.

[P-R06]G.Pacienza, E. Rousseau. On the logarithmic Kobayashi conjecture. arXiv:math/0603712

[R72] M. Raynaud. Flat modules in algebraic geometry. Comp. Math. 24 (1972), 11-31.

[R08] E. Rousseau. Hyperbolicity of geometric orbifolds. arXiv:08091356.

[Sak 74] F. Sakai. Degeneracy of holomorphic maps with ramifications. Inv. Math. 28 (1974), 213-229.

[SB 92] N. Shepherd-Barron. Miyaoka's theorem on the semi-negativity of $T_X$. Ast\'erisque 211, 103-114 (1992)

[Siu02]Y.T.Siu.  Extension of twisted pluricanonical sections  with plurisubharmonic weights. Complex Geometry (G\" ottingen 2000),223-277, Springer (2002). 

[Tak00]S. Takayama. Simple connectedness of weak Fano varieties. J. Alg. Geom. 9 (2000), 403-407.

[U75] K. Ueno. Classification theory of complex analytic manifolds. LNM 439 (1975)

[Vie83] E. Viehweg. Weak positivity and the additivity f the Kodaira dimension for certain fibre spaces. Ad. Studies in Pure Math. 1 (1983), 329-353.

[V-Z02] E. Viehweg-K. Zuo. Base spaces of non-isotrivial families of smooth minimal models. Complex geometry (G\" ottingen 2000), 279-328. Springer Verlag 2002.


\

{ F.~Campana

D\'epartement de Math\'ematiques

Universit\'e Nancy 1

BP 239

F-54506 Vandoeuvre-l\`es-Nancy Cedex}

campana@iecn.u-nancy.fr

 \end{document}